\documentclass[english]{smfart}
\usepackage{smfthm,smfenum,bull}
\usepackage{amsfonts,amssymb,amscd,amsmath,latexsym}
\usepackage{enumerate}
\usepackage[english]{babel}
\usepackage[dvips]{graphicx}

\theoremstyle{plain}

\def\surl#1_#2{\mathrel{\mathop{\kern 0pt #1}\limits_{#2}}}

\author{Franck LESIEUR}
\address{LMNO\\
         Universit\'e de Caen\\
         BP 5186\\
         F-14032 Caen Cedex}
\email{franck.lesieur@math.unicaen.fr}
\urladdr{http://www.math.unicaen.fr/~lesieur}
\title[MEASURED QUANTUM GROUPOIDS]{Measured quantum groupoids}
\alttitle{Groupo\"{\i}des quantiques mesur\'es}

\begin{document}
\frontmatter

\begin{abstract}
In this article, part of the author's thesis \cite{Les}, we
propose a definition for measured quantum groupoid. The aim is the
construction of objects with duality including both quantum
groups and groupoids. We base ourselves on J. Kustermans and S.
Vaes' works about locally compact quantum groups that we
generalize thanks to formalism introduced by M. Enock and J.M.
Vallin in the case of inclusion of von Neumann algebras. From a
structure of Hopf-bimodule with left and right invariant
operator-valued weights, we define a fundamental
pseudo-multiplicative unitary. We introduce the notion of
quasi-invariant weight on the basis and, then, we construct an
antipode with polar decomposition, a coinvolution, a scaling
group, a modulus and a scaling operator. This theory is
illustrated with different examples. Duality of measured quantum
groupoids will be discussed in a forthcoming article.
\end{abstract}

\begin{altabstract}
Dans cet article, extrait d'une partie de la th\`{e}se \cite{Les}
de l'auteur, on propose une d\'efinition des groupo\"{\i}des
quantiques mesur\'es. L'objectif est la construction d'objets,
munis d'une dualit\'e, qui englobent \`a la fois les
groupo\"{\i}des et les groupes quantiques. On s'appuie sur les
travaux de J. Kustermans et S. Vaes concernant les groupes
quantiques localement compacts qu'on g\'en\'eralise gr\^{a}ce au
formalisme introduit par M. Enock et J.M. Vallin \`a propos des
inclusions d'alg\`ebres de von Neumann. \`{A} partir d'un
bimodule de Hopf muni de poids op\'eratoriels invariants \`a
gauche et \`a droite, on d\'efinit un unitaire
pseudo-multiplicatif fondamental. On introduit la notion de poids
quasi-invariant sur la base et on construit une antipode avec
d\'ecomposition polaire, une coinvolution, un groupe d'\'echelle,
un module et un op\'erateur d'\'echelle. Cette th\'eorie est
illustr\'ee par diff\'erents exemples. La dualit\'e de ces objets
sera discut\'ee dans un prochain article.
\end{altabstract}

\altkeywords{Groupo\"{\i}des
quantiques ; antipode ; unitaire pseudo-multiplicatif}
\thanks{The author is mostly indebted to Michel Enock, Stefaan Vaes,
Leonid Va\u{\i}nerman and Jean-Michel Vallin for many fruitful
conversations.}

\maketitle
\mainmatter

\section{Introduction}

\subsection{Historique}
Theory of quantum groups has lot of developments in operator
algebras setting. Many contributions are given by \cite{KaV},
\cite{W6}, \cite{ES}, \cite{MN}, \cite{BS}, \cite{W1}, \cite{W},
\cite{VD2}, \cite{KV1}. In particular, J. Kustermans and S. Vaes'
work is crucial: in \cite{KV1}, they propose a simple definition
for locally compact quantum groups which gathers all known
examples (locally compact groups, quantum compacts groupe
\cite{W1}, quantum group $ax+b$ \cite{W3}, \cite{W4}, Woronowicz'
algebra \cite{MN}...) and they find a general framework for
duality of theses objects. The very few number of axioms gives
the theory a high manageability which is proved with recent
developments in many directions (actions of locally compact
quantum groups \cite{Vae2}, induced co-representations
\cite{Ku2}, cocycle bi-crossed products \cite{VaV}). They
complete their work with a theory of locally compact quantum
groups in the von Neumann setting \cite{KV2}.

In geometry, groups are rather defined by their actions.
Groupoids category contains groups, group actions and equivalence
relation. It is used by G.W Mackey and P. Hahn (\cite{Ma},
\cite{Ha1} and \cite{Ha2}), in a measured version, to link theory
of groups and ergodic theory. Locally compact groupoids and the
operator theory point of view are introduced and studied by J.
Renault in \cite{R1} and \cite{R2}. It covers many interesting
examples in differential geometry \cite{C2} e.g holonomy groupoid
of a foliation.

In \cite{V1}, J.M Vallin introduces the notion of Hopf bimodule
from which he is able to prove a duality for groupoids. Then, a
natural question is to construct a category, containing quantum
groups and groupoids, with a duality theory.

In the quantum group case, duality is essentially based on a
multiplicative unitary \cite{BS}. To generalize the notion up to
the groupoid case, J.M Vallin introduces pseudo-multiplicative
unitaries. In \cite{V2}, he exhibits such an object coming from
Hopf bimodule structures for groupoids. Technically speaking,
Connes-Sauvageot's theory of relative tensor products is
intensively used.

In the case of depth $2$ inclusions of von Neumann algebras, M.
Enock and J.M Vallin, and then, M. Enock underline two "quantum
groupoids" in duality. They also use Hopf bimodules and
pseudo-multiplicative unitaries. At this stage, a non trivial
modular theory on the basis (the equivalent for units of a
groupoid) is revealed to be necessary and a simple generalization
of axioms quantum groups is not sufficient to construct quantum
groupoids category: we have to add an axiom on the basis
\cite{E1} i.e we use a special weight to do the construction. The
results are improved in \cite{E3}.

In \cite{E2}, M. Enock studies in detail pseudo-multiplicative
unitaries and introduces an analogous notion of S. Baaj and G.
Skandalis' regularity. In quantum groups, the fundamental
multiplicative unitary is weakly regular and manageable in the
sense of Woronowicz. Such properties have to be satisfied in
quantum groupoids. Moreover, M. Enock defines and studies compact
(resp. discrete) quantum groupoids which have to enter into the
general theory.

Lot of works have been led about quantum groupoids but
essentially in finite dimension. We have to quote weak Hopf
$\text{C}^*$-algebras introduced by G. B\"{o}hm, F. Nill and K.
Szlach\'{a}nyi \cite{BNSz}, \cite{BSz}, and then studied by F.
Nill and L. Va\u{\i}nerman \cite{Ni}, \cite{N}, \cite{NV1},
\cite{NV2}. J.M Vallin develops a quantum groupoids theory in
finite dimension thanks to multiplicative partial isometries
\cite{V3}, \cite{V4}. He proves that his theory coincide exactly
with weak Hopf $\text{C}^*$-algebras.

\subsection{Aims and Methods}
In this article, we propose a definition for measured quantum
groupoids in any dimensions. "Measured" means we are in the von
Neumann setting and we assume existence of the analogous of a
measure. We use a similar approach as J. Kustermans and S. Vaes'
theory with the formalism of Hopf bimodules and
pseudo-multiplicative unitaries. We develop the theory by
constructing all fundamental objects and we give some examples.
In a forthcoming article, we will study duality within the
category.

\subsection{Study plan}
After brief recalls about tools and technical points, we define
objects we will use. We start by associating a
pseudo-multiplicative unitary to every Hopf bimodule with
invariant operator-valued weights. This unitary gathers all
informations on the structure so that we can re-construct von
Neumann algebra and co-product. Then a measured quantum groupoid
will be a Hopf bimodule with invariant operator-valued weights
which are "adapted" in a certain sense. This hypothesis
corresponds to the choice of a special weight on the basis to do
the constructions.

Thanks to this axiom, we are able to construct fundamental
objects of the structure: first of all, the antipode $S$, the
polar decomposition of which is given by a co-involution $R$ and
a one-parameter group of automorphisms called scaling group
$\tau$. In particular, we show that $S,R$ and $\tau$ are
independent of operator-valued weights. Also, we introduce a
modulus, which corresponds to modulus of groupoids, and a scaling
operator, affiliated to the hyper-center of the Hopf bimodule.
They come from Radon-Nikodym's cocycle of right invariant
operator-valued weight with respect to left invariant one thanks
to proposition 5.2 of \cite{Vae}. So, it is the existence of a
suitable weight on the basis $N$ which allows us to construct
modulus like in the groupoid case with a quasi-invariant measure
on $G^{\{0\}}$. Scaling operator is the object which corresponds
to scaling factor in locally compact quantum groups. We are also
able to prove uniqueness of invariant operator-valued weight up to
an element of basis center. Finally, we prove a "manageability"
property of the fundamental pseudo-multiplicative unitary.

We have a lot of examples for locally compact quantum groups
thanks to Woronowicz \cite{W2}, \cite{W3}, \cite{W4}, \cite{W5}
and the cocycle bi-crossed product due to S. Vaes and L.
Va\u{\i}nerman \cite{VaV}. Theory of measured quantum groupoids
has also a lot of examples: groupoids, weak Hopf
$\text{C}^*$-algebras, quantum groups, quantum groupoids of
compact (resp. discrete) type\ldots which are characterized in
the general theory. Depth $2$ inclusions of von Neumann algebras
with semi-finite basis also enter into our setting. Finally, we
state stability of the category by direct sum (which reflects the
stability of groupoids under disjoint unions), finite tensor
product and direct integrals. Then, we are able to construct new
examples: in particular we can exhibit quantum groupoids with non
scalar scaling operator.

\section{Recalls}
\subsection{Weights and operator-valued weights \cite{St}, \cite{T}}

Let $N$ be a von Neumann and $\psi$ a normal, semi-finite
faithful (n.s.f) weight on $N$; we denote by ${\mathcal
N}_{\psi}$, ${\mathcal M}_{\psi}$, $H_{\psi}$, $\pi_{\psi}$,
$\Lambda_{\psi}$, $J_{\psi}$, $\Delta_{\psi}\ldots$ canonical
objects of Tomita's theory with respect to (w.r.t) $\psi$.

\begin{defi}
Let denote by ${\mathcal T}_{\psi}$ {\bf Tomita's algebra} w.r.t
$\psi$ defined by:
$$\{ x \in {\mathcal N}_{\psi } \cap {\mathcal N}_{\psi }^*|\ x \text{
analytic w.r.t }\sigma^{\psi}\text{ such that } \sigma_z^{\psi}(x)
\in {\mathcal N}_{\psi } \cap {\mathcal N}_{\psi }^* \text{ for
all } z \in {\mathbb C} \}$$
\end{defi}

By (\cite{St}, 2.12), we have the following approximating result:

\begin{lemm}\label{analyse}
For all $x\in {\mathcal N}_{\psi}$, there exists a sequence
$(x_n)_{n\in \mathbb{N}}$ of ${\mathcal T}_{\psi}$ such that:

\begin{center}
\begin{minipage}{13cm}
\begin{enumerate}[i)]
\item $||x_n||\leq ||x||$ for all $n\in\mathbb{N}$;
\item $(x_n)_{n\in\mathbb{N}}$ converges to $x$ in the strong topology;
\item $(\Lambda_{\psi}(x_n))_{n\in\mathbb{N}}$ converges to
$\Lambda_{\psi}(x)$ in the norm topology of $H_{\psi}$.
\end{enumerate}
\end{minipage}
\end{center}
Moreover, if $x\in {\mathcal N}_{\psi}\cap {\mathcal
N}_{\psi}^*$, then we have:

\begin{center}
\begin{minipage}{13cm}
\begin{enumerate}[iv)]
\item $(x_n)_{n\in\mathbb{N}}$ converges to $x$ in the *-strong topology;
\item $(\Lambda_{\psi}(x_n^*))_{n\in\mathbb{N}}$ converges to
$\Lambda_{\psi}(x^*)$ in the norm topology of $H_{\psi}$.
\end{enumerate}
\end{minipage}
\end{center}
\end{lemm}

Let $N\subset M$ be an inclusion of von Neumann algebras and $T$
a normal, semi-finite, faithful (n.s.f) operator-valued weight
from $M$ to $N$. We put:
$${\mathcal N}_T=\left\lbrace x\in M\,
/\, T(x^*x)\in N^+ \right\rbrace  \text{ and } {\mathcal
M}_T={\mathcal N}_T^*{\mathcal N}_T$$ We can define a n.s.f
weight $\psi\circ T$ on $M$ in a natural way.  Let us recall
theorem 10.6 of \cite{EN}:
\begin{prop}\label{preintro}
Let $N\subset M$ be an inclusion of von Neumann algebras and $T$
be a normal, semi-finite, faithful (n.s.f) operator-valued weight
from $M$ to $N$ and $\psi$ a n.s.f weight on $N$. Then we have:

\begin{enumerate}[i)]
\item for all
$x\in {\mathcal N}_T$ and $a\in {\mathcal N}_{\psi}$, $xa$
belongs to ${\mathcal N}_T\cap {\mathcal N}_{\psi\circ T}$, there
exists $\Lambda_T(x)\in\text{Hom}_{N^o}(H_{\psi},H_{\psi\circ
T})$ such that: $$\Lambda_T(x)\Lambda_{\psi}(a)=\Lambda_{\psi\circ
T}(xa)$$ and $\Lambda_T$ is a morphism of $M-N$-bimodules from
${\mathcal N}_T$ to $\text{Hom}_{N^o}(H_{\psi},H_{\psi\circ T})$;
\item ${\mathcal N}_T\cap {\mathcal N}_{\psi\circ T}$ is a weakly
dense ideal of $M$ and $\Lambda_{\psi\circ T}({\mathcal N}_T\cap
{\mathcal N}_{\psi\circ T})$ is dense in $H_{\psi\circ T}$,
$\Lambda_{\psi\circ T}({\mathcal N}_T\cap {\mathcal N}_{\psi\circ
T}\cap {\mathcal N}^*_T\cap {\mathcal N}^*_{\psi\circ T})$ is a
core for $\Delta_{\psi\circ T}^{1/2}$ and $\Lambda_T({\mathcal
N}_T)$ is dense in $\text{Hom}_{N^o}(H_{\psi},H_{\psi\circ T})$
for the s-topology defined by (\cite{BDH}, 1.3);
\item for all $x\in {\mathcal N}_T$ and
$z\in {\mathcal N}_T\cap {\mathcal N}_{\psi\circ T}$, $T(x^*z)$
belongs to ${\mathcal N}_{\psi}$ and:
$$\Lambda_T(x)^*\Lambda_{\psi\circ T}(z)=\Lambda_{\psi}(T(x^*z))$$
\item for all $x,y\in {\mathcal
N}_T$: $$\Lambda_T(y)^*\Lambda_T(x)=\pi_{\psi}(T(x^*y)) \text{
and } ||\Lambda_T(x)||=||T(x^*x)||^{1/2}$$ and $\Lambda_T$ is
injective.
\end{enumerate}
\end{prop}

Let us also recall lemma 10.12 of \cite{EN}:

\begin{prop}
Let $N\subseteq M$ be an inclusion of von Neumann algebras, $T$ a
n.s.f operator-valued weight from $M$ to $N$, $\psi$ a n.s.f
weight on $N$ and $x\in {\mathcal M}_T\cap {\mathcal
M}_{\psi\circ T}$. If we put:
$$x_n=\sqrt\frac{n}{\pi}\int^{+\infty}_{-\infty}\!
e^{-nt^2}\sigma_t^{\psi\circ T}(x)\ dt$$ then $x_n$ belongs to
${\mathcal M}_T\cap {\mathcal M}_{\psi\circ T}$ and is analytic
w.r.t $\psi\circ T$. The sequence converges to $x$ and is bounded
by $||x||$. Moreover, $(\Lambda_{\psi\circ
T}(x_n))_{n\in\mathbb{N}}$ converges to $\Lambda_{\psi\circ T}(x)$
and  $\sigma_z^{\psi\circ T}(x_n)\in {\mathcal M}_T\cap {\mathcal
M}_{\psi\circ T}$ for all $z\in\mathbb{C}$.
\end{prop}

\begin{defi}
The set of $x\in {\mathcal N}_{\Phi } \cap {\mathcal N}_{\Phi }^*
\cap {\mathcal N}_T \cap {\mathcal N}_T^*$, analytic w.r.t
$\sigma^{\Phi}$ such that $\sigma_z^{\Phi}(x)\in {\mathcal
N}_{\Phi } \cap {\mathcal N}_{\Phi }^* \cap {\mathcal N}_T \cap
{\mathcal N}_T^*$ for all $z\in {\mathbb C}$ is denoted by
${\mathcal T}_{\Phi}$ and is called {\bf Tomita's algebra} w.r.t
$\psi\circ T=\Phi$ and $T$.
\end{defi}
Lemma \ref{analyse} is still satisfied with Tomita's algebra w.r.t
$\Phi$ and $T$.

\subsection{Spatial theory \cite{C1},
\cite{S2}, \cite{T}}\label{intre}

Let $\alpha$ be a normal, non-degenerated representation of $N$
on a Hilbert space $H$. So, $H$ becomes a left $N$-module and we
write $\ _{\alpha}H$.

\begin{defi}\cite{C1}
An element $\xi$ of $\ _{\alpha}H$ is said to be bounded w.r.t
$\psi$ if there exists $C\in\mathbb{R}^+$ such that, for all $y\in
{\mathcal N}_{\psi}$, we have $||\alpha(y)\xi||\leq C||
\Lambda_{\psi}(y) ||$. The set of {\bf bounded elements} w.r.t
$\psi$ is denoted by $D(_{\alpha}H,\psi)$.
\end{defi}

By \cite{C1} (lemma 2), $D(_{\alpha}H,\psi)$ is dense in $H$ and
$\alpha(N)'$-stable. An element $\xi$ of $D(_{\alpha}H,\psi)$
gives rise to a bounded operator $R^{\alpha,\psi}(\xi)$ of
$Hom_N(H_{\psi},H)$ such that, for all $y\in {\mathcal N}_{\psi}$:
$$R^{\alpha,\psi}(\xi)\Lambda_{\psi}(y)=\alpha(y)\xi$$
For all $\xi,\eta \in D(_{\alpha}H,\psi)$, we put:
$$\theta^{\alpha,\psi}(\xi,\eta)=R^{\alpha,\psi}(\xi)R^{\alpha,\psi}(\eta)^*\text{ and }
<\xi,\eta>_{\alpha,\psi}=R^{\alpha,\psi}(\eta)^*R^{\alpha,\psi}(\xi)^*$$
By \cite{C1} (lemma 2), the linear span of
$\theta^{\alpha,\psi}(\xi,\eta)$ is a weakly dense ideal of
$\alpha(N)'$. $<\xi,\eta>_{\alpha,\psi}$ belongs to
$\pi_{\psi}(N)'=J_{\psi}\pi_{\psi}(N)J_{\psi}$ which is
identified with the opposite von Neumann algebra $N^o$. The
linear span of $<\xi,\eta>_{\alpha,\psi}$ is weakly dense in
$N^o$.

\noindent By \cite{C1} (proposition 3), there exists a net
$(\eta_i)_{i\in I}$ of $D(_{\alpha}H,\psi)$ such that:
$$\sum_{i\in I}\theta^{\alpha,\psi}(\eta_i,\eta_i)=1$$ Such a net is called a
$(N,\psi)$-\textbf{basis} of $\ _{\alpha}H$. By \cite{EN}
(proposition 2.2), we can choose $\eta_i$ such that
$R^{\alpha,\psi}(\eta_i)$ is a partial isometry with two-by-two
orthogonal final supports and such that
$<\eta_i,\eta_j>_{\alpha,\psi}=0$ unless $i=j$. In the following,
we assume these properties hold for all $(N,\psi)$-basis of $\
_{\alpha}H$.

Now, let $\beta$ be a normal, non-degenerated anti-representation
from $N$ on $H$. So $H$ becomes a right $N$-module and we write
$H_{\beta}$. But $\beta$ is also a representation of $N^o$. If
$\psi^o$ is the n.s.f weight on $N^o$ coming from $\psi$ then
${\mathcal N}_{\psi^o}={\mathcal N}_{\psi}^*$ and we identify
$H_{\psi^o}$ with $H_{\psi}$ thanks to:
$$(\Lambda_{\psi^o}(x^*)\mapsto J_{\psi}\Lambda_{\psi}(x))$$

\begin{defi}\cite{C1}
An element $\xi$ of $H_{\beta}$ is said to be bounded w.r.t
$\psi^o$ if there exists $C\in\mathbb{R}^+$ such that, for all
$y\in {\mathcal N}_{\psi}$, we have $||\beta(y^*)\xi||\leq C||
\Lambda_{\psi}(y) ||$. The set of {\bf bounded elements} w.r.t
$\psi^o$ is denoted by $D(H_{\beta},\psi^o)$.
\end{defi}

$D(_{\alpha}H,\psi)$ is dense in $H$ and $\beta(N)'$-stable. An
element $\xi$ of $D(H_{\beta},\psi^o)$ gives rise to a bounded
operator $R^{\beta,\psi^o}(\xi)$ of $Hom_{N^o}(H_{\psi},H)$ such
that, for all $y\in {\mathcal N}_{\psi}$:
$$R^{\beta,\psi^o}(\xi)\Lambda_{\psi}(y)=\beta(y^*)\xi$$
For all $\xi,\eta \in D(H_{\beta},\psi^o)$, we put:
$$\theta^{\beta,\psi^o}(\xi,\eta)=R^{\beta,\psi^o}(\xi)R^{\beta,\psi^o}(\eta)^*\text{ and }
<\xi,\eta>_{\beta,\psi^o}=R^{\beta,\psi^o}(\eta)^*R^{\beta,\psi^o}(\xi)^*$$
The linear span of $\theta^{\beta,\psi^o}(\xi,\eta)$ is a weakly
dense ideal of $\beta(N)'$. $<\xi,\eta>_{\beta,\psi^o}$ belongs to
$\pi_{\psi}(N)$ which is identified with $N$. The linear span of
$<\xi,\eta>_{\beta,\psi^o}$ is weakly dense in $N$. In fact, we
know that $<\xi,\eta>_{\beta,\psi^o}\in {\mathcal M}_{\psi}$ by
\cite{C1} (lemma 4) and by \cite{S2} (lemma 1.5), we have
$$\Lambda_{\psi}(<\xi,\eta>_{\beta,\psi^o})=R^{\beta,\psi^o}(\eta)^*\xi$$
A net $(\xi_i)_{i\in I}$ of $\psi^o$-bounded elements of is said
to be a $(N^o,\psi^o)$-basis of $H_{\beta}$ if:
$$\sum_{i\in I}\theta^{\beta,\psi^o}(\xi_i,\xi_i)=1$$
and if $\xi_i$ such that $R^{\beta,\psi^o}(\xi_i)$ is a partial
isometry with two-by-two orthogonal final supports and such that
$<\xi_i,\xi_j>_{\alpha,\psi}=0$ unless $i=j$. Therefore, we have:
$$R^{\beta,\psi^o}(\xi_i)=\theta^{\beta,\psi^o}(\xi_i,\xi_i)R^{\beta,\psi^o}(\xi_i)=
R^{\beta,\psi^o}(\xi_i)<\xi_i,\xi_i>_{\beta,\psi^o}$$ and, for
all $\xi\in D(H_{\beta},\psi^o)$:
$$\xi=\sum_{i\in I}
R^{\beta,\psi^o}(\xi_i)\Lambda_{\psi}(<\xi,\xi_i>_{\beta,\psi^o})$$

\begin{prop}(\cite{E2}, proposition 2.10)
Let $N\subseteq M$ be an inclusion of von Neumann algebras and
$T$ be a n.s.f operator-valued weight from $M$ to $N$. There
exists a net $(e_i)_{i\in I}$ of ${\mathcal N}_{T}\cap {\mathcal
N}_{T}^*\cap {\mathcal N}_{\psi\circ T}\cap {\mathcal
N}_{\psi\circ T}^*$ such that $\Lambda_{T}(e_i)$ is a partial
isometry, $T(e_j^*e_i)=0$ unless $i=j$ and with orthogonal final
supports of sum $1$. Moreover, we have $e_i=e_iT(e_i^*e_i)$ for
all $i\in I$, and, for all $x\in {\mathcal N}_{T}$:
$$\Lambda_{T}(x)=\sum_{i\in I}\Lambda_{T}(e_i)T(e_i^*x)\quad\text{ and }\quad
x=\sum_{i\in I} e_iT(e_i^*x)$$ in the weak topology. Such a net
is called a basis for $(T,\psi^o)$. Finally, the net
$(\Lambda_{\psi\circ T}(e_i))_{i\in I}$ is a $(N^o,\psi^o)$-basis
of $(H_{\psi\circ T})_s$ where $s$ is the anti-representation
which sends $y\in N$ to $J_{\psi\circ T}y^*J_{\psi\circ
T}$.\label{partis}
\end{prop}

\subsection{Relative tensor product \cite{C1}, \cite{S2}, \cite{T}}
Let $H$ and $K$ be Hilbert space. Let $\alpha$ (resp. $\beta$) be
a normal and non-degenerated (resp. anti-) representation of $N$
on $K$ (resp. $H$). Let $\psi$ be a n.s.f weight on $N$. Following
\cite{S2}, we put on $D(H_{\beta},\psi^o)\odot K$ a scalar
product defined by:
$$(\xi_1 \odot \eta_1|\xi_2 \odot \eta_2)=(\alpha(<\xi_1,\xi_2>_{\beta,\psi^o})\eta_1|\eta_2)$$
for all $\xi_1,\xi_2\in D(H_{\beta},\psi^o)$ and
$\eta_1,\eta_2\in K$. We have identified $\pi_{\psi}(N)$ with $N$.

\begin{defi}
The completion of $D(H_{\beta},\psi^o)\odot K$ is called {\bf
relative tensor product} and is denoted by $H\surl{\ _{\beta}
\otimes_{\alpha}}_{\ \psi}K$.
\end{defi}

The image of $\xi\odot\eta$ in $H\surl{\ _{\beta}
\otimes_{\alpha}}_{\ \psi} K$ is denoted by $\xi\surl{\ _{\beta}
\otimes_{\alpha}}_{\ \psi}\eta$. One should bear in mind that, if
we start from another n.s.f weight $\psi'$ on $N$, we get another
Hilbert space which is canonically isomorphic to $H\surl{\
_{\beta} \otimes_{\alpha}}_{\ \psi}K$ by (\cite{S2}, proposition
2.6). However this isomorphism does not send $\xi\surl{\ _{\beta}
\otimes_{\alpha}}_{\ \psi}\eta$ on $\xi\surl{\ _{\beta}
\otimes_{\alpha}}_{\ \ \psi '}\eta$.

By \cite{S2} (definition 2.1), relative tensor product can be
defined from the scalar product:
$$(\xi_1 \odot \eta_1| \xi_2 \odot \eta_2)=(\beta(<\eta_1,\eta_2>_{\alpha,\psi})\xi_1|\xi_2)$$
for all $\xi_1,\xi_2 \in H$ and $\eta_1,\eta_2\in
D(_{\alpha}K,\psi)$ that's why we can define a one-to-one flip
from $H\surl{\ _{\beta} \otimes_{\alpha}}_{\ \psi}K$ onto $K
\surl{\ _{\alpha} \otimes_{\beta}}_{\ \ \psi^o}H$ such that:
$$\sigma_{\psi}(\xi\surl{\ _{\beta} \otimes_{\alpha}}_{\ \psi}\eta)=\eta \surl{\
_{\alpha} \otimes_{\beta}}_{\ \ \psi^o} \xi$$ for all $\xi\in
D(H_{\beta},\psi)$ (resp. $\xi\in H$) and $\eta\in K$ (resp.
$\eta\in D(_{\alpha}K,\psi)$). The flip gives rise at the
operator level to $\varsigma_{\psi}$ from $\mathcal{L}(H \surl{\
_{\beta} \otimes_{\alpha}}_{\ \psi} K)$ onto $\mathcal{L}(K
\surl{\ _{\alpha} \otimes_{\beta}}_{\ \ \psi^o} H)$ such that:
$$\varsigma_{\psi}(X)=\sigma_{\psi}X\sigma_{\psi}^*$$
Canonical isomorphisms of change of weights send
$\varsigma_{\psi}$ on $\varsigma_{\psi '}$ so that we write
$\varsigma_N$ without any reference to the weight on $N$.

For all $\xi\in D(H_{\beta},\psi^o)$ and $\eta\in
D(_{\alpha}K,\psi)$, we define bounded operators:

$$
\begin{aligned}
\lambda^{\beta,\alpha}_{\xi}: K &\rightarrow H\surl{\ _{\beta}
\otimes_{\alpha}}_{\ \psi} K& \text{ and }\quad
\rho^{\beta,\alpha}_{\eta}: H &\rightarrow H\surl{\ _{\beta}
\otimes_{\alpha}}_{\ \psi} K \\
\eta &\mapsto \xi\surl{\ _{\beta} \otimes_{\alpha}}_{\ \psi}\eta&
\xi &\mapsto \xi\surl{\ _{\beta} \otimes_{\alpha}}_{\ \psi}\eta
\end{aligned}$$
Then, we have:
$$(\lambda^{\beta,\alpha}_{\xi})^*\lambda^{\beta,\alpha}_{\xi}=
\alpha(<\xi,\xi>_{\beta,\psi^o})\text{ and }
(\rho^{\beta,\alpha}_{\eta})^*\rho^{\beta,\alpha}_{\eta}=
\beta(<\eta,\eta>_{\alpha,\psi})$$

By \cite{S2} (remark 2.2), we know that $D(_{\alpha}K,\psi)$ is
$\alpha(\sigma_{-i/2}^{\psi}({\mathcal
D}(\sigma_{-i/2}^{\psi})))$-stable and for all $\xi\in H$,
$\eta\in D(_{\alpha}K,\psi)$ and $y\in {\mathcal
D}(\sigma_{-i/2}^{\psi})$, we have:
$$\beta(y)\xi\surl{\ _{\beta} \otimes_{\alpha}}_{\ \psi}\eta=
\xi\surl{\ _{\beta} \otimes_{\alpha}}_{\
\psi}\alpha(\sigma_{-i/2}^{\psi}(y))\eta$$

\begin{lemm}\label{aidintre}
If $\xi'\surl{\ _{\beta} \otimes_{\alpha}}_{\ \psi}\eta=0$ for
all $\xi'\in D(H_{\beta},\psi^o)$ then $\eta=0$.
\end{lemm}

\begin{proof}
For all $\xi,\xi'\in D(H_{\beta},\psi^o)$, we have:
$$\alpha(<\xi',\xi>_{\beta,\psi^o})\eta=(\lambda_{\xi}^{\beta,\alpha})^*
\lambda_{\xi'}^{\beta,\alpha}\eta=(\lambda_{\xi}^{\beta,\alpha})^*(\xi'
\surl{\ _{\beta} \otimes_{\alpha}}_{\ \psi} \eta)=0$$ Since the
linear span of $<\xi',\xi>_{\beta,\psi^o}$ is dense in $N$, we get
$\eta=0$.
\end{proof}

\begin{prop}\label{preli}
Assume $H\neq\{0\}$. Let $K'$ be a closed subspace of $K$ such
that $\alpha(N)K'\subseteq K'$. Then:
$$H\surl{\ _{\beta}
\otimes_{\alpha}}_{\ \psi} K=H\surl{\ _{\beta}
\otimes_{\alpha}}_{\ \psi} K'\quad\Rightarrow\quad K=K'$$
\end{prop}

\begin{proof}
Let $\eta \in K'^{\bot}$. For all $\xi,\xi'\in
D(H_{\beta},\psi^o)$ and $k \in K'$, we have:
$$(\xi \surl{\
_{\beta} \otimes_{\alpha}}_{\ \psi} k|\xi' \surl{\ _{\beta}
\otimes_{\alpha}}_{\ \psi}
\eta)=(\alpha(<\xi,\xi'>_{\beta,\nu^o})k|\eta)=0$$ Therefore, for
all $\xi'\in D(H_{\beta},\psi^o)$, we have:
$$\xi'\surl{\ _{\beta}
\otimes_{\alpha}}_{\ \psi} \eta \in (H\surl{\ _{\beta}
\otimes_{\alpha}}_{\ \psi} K')^{\bot}=(H\surl{\ _{\beta}
\otimes_{\alpha}}_{\ \psi} K)^{\bot}=\{0\}$$ By the previous
lemma, we get $\eta=0$ and $K=K'$.
\end{proof}

Let $x\in\beta(N)'\cap\mathcal{L}(H)$ and $y\in\alpha(N)'\cap
\mathcal{L}(K)$. By \cite{S1}, 2.3 and 2.6, we can naturally
define an operator $x\surl{\ _{\beta}\otimes_{\alpha}}_{\ \psi}y$
on $H\surl{\ _{\beta}\otimes_{\alpha}}_{\ \psi} K$. Canonical
isomorphism of change of weights sends $x\surl{\
_{\beta}\otimes_{\alpha}}_{\ \psi} y$ on $x\surl{\
_{\beta}\otimes_{\alpha}}_{\ \ \psi '}y$ so that we write
$x\surl{\ _{\beta}\otimes_{\alpha}}_{\ N} y$ without any
reference to the weight.

Let $P$ be a von Neumann algebra and $\epsilon$ a normal and
non-degenerated anti-representation of $P$ on $K$ such that
$\epsilon(P)'\subseteq \alpha(N)$. $K$ is equipped with a
$N-P$-bimodule structure denoted by $\ _{\alpha}K_{\epsilon}$.
For all $y\in P$, $1_H \surl{\ _{\beta}\otimes_{\alpha}}_{\ \psi}
\epsilon(y)$ is an operator on $H\surl{\
_{\beta}\otimes_{\alpha}}_{\ \psi}K$ so that we define a
representation of $P$ on $H\surl{\ _{\beta}\otimes_{\alpha}}_{\
\psi} K$ still denoted by $\epsilon$. If $H$ is a $Q-N$-bimodule,
then $H\surl{\ _{\beta}\otimes_{\alpha}}_{\ \psi} K$ becomes a
$Q-P$-bimodule (Connes' fusion of bimodules). If $\nu$ is a n.s.f
weight on $P$ and $\ _{\zeta}L$ a left $P$-module. It is possible
to define two Hilbert spaces $(H\surl{\
_{\beta}\otimes_{\alpha}}_{\ \psi} K) \surl{\
_{\epsilon}\otimes_{\zeta}}_{\ \nu} L$ and $H\surl{\
_{\beta}\otimes_{\alpha}}_{\ \psi} (K \surl{\
_{\epsilon}\otimes_{\zeta}}_{\ \nu} L)$. These two
$\beta(N)'-\zeta(P)'^o$-bimodules are isomorphic. (The proof of
\cite{V1}, lemme 2.1.3, in the case of commutative $N=P$ is still
valid). We speak about associativity of relative tensor product
and we write $H\surl{\ _{\beta}\otimes_{\alpha}}_{\ \psi} K
\surl{\ _{\epsilon}\otimes_{\zeta}}_{\ \nu} L$ without
parenthesis.

We identify $H_{\psi} \surl{\ _{\beta} \otimes_{\alpha}}_{\ \psi}
K$ and $K$ as left $N$-modules by $\Lambda_{\psi}(y) \surl{\
_{\beta} \otimes_{\alpha}}_{\ \psi}\eta\mapsto\alpha(y)\eta$ for
all $y\in {\mathcal N}_{\psi}$. By \cite{EN}, 3.10, we have:
$$\lambda^{\beta,\alpha}_{\xi}=R^{\beta,\psi^o}(\xi) \surl{\
_{\beta} \otimes_{\alpha}}_{\ \psi} 1_K$$

We recall proposition 2.3 of \cite{E2}:

\begin{prop}\label{decomposition}
Let $(\xi_i)_{i\in I}$ be a $(N^o,\psi^o)$-basis of $H_{\beta}$.
Then:
\begin{enumerate}[i)]
\item for all $\xi\in D(H_{\beta},\psi^o)$ and $\eta\in K$, we have:
$$\xi\surl{\
_{\beta}\otimes_{\alpha}}_{\ \psi} \eta=\sum_{i\in I}\xi_i\surl{\
_{\beta}\otimes_{\alpha}}_{\
\psi}\alpha(<\xi,\xi_i>_{\beta,\psi^o})\eta$$
\item we have the following decomposition:
$$ H\surl{\
_{\beta}\otimes_{\alpha}}_{\ \psi}K=\bigoplus_{i\in I}(\xi_i
\surl{\ _{\beta}\otimes_{\alpha}}_{\ \psi}
\alpha(<\xi_i,\xi_i>_{\beta,\psi^o})K)$$
\end{enumerate}
\end{prop}

To end the paragraph, we detail finite dimension case. We assume
that $N$, $H$ and $K$ are of finite dimensions. $H\surl{\
_{\beta}\otimes_{\alpha}}_{\ \psi} K$ can be identified with a
subspace of $H \otimes K$. We denote by $\text{Tr}$ the normalized
canonical trace on $K$ and $\tau=\text{Tr}\circ\alpha$. There
exist a projection $e_{\beta,\alpha}\in\beta(N)\otimes\alpha(N)$
and $n_o\in Z(N)^+$ such that
$(id\otimes\text{Tr})(e_{\beta,\alpha})=\beta(n_o)$. Let $d$ be
the Radon-Nikodym derivative of $\psi$ w.r.t $\tau$. By \cite{EV},
2.4, and proposition 2.7 of \cite{S2}, for all $\xi,\eta\in H$:
$$I_{\beta,\alpha}^{\psi}:\xi\surl{\ _{\beta}\otimes_{\alpha}}_{\ \psi}\eta
\mapsto\xi\surl{\ _{\beta}\otimes_{\alpha}}_{\
\tau}\alpha(d)^{1/2}\eta \mapsto
e_{\beta,\alpha}(\beta(n_o)^{-1/2}\xi\otimes\alpha(d)^{1/2}\eta)$$
defines an isometric isomorphism of
$\beta(N)'-\alpha(N)'^o$-bimodules from $H\surl{\
_{\beta}\otimes_{\alpha}}_{\ \psi}K$ onto a subspace of $H\otimes
K$, the final support of which is $e_{\beta,\alpha}$.\label{iden}

\subsection{Fiber product \cite{V1}, \cite{EV}}
We use previous notations. Let $M_1$ (resp. $M_2$) be a von
Neumann algebra on $H$ (resp. $K$) such that $\beta(N)\subseteq
M_1$ (resp. $\alpha(N)\subseteq M_2$). We denote by $M_1'\surl{\
_{\beta}\otimes_{\alpha}}_{\ N} M_2'$ the von Neumann algebra
generated by $x\surl{\ _{\beta}\otimes_{\alpha}}_{\ N}y$ with
$x\in M_1'$ and $y\in M_2'$.

\begin{defi}
The commutant of $M_1'\surl{\ _{\beta}\otimes_{\alpha}}_{\ N}
M_2'$ in $\mathcal{L}(H\surl{\ _{\beta}\otimes_{\alpha}}_{\
\psi}K)$ is denoted by $M_1\surl{\ _{\beta}\star_{\alpha}}_{\ N}
M_2$ and is called \textbf{fiber product}.
\end{defi}

If $P_1$ and $ P_2$ are von Neumann algebras like $M_1$ and $M_2$,
we have:
$$
\begin{aligned}
i) &\quad (M_1\surl{\ _{\beta}\star_{\alpha}}_{\ N} M_2)\cap
(P_1\surl{\ _{\beta}\star_{\alpha}}_{\ N} P_2)=(M_1\cap
P_1)\surl{\ _{\beta}\star_{\alpha}}_{\ N} (M_2\cap P_2)\\
ii)&\quad \varsigma_N(M_1\surl{\ _{\beta}\star_{\alpha}}_{\ N}
M_2)=M_2\surl{\ _{\alpha}\star_{\beta}}_{\ \ N^o} M_1 \\
iii) &\quad (M_1\cap\beta(N)')\surl{\
_{\beta}\otimes_{\alpha}}_{\ N} (M_2\cap\alpha(N)')\subseteq
M_1\surl{\
_{\beta}\star_{\alpha}}_{\ N} M_2 \\
iv) &\quad M_1\surl{\ _{\beta}\star_{\alpha}}_{\ N}
\alpha(N)=(M_1\cap\beta(N)')\surl{\ _{\beta}\otimes_{\alpha}}_{\
N} 1
\end{aligned}$$
More generally, if $\beta$ (resp. $\alpha$) is a normal,
non-degenerated *-anti-homomorphism (resp. homomorphism) from $N$
to a von Neumann algebra $M_1$ (resp. $M_2$), it is possible to
define a von Neumann algebra $M_1\surl{\
_{\beta}\star_{\alpha}}_{\ N} M_2$ without any reference to a
specific Hilbert space. If $P_1$, $P_2$, $\alpha'$ and $\beta'$
are like $M_1$, $M_2$, $\alpha$ and $\beta$ and if $\Phi$ (resp.
$\Psi$) is a normal *-homomorphism from $M_1$ (resp. $M_2$) to
$P_1$ (resp. $P_2$) such that $\Phi\circ\beta=\beta'$ (resp.
$\Psi\circ\alpha=\alpha'$), then we define a normal
*-homomorphism by \cite{S1}, 1.2.4: $$\Phi\surl{\
_{\beta}\star_{\alpha}}_{\ N}\Psi: M_1\surl{\
_{\beta}\star_{\alpha}}_{\ N}M_2\rightarrow P_1\surl{\
_{\beta'}\star_{\alpha'}}_{\ N}P_2$$

Assume $\ _{\alpha}K_{\epsilon}$ is a $N-P^o$-bimodule and $\
_{\zeta}L$ a left $P$-module. If $\alpha(N)\subseteq M_2$,
$\epsilon(P)\subseteq M_2$ and if $\zeta(P)\subseteq M_3$ where
$M_3$ is a von Neumann algebrasur on $L$, then we can construct
$M_1\surl{\ _{\beta}\star_{\alpha}}_{\ N}(M_2\surl{\
_{\epsilon}\star_{\zeta}}_{\ N}M_3)$ and $(M_1\surl{\
_{\beta}\star_{\alpha}}_{\ N}M_2)\surl{\
_{\epsilon}\star_{\zeta}}_{\ N}M_3$. Associativity of relative
tensor product induces an isomorphism between these fiber products
and we write $M_1\surl{\ _{\beta}\star_{\alpha}}_{\ N}M_2\surl{\
_{\epsilon}\star_{\zeta}}_{\ N}M_3$ without parenthesis.

Finally, if $M_1$ and $M_2$ are of finite dimensions, then we
have:
$$M_1'\surl{\ _{\beta}\otimes_{\alpha}}_{\
N}M_2'=(I_{\beta,\alpha}^{\psi})^*(M_1'\otimes
M_2')I_{\beta,\alpha}^{\psi}\text{ and } M_1\surl{\
_{\beta}\star_{\alpha}}_{\
N}M_2=(I_{\beta,\alpha}^{\psi})^*(M_1\otimes
M_2)I_{\beta,\alpha}^{\psi}$$ Therefore the fiber product can be
identified with a reduction of $M_1\otimes M_2$ by
$e_{\beta,\alpha}$ by \cite{EV}, 2.4.

\subsection{Slice map \cite{E1}}
\subsubsection{For normal forms}
Let $A\in M_1\surl{\ _{\beta}\star_{\alpha}}_{\ N}M_2$ and
$\xi_1,\xi_2\in D(H_{\beta},\psi^o)$. We define an element of
$M_2$ by:
$$(\omega_{\xi_1,\xi_2}\surl{\ _{\beta}\star_{\alpha}}_{\ \psi}
id)(A)=(\lambda^{\beta,\alpha}_{\xi_2})^*A\lambda^{\beta,\alpha}_{\xi_1}$$
so that we have $((\omega_{\xi_1,\xi_2}\surl{\
_{\beta}\star_{\alpha}}_{\ \psi}
id)(A)\eta_1|\eta_2)=(A(\xi_1\surl{\ _{\beta}\otimes_{\alpha}}_{\
\psi} \eta_1)|\xi_2\surl{\ _{\beta}\otimes_{\alpha}}_{\
\psi}\eta_2)$ for all $\eta_1,\eta_2\in K$. Also, we define an
operator of $M_1$ by:
$$(id\surl{\ _{\beta}\star_{\alpha}}_{\
\psi}\omega_{\eta_1,\eta_2})(A)=
(\rho^{\beta,\alpha}_{\eta_2})^*A\rho^{\beta,\alpha}_{\eta_1}$$
for all $\eta_1,\eta_2\in D(_{\alpha}K,\psi)$. We have a Fubini's
formula:
$$\omega_{\eta_1,\eta_2}((\omega_{\xi_1,\xi_2}\surl{\ _{\beta}\star_{\alpha}}_{\ \psi}
id)(A))=\omega_{\xi_1,\xi_2}((id\surl{\
_{\beta}\star_{\alpha}}_{\ \psi}\omega_{\eta_1,\eta_2})(A))$$ for
all $\xi_1,\xi_2\in D(H_{\beta},\psi^o)$ and $\eta_1,\eta_2\in
D(_{\alpha}K,\psi)$.

Equivalently, by (\cite{E1}, proposition 3.3), for all
$\omega_1\in M_{1*}^+$ and $k_1\in\mathbb{R}^+$ such that
$\omega_1\circ\beta\leq k_1\psi$ and for all $\omega_2\in
M_{2*}^+$ and $k_2\in\mathbb{R}^+$ such $\omega_2\circ\alpha\leq
k_2\psi$, we have:
$$\omega_2((\omega_1\surl{\ _{\beta}\star_{\alpha}}_{\ \psi}
id)(A))=\omega_1((id\surl{\ _{\beta}\star_{\alpha}}_{\
\psi}\omega_2)(A))$$

\subsubsection{For conditional expectations}
If $P_2$ is a von Neumann algebra such that $\alpha(N)\subseteq
P_2 \subseteq M_2$ and if $E$ is a normal, faithful conditional
expectation from $M_2$ onto $P_2$, we can define a normal,
faithful conditional expectation $(id\surl{\
_{\beta}\star_{\alpha}}_{\ N}E)$ from $M_1\surl{\
_{\beta}\star_{\alpha}}_{\ N} M_2$ onto $M_1\surl{\
_{\beta}\star_{\alpha}}_{\ N} P_2$ such that:
$$(\omega\surl{\ _{\beta}\star_{\alpha}}_{\
\psi}id) (id\surl{\ _{\beta}\star_{\alpha}}_{\
N}E)(A)=E((\omega\surl{\ _{\beta}\star_{\alpha}}_{\
\psi}id)(A))$$ for all $A\in M_1\surl{\
_{\beta}\star_{\alpha}}_{\ N} M_2$, $\omega\in M_{1*}^+$ and
$k_1\in\mathbb{R}^+$ such that $\omega\circ\beta\leq k_1\psi$.

\subsubsection{For weights}
If $\phi_1$ is n.s.f weight on $M_1$ and if $A$ is a positive
element of $M_1\surl{\ _{\beta}\star_{\alpha}}_{\ N} M_2$, we can
define an element of the extended positive part of $M_2$, denoted
by $(\phi_1\surl{\ _{\beta}\star_{\alpha}}_{\ \psi}id)(A)$, such
that, for all $\eta\in D(_{\alpha}L^2(M_2),\psi)$, we have:
$$||((\phi_1\surl{\ _{\beta}\star_{\alpha}}_{\ \psi}
id)(A))^{1/2}\eta||^2 =\phi_1((id\surl{\
_{\beta}\star_{\alpha}}_{\ \psi}\omega_{\eta})(A))$$ Moreover, if
$\phi_2$ is a n.s.f weight on $M_2$, we have:
$$\phi_2((\phi_1\surl{\ _{\beta}\star_{\alpha}}_{\ \psi} id)(A))=
\phi_1((id\surl{\ _{\beta}\star_{\alpha}}_{\ \psi}\phi_2)(A))$$

Let $(\omega_i)_{i\in I}$ be an increasing net of normal forms
such that $\phi_1=\sup_{i\in I}\omega_i$. Then we have
$(\phi_1\surl{\ _{\beta}\star_{\alpha}}_{\ \psi}
id)(A)=\sup_i(\omega_i\surl{\ _{\beta}\star_{\alpha}}_{\ \psi}
id)(A)$.

\subsubsection{For operator-valued weights}
Let $P_1$ be a von Neumann algebra such that $\beta(N)\subseteq
P_1\subseteq M_1$ and $\Phi_i$ $(i=1,2)$ be operator-valued n.s.f
weights from $M_i$ to $P_i$. By \cite{E1}, for all positive
operator $A\in M_1\surl{\ _{\beta}\star_{\alpha}}_{\ N} M_2$,
there exists an element $(\Phi_1\surl{\ _{\beta}\star_{\alpha}}_{\
N} id)(A)$ belonging to $P_1\surl{\ _{\beta}\star_{\alpha}}_{\ N}
M_2$ such that, for all $\xi\in L^2(P_1)$ and $\eta \in
D(_{\alpha}K,\psi)$, we have:
$$||((\Phi_1\surl{\
_{\beta}\star_{\alpha}}_{\ N} id)(A))^{1/2}(\xi\surl{\
_{\beta}\otimes_{\alpha}}_{\ \psi}\eta)||^2=||[\Phi_1((id\surl{\
_{\beta}\star_{\alpha}}_{\
\psi}\omega_{\eta,\eta})(A))]^{1/2}\xi||^2$$ If $\phi_1$ is a
n.s.f weight on $P_1$, we have:
$$(\phi_1\circ\Phi_1\surl{\
_{\beta}\star_{\alpha}}_{\ N} id)(A)=(\phi_1\surl{\
_{\beta}\star_{\alpha}}_{\ \psi} id)(\Phi_1\surl{\
_{\beta}\star_{\alpha}}_{\ N} id)(A)$$ Also, we define an element
$(id\surl{\ _{\beta}\star_{\alpha}}_{\ N} \Phi_2)(A)$ of the
extended positive part of $M_1\surl{\ _{\beta}\star_{\alpha}}_{\
N} P_2$ and we have:
$$(id\surl{\
_{\beta}\star_{\alpha}}_{\ N} \Phi_2)((\Phi_1\surl{\
_{\beta}\star_{\alpha}}_{\ N} id)(A))=(\Phi_1\surl{\
_{\beta}\star_{\alpha}}_{\ N} id)((id\surl{\
_{\beta}\star_{\alpha}}_{\ N} \Phi_2)(A))$$

\begin{rema}
We have seen that we can identify $M_1\surl{\
_{\beta}\star_{\alpha}}_{\ N} \alpha(N)$ with $M_1\cap\beta(N)'$.
Then, it is easy to check that the slice map $id\surl{\
_{\beta}\star_{\alpha}}_{\ \psi}\psi\circ\alpha^{-1}$ (if
$\alpha$ is injective) is just the injection of $M_1\surl{\
_{\beta}\star_{\alpha}}_{\ N} \alpha(N)$ into $M_1$. Also we see
on that example that, if $\phi_1$ is a n.s.f weight on $M_1$, then
$\phi_1\surl{\ _{\beta}\star_{\alpha}}_{\ N}id$ (which is equal to
$\phi_{1|M_1\cap\beta(N)'}$) needs not to be semi-finite.
\end{rema}

\section{Fundamental pseudo-multiplicative unitary}
In this section, we construct a fundamental pseudo-multiplicative
unitary from a Hopf bimodule with a left invariant
operator-valued weight and a right invariant operator-valued
weight. Let $N$ and $M$ be von Neumann algebras, $\alpha$ (resp.
$\beta$) be a faithful, non-degenerate, normal (resp. anti-)
representation from $N$ to $M$. We suppose that $\alpha(N)
\subseteq \beta(N)'$.

\subsection{Definitions}

\begin{defi}
A quintuplet $(N,M,\alpha,\beta,\Gamma)$ is said to be a {\bf
Hopf bimodule} of basis $N$ if $\Gamma$ is a normal *-homomorphism
from $M$ into $M \surl{\nonumber_{\beta}\star_{\alpha}}_{N} M$
such that, for all $n,m\in N$, we have:

\begin{center}
\begin{minipage}{10cm}
\begin{enumerate}[i)]
\item $\Gamma(\alpha(n)\beta(m))=\alpha(n)\surl{\nonumber_{\beta} \otimes_{\alpha}}_{\ N}
  \beta(m)$
\item $\Gamma$ is co-associative: $(\Gamma \surl{\nonumber _{\beta}
\star_{\alpha}}_{\ N} id)\circ\Gamma=(id \surl{\ _{\beta}
\star_{\alpha}}_{\ N}\Gamma)\circ\Gamma$
\end{enumerate}
\end{minipage}
\end{center}

\end{defi}

One should notice that property i) is necessary in order to write
down the formula given in ii).
$(N^o,M,\beta,\alpha,\varsigma_N\circ\Gamma)$ is a Hopf bimodule
called opposite Hopf bimodule. If $N$ is commutative,
$\alpha=\beta$ and $\Gamma=\varsigma_N\circ\Gamma$, then
$(N,M,\alpha,\alpha,\Gamma)$ is equal to its opposite: we shall
speak about a symmetric Hopf bimodule.

\begin{defi} Let
$(N,M,\alpha,\beta,\Gamma)$ be a Hopf bimodule. A normal,
semi-finite, faithful operator-valued weight from $M$ to
$\alpha(N)$ is said to be {\bf left invariant} if:
$$(id\surl{\ _{\beta}\star_{\alpha}}_{\ N} T_L)\Gamma(x)=T_L(x)\surl{\ _{\beta}
\otimes_{\alpha}}_{\ N} 1\quad\quad\text{for all }x\in
\mathcal{M}_{T_L}^+$$ In the same way, a normal, semi-finite,
faithful operator-valued weight from $M$ to $\beta(N)$ is said to
be {\bf right invariant} if:
$$(T_R \surl{\
_{\beta}
  \star_{\alpha}}_{\ N} id)\Gamma(x)=1 \surl{\ _{\beta}
  \otimes_{\alpha}}_{\ N}T_R(x)\quad\quad\text{for all }x\in \mathcal{M}_{T_R}^+$$
\end{defi}

We give several examples in the last section. In this section,
$(N,M,\alpha,\beta,\Gamma)$ is a Hopf bimodule with a left
operator-valued weight $T_L$ and a right operator-valued weight
$T_R$.

\begin{defi}
A *-anti-automorphism $R$ of $M$ is said to be a {\bf
co-involution} if $R\circ\alpha=\beta$, $R^2=id$ and
$\varsigma_{N^o}\circ(R\surl{\ _{\beta} \star_{\alpha}}_{\ N}
R)\circ\Gamma=\Gamma \circ R$.
\end{defi}

\begin{rema}
With the previous notations, let us notice that $R\circ T_L\circ
R$ is a right invariant operator-valued weight from $M$ to
$\beta(N)$. Also, let us say that $R$ is an anti-isomorphism of
Hopf bimodule from the bimodule and its symmetric.
\end{rema}

Let $\mu$ be a normal, semi-finite, faithful weight of $N$. We
put:
$$\Phi =\mu\circ\alpha^{-1}\circ T_L \text{ and } \Psi
=\mu\circ\beta^{-1}\circ T_R$$ so that, for all $x\in M^+$, we
have:
$$(id \surl{\ _{\beta}
  \star_{\alpha}}_{\ \mu} \Phi)\Gamma(x)=T_L(x) \text{ and }
(\Psi \surl{\ _{\beta} \star_{\alpha}}_{\ \mu}
id)\Gamma(x)=T_R(x)$$

If $H$ denote a Hilbert space on which $M$ acts, then $N$ acts on
$H$, also, by way of $\alpha$ and $\beta$. We shall denote again
$\alpha$ (resp. $\beta$) for (resp. anti-) the representation of
$N$ on $H$.

\subsection{Construction of the fundamental isometry}

\begin{defi}
Let define $\hat{\beta}$ and $\hat{\alpha}$ by:

$$
\begin{aligned}
\hat{\beta}: N &\rightarrow \mathcal{L}(H_{\Phi}) &\quad\text{ and }\quad\quad\hat{\alpha}: N &\rightarrow \mathcal{L}(H_{\Psi}) \\
x &\mapsto J_{\Phi}\alpha(x^*)J_{\Phi}& x &\mapsto J_{\Psi}\beta
(x^*)J_{\Psi}
\end{aligned}$$

Then $\hat{\beta}$ (resp. $\hat{\alpha}$) is a normal,
non-degenerate and faithful anti-representation (resp.
representation) from $N$ to $\mathcal{L}(H_{\Phi})$ (resp.
$\mathcal{L}(H_{\Psi})$). \label{defutil}
\end{defi}

\begin{prop} \label{prem}
We have $\Lambda_{\Phi}({\mathcal N}_{T_L} \cap {\mathcal
N}_{\Phi}) \subseteq D((H_{\Phi})_{\hat{\beta}},\mu^o)$ and for
all $a\in {\mathcal N}_{T_L} \cap {\mathcal N}_{\Phi}$, we have:
$$R^{\hat{\beta},\mu^o}(\Lambda_{\Phi}(a))=\Lambda_{T_L}(a)$$ Also,
we have $\Lambda_{\Psi}({\mathcal N}_{T_R} \cap {\mathcal
N}_{\Psi}) \subseteq D(_{\hat{\alpha}}(H_{\Psi}),\mu)$ and for all
$b\in {\mathcal N}_{T_R} \cap {\mathcal N}_{\Psi}$, then:
$$R^{\hat{\alpha},\mu}(\Lambda_{\Psi}(b))=\Lambda_{T_R}(b)$$
\end{prop}

\begin{rema}
We identify $H_{\mu}$ with $H_{\mu\circ\alpha^{-1}}$ and $H_{\mu}$
with $H_{\mu\circ\beta^{-1}}$.
\end{rema}

\begin{proof}
Let $y\in {\mathcal N}_{\mu}$ analytic w.r.t $\mu$. We have:
$$
\begin{aligned}
\hat{\beta}(y^*)\Lambda_{\Phi}(a)&=\Lambda_{\Phi}(a\sigma_{-i/2}^{\Phi}(\alpha(y^*)))
=\Lambda_{\Phi}(a\sigma_{-i/2}^{\mu\circ\alpha^{-1}}(\alpha(y^*)))\\
&=\Lambda_{\Phi}(a\alpha(\sigma_{-i/2}^{\mu}(y^*)))
=\Lambda_{T_L}(a)\Lambda_{\mu}(\sigma_{-i/2}^{\mu}(y^*))
=\Lambda_{T_L}(a)J_{\mu}\Lambda_{\mu}(y)
\end{aligned}$$
Thanks to lemma \ref{analyse}, we get
$\hat{\beta}(y^*)\Lambda_{\Phi}(a)
=\Lambda_{T_L}(a)J_{\mu}\Lambda_{\mu}(y)$, for all $y\in {\mathcal
N}_{\mu}$, which gives the first part of the proposition. The end
of the proof is very similar.
\end{proof}

\begin{prop}\label{evi}
We have $J_{\Phi}D((H_{\Phi})_{\hat{\beta}},
\mu^o)=D(_{\alpha}(H_{\Phi}),\mu)$ and for all $\eta\in
D((H_{\Phi})_{\hat{\beta}}, \mu^o)$, we have:
$$R^{\alpha,\mu}(J_{\Phi}\eta)=J_{\Phi}R^{\hat{\beta},\mu^o}(\eta)J_{\mu}$$
Also, we have
$J_{\Psi}D(_{\hat{\alpha}}(H_{\Psi}),\mu)=D((H_{\Phi})_{\beta},\mu^o)$
and for all $\xi\in D((H_{\Phi})_{\beta},\mu^o)$, we have:
$$R^{\beta,\mu^o}(J_{\Psi}\xi)=J_{\Psi}R^{\hat{\alpha},\mu}(\xi)J_{\mu}$$
\end{prop}

\begin{proof}
Straightforward.
\end{proof}

\begin{coro}
We have $\Lambda_{\Phi}({\mathcal T}_{\Phi,T_L}) \subseteq
D((H_{\Phi})_{\hat{\beta}}, \mu^o) \cap
D(_{\alpha}(H_{\Phi}),\mu)$ and $\Lambda_{\Psi}({\mathcal
T}_{\Psi,T_R}) \subseteq D(_{\hat{\alpha}}(H_{\Psi}),\mu) \cap
D((H_{\Psi})_{\beta},\mu^o)$.
\end{coro}

\begin{proof}
This is a corollary of the two previous propositions.
\end{proof}

\begin{rema}
The invariance of operator-valued weights does not play a part in
the previous propositions.
\end{rema}

\begin{prop}\label{semi}
We have $(\omega_{v,\xi}\surl{\ _{\beta} \star_{\alpha}}_{\ \mu}
id)(\Gamma(a)) \in {\mathcal N}_{T_L}\cap {\mathcal N}_{\Phi}$
for all elements $a\in {\mathcal N}_{T_L}\cap {\mathcal
N}_{\Phi}$ and $v,\xi\in D(H_{\beta},\mu^o)$.
\end{prop}

\begin{proof} By definition of the slice maps, we have:
$$\begin{aligned}
(\omega_{v,\xi}\surl{\ _{\beta} \star_{\alpha}}_{\ \mu}
  id)(\Gamma(a))^*(\omega_{v,\xi}\surl{\ _{\beta}
  \star_{\alpha}}_{\ \mu} id)(\Gamma(a))
&=(\lambda^{\beta,
  \alpha}_v)^*\Gamma(a^*)\lambda^{\beta,
  \alpha}_{\xi}(\lambda^{\beta,
  \alpha}_{\xi})^*\Gamma(a)\lambda^{\beta, \alpha}_v\\
&\leq \|\lambda^{\beta, \alpha}_{\xi}\|^2(\omega_{v,v}\surl{\
_{\beta}
  \star_{\alpha}}_{\ \mu} id)(\Gamma(a^*a))\\
&\leq \|R^{\beta, \mu^o}(\xi)\|^2(\omega_{v,v}\surl{\ _{\beta}
  \star_{\alpha}}_{\ \mu} id)(\Gamma(a^*a))
  \end{aligned}$$
Then, on one hand, we get, thanks to left invariance of $T_L$:
$$\begin{aligned}
&\ \quad {T_L}((\omega_{v,\xi}\surl{\ _{\beta} \star_{\alpha}}_{\
\mu}
  id)(\Gamma(a))^*(\omega_{v,\xi}\surl{\ _{\beta}
  \star_{\alpha}}_{\ \mu} id)(\Gamma(a)))\\
&\leq \|R^{\beta, \mu^o}(\xi)\|^2{T_L}((\omega_{v,v}\surl{\
_{\beta}\star_{\alpha}}_{\ \mu} id)(\Gamma(a^*a)))\\
&=\|R^{\beta, \mu^o}(\xi)\|^2 (\omega_{v,v}\surl{\ _{\beta}
  \star_{\alpha}}_{\ \mu} id)(id \surl{\ _{\beta}
  \star_{\alpha}}_{\ \mu} {T_L})(\Gamma(a^*a))\\
&\leq \|R^{\beta,\mu^o}(\xi)\|^2(\lambda^{\beta,\alpha}_{v})^*
({T_L}(a^*a)\surl{\ _{\beta}\otimes_{\alpha}}_{\ \mu} 1)
\lambda^{\beta,\alpha}_{v}\\
&\leq\|R^{\beta,\mu^o}(\xi)\|^2||{T_L}(a^*a)||
\|\alpha(<v,v>_{\beta,\mu^o})\|1\\
&\leq \|R^{\beta,\mu^o}(\xi)\|^2\|{T_L}(a^*a)\|\|R^{\beta,
\mu^o}(v)\|^21 \end{aligned}$$ So, we get that
$(\omega_{v,\xi}\surl{\ _{\beta}
  \star_{\alpha}}_{\ \mu} id)(\Gamma(a)) \in {\mathcal N}_{T_L}$.
On the other hand, thanks to left invariance of $T_L$, we know
that:
$$\Phi(((\omega_{v,\xi}\surl{\ _{\beta}
  \star_{\alpha}}_{\ \mu} id)(\Gamma(a)))^*(\omega_{v,\xi}\surl{\
  _{\beta} \star_{\alpha}}_{\ \mu} id)(\Gamma(a)))$$ is less or
equal to:
$$
\begin{aligned}
&\ \quad\|R^{\beta,\mu^o}(\xi)\|^2 \Phi((\omega_{v,v}\surl{\
_{\beta}\star_{\alpha}}_{\ \mu} id)(\Gamma(a^*a)))\\
&=\|R^{\beta, \mu^o}(\xi)\|^2\omega_{v,v}((id \surl{\ _{\beta}
  \star_{\alpha}}_{\ \mu}\Phi)(\Gamma(a^*a)))\\
&=\|R^{\beta,
\mu^o}(\xi)\|^2({T_L}(a^*a)v|v)\leq\|R^{\beta,\mu^o}(\xi)\|^2\|{T_L}(a^*a)\|\|v\|^2
< + \infty
\end{aligned}$$  So, we get that
$(\omega_{v,\xi}\surl{\ _{\beta} \star_{\alpha}}_{\ \mu}
id)(\Gamma(a)) \in {\mathcal N}_{\Phi}$.
\end{proof}

\begin{prop}\label{rapide}
For all $v,w\in H$ and $a,b\in {\mathcal N}_{\Phi}\cap {\mathcal
N}_{T_L}$, we have:
$$(v\surl{\ _{\alpha}
\otimes_{\hat{\beta}}}_{\ \ \mu^o}\Lambda_{\Phi}(a)|w\surl{\
_{\alpha} \otimes_{\hat{\beta}}}_{\ \
\mu^o}\Lambda_{\Phi}(b))=(T_L(b^*a)v|w)$$

For all $v,w\in H$ and $c,d\in {\mathcal N}_{\Psi}\cap {\mathcal
N}_{T_R}$, we have:
$$(\Lambda_{\Psi}(c)\surl{\ _{\hat{\alpha}} \otimes_{\beta}}_{\ \
\mu^o}v|\Lambda_{\Psi}(d)\surl{\ _{\hat{\alpha}}
\otimes_{\beta}}_{\ \ \mu^o}w)=(T_R(d^*c)v|w)$$
\end{prop}

\begin{proof}
Using \ref{prem} and \ref{preintro}, we get that:
$$\begin{aligned}
(v \surl{\ _{\alpha} \otimes_{\hat{\beta}}}_{\ \ \mu^o}
\Lambda_{\Phi}(a)|w\surl{\ _{\alpha} \otimes_{\hat{\beta}}}_{\ \
\mu^o}\Lambda_{\Phi}(b))
&=(\alpha(<\Lambda_{\Phi}(a),\Lambda_{\Phi}(b)>_{\hat{\beta},\mu^o})v|w)\\
&=(\alpha(\Lambda_{T_L}(b)^*\Lambda_{T_L}(a))v|w)\\
&=(\alpha(\pi_{\mu}(\alpha^{-1}({T_L}(b^*a))))v|w)
\end{aligned}$$
which gives the result after the identification of $\pi_{\mu}(N)$
with $N$. The second point is very similar.
\end{proof}

\begin{lemm}\label{inde}
Let $a\in {\mathcal N}_{\Phi}\cap {\mathcal N}_{T_L}$ and $v\in
D(H_{\beta},\mu^o)$. The following sum: $$\sum_{i\in I} \xi_{i}
\surl{\ _{\beta} \otimes_{\alpha}}_{\ \mu}
  \Lambda_{\Phi} ((\omega_{v,\xi_i} \surl{\
      _{\beta} \star_{\alpha}}_{\ \mu} id)(\Gamma(a)))$$ converges
in $H \surl{\ _{\beta} \otimes_{\alpha}}_{\ \mu} H_{\Phi}$ for
all $(N^o,\mu^o)$-basis $(\xi_i)_{i\in I}$ of $H_{\beta}$ and it
does not depend on the $(N^o,\mu^o)$-basis of $H_{\beta}$.
\end{lemm}

\begin{proof}
By \ref{semi}, we have $(\omega_{v,\xi_i} \surl{\ _{\beta}
\star_{\alpha}}_{\ \mu} id)(\Gamma(a))\in {\mathcal N}_{\Phi}\cap
{\mathcal N}_{T_L}$ for all $i\in I$, and the vectors $\xi_{i}
\surl{\ _{\beta} \otimes_{\alpha}}_{\ \mu} \Lambda_{\Phi}
((\omega_{v,\xi_i} \surl{\ _{\beta} \star_{\alpha}}_{\ \mu}
id)(\Gamma(a)))$ are two-by-two orthogonal. Normality and left
invariance of $\Phi$ imply:
$$
\begin{aligned}
&\quad\sum_{i
    \in I} ||\xi_{i} \surl{\ _{\beta} \otimes_{\alpha}}_{\ \mu}
  \Lambda_{\Phi} ((\omega_{v,\xi_i} \surl{\
      _{\beta} \star_{\alpha}}_{\ \mu} id)(\Gamma(a)))||^2\\
&=\sum_{i\in I} (\alpha(<\xi_i,\xi_i>_{\beta,\mu^o})
\Lambda_{\Phi} ((\omega_{v,\xi_i} \surl{\
      _{\beta} \star_{\alpha}}_{\ \mu} id)(\Gamma(a)))|
\Lambda_{\Phi} ((\omega_{v,\xi_i} \surl{\
      _{\beta} \star_{\alpha}}_{\ \mu} id)(\Gamma(a))))\\
&=\Phi((\lambda_v^{\beta,\alpha})^*\Gamma(a^*)[\sum_{i\in I}
\lambda_{\xi_i}^{\beta,\alpha}(\lambda_{\xi_i}^{\beta,\alpha})^*
\lambda_{\xi_i}^{\beta,\alpha}(\lambda_{\xi_i}^{\beta,\alpha})^*]
\Gamma(a)\lambda_v^{\beta,\alpha})\\
&=\Phi((\omega_{v,v}\surl{\ _{\beta} \star_{\alpha}}_{\ \mu}
id)(\Gamma(a^*a)))=((id \surl{\ _{\beta} \star_{\alpha}}_{\ \mu}
\Phi)(\Gamma(a^*a))v|v)=(T_L(a^*a)v|v)<\infty
\end{aligned}$$ We deduce that the sum $\sum_{i\in I} \xi_{i} \surl{\ _{\beta}
\otimes_{\alpha}}_{\ \mu}\Lambda_{\Phi} ((\omega_{v,\xi_i}
\surl{\ _{\beta} \star_{\alpha}}_{\ \mu} id)(\Gamma(a)))$
converges in $H \surl{\ _{\beta} \otimes_{\alpha}}_{\ \mu}
H_{\Phi}$. To prove that the sum does not depend on the $(N^o,
\mu^{o})$-basis, we compute for all $b \in {\mathcal N}_{T_L} \cap
{\mathcal N}_{\Phi}$ and $w \in D(H_{\beta},\mu^o)$:
$$
\begin{aligned}
&\quad(\sum_{i\in I} \xi_{i} \surl{\ _{\beta}
\otimes_{\alpha}}_{\ \mu} \Lambda_{\Phi} ((\omega_{v,\xi_i}
\surl{\ _{\beta} \star_{\alpha}}_{\ \mu} id)(\Gamma(a)))) | w
\surl{\ _{\beta} \otimes_{\alpha}}_{\ \mu} \Lambda_{\Phi}(b))\\
&=\sum_{i \in I}
(\alpha(<\xi_i,w>_{\beta,\mu^o})\Lambda_{\Phi}((\omega_{v,\xi_i}
\surl{\ _{\beta}\star_{\alpha}}_{\ \mu}id)
(\Gamma(a)))|\Lambda_{\Phi}(b))\\
&=\sum_{i\in I}
\Phi(b^*\alpha(<\xi_i,w>_{\beta,\mu^o})(\omega_{v,\xi_i} \surl{\
_{\beta} \star_{\alpha}}_{\ \mu} id)(\Gamma(a)))\\
&=\Phi(b^*\lambda_w^{\beta,\alpha}[\sum_{i \in I}
\lambda^{\beta,\alpha}_{\xi_i}(\lambda^{\beta,\alpha}_{\xi_i})^*]
\Gamma(a)\lambda_v^{\beta,\alpha})=\Phi(b^*(\omega_{v,w}\surl{\
_{\beta} \star_{\alpha}}_{\ \mu} id)(\Gamma(a))).
\end{aligned}$$ As $D(H_{\beta}, \mu^{o})
\odot \Lambda_{\Phi}({\mathcal N}_{T_L} \cap {\mathcal
N}_{\Phi})$ is dense in $H \surl{\ _{\beta} \otimes_{\alpha}}_{\
\mu} H_{\Phi}$ and the last expression is independent of the
$(N^o, \mu^{o})$-basis, we can conclude.
\end{proof}

\begin{theo}\label{isom}
Let $H$ be a Hilbert space on which $M$ acts. There exists a
unique isometry $U_H$, called \textbf{fundamental isometry}, from
$ H \surl{\ _{\alpha} \otimes_{\hat{\beta}}}_{\ \ \mu^o} H_{\Phi}$
to $H \surl{\ _{\beta} \otimes_{\alpha}}_{\ \mu} H_{\Phi}$ such
that, for all $(N^o, \mu^{o})$-basis $(\xi_i)_{i \in I}$ of
$H_{\beta}$, $a \in {\mathcal N}_{T_L} \cap {\mathcal N}_{\Phi}$
and $v \in D(H_{\beta}, \mu^{o})$:
$$U_H(v \surl{\ _{\alpha} \otimes_{\hat{\beta}}}_{\ \ \mu^o}
\Lambda_{\Phi}(a))=\sum_{i\in I} \xi_{i} \surl{\ _{\beta}
\otimes_{\alpha}}_{\ \mu}\Lambda_{\Phi} ((\omega_{v,\xi_i}
\surl{\ _{\beta} \star_{\alpha}}_{\ \mu} id)(\Gamma(a))))$$
\end{theo}

\begin{proof}
By \ref{inde}, we can define the following application:
$$\begin{aligned}
\tilde{U}:\  D(H_{\beta}, \mu^{o}) \times
\Lambda_{\Phi}({\mathcal N}_{T} \cap {\mathcal N}_{\Phi})
&\rightarrow H \surl{\ _{\beta} \otimes_{\alpha}}_{\ \mu}
H_{\Phi} \\
(v,\Lambda_{\Phi}(a)) & \mapsto \sum_{i \in I} \xi_{i} \surl{\
_{\beta} \otimes_{\alpha}}_{\ \mu} \Lambda_{\Phi}
((\omega_{v,\xi_i} \surl{\ _{\beta} \star_{\alpha}}_{\ \mu}
id)(\Gamma(a))))
\end{aligned}$$
Let $b \in {\mathcal N}_{T_L} \cap {\mathcal N}_{\Phi}$ and $w \in
D(H_{\beta},\mu^o)$. Then, by normality and left invariance of
$\Phi$, we have:
$$
\begin{aligned}
&\ \quad(\tilde{U}(v,\Lambda_{\Phi}(a))|\tilde{U}(w,\Lambda_{\Phi}(b)))\\
&=\sum_{i,j\in I}
(\alpha(<\xi_i,\xi_j>_{\beta,\mu^o})\Lambda_{\Phi}((\omega_{v,\xi_i}
\surl{\ _{\beta} \star_{\alpha}}_{\ \mu} id)(\Gamma(a)))
|\Lambda_{\Phi}((\omega_{w,\xi_i} \surl{\ _{\beta}
\star_{\alpha}}_{\ \mu} id)(\Gamma(b))))\\
&=\sum_{i\in I}
(\Lambda_{\Phi}(\alpha(<\xi_i,\xi_i>_{\beta,\mu^o})
(\omega_{v,\xi_i} \surl{\ _{\beta} \star_{\alpha}}_{\ \mu}
id)(\Gamma(a))) | \Lambda_{\Phi}((\omega_{w,\xi_i} \surl{\
_{\beta} \star_{\alpha}}_{\ \mu} id)(\Gamma(b))))\\
&=\sum_{i \in I}\Phi((\lambda_w^{\beta,\alpha})^*\Gamma(b^*)
\lambda_{\xi_i}^{\beta,\alpha}\alpha(<\xi_i,\xi_i>_{\beta,\mu^o})
(\lambda_{\xi_i}^{\beta,\alpha})^*\Gamma(a)\lambda_v^{\beta,\alpha})\\
&=\Phi ((\lambda_w^{\beta,\alpha})^*\Gamma(b^*)[\sum_{i \in I}
\lambda_{\xi_i}^{\beta,\alpha}(\lambda_{\xi_i}^{\beta,\alpha})^*
\lambda_{\xi_i}^{\beta,\alpha}(\lambda_{\xi_i}^{\beta,\alpha})^*]
\Gamma(a)\lambda_v^{\beta,\alpha}) \end{aligned}$$ Then,
properties of $(N^o, \mu^{o})$-basis $(\xi_i)_{i \in I}$ of
$H_{\beta}$ imply that:
$$\begin{aligned}
\Phi((\omega_{v,w}\surl{\ _{\beta}\star_{\alpha}}_{\ \mu}
id)(\Gamma(b^*a)))&=\omega_{v,w}((id \surl{\ _{\beta}
\star_{\alpha}}_{\ \mu}\Phi)(\Gamma(b^*a)))\\
&=\omega_{v,w}(T_L(b^*a))=(T_L(b^*a)v|w)
\end{aligned}$$
By \ref{rapide}, we get:
$$(\tilde{U}((v,\Lambda_{\Phi}(a))|\tilde{U}((w,\Lambda_{\Phi}(b))))
=(v\surl{\ _{\alpha}\otimes_{\hat{\beta}}}_{\ \ \mu^o}
\Lambda_{\Phi}(a)|w\surl{\ _{\alpha} \otimes_{\hat{\beta}}}_{\ \
\mu^o} \Lambda_{\Phi}(b))$$ so that, from $\tilde{U}$, we can
easily define a suitable application $U_H$ which is independent
of the $(N^o,\mu^{o})$-basis by \ref{inde}.
\end{proof}

One can define a right version of $U_H$ from the right invariant
weight:

\begin{theo}
Let $H$ be a Hilbert space on which $M$ acts. There exists a
unique isometry $U'_H$ from $H_{\Psi} \surl{\ _{\hat{\alpha}}
\otimes_{\beta}}_{\ \ \mu^o} H$ to $H_{\Psi} \surl{\ _{\beta}
\otimes_{\alpha}}_{\ \mu} H$ such that, for all $(N,\mu)$-basis
$(\eta_i)_{i \in I}$ of $\ _{\alpha}H$, $a \in {\mathcal N}_{T_R}
\cap {\mathcal N}_{\Psi}$ and $v \in D(_{\alpha}H, \mu)$:
$$U'_H(\Lambda_{\Psi}(a) \surl{\ _{\hat{\alpha}} \otimes_{\beta}}_{\ \ \mu^o} v
)=\sum_{i\in I} \Lambda_{\Psi} ((id\surl{\
_{\beta}\star_{\alpha}}_{\ \mu} \omega_{v,\eta_i})(\Gamma(a)))
\surl{\ _{\beta} \otimes_{\alpha}}_{\ \mu} \eta_{i}$$
\end{theo}

\subsection{Relations between the fundamental isometry and the co-product}

\begin{prop}\label{raccourci}
We have $(1\surl{\ _{\beta} \otimes_{\alpha}}_{\
N}J_{\Phi}eJ_{\Phi})U_H
\rho^{\alpha,\hat{\beta}}_{\Lambda_{\Phi}(x)}=\Gamma(x)
\rho^{\beta,\alpha}_{J_{\Phi}\Lambda_{\Phi}(e)}$ for all $e,x\in
{\mathcal N}_{\Phi}\cap {\mathcal N}_{T_L}$ and
$(J_{\Psi}fJ_{\Psi}\surl{\ _{\beta} \otimes_{\alpha}}_{\ N}1)U'_H
\lambda^{\hat{\alpha},\beta}_{\Lambda_{\Psi}(y)}=\Gamma(y)
\lambda^{\beta,\alpha}_{J_{\Psi}\Lambda_{\Psi}(f)}$\! for all
$f,y\in {\mathcal N}_{\Psi}\cap {\mathcal N}_{T_R}$.
\end{prop}

\begin{proof}
Let $v\in D(H_{\beta},\mu^o)$ and $(\xi_i)_{i \in I}$ a
$(N^o,\mu^o)$-basis of $H_{\beta}$. We have:

$$
\begin{aligned}
&\ \quad (1\surl{\ _{\beta} \otimes_{\alpha}}_{\
N}J_{\Phi}eJ_{\Phi})U_H(v\surl{\ _{\alpha}
\otimes_{\hat{\beta}}}_{\ \ \mu^o} \Lambda_{\Phi}(x))\\
&=\sum_{i\in I}\xi_i\surl{\ _{\beta} \otimes_{\alpha}}_{\ \mu}
J_{\Phi}eJ_{\Phi}\Lambda_{\Phi}((\omega_{v,\xi_i}\surl{\ _{\beta}
\star_{\alpha}}_{\ \mu} id) (\Gamma(x)))\\
&=\sum_{i\in I}\xi_i\surl{\ _{\beta} \otimes_{\alpha}}_{\ \mu}
(\omega_{v,\xi_i}\surl{\ _{\beta} \star_{\alpha}}_{\ \mu} id)
(\Gamma(x))J_{\Phi}\Lambda_{\Phi}(e)=\Gamma(x)(v \surl{\ _{\beta}
\otimes_{\alpha}}_{\ \mu} J_{\Phi}\Lambda_{\Phi}(e))
\end{aligned}$$

By \ref{prem} and \ref{evi}, we have $\Lambda_{\Phi}(x)\in
D((H_{\Phi})_{\hat{\beta}},\mu^o)$ and
$J_{\Phi}\Lambda_{\Phi}(e)\in D( _{\alpha}(H_{\Phi}),\mu)$ so
that each term of the previous equality is continuous in $v$.
Density of $D(H_{\beta},\mu^o)$ in $H$ finishes the proof. The
last part is very similar.
\end{proof}

\begin{prop}\label{rap}
For all $v,w\in D(H_{\beta},\mu^o)$ and $a\in {\mathcal
N}_{\Phi}\cap {\mathcal N}_{T_L}$, we have:
$$(\lambda_w^{\beta,\alpha})^*U_H(v\surl{\
_{\alpha}\otimes_{\hat{\beta}}}_{\ \ \mu^o}\Lambda_{\Phi}(a))=
\Lambda_{\Phi}((\omega_{v,w}\surl{\ _{\beta}\star_{\alpha}}_{\
\mu} id)(\Gamma(a)))$$ Also, for all $v',w'\in D(_{\alpha}H,\mu)$
and $b\in {\mathcal N}_{\Psi}\cap {\mathcal N}_{T_R}$, we have:
$$(\rho_{w'}^{\beta,\alpha})^*U'_H(\Lambda_{\Psi}(b)\surl{\
_{\alpha}\otimes_{\hat{\beta}}}_{\ \ \mu^o}v')=
\Lambda_{\Psi}((id\surl{\ _{\beta}\star_{\alpha}}_{\ \mu}
\omega_{v',w'})(\Gamma(b)))$$
\end{prop}

\begin{proof}
Let $e\in {\mathcal N}_{\Phi}\cap {\mathcal N}_{T_L}$. By
\ref{raccourci}, we can compute:

$$
\begin{aligned}
J_{\Phi}eJ_{\Phi}(\lambda_w^{\beta,\alpha})^*U_H(v\surl{\
_{\alpha}\otimes_{\hat{\beta}}}_{\ \ \mu^o}\Lambda_{\Phi}(a)) &=
(\lambda_w^{\beta,\alpha})^*(1\surl{\ _{\beta}\otimes_{\alpha}}_{\
N}J_{\Phi}eJ_{\Phi})U_H\rho_{\Lambda_{\Phi}(a)}^{\alpha,\hat{\beta}}v\\
&=(\lambda_w^{\beta,\alpha})^*\Gamma(a)\rho^{\beta,\alpha}_{J_{\Phi}\Lambda_{\Phi}(e)}v\\
&=(\omega_{v,w}\surl{\ _{\beta}\star_{\alpha}}_{\ \mu}
id)(\Gamma(a))J_{\Phi}\Lambda_{\Phi}(e)\\
&=J_{\Phi}eJ_{\Phi}\Lambda_{\Phi}((\omega_{v,w}\surl{\
_{\beta}\star_{\alpha}}_{\ \mu} id)(\Gamma(a)))
\end{aligned}$$
Density of ${\mathcal N}_{\Phi}\cap {\mathcal N}_{T_L}$ in $N$
finishes the proof. The second part is very similar.
\end{proof}

\begin{coro}\label{lienGV}
For all $a \in {\mathcal N}_{T_L}\cap {\mathcal N}_{\Phi}$, $v\in
D( _{\alpha}H,\mu) \cap D(H_{\beta},\mu^o)$ and $w\in
D(H_{\beta},\mu^o)$, we have: $$(\omega_{v,w} *
id)(U_H)\Lambda_{\Phi}(a)= \Lambda_{\Phi}((\omega_{v,w} \surl{\
_{\beta} \star_{\alpha}}_{\ \mu} id)(\Gamma(a)))$$ where we denote
by $(\omega_{v,w}*id)(U_H)$ the operator
$(\lambda_w^{\beta,\alpha})^*U_H\lambda_v^{\alpha,\hat{\beta}}$ of
$\mathcal{L}(H_{\Phi})$.
\end{coro}

\begin{proof}
Straightforward.
\end{proof}

\begin{coro}\label{corres}
For all $e,x\in {\mathcal N}_{\Phi}\cap {\mathcal N}_{T_L}$ and
$\eta\in D(_{\alpha}H_{\Phi},\mu^o)$, we have: $$(id\surl{\
_{\beta} \star_{\alpha}}_{\ \mu}
\omega_{J_{\Phi}\Lambda_{\Phi}(e),\eta} )(\Gamma(x))=(id*
\omega_{\Lambda_{\Phi}(x),J_{\Phi}e^*J_{\Phi}\eta})(U_H)$$ Also,
for all $f,y\in {\mathcal N}_{\Psi}\cap {\mathcal N}_{T_R}$ and
$\xi\in D((H_{\Psi})_{\beta},\mu^o)$, we have:
$$(\omega_{J_{\Psi}\Lambda_{\Psi}(f),\xi} \surl{\ _{\beta}
\star_{\alpha}}_{\
\mu}id)(\Gamma(y))=(\omega_{\Lambda_{\Psi}(y),J_{\Psi}f^*J_{\Psi}\xi}*
id)(U'_H)$$
\end{coro}

\begin{proof}
Straightforward by \ref{raccourci}.
\end{proof}

\begin{coro}\label{switch}
For all $a,b \in {\mathcal N}_{\Psi}\cap {\mathcal N}_{T_R} \cap
{\mathcal N}_{\Psi}^*\cap {\mathcal N}_{T_R}^*$, we have:
$$(\omega_{\Lambda_{\Psi}(a),J_{\Psi}\Lambda_{\Psi}(b)}*id)(U'_H)^*
=(\omega_{\Lambda_{\Psi}(a^*),J_{\Psi}\Lambda_{\Psi}(b^*)}*id)(U'_H)$$
\end{coro}

\begin{proof}
By \ref{corres}, we have for all $e\in {\mathcal N}_{\Psi}\cap
{\mathcal N}_{T_R}$:
$$
\begin{aligned}
(\omega_{\Lambda_{\Psi}(a),J_{\Psi}\Lambda_{\Psi}(e^*b)}*
id)(U'_H)^* &=
(\omega_{J_{\Psi}\Lambda_{\Psi}(e),J_{\Psi}\Lambda_{\Psi}(b)}
\surl{\ _{\beta} \star_{\alpha}}_{\ \mu}
id)(\Gamma(a))^* \\
&= (\omega_{J_{\Psi}\Lambda_{\Psi}(b),J_{\Psi}\Lambda_{\Psi}(e)}
\surl{\ _{\beta} \star_{\alpha}}_{\ \mu}
id)(\Gamma(a^*)) \\
&= (\omega_{\Lambda_{\Psi}(a^*),J_{\Psi}\Lambda_{\Psi}(b^*e)}*
id)(U'_H).
\end{aligned}$$

Let $(u_k)_{k\in K}$ be a family in ${\mathcal N}_{\Psi} \cap
{\mathcal N}_{\Psi}^*$ such that $u_k\rightarrow 1$ in the
*-strong topology. We denote:
$$e_k=\frac{1}{\sqrt{\pi}} \int\! e^{-t^2}\sigma_{t}^{\Psi}(u_k)\ dt$$
For all $k\in K$, $e_k$ and $\sigma_{-i/2}^{\Psi}(e^*_k)$ are
bounded and belong to ${\mathcal N}_{\Psi}$ and converge to $1$
in the *-strong topology so that $J_{\Psi}\Lambda_{\Psi}(b^*e_k)
=\sigma_{-i/2}^{\Psi}(e^*_k)J_{\Psi}\Lambda_{\Psi}(b^*)$ converge
to $J_{\Psi}\Lambda_{\Psi}(b^*)$ in norm of $H_{\Psi}$. Let
$\xi,\eta \in D( _{\alpha}H,\mu)$ and we compute:
$$
\begin{aligned}
((\omega_{\Lambda_{\Psi}(a),J_{\Psi}\Lambda_{\Psi}(b)}*
id)(U'_H)^*\xi|\eta)&=(J_{\Psi}\Lambda_{\Psi}(b) \surl{\ _{\beta}
\otimes_{\alpha}}_{\ \mu}\xi|U'_H(\Lambda_{\Psi}(a)\surl{\
_{\hat{\alpha}} \otimes_{\beta}}_{\ \ \mu^o} \eta))\\
&=\lim_{k\in K}(J_{\Psi}\Lambda_{\Psi}(e_k^*b) \surl{\ _{\beta}
\otimes_{\alpha}}_{\ \mu}\xi|U'_H(\Lambda_{\Psi}(a)\surl{\
_{\hat{\alpha}} \otimes_{\beta}}_{\ \ \mu^o} \eta))\\
&=\lim_{k\in
K}((\omega_{\Lambda_{\Psi}(a),J_{\Psi}\Lambda_{\Psi}(e_k^*b)}*
id)(U'_H)^*\xi|\eta)\end{aligned}$$ By the previous computation,
this last expression is equal to:
$$
\begin{aligned}
&\ \quad\lim_{k\in
K}((\omega_{\Lambda_{\Psi}(a^*),J_{\Psi}\Lambda_{\Psi}(b^*e_k)}*
id)(U'_H)\xi|\eta)\\
&=\lim_{k\in K}(U'_H(\Lambda_{\Psi}(a)\surl{\ _{\hat{\alpha}}
\otimes_{\beta}}_{\ \ \mu^o} \xi)|J_{\Psi}\Lambda_{\Psi}(b^*e_k)
\surl{\ _{\beta} \otimes_{\alpha}}_{\ \mu}\eta)\\
&=(U'_H(\Lambda_{\Psi}(a^*)\surl{\ _{\hat{\alpha}}
\otimes_{\beta}}_{\ \ \mu^o} \xi)|J_{\Psi}\Lambda_{\Psi}(b^*)
\surl{\ _{\beta}\otimes_{\alpha}}_{\
\mu}\eta)=((\omega_{\Lambda_{\Psi}(a^*),J_{\Psi}\Lambda_{\Psi}(b^*)}*
id)(U'_H)\xi|\eta)\\
\end{aligned}$$
By density of $D( _{\alpha}H,\mu)$ in $H$, the result holds.
\end{proof}

\subsection{Commutation relations}\label{rcom}
In this section, we verify commutation relations which are
necessary for $U_H$ to be a pseudo-multiplicative unitary and we
establish a link between $U_H$ and $\Gamma$. We also have similar
formulas for $U'_H$.

\begin{lemm}\label{base}
Let $\xi\in D(H_{\beta},\mu^o)$ and $\eta\in D( _{\alpha}H,\mu)$.

\begin{center}
\begin{minipage}{11cm}
\begin{enumerate}[i)]
\item For all $a \in \alpha(N)'$, we have
$\lambda_\xi^{\beta,\alpha}\circ a=(1 \surl{\ _{\beta}
\otimes_{\alpha}}_{\ N} a)\lambda_{\xi}^{\beta,\alpha}$.
\item For all $b \in \beta(N)'$, we have
$\lambda_{b\xi}^{\beta,\alpha}=(b \surl{\ _{\beta}
  \otimes_{\alpha}}_{\ N} 1)\lambda_{\xi}^{\beta,\alpha}$.
\item For all $x \in {\mathcal D}(\sigma_{-i/2}^{\mu})$, we have
$\lambda_{\beta(x)\xi}^{\beta,\alpha}=\lambda_{\xi}^{\beta,\alpha}\circ
\alpha(\sigma_{-i/2}^{\mu}(x))$.
\item For all $x \in {\mathcal D}(\sigma_{i/2}^{\mu})$, we have
$\rho_{\alpha(x)\eta}^{\beta,\alpha}=\rho_{\eta}^{\beta,\alpha}\circ
\beta(\sigma_{i/2}^{\mu}(x))$.
\end{enumerate}
\end{minipage}
\end{center}
\end{lemm}

\begin{proof}
Straightforward.
\end{proof}

We recall that $\alpha(N)$ and $\beta(N)$ commute with
$\hat{\beta}(N)'$.

\begin{prop}\label{comm}
For all $n\in N$, we have:

\begin{center}
\begin{minipage}{8cm}
\begin{enumerate}[i)]
\item $U_H(1\surl{\ _{\alpha}\otimes_{\hat{\beta}}}_{\ \ N^o}
  \alpha(n))=(\alpha(n)\surl{\ _{\beta}\otimes_{\alpha}}_{\ N}
  1)U_H$;
\item $U_H(1\surl{\ _{\alpha}\otimes_{\hat{\beta}}}_{\ \ N^o}
  \beta(n))=(1\surl{\ _{\beta}\otimes_{\alpha}}_{\ N}
  \beta(n))U_H$;
\item $U_H(\beta(n)\surl{\ _{\alpha}\otimes_{\hat{\beta}}}_{\ \ N^o}
  1)=(1\surl{\ _{\beta}\otimes_{\alpha}}_{\ N}
  \hat{\beta}(n))U_H$.
\end{enumerate}
\end{minipage}
\end{center}

\end{prop}

\begin{proof}
By \ref{raccourci}, we can compute for all $n\in N $ and $e,x \in
{\mathcal N}_{T_L} \cap {\mathcal N}_{\Phi}$:
$$\begin{aligned}
(\alpha(n) \surl{\ _{\beta} \otimes_{\alpha}}_{\
N}J_{\Phi}eJ_{\Phi})U_H\rho^{\alpha,\hat{\beta}}_{\Lambda_{\Phi}(x)}
&=(\alpha(n) \surl{\ _{\beta} \otimes_{\alpha}}_{\
N}1)\Gamma(x)\rho_{J_{\Phi}\Lambda_{\Phi}(e)}^{\beta,\alpha}\\
&=\Gamma(\alpha(n)x)\rho_{J_{\Phi}\Lambda_{\Phi}(e)}^{\beta,\alpha}\\
&=(1\surl{\ _{\beta} \otimes_{\alpha}}_{\
N}J_{\Phi}eJ_{\Phi})U_H\rho^{\alpha,\hat{\beta}}_{\Lambda_{\Phi}(\alpha(n)x)}\\
&=(1\surl{\ _{\beta} \otimes_{\alpha}}_{\
N}J_{\Phi}eJ_{\Phi})U_H(1 \surl{\ _{\alpha}
\otimes_{\hat{\beta}}}_{\ \ N^o}
\alpha(n))\rho^{\alpha,\hat{\beta}}_{\Lambda_{\Phi}(x)}
\end{aligned}$$
Usual arguments of density imply the first equality. The second
one can be proved in a very similar way. By \ref{raccourci} and
\ref{base}, we can compute for all $n\in {\mathcal T}_{\mu}$ and
$e,x \in {\mathcal N}_{T_L} \cap {\mathcal N}_{\Phi}$:
$$
\begin{aligned}
(1\surl{\ _{\beta} \otimes_{\alpha}}_{\
N}J_{\Phi}eJ_{\Phi}\hat{\beta}(n))U_H\rho^{\alpha,\hat{\beta}}_{\Lambda_{\Phi}(x)}
&=\Gamma(x)\rho_{J_{\Phi}\Lambda_{\Phi}(e\alpha(n^*))}^{\beta,\alpha}\\
&=\Gamma(x)\rho_{\alpha(\sigma_{-i/2}^{\mu}(n))J_{\Phi}\Lambda_{\Phi}(e)}^{\beta,\alpha}\\
&=\Gamma(x)\rho_{J_{\Phi}\Lambda_{\Phi}(e)}^{\beta,\alpha}\beta(n)\\
&=(1\surl{\ _{\beta} \otimes_{\alpha}}_{\
N}J_{\Phi}eJ_{\Phi})U_H\rho^{\alpha,\hat{\beta}}_{\Lambda_{\Phi}(x)}\beta(n)\\
&=(1\surl{\ _{\beta} \otimes_{\alpha}}_{\
N}J_{\Phi}eJ_{\Phi})U_H(\beta(n)\surl{\ _{\alpha}
\otimes_{\hat{\beta}}}_{\ \ N^o}
1)\rho^{\alpha,\hat{\beta}}_{\Lambda_{\Phi}(x)}
\end{aligned}$$
Density of ${\mathcal T}_{\mu}$ in $N$ and normality of $\beta$
and $\hat{\beta}$ finish the proof.
\end{proof}

\begin{prop}\label{appartenance}
For all $x\in M'\cap\mathcal{L}(H)$, we have: $$U_H(x \surl{\
_{\alpha} \otimes_{\hat{\beta}}}_{\ \ N^o} 1)= (x \surl{\ _{\beta}
\otimes_{\alpha}}_{\ N} 1)U_H$$
\end{prop}

\begin{proof}
For all $e,y \in {\mathcal N}_{T_L} \cap {\mathcal N}_{\Phi}$ and
$x\in M'\cap\mathcal{L}(H)\subseteq
\alpha(N)'\cap\beta(N)'\cap\mathcal{L}(H)$, we have by
\ref{raccourci}:
$$\begin{aligned}
(x\surl{\ _{\beta} \otimes_{\alpha}}_{\
N}J_{\Phi}eJ_{\Phi})U_H\rho^{\alpha,\hat{\beta}}_{\Lambda_{\Phi}(y)}
&=(x\surl{\ _{\beta} \otimes_{\alpha}}_{\
N}1)\Gamma(y)\rho_{J_{\Phi}\Lambda_{\Phi}(e)}^{\beta,\alpha}\\
&=\Gamma(y)\rho_{J_{\Phi}\Lambda_{\Phi}(e)}^{\beta,\alpha}x\\
&=(1\surl{\ _{\beta} \otimes_{\alpha}}_{\
N}J_{\Phi}eJ_{\Phi})U_H\rho^{\alpha,\hat{\beta}}_{\Lambda_{\Phi}(y)}x\\
&=(1\surl{\ _{\beta} \otimes_{\alpha}}_{\
N}J_{\Phi}eJ_{\Phi})U_H(x \surl{\ _{\alpha}
\otimes_{\hat{\beta}}}_{\ \ N^o}
1)\rho^{\alpha,\hat{\beta}}_{\Lambda_{\Phi}(y)}
\end{aligned}$$ Usual arguments of density imply the result.
\end{proof}

\begin{coro}
For all $n \in N$, we have:

\begin{center}
\begin{minipage}{8cm}
\begin{enumerate}[i)]
\item $U_{H_{\Phi}}(\hat{\beta}(n) \surl{\ _{\alpha}
\otimes_{\hat{\beta}}}_{\ \ N^o} 1)=(\hat{\beta}(n) \surl{\
_{\beta} \otimes_{\alpha}}_{\ N} 1)U_{H_{\Phi}}$
\item $U_{H_{\Psi}}(\hat{\alpha}(n) \surl{\ _{\alpha}
\otimes_{\hat{\beta}}}_{\ \ N^o} 1)=(\hat{\alpha}(n) \surl{\
_{\beta} \otimes_{\alpha}}_{\ N} 1)U_{H_{\Psi}}$
\end{enumerate}
\end{minipage}
\end{center}

\end{coro}

\begin{prop}\label{impl}
We have $\Gamma(m)U_H=U_H(1 \surl{\ _{\alpha}
\otimes_{\hat{\beta}}}_{\ \ N^o} m)$ for all $m \in M$.
\end{prop}

\begin{proof}
By \ref{raccourci}, we can compute for all $e,x\in {\mathcal
N}_{T_L} \cap {\mathcal N}_{\Phi}$:
$$\begin{aligned}
(1\surl{\ _{\beta} \otimes_{\alpha}}_{\
N}J_{\Phi}eJ_{\Phi})\Gamma(m)U_H\rho^{\alpha,\hat{\beta}}_{\Lambda_{\Phi}(x)}
&=\Gamma(m)(1\surl{\ _{\beta} \otimes_{\alpha}}_{\
N}J_{\Phi}eJ_{\Phi})U_H\rho^{\alpha,\hat{\beta}}_{\Lambda_{\Phi}(x)}\\
&=\Gamma(mx)\rho_{J_{\Phi}\Lambda_{\Phi}(e)}^{\beta,\alpha}\\
&=(1\surl{\ _{\beta} \otimes_{\alpha}}_{\
N}J_{\Phi}eJ_{\Phi})U_H\rho^{\alpha,\hat{\beta}}_{\Lambda_{\Phi}(mx)}\\
&=(1\surl{\ _{\beta} \otimes_{\alpha}}_{\
N}J_{\Phi}eJ_{\Phi})U_H(1\surl{\ _{\alpha}
\otimes_{\hat{\beta}}}_{\ \ N^o}
m)\rho^{\alpha,\hat{\beta}}_{\Lambda_{\Phi}(x)}
\end{aligned}$$ Usual arguments of density imply the result.
\end{proof}

\subsection{Unitarity of the fundamental isometry}

To prove unitary of $U_H$ (resp. $U'_H$), we establish a
reciprocity law where both left and right operator-valued weights
are at stake.

\subsubsection{First technical result}

We establish results needed for \ref{reciprocite}. In the
following proposition, we compute some functions $\theta$ defined
in section \ref{intre}.

\begin{prop}\label{inov}
We have for all $c\in {\mathcal N}_{\Psi}\cap {\mathcal N}_{T_R}$,
$m\in ({\mathcal N}_{\Psi}\cap {\mathcal N}_{T_R})^*$ and $v\in
D(H_{\beta},\mu^o)$:
$$\theta^{\beta,\mu^o}(v,J_{\Psi}\Lambda_{\Psi}(c))m=(\lambda_{\Lambda_{\Psi}
(m^*)}^{\hat{\alpha},\beta})^*\rho_v^{\hat{\alpha},\beta}J_{\Psi}c^*J_{\Psi}$$
\end{prop}

\begin{proof}
Let $x\in {\mathcal N}_{\Psi}\cap {\mathcal N}_{T_R}$. On one
hand, we get by \ref{prem} and \ref{evi}:
$$
\begin{aligned}
\theta^{\beta,\mu^o}(v,J_{\Psi}\Lambda_{\Psi}(c))m\Lambda_{\Psi}(x)
&=
R^{\beta,\mu^o}(v)R^{\beta,\mu^o}(J_{\Psi}\Lambda_{\Psi}(c))^*\Lambda_{\Psi}(mx)\\
&=
R^{\beta,\mu^o}(v)J_{\mu}\Lambda_{T_R}(c)^*J_{\Psi}\Lambda_{\Psi}(mx).
\end{aligned}$$

On the other hand, if $c\in {\mathcal T}_{\Psi,T_R}$, then we
have by \ref{rapide}:
$$
\begin{aligned}
(\lambda_{\Lambda_{\Psi}(m^*)}^{\hat{\alpha},\beta})^*\rho_v^{\hat{\alpha},\beta}
J_{\Psi}c^*J_{\Psi}\Lambda_{\Psi}(x)
&=(\lambda_{\Lambda_{\Psi}(m^*)}^{\hat{\alpha},\beta})^*
(J_{\Psi}c^*J_{\Psi}\Lambda_{\Psi}(x) \surl{\
_{\hat{\alpha}}\otimes_{\beta}}_{\ \mu^o} v) \\
&=T_R(mx\sigma_{-i/2}^{\Psi}(c))v\\
&=R^{\beta,\mu^o}(v)J_{\mu}
\Lambda_{\mu}(\beta^{-1}(T_R(\sigma_{i/2}^{\Psi}(c^*)x^*m^*)))\\
&=R^{\beta,\mu^o}(v)J_{\mu}\Lambda_{T_R}(c)^*J_{\Psi}\Lambda_{\Psi}(mx)
\end{aligned}$$
We obtain:
$$(\lambda_{\Lambda_{\Psi}(m^*)}^{\hat{\alpha},\beta})^*\rho_v^{\hat{\alpha},\beta}
J_{\Psi}c^*J_{\Psi}\Lambda_{\Psi}(x)=R^{\beta,\mu^o}(v)J_{\mu}\Lambda_{T_R}(c)^*J_{\Psi}\Lambda_{\Psi}(mx)$$
for all $c\in {\mathcal N}_{\Psi}\cap {\mathcal N}_{T_R}$ by
normality which finishes the proof.
\end{proof}

\begin{coro}\label{prepa}
Let $a\in ({\mathcal N}_{\Psi}\cap {\mathcal N}_{T_R})^*
({\mathcal N}_{\Phi}\cap {\mathcal N}_{T_L})$. If $c\in {\mathcal
N}_{\Psi}\cap {\mathcal N}_{T_R}$, $e\in {\mathcal N}_{\Phi}\cap
{\mathcal N}_{T_L}$ and $\xi\in H_{\Psi},\eta \in D(_{\alpha}
(H_{\Phi}),\mu)$, $u\in H$, $v\in D(H_{\beta},\mu^o)$, then we
have:
$$(v\surl{\ _{\beta}\otimes_{\alpha}}_{\ \mu}
(\lambda^{\beta,\alpha}_{J_{\Psi}\Lambda_{\Psi}(c)})^*
U_{H_{\Psi}}(\xi\surl{\ _{\alpha}\otimes_{\hat{\beta}}}_{\ \
\mu^o}\Lambda_{\Phi}(a))|u\surl{\ _{\beta}\otimes_{\alpha}}_{\
\mu} J_{\Phi}e^*J_{\Phi}\eta)$$
$$=(J_{\Psi}c^*J_{\Psi}\xi\surl{\ _{\hat{\alpha}}\otimes_{\beta}}_{\ \ \mu^o} v|
\Lambda_{\Psi}((id\surl{\ _{\beta}\star_{\alpha}}_{\ \mu}
\omega_{\eta, J_{\Phi}\Lambda_{\Phi}(e)})(\Gamma(a^*)))\surl{\
_{\hat{\alpha}}\otimes_{\beta}}_{\ \ \mu^o}u)$$
\end{coro}

\begin{proof}
By \ref{raccourci} and \ref{inov}, we can compute:
$$\begin{aligned}
&\ \quad(v\surl{\ _{\beta}
  \otimes_{\alpha}}_{\ \mu}
  (\lambda^{\beta,\alpha}_{J_{\Psi}\Lambda_{\Psi}(c)})^*
  U_{H_{\Psi}}(\xi\surl{\ _{\alpha}\otimes_{\hat{\beta}}}_{\ \ \mu^o}\Lambda_{\Phi}(a))
  |u\surl{\ _{\beta}
  \otimes_{\alpha}}_{\ \mu} J_{\Phi}e^*J_{\Phi}\eta)\\
&=((\rho^{\beta,\alpha}_{\eta})^*
\lambda_v^{\beta,\alpha}(\lambda^{\beta,\alpha}_{J_{\Psi}\Lambda_{\Psi}(c)})^*
(1\surl{\ _{\beta}
  \otimes_{\alpha}}_{\ N}J_{\Phi}e^*J_{\Phi})
  U_{H_{\Psi}}(\xi\surl{\ _{\alpha}\otimes_{\hat{\beta}}}_{\ \
  \mu^o}\Lambda_{\Phi}(a))|u)\\
&=((\rho^{\beta,\alpha}_{\eta})^*
\lambda_v^{\beta,\alpha}(\lambda^{\beta,\alpha}_{J_{\Psi}\Lambda_{\Psi}(c)})^*
\Gamma(a)\rho^{\beta,\alpha}_{J_{\Phi}\Lambda_{\Phi}(e)}\xi|u)\\
&=\theta^{\beta,\mu^o}(v,J_{\Psi}\Lambda_{\Psi}(c))(\rho^{\beta,\alpha}_{\eta})^*
\Gamma(a)\rho^{\beta,\alpha}_{J_{\Phi}\Lambda_{\Phi}(e)}\xi|u)\\
&=((\lambda^{\hat{\alpha},\beta}_{\Lambda_{\Psi}((id\surl{\
_{\beta} \star_{\alpha}}_{\ \mu} \omega_{\eta,
J_{\Phi}\Lambda_{\Phi}(e))})(\Gamma(a^*))})^*\rho^{\hat{\alpha},\beta}_v
J_{\Psi}c^*J_{\Psi}\xi|u)\\
&=(J_{\Psi}c^*J_{\Psi}\xi\surl{\ _{\hat{\alpha}}
  \otimes_{\beta}}_{\ \ \mu^o} v|\Lambda_{\Psi}((id\surl{\
_{\beta} \star_{\alpha}}_{\ \mu} \omega_{\eta,
J_{\Phi}\Lambda_{\Phi}(e))})(\Gamma(a^*)))\surl{\ _{\hat{\alpha}}
  \otimes_{\beta}}_{\ \ \mu^o} u)
\end{aligned}$$

\end{proof}

\subsubsection{Second technical result}

In this section, results only depend on \ref{raccourci} and
co-product relation but not on the previous technical result. Let
${\mathcal H}$ be an other Hilbert space on which $M$ acts.

\begin{lemm}\label{simple}
Let $a,e\in {\mathcal N}_{\Phi}\cap {\mathcal N}_{T_L}$, $\xi\in
D({\mathcal H}_{\beta},\mu^o)$, $\eta\in D(_{\alpha}H,\mu)$, and
$\zeta\in {\mathcal H}$. We have:
$$(1\surl{\ _{\beta} \otimes_{\alpha}}_{\ N}J_{\Phi}eJ_{\Phi})
U_H(\eta\surl{\ _{\alpha} \otimes_{\hat{\beta}}}_{\ \ \mu^o}
[(\lambda_{\xi}^ {\beta,\alpha})^*U_{\mathcal H}(\zeta\surl{\
_{\alpha} \otimes_{\hat{\beta}}}_{\ \ \mu^o}\Lambda_{\Phi}
(a))])$$
$$=(\lambda_{\xi}^ {\beta,\alpha}\surl{\ _{\beta} \otimes_{\alpha}}_{\
N}1)^*(id \surl{\ _{\beta}
  \star_{\alpha}}_{\ N} \Gamma)(\Gamma(a))(\zeta
  \surl{\ _{\beta} \otimes_{\alpha}}_{\ \mu}
\eta \surl{\ _{\beta} \otimes_{\alpha}}_{\ \mu}
J_{\Phi}\Lambda_{\Phi}(e))$$
\end{lemm}

\begin{proof}
First let assume $\zeta\in D({\mathcal H}_{\beta},\mu^o)$. By
\ref{rap} and \ref{raccourci}, we can compute:
$$
\begin{aligned} &\quad(1\surl{\ _{\beta} \otimes_{\alpha}}_{\
N}J_{\Phi}eJ_{\Phi}) U_H(\eta\surl{\ _{\alpha}
\otimes_{\hat{\beta}}}_{\ \ \mu^o} [(\lambda_{\xi}^
{\beta,\alpha})^*U_{\mathcal H}(\zeta\surl{\ _{\alpha}
\otimes_{\hat{\beta}}}_{\ \ \mu^o}\Lambda_{\Phi} (a))])\\
&=(1\surl{\ _{\beta} \otimes_{\alpha}}_{\ N}J_{\Phi}eJ_{\Phi})
  U_H(\eta\surl{\ _{\alpha} \otimes_{\hat{\beta}}}_{\ \ \mu^o} \Lambda_{\Phi}
((\omega_{\zeta,\xi}\surl{\ _{\beta}
  \star_{\alpha}}_{\ \mu} id)(\Gamma(a)))\\
&=\Gamma((\omega_{\zeta,\xi}\surl{\ _{\beta}
  \star_{\alpha}}_{\ \mu} id)(\Gamma(a)))
  (\eta \surl{\ _{\beta} \otimes_{\alpha}}_{\ \mu}
  J_{\Phi}\Lambda_{\Phi}(e))\\
&=(\lambda_{\xi}^ {\beta,\alpha}\surl{\ _{\beta}
\otimes_{\alpha}}_{\ N}1)^*(id \surl{\ _{\beta}
  \star_{\alpha}}_{\ N} \Gamma)(\Gamma(a))(\zeta
  \surl{\ _{\beta} \otimes_{\alpha}}_{\ \mu}
\eta \surl{\ _{\beta} \otimes_{\alpha}}_{\ \mu}
J_{\Phi}\Lambda_{\Phi}(e))
\end{aligned}$$
So, we get the result for all $\zeta\in D({\mathcal
H}_{\beta},\mu^o)$. The first term of the equality is continuous
in $\zeta$ because $\eta\in D(_{\alpha}H,\mu)$ and
$\Lambda_{\Phi}(a)\in D((H_{\Phi})_{\hat{\beta}},\mu^o)$. Also,
since $\eta\in D(_{\alpha}H,\mu)$ and $\Lambda_{\Phi}(a)\in
D((H_{\Phi})_{\hat{\beta}},\mu^o)$, the last term of the equality
is continuous in $\zeta$. Density of $D({\mathcal
H}_{\beta},\mu^o)$ in ${\mathcal H}$ finishes the proof.
\end{proof}

\begin{lemm}\label{simple2}
The sum $\sum_{i\in I} \eta_i\surl{\ _{\alpha}
\otimes_{\hat{\beta}}}_{\ \ \mu^o}
[(\lambda_{\xi}^{\beta,\alpha})^*U_{\mathcal
H}((\rho_{\eta_i}^{\beta,\alpha})^*\Xi\surl{\ _{\alpha}
\otimes_{\hat{\beta}}}_{\ \ \mu^o} \Lambda_{\Phi}(a))]$ converges
for all $\xi\in D({\mathcal H}_{\beta},\mu^o)$, $\Xi\in {\mathcal
H}\surl{\ _{\beta}\otimes_{\alpha}}_{\ \mu} H$, $a\in {\mathcal
N}_{\Phi}\cap {\mathcal N}_{T_L}$ and $(N,\mu)$-basis
$(\eta_i)_{i\in I}$ of $\ _{\alpha}H$.
\end{lemm}

\begin{proof}
First, observe that $\eta_i\surl{\
_{\alpha}\otimes_{\hat{\beta}}}_{\ \ \mu^o}
[(\lambda_{\xi}^{\beta,\alpha})^*U_{\mathcal
H}((\rho_{\eta_i}^{\beta,\alpha})^*\Xi\surl{\ _{\hat{\alpha}}
\otimes_{\beta}}_{\ \ \mu^o} \Lambda_{\Phi}(a))]$ are orthogonal.
To compute, we put:
$\Omega_i=\rho_{\eta_i}^{\beta,\alpha})^*\Xi\surl{\ _{\alpha}
  \otimes_{\hat{\beta}}}_{\ \ \mu^o} \Lambda_{\Phi}(a)$. By \ref{base} and \ref{comm}, we have:
$$
\begin{aligned}
&\ \quad ||\eta_i\surl{\ _{\alpha}
  \otimes_{\hat{\beta}}}_{\ \ \mu^o} [(\lambda_{\xi}^{\beta,\alpha})^*
  U_{\mathcal H}(\Omega_i)]||^2\\
&=(\hat{\beta}(<\eta_i,\eta_i>_{\alpha,\mu})
(\lambda_{\xi}^{\beta,\alpha})^*
  U_{\mathcal H}(\Omega_i)|(\lambda_{\xi}^{\beta,\alpha})^*
  U_{\mathcal H}(\Omega_i))\\
&=((\lambda_{\xi}^{\beta,\alpha})^*(1\surl{\ _{\beta}
  \otimes_{\alpha}}_{\ \mu}\hat{\beta}(<\eta_i,\eta_i>_{\alpha,\mu}))
  U_{\mathcal H}(\Omega_i)|(\lambda_{\xi}^{\beta,\alpha})^*
  U_{\mathcal H}(\Omega_i))\\
&=((\lambda_{\xi}^{\beta,\alpha})^*
  U_{\mathcal H}(\beta(<\eta_i,\eta_i>_{\alpha,\mu})
  (\Omega_i)|(\lambda_{\xi}^{\beta,\alpha})^*U_{\mathcal H}(\Omega_i))\\
&=(\lambda_{\xi}^{\beta,\alpha}(\lambda_{\xi}^{\beta,\alpha})^*
  U_{\mathcal H}(\Omega_i)|U_{\mathcal H}(\Omega_i))
\end{aligned}$$
By \ref{rapide}, it follows that we have, for all $i\in I$:
$$
\begin{aligned}
&\ \quad ||\eta_i\surl{\ _{\alpha}
  \otimes_{\hat{\beta}}}_{\ \ \mu^o} [(\lambda_{\xi}^{\beta,\alpha})^*
  U_{\mathcal H}((\rho_{\eta_i}^{\beta,\alpha})^*\Xi\surl{\ _{\alpha}
  \otimes_{\hat{\beta}}}_{\ \ \mu^o} \Lambda_{\Phi}(a))]||^2\\
&\leq ||R^{\beta,\alpha}(\xi)||^2
  ((\rho_{\eta_i}^{\beta,\alpha})^*\Xi\surl{\ _{\alpha}
  \otimes_{\hat{\beta}}}_{\ \ \mu^o} \Lambda_{\Phi}(a)|
  (\rho_{\eta_i}^{\beta,\alpha})^*\Xi\surl{\ _{\alpha}
  \otimes_{\hat{\beta}}}_{\ \ \mu^o} \Lambda_{\Phi}(a))\\
&\leq ||R^{\beta,\alpha}(\xi)||^2(T_L(a^*a)
  (\rho_{\eta_i}^{\beta,\alpha})^*\Xi
  |(\rho_{\eta_i}^{\beta,\alpha})^*\Xi)\\
&\leq ||R^{\beta,\alpha}(\xi)||^2||T(a^*a)||(
  (\rho_{\eta_i}^{\beta,\alpha})^*\Xi
  |(\rho_{\eta_i}^{\beta,\alpha})^*\Xi)
\end{aligned}$$
So, we can sum over $i\in I$ to get that:
$$\sum_{i\in I}||\eta_i\surl{\ _{\alpha}\otimes_{\hat{\beta}}}_{\ \ \mu^o}
[(\lambda_{\xi}^{\beta,\alpha})^*U_{\mathcal
H}((\rho_{\eta_i}^{\beta,\alpha})^*\Xi\surl{\ _{\alpha}
\otimes_{\hat{\beta}}}_{\ \ \mu^o} \Lambda_{\Phi}(a))]||^2$$ is
less or equal to
$||R^{\beta,\alpha}(\xi)||^2||T(a^*a)||||\Xi||^2<\infty$. That's
why the sum converges.
\end{proof}

\begin{prop}\label{techinter}
Let $a,e\in {\mathcal N}_{\Phi}\cap {\mathcal N}_{T_L}$, $\Xi\in
{\mathcal H}\surl{\ _{\beta}\otimes_{\alpha}}_{\ \mu} H$, $\xi\in
D({\mathcal H}_{\beta},\mu^o)$, $\eta\in
D(_{\alpha}(H_{\Phi}),\mu)$ and $(\eta_i)_{i\in I}$ a
$(N,\mu)$-basis of $\ _{\alpha}H$. We have:
$$(\rho_{J_{\Phi}eJ_{\Phi}\eta}^{\beta,\alpha})^*U_H(\sum_{i\in I} \eta_i\surl{\ _{\alpha}
  \otimes_{\hat{\beta}}}_{\ \ \mu^o} [(\lambda_{\xi}^{\beta,\alpha})^*
  U_{\mathcal H}((\rho_{\eta_i}^{\beta,\alpha})^*\Xi\surl{\ _{\alpha}
  \otimes_{\hat{\beta}}}_{\ \ \mu^o} \Lambda_{\Phi}(a))])$$
$$=(\lambda_{\xi}^{\beta,\alpha})^*\Gamma((id\surl{\ _{\beta}
  \star_{\alpha}}_{\ \mu}\omega_{J_{\Phi}
  \Lambda_{\Phi}(e),\eta})(\Gamma(a)))\Xi$$
\end{prop}

\begin{proof}
The existence of the first term comes from the previous lemma. By
\ref{simple} and the co-product relation, we can compute:
$$
\begin{aligned}
&\quad\sum_{i\in I}(\rho_{\eta}^{\beta,\alpha})^*(1\surl{\
_{\beta}
  \otimes_{\alpha}}_{\ N}
J_{\Phi}eJ_{\Phi})U_H( \eta_i\surl{\ _{\alpha}
  \otimes_{\hat{\beta}}}_{\ \ \mu^o} [(\lambda_{\xi}^{\beta,\alpha})^*
  U_{\mathcal H}((\rho_{\eta_i}^{\beta,\alpha})^*\Xi\surl{\ _{\alpha}
  \otimes_{\hat{\beta}}}_{\ \ \mu^o} \Lambda_{\Phi}(a))])\\
&=\sum_{i\in I}(\rho_{\eta}^{\beta,\alpha})^*
  (\lambda_{\xi}^ {\beta,\alpha}\surl{\ _{\beta} \otimes_{\alpha}}_{\
  N}1)^*(id \surl{\ _{\beta}
  \star_{\alpha}}_{\ N} \Gamma)(\Gamma(a))((\rho_{\eta_i}^{\beta,\alpha})^*
  \Xi \surl{\ _{\beta} \otimes_{\alpha}}_{\ \mu}
  \eta_i \surl{\ _{\beta} \otimes_{\alpha}}_{\ \mu}
  J_{\Phi}\Lambda_{\Phi}(e))\\
&=(\rho_{\eta}^{\beta,\alpha})^*
  (\lambda_{\xi}^ {\beta,\alpha}\surl{\ _{\beta} \otimes_{\alpha}}_{\
  N}1)^*(\Gamma\surl{\ _{\beta}
  \star_{\alpha}}_{\ N} id)(\Gamma(a))([\sum_{i\in I}
  \rho_{\eta_i}^{\beta,\alpha}
  (\rho_{\eta_i}^{\beta,\alpha})^*]\Xi\surl{\ _{\beta}
  \otimes_{\alpha}}_{\ \mu} J_{\Phi}\Lambda_{\Phi}(e)))\\
&=(\rho_{\eta}^{\beta,\alpha})^*
  (\lambda_{\xi}^ {\beta,\alpha}\surl{\ _{\beta} \otimes_{\alpha}}_{\
  N}1)^*(\Gamma\surl{\ _{\beta}
  \star_{\alpha}}_{\ N} id)(\Gamma(a))(\Xi\surl{\ _{\beta}
  \otimes_{\alpha}}_{\ \mu} J_{\Phi}\Lambda_{\Phi}(e)))\\
&=(\lambda_{\xi}^ {\beta,\alpha})^*
  (1\surl{\ _{\beta} \otimes_{\alpha}}_{\
  N}\rho_{\eta}^{\beta,\alpha})^*(\Gamma\surl{\ _{\beta}
  \star_{\alpha}}_{\ N} id)(\Gamma(a))(\Xi\surl{\ _{\beta}
  \otimes_{\alpha}}_{\ \mu} J_{\Phi}\Lambda_{\Phi}(e)))\\
&=(\lambda_{\xi}^{\beta,\alpha})^*\Gamma((id\surl{\ _{\beta}
  \star_{\alpha}}_{\ \mu}
  \omega_{J_{\Phi}\Lambda_{\Phi}(e),
  \eta})(\Gamma(a)))\Xi
\end{aligned}$$

\end{proof}

With results of the two last sections in hand, we can prove now a
reciprocity law where ${\mathcal H}$ will be equal to $H_{\Psi}$.

\subsubsection{Reciprocity law}\label{reciprocite}

For all monotone increasing net $(e_k)_{k\in K}$ in ${\mathcal
N}_{\Psi}\cap {\mathcal N}_{T_R}$ of limit equal to $1$, the
following  $(\omega_{J_{\Psi}\Lambda_{\Psi}(e_k)})_{k\in K}$ is
monotone increasing and converges to $\Psi$. So, for all $x\in
{\mathcal N}_{\Psi}\cap {\mathcal N}_{T_R}$,
$(\omega_{J_{\Psi}\Lambda_{\Psi}(e_k)}\surl{\ _{\beta}
\star_{\alpha}}_{\ \mu} id)(\Gamma(x))$ converges to $(\Psi\surl{\
_{\beta}\star_{\alpha}}_{\ \mu} id)(\Gamma(x))$ in the weak
topology. We denote $\zeta_k=J_{\Psi}\Lambda_{\Psi}(e^*_ke_k)\in
D((H_{\Psi})_{\beta},\mu^o)$ for all $k\in K$.

\begin{prop}
For all $a\in ({\mathcal N}_{\Psi}\cap {\mathcal N}_{T_R})^*
({\mathcal N}_{\Phi}\cap {\mathcal N}_{T_L}))$, $e\in {\mathcal
N}_{\Phi}\cap {\mathcal N}_{T_L},b\in {\mathcal N}_{\Psi}\cap
{\mathcal N}_{T_R}, c\in {\mathcal T}_{\Psi,T_R}$, $v \in
D(H_{\beta},\mu^o), \eta \in D(_{\alpha} (H_{\Phi}),\mu)$ and
$(N,\mu)$-basis of $\ _{\alpha}H$, $(\eta_i)_{i\in I}$ , we have
that the image of:
$$\sum_{i\in I} \eta_i\!\!\surl{\ _{\alpha}\otimes_{\hat{\beta}}}_{\ \ \mu^o}
[(\lambda_{\zeta_k}^{\beta,\alpha})^*U_{H_{\Psi}}([(\rho_{\eta_i}^{\beta,\alpha})^*
U'_H(J_{\Psi}c^*J_{\Psi}\Lambda_{\Psi}(b)\!\!\surl{\
_{\hat{\alpha}}\otimes_{\beta}}_{\ \ \mu^o}v)]\!\!\surl{\
_{\alpha}\otimes_{\hat{\beta}}}_{\ \ \mu^o} \Lambda_{\Phi}(a))]$$
by $(\rho^{\beta,\alpha}_{J_{\Phi}e^*J_{\Phi}\eta})^*U_H$
converges, in the weak topology, to:
$$(\rho^{\beta,\alpha}_{J_{\Phi}e^*J_{\Phi}\eta})^*
(v\surl{\ _{\beta}\otimes_{\alpha}}_{\
\mu}(\lambda^{\beta,\alpha}_{J_{\Psi}\Lambda_{\Psi}(c)})^*
U_{H_{\Psi}}(\Lambda_{\Psi}(b)\surl{\
_{\alpha}\otimes_{\hat{\beta}}}_{\ \ \mu^o}\Lambda_{\Phi}(a)))$$
\end{prop}

\begin{proof}
Let $u\in H$. We compute the value of the scalar product of:
$$U_H(\sum_{i\in I} \eta_i\surl{\ _{\alpha}
  \otimes_{\hat{\beta}}}_{\ \ \mu^o} [(\lambda_{\zeta_k}^{\beta,\alpha})^*
  U_{H_{\Psi}}([(\rho_{\eta_i}^{\beta,\alpha})^*U'_H(\Lambda_{\Psi}(bc)
  \surl{\ _{\hat{\alpha}}
  \otimes_{\beta}}_{\ \ \mu^o}v)]\surl{\ _{\alpha}
  \otimes_{\hat{\beta}}}_{\ \ \mu^o} \Lambda_{\Phi}(a))])$$ by $u\surl{\ _{\beta}
  \otimes_{\alpha}}_{\ \mu} J_{\Phi}e^*J_{\Phi}\eta$. By \ref{techinter},
we get that it is equal to:
$$(\Gamma((id\surl{\ _{\beta}
  \star_{\alpha}}_{\ \mu}
  \omega_{J_{\Phi}\Lambda_{\Phi}(e),
  \eta})(\Gamma(a)))U'_H(\Lambda_{\Psi}(bc)
  \surl{\ _{\hat{\alpha}}
  \otimes_{\beta}}_{\ \ \mu^o}v)|\zeta_k\surl{\ _{\beta}
  \otimes_{\alpha}}_{\ \mu} u)$$
By the right version of \ref{impl}, this is equal to:
$$(U'_H(\Lambda_{\Psi}((id\surl{\ _{\beta}
  \star_{\alpha}}_{\ \mu}
  \omega_{J_{\Phi}\Lambda_{\Phi}(e),
  \eta})(\Gamma(a))bc)
  \surl{\ _{\hat{\alpha}}
  \otimes_{\beta}}_{\ \ \mu^o}v)|\zeta_k\surl{\ _{\beta}
  \otimes_{\alpha}}_{\ \mu} u)$$
By \ref{raccourci}, we obtain:
$$((\omega_{J_{\Psi}\Lambda_{\Psi}(e_k)}\surl{\ _{\beta}
  \star_{\alpha}}_{\ \mu}id)(\Gamma((id\surl{\ _{\beta}
  \star_{\alpha}}_{\ \mu}
  \omega_{J_{\Phi}\Lambda_{\Phi}(e),
  \eta})(\Gamma(a))bc))v|u)$$ which converges to:
$$((\Psi\surl{\ _{\beta}
  \star_{\alpha}}_{\ \mu} id)(\Gamma((id\surl{\ _{\beta}
  \star_{\alpha}}_{\ \mu}
  \omega_{J_{\Phi}\Lambda_{\Phi}(e),
  \eta})(\Gamma(a)))bc)v|u)$$
Now, by right invariance of $T_R$, \ref{rapide} and \ref{prepa},
we can compute this last expression:
$$
\begin{aligned}
&\ \quad ((\Psi\surl{\ _{\beta}
  \star_{\alpha}}_{\ \mu} id)(\Gamma((id\surl{\ _{\beta}
  \star_{\alpha}}_{\ \mu}
  \omega_{J_{\Phi}\Lambda_{\Phi}(e),
  \eta})(\Gamma(a)))bc)v|u)\\
&=(T_R((id\surl{\ _{\beta} \star_{\alpha}}_{\ \mu}
\omega_{J_{\Phi}\Lambda_{\Phi}(e),\eta})(\Gamma(a))bc)v|u)\\
&=(\Lambda_{\Psi}(bc) \surl{\ _{\hat{\alpha}}
  \otimes_{\beta}}_{\ \ \mu^o} v|\Lambda_{\Psi}((id\surl{\
_{\beta} \star_{\alpha}}_{\ \mu} \omega_{\eta,
J_{\Phi}\Lambda_{\Phi}(e)})(\Gamma(a^*)))) \surl{\
_{\hat{\alpha}}\otimes_{\beta}}_{\ \ \mu^o} u)\\
&=(v\surl{\ _{\beta}
  \otimes_{\alpha}}_{\ \mu}
  (\lambda^{\beta,\alpha}_{\Lambda_{\Psi}(\sigma_{-i}^{\Psi}(c^*))})^*
  U_{H_{\Psi}}(\Lambda_{\Psi}(b)\surl{\ _{\alpha}\otimes_{\hat{\beta}}}_{\ \ \mu^o}\Lambda_{\Phi}(a))
  |u\surl{\ _{\beta}
  \otimes_{\alpha}}_{\ \mu} J_{\Phi}e^*J_{\Phi}\eta)
\end{aligned}$$
which finishes the proof.
\end{proof}

Let $(\eta_i)_{i\in I}$be a $(N,\mu)$-basis of $\ _{\alpha}H$.
For all finite subset $J$ of $I$, we denote by $P_J$ the
projection $\sum_{i\in J}\theta^{\alpha,\mu}(\eta_i,\eta_i)\in
\alpha(N)'$ so that:
$$\sum_{i\in J}\rho_{\eta_i}^{\beta,\alpha}(\rho_{\eta_i}^{\beta,\alpha})^*=
1\surl{\ _{\beta} \otimes_{\alpha}}_{\ N} P_J$$ For all $e\in
{\mathcal N}_{\Phi}\cap {\mathcal N}_{T_L}$, we also denote by
$P_J^e$:
$$1\surl{\ _{\beta} \otimes_{\alpha}}_{\ N}
J_{\Phi}e^*J_{\Phi}P_JJ_{\Phi}eJ_{\Phi}=\sum_{i\in
J}\rho_{J_{\Phi}e^*J_{\Phi}\eta_i}^{\beta,\alpha}
(\rho_{J_{\Phi}e^*J_{\Phi}\eta_i}^{\beta,\alpha})^*$$

\begin{coro}\label{mef}
For all $a\in ({\mathcal N}_{\Psi}\cap {\mathcal N}_{T_R})^*
({\mathcal N}_{\Phi}\cap {\mathcal N}_{T_L})$, $b\in {\mathcal
N}_{\Psi}\cap {\mathcal N}_{T_R}$, and $c\in {\mathcal
T}_{\Psi,T_R}$, $v\in D(H_{\beta},\mu^o)$, $e\in {\mathcal
N}_{\Phi}\cap {\mathcal N}_{T_L}$ and $J$ finite subset of $I$, we
have:
$$P_J^eU_H(\sum_{i\in I} \eta_i\!\!\surl{\ _{\alpha}
  \otimes_{\hat{\beta}}}_{\ \ \mu^o}\! [(\lambda_{\zeta_k}^{\beta,\alpha})^*
  U_{H_{\Psi}}([(\rho_{\eta_i}^{\beta,\alpha})^*U'_H(J_{\Psi}c^*J_{\Psi}\Lambda_{\Psi}(b)\surl{\ _{\hat{\alpha}}
  \otimes_{\beta}}_{\ \ \mu^o}v)]\!\!\surl{\ _{\alpha}
  \otimes_{\hat{\beta}}}_{\ \ \mu^o}\! \Lambda_{\Phi}(a))])$$
converges, in the weak topology, to:
$$P_J^e(v\surl{\ _{\beta}\otimes_{\alpha}}_{\ \mu}
(\lambda^{\beta,\alpha}_{J_{\Psi}\Lambda_{\Psi}(c)})^*
U_{H_{\Psi}}(\Lambda_{\Psi}(b)\surl{\
_{\alpha}\otimes_{\hat{\beta}}}_{\ \ \mu^o}\Lambda_{\Phi}(a)))$$
\end{coro}

\begin{proof}
We apply to the reciprocity law
$\rho_{J_{\Phi}e^*J_{\Phi}\eta}^{\beta,\alpha}$ which is a
continuous linear operator of $H$ in $H\surl{\ _{\beta}
\otimes_{\alpha}}_{\ \mu}H_{\Phi}$, and also a continuous linear
operator of $H$ with weak topology in $H\surl{\ _{\beta}
\otimes_{\alpha}}_{\ \mu}H_{\Phi}$ with weak topology. Then, we
take finite sums.
\end{proof}

Until the end of the section, we denote by ${\mathcal H}_{\Phi}$
the closed linear span in $H_{\Phi}$ of
$(\lambda_w^{\beta,\alpha})^*U_{H_{\Psi}}(v\surl{\ _{\alpha}
\otimes_{\hat{\beta}}}_{\ \ \mu^o} \Lambda_{\Phi}(a))$ where
$v\in H_{\Psi}$, $w\in J_{\Psi}\Lambda_{\Psi}({\mathcal
N}_{\Psi}\cap {\mathcal N}_{T_R})$, and $a\in ({\mathcal
N}_{\Psi}\cap {\mathcal N}_{T_R})^*{\mathcal N}_{\Phi}\cap
{\mathcal N}_{T_L}$. By the third relation of lemma \ref{base}
(resp. proposition \ref{comm}), $\alpha$ (resp. $\hat{\beta}$) is
a non-degenerated (resp. anti-) representation of $N$ on
${\mathcal H}_{\Phi}$.

\begin{lemm}
Let $a\in ({\mathcal N}_{\Psi}\cap {\mathcal N}_{T_R})^*({\mathcal
N}_{\Phi}\cap {\mathcal N}_T)$, $b\in {\mathcal N}_{\Psi}\cap
{\mathcal N}_{T_R}$, $c\in {\mathcal T}_{\Psi,T_R}$, $v\in
D(H_{\beta},\mu^o)$ and $(\eta_i)_{i\in I}$ a $(N,\mu)$-basis of
$\ _{\alpha}H$. We put, for all $k\in K$:
$$\Xi_k=(\sum_{i\in I} \eta_i\surl{\ _{\alpha}
\otimes_{\hat{\beta}}}_{\ \ \mu^o}
[(\lambda_{\zeta_k}^{\beta,\alpha})^*
U_{H_{\Psi}}([(\rho_{\eta_i}^{\beta,\alpha})^*U'_H(J_{\Psi}c^*J_{\Psi}
\Lambda_{\Psi}(b)\surl{\ _{\hat{\alpha}}\otimes_{\beta}}_{\ \
\mu^o}v)]\surl{\ _{\alpha}\otimes_{\hat{\beta}}}_{\ \ \mu^o}
\Lambda_{\Phi}(a))]$$ Then the net $(\Xi_k)_{k\in K}$ is bounded.
\end{lemm}

\begin{proof}
Let $\Xi=v\surl{\ _{\beta}\otimes_{\alpha}}_{\ \mu}
(\lambda^{\beta,\alpha}_{J_{\Psi}\Lambda_{\Psi}(c)})^*
U_{H_{\Psi}}(\Lambda_{\Psi}(b)\surl{\
_{\alpha}\otimes_{\hat{\beta}}}_{\ \ \mu^o}\Lambda_{\Phi}(a))$.
By the previous corollary, we know that $P_J^eU_H\Xi_k$ weakly
converges to $P_J^e\Xi$, so that:
$$\lim_{J,||e||\leq 1}\lim_k P_J^eU_H\Xi_k =\Xi$$

Consequently, there exists $C\in\mathbb{R}^+$ such that:
$$\sup_{J,||e||\leq 1}\sup_k ||P_J^eU_H\Xi_k||\leq C$$
and, the interversion of the supremum gives:
$$C\geq\sup_k \sup_{J,||e||\leq 1}||P_J^eU_H\Xi_k||=\sup_k||U_H\Xi_k||=\sup_k||\Xi_k||$$
\end{proof}

\begin{coro}\label{precau}
For all $a\in ({\mathcal N}_{\Psi}\cap {\mathcal
N}_{T_R})^*({\mathcal N}_{\Phi}\cap {\mathcal N}_T)$, $b\in
{\mathcal N}_{\Psi}\cap {\mathcal N}_{T_R}$, $c\in {\mathcal
T}_{\Psi,T_R}$, $v\in D(H_{\beta},\mu^o)$ and $(\eta_i)_{i\in I}$
a $(N,\mu)$-basis of $\ _{\alpha}H$, we put: $$\Xi_k=(\sum_{i\in
I} \eta_i\surl{\ _{\alpha} \otimes_{\hat{\beta}}}_{\ \ \mu^o}
[(\lambda_{\zeta_k}^{\beta,\alpha})^*U_{H_{\Psi}}([(\rho_{\eta_i}^{\beta,\alpha})^*
U'_H(J_{\Psi}c^*J_{\Psi}\Lambda_{\Psi}(b)\surl{\ _{\hat{\alpha}}
\otimes_{\beta}}_{\ \ \mu^o}v)]\surl{\ _{\alpha}
\otimes_{\hat{\beta}}}_{\ \ \mu^o} \Lambda_{\Phi}(a))]$$ for all
$k\in K$, and:
$$\Xi=v\surl{\ _{\beta}\otimes_{\alpha}}_{\ \mu}
(\lambda^{\beta,\alpha}_{J_{\Psi}\Lambda_{\Psi}(c)})^*
U_{H_{\Psi}}(\Lambda_{\Psi}(b)\surl{\
_{\alpha}\otimes_{\hat{\beta}}}_{\ \ \mu^o}\Lambda_{\Phi}(a))$$
Then $U_H\Xi_k$ converges to $\Xi$ in the weak topology.
\end{coro}

\begin{proof}
Let $\Theta\in H\surl{\ _{\beta} \otimes_{\alpha}}_{\
\mu}H_{\Phi}$ and $\epsilon>0$. Then, there exists $e\in {\mathcal
N}_{\Phi}\cap {\mathcal N}_{T_L}$ of norm less than equal to $1$
and a finite subset $J$ of $I$ such that
$||(1-P_J^e)\Theta||\leq\epsilon$. By \ref{mef}, there also exists
$k_0$ such that $|(P_J^eU_H\Xi_k-P_J^e\Xi|\Theta)|\leq\epsilon$
for all $k\geq k_0$. Then, we get:
$$
\begin{aligned}
&\ \quad |(U_H\Xi_k-\Xi|\Theta)|\\
&\leq|(U_H\Xi_k-P_J^eU_H\Xi_k|\Theta)|+|(P_J^eU_H\Xi_k-P_J^e\Xi|\Theta)|+|(P_J^e\Xi-\Xi|\Theta)|\\
&\leq|(U_H\Xi_k|(1-P_J^e)\Theta)|+\epsilon+|(\Xi|(1-P_J^e)\Theta)|\\
&\leq|(U_H\Xi_k|(1-P_J^e)\Theta)|+\epsilon+|(\Xi|(1-P_J^e)\Theta)|
\leq (sup_{k\in k}||\Xi_k||+||\Xi||+1)\epsilon
\end{aligned}$$

\end{proof}

\begin{coro}
We have the following inclusion:
$$H\surl{\ _{\beta}\otimes_{\alpha}}_{\ \mu} {\mathcal H}_{\Phi} \subseteq
U_H(H\surl{\ _{\alpha}\otimes_{\hat{\beta}}}_{\ \ \mu^o}
{\mathcal H}_{\Phi})$$
\end{coro}

\begin{proof}
By the previous corollary, we know that $\Xi$ belongs to the weak
closure of $U_H(H\surl{\ _{\alpha}\otimes_{\hat{\beta}}}_{\ \
\mu^o} {\mathcal H}_{\Phi})$ which is also the norm closure. Now,
$U_H$ is an isometry, that's why $U_H(H\surl{\
_{\alpha}\otimes_{\hat{\beta}}}_{\ \ \mu^o} {\mathcal H}_{\Phi})$
is equal to $U_H(H\surl{\ _{\alpha}\otimes_{\hat{\beta}}}_{\ \
\mu^o} {\mathcal H}_{\Phi})$.
\end{proof}

\begin{theo}
$U_H: H \surl{\ _{\alpha} \otimes_{\hat{\beta}}}_{\ \ \mu^o}
H_{\Phi} \rightarrow H\surl{\ _{\beta} \otimes_{\alpha}}_{\ \mu}
H_{\Phi}$ is a unitary.
\end{theo}

\begin{proof}
By the previous corollary, we have:
\begin{equation}\label{inc}
H\surl{\ _{\beta}\otimes_{\alpha}}_{\ \mu} {\mathcal H}_{\Phi}
\subseteq U_H(H\surl{\ _{\alpha}\otimes_{\hat{\beta}}}_{\ \
\mu^o}{\mathcal H}_{\Phi}) \subseteq U_H(H \surl{\
_{\alpha}\otimes_{\hat{\beta}}}_{\ \ \mu^o} H_{\Phi}) \subseteq H
\surl{\ _{\beta} \otimes_{\alpha}}_{\ \mu} H_{\Phi}.
\end{equation}

Also, using a $(N^o,\mu^o)$-basis, we have, for all $v\in
H_{\Psi}$ and $a\in {\mathcal N}_{T_L}\cap {\mathcal N}_{\Phi}$:
$$U_{H_{\Psi}}(v\surl{\ _{\alpha} \otimes_{\hat{\beta}}}_{\ \
\mu^o}\Lambda_{\phi}(a))=\sum_i\xi_i\surl{\ _{\beta}
\otimes_{\alpha}}_{\
\mu}(\lambda_{\xi_i}^{\beta,\alpha})^*U_{H_{\Psi}}(v\surl{\
_{\alpha} \otimes_{\hat{\beta}}}_{\ \ \mu^o}\Lambda_{\phi}(a))$$
so that $U_{H_{\Psi}}(H_{\Psi}\surl{\ _{\alpha}
\otimes_{\hat{\beta}}}_{\ \ \mu^o}H_{\Phi}) \subseteq
H_{\Psi}\surl{\ _{\beta} \otimes_{\alpha}}_{\ \mu} {\mathcal
H}_{\Phi}$. The reverse inclusion is the relation \eqref{inc}
applied to $H_{\Psi}$. Consequently, we get:
$$U_{H_{\Psi}}(H_{\Psi}\surl{\ _{\alpha}
\otimes_{\hat{\beta}}}_{\ \
\mu^o}H_{\Phi})=U_{H_{\Psi}}(H_{\Psi}\surl{\ _{\alpha}
\otimes_{\hat{\beta}}}_{\ \ \mu^o}{\mathcal H}_{\Phi})$$ Since
$U_{H_{\Psi}}$ is an isometry, $H_{\Psi}\surl{\ _{\alpha}
\otimes_{\hat{\beta}}}_{\ \ \mu^o}H_{\Phi}=H_{\Psi}\surl{\
_{\alpha} \otimes_{\hat{\beta}}}_{\ \ \mu^o}{\mathcal H}_{\Phi}$
and, so ${\mathcal H}_{\Phi}=H_{\Phi}$. Finally, by inclusion
\eqref{inc}, we obtain $U_H(H \surl{\ _{\hat{\alpha}}
\otimes_{\beta}}_{\ \ \mu^o} H_{\Phi})=H \surl{\ _{\beta}
\otimes_{\alpha}}_{\ \mu} H_{\Phi}$.
\end{proof}


\begin{coro}\label{dense} If $[F]$ denote the linear span of a
subset $F$ of a vector space $E$, we have:

$$
\begin{aligned}
H_{\Phi}&=[\Lambda_{\Phi}((\omega_{v,w}\surl{\ _{\beta}
  \otimes_{\alpha}}_{\ \mu}id)(\Gamma(a)))|v,w\in
D(H_{\beta},\mu^o), a\in {\mathcal N}_{\Phi}\cap {\mathcal N}_{T_L}]\\
&=[(\lambda_w^{\beta,\alpha})^*U_H(v\surl{\ _{\alpha}
\otimes_{\hat{\beta}}}_{\ \ \mu^o} \Lambda_{\Phi}(a))|v\in H,w\in
D(H_{\beta},\mu^o), a\in {\mathcal N}_{\Phi}\cap {\mathcal N}_{T_L}]\\
&=[(\omega_{v,w}*id)(U_H)\xi|v\in D(_{\alpha}H,\mu),w\in
D(H_{\beta},\mu^o),\xi\in H_{\Phi}]
\end{aligned}$$
\end{coro}

\begin{proof}
The second equality comes from \ref{lienGV}. The last one is
clear. It's sufficient to prove that the last subspace is equal to
$H_{\Phi}$. Let $\eta\in H_{\Phi}$ in the orthogonal of:
$$[(\omega_{v,w}*id)(U_H)\xi|v\in D(_{\alpha}H,\mu),w\in
D(H_{\beta},\mu^o),\xi\in H_{\Phi}]$$ Then, for all $v\in
D(_{\alpha}H,\mu),w\in D(H_{\beta},\mu^o)$ and $\xi\in H_{\Phi}$,
we have:
$$(U_H(v\surl{\ _{\alpha}
\otimes_{\hat{\beta}}}_{\ \ \mu^o}\xi)|w\surl{\ _{\beta}
\otimes_{\alpha}}_{\
\mu}\eta)=((\omega_{v,w}*id)(U_H)\xi|\eta)=0$$ Since $U_H$ is a
unitary, $w\surl{\ _{\beta} \otimes_{\alpha}}_{\ \mu}\eta=0$ for
all $w\in D(H_{\beta},\mu^o)$ from which we easily deduce that
$\eta=0$ (by \ref{aidintre} for example).
\end{proof}

\begin{coro}\label{implementation}
We have $\Gamma(m)=U_H(1 \surl{\ _{\alpha}
\otimes_{\hat{\beta}}}_{\ \ N^o} m)U^*_H$ for all $m\in M$.
\end{coro}

\begin{proof}
Straightforward thanks to unitary of $U_H$ and \ref{impl}.
\end{proof}

\subsection{Pseudo-multiplicativity}
Let put $W=U^*_{H_{\Phi}}$. We have already proved commutation
relations of section \ref{rcom} and, now the aim is to prove that
$W$ is a pseudo-multiplicative unitary in the sense of M. Enock
and J.M Vallin (\cite{EV}, definition 5.6):

\begin{defi}
We call {\bf pseudo-multiplicative unitary} over $N$ w.r.t
$\alpha,\hat{\beta},\beta$ each unitary $V$ from $H \surl{\
_{\beta} \otimes_{\alpha}}_{\ \mu} H$ onto $H \surl{\ _{\alpha}
\otimes_{\hat{\beta}}}_{\ \ \mu^o} H$ which satisfies the
following commutation relations, for all $n,m\in N$:
$$(\beta(n)\surl{\ _{\alpha} \otimes_{\hat{\beta}}}_{\ \ N^o}
\alpha(m))V=V(\alpha(m) \surl{\ _{\beta} \otimes_{\alpha}}_{\
N}\hat{\beta}(n))$$ and $$(\hat{\beta}(n)\surl{\ _{\alpha}
\otimes_{\hat{\beta}}}_{\ \ N^o}\beta(m))V=
V(\hat{\beta}(n)\surl{\ _{\beta} \otimes_{\alpha}}_{\ N}
\beta(m))$$ and the formula:
$$(V \surl{\ _{\alpha}
\otimes_{\hat{\beta}}}_{\ \ N^o} 1)(\sigma_{\mu^o}\surl{\
_{\alpha}\otimes_{\hat{\beta}}}_{\ \ N^o} 1)(1\surl{\ _{\alpha}
\otimes_{\hat{\beta}}}_{\ \ N^o} V)\sigma_{2\mu}(1\surl{\ _{\beta}
\otimes_{\alpha}}_{\ N} \sigma_{\mu^o})(1\surl{\ _{\beta}
\otimes_{\alpha}}_{\ N} V)=$$
$$(1\surl{\ _{\alpha} \otimes_{\hat{\beta}}}_{\ \ N^o} V)
(V\surl{\ _{\beta} \otimes_{\alpha}}_{\ N}1)$$ where the first
$\sigma_{\mu^o}$ is the flip from $H \surl{\ _{\alpha}
\otimes_{\hat{\beta}}}_{\ \ \mu^o} H$ onto $H \surl{\
_{\hat{\beta}} \otimes_{\alpha}}_{\ \mu} H$, the second is the
flip from $H \surl{\ _{\alpha} \otimes_{\beta}}_{\ \ \mu^o} H$
onto $H \surl{\ _{\beta} \otimes_{\alpha}}_{\ \mu} H$ and
$\sigma_{2\mu}$ is the flip from $H \surl{\ _{\beta}
\otimes_{\alpha}}_{\ \mu} H \surl{\ _{\hat{\beta}}
\otimes_{\alpha}}_{\ \mu} H$ onto $H \surl{\ _{\alpha}
\otimes_{\hat{\beta}}}_{\ \ \mu^o} (H \surl{\ _{\beta}
\otimes_{\alpha}}_{\ \mu} H)$. This last flip turns around the
second tensor product. Moreover, parenthesis underline the fact
that the representation acts on the furthest leg.
\end{defi}

We recall, following (\cite{E2}, 3.5), if we use an other n.s.f
weight for the construction of relative tensor product, then
canonical isomorphisms of bimodules change the
pseudo-multiplicative unitary into another pseudo-multiplicative
unitary. The pentagonal relation comes essentially from the
co-product relation. So, we compute $(id\surl{\
_{\beta}\star_{\alpha}}_{\ N}\Gamma)\circ\Gamma $ and
$(\Gamma\surl{\ _{\beta}\star_{\alpha}}_{\ N}id)\circ\Gamma $ in
terms of $U_H$ with the following propositions \ref{igamma} and
\ref{gammai}. Until the end of the section, ${\mathcal H}$ is an
other Hilbert space on which $M$ acts.

\begin{lemm}
We have, for all $\xi_1\in D(_{\alpha}{\mathcal H},\mu)$ and
$\xi_2'\in D(H_{\beta},\mu^o)$:
$$\lambda_{\xi_1}^{\alpha,\hat{\beta}}(\lambda_{\xi_2'}^{\beta,\alpha})^*
=(\lambda_{\xi_2'}^{\beta,\alpha})^*\sigma_{2\mu^o}(1\surl{\
_{\alpha}\otimes_{\hat{\beta}}}_{\ \
N^o}\sigma_{\mu})\lambda_{\xi_1}^{\alpha,\hat{\beta}}$$ and:
$$U_{\mathcal
H}\lambda_{\xi_1}^{\alpha,\hat{\beta}}(\lambda_{\xi_2'}^{\beta,\alpha})^*U_H
=(\lambda_{\xi_2'}^{\beta,\alpha})^*(1\surl{\ _{\beta}
\otimes_{\alpha}}_{\ \ N}U_{\mathcal H})\sigma_{2\mu^o}(1\surl{\
_{\alpha}\otimes_{\hat{\beta}}}_{\ \ N^o}\sigma_{\mu})(1\surl{\
_{\alpha}\otimes_{\beta}}_{\ \
N^o}U_H)\lambda_{\xi_1}^{\alpha,\beta}$$
\end{lemm}

\begin{proof}
The first equality is easy to verify and the second one comes from
the first one.
\end{proof}

\begin{prop} The two following equations hold:
\begin{enumerate}[i)]\label{calcule}
\item for all $\xi_1\in D(
_{\alpha}{\mathcal H},\mu),\xi'_1\in D( _{\alpha}H,\mu),\xi_2\in
D({\mathcal H}_{\beta},\mu^o),\xi'_2\in D(H_{\beta},\mu^o)$ and
$\eta_1,\eta_2\in H_{\Phi}$, the scalar product of:
$$(1\surl{\ _{\beta}\otimes_{\alpha}}_{\ N}U_{\mathcal H})\sigma_{2\mu^o}(1\surl{\
_{\alpha}\otimes_{\hat{\beta}}}_{\ \ N^o}\sigma_{\mu})(1\surl{\
_{\alpha}\otimes_{\beta}}_{\ \ N^o}U_H)(\sigma_{\mu}\surl{\
_{\alpha}\otimes_{\hat{\beta}}}_{\ \ N^o}1)([\xi'_1\surl{\
_{\beta}\otimes_{\alpha}}_{\ \ \mu}\xi_1]\surl{\
_{\alpha}\otimes_{\hat{\beta}}}_{\ \ \mu^o} \eta_1)$$ by
$\xi'_2\surl{\ _{\beta}\otimes_{\alpha}}_{\ \mu}\xi_2\surl{\
_{\beta}\otimes_{\alpha}}_{\ \mu}\eta_2$ is equal to
$((\omega_{\xi_1,\xi_2}*id)(U_{\mathcal
H})(\omega_{\xi'_1,\xi'_2}*id)(U_H)\eta_1|\eta_2)$.
\item for all $a\in {\mathcal N}_{\Phi}\cap {\mathcal N}_{T_L}$,
$\xi_1\in {\mathcal H}$ and $\xi'_1,\xi'_2\in D(H_{\beta},\mu^o)$,
the value of:
$$(\lambda_{\xi_2'}^{\beta,\alpha})^*(1\surl{\ _{\beta}
\otimes_{\alpha}}_{\ \ N}U_{\mathcal H})\sigma_{2\mu^o}(1\surl{\
_{\alpha} \otimes_{\hat{\beta}}}_{\ \ N^o}\sigma_{\mu})(1\surl{\
_{\alpha}\otimes_{\beta}}_{\ \ N^o}U_H)(\sigma_{\mu}\surl{\
_{\alpha}\otimes_{\hat{\beta}}}_{\ \ N^o}1)$$ on $[\xi'_1\surl{\
_{\beta}\otimes_{\alpha}}_{\ \ \mu}\xi_1]\surl{\ _{\alpha}
\otimes_{\hat{\beta}}}_{\ \ \mu^o}\Lambda_{\Phi}(a)$ is equal to:
$$U_{\mathcal H}(\xi_1\surl{\ _{\alpha}\otimes_{\hat{\beta}}}_{\ \
N^o}\Lambda_{\Phi}((\omega_{\xi'_1,\xi'_2}\surl{\ _{\beta}
\star_{\alpha}}_{\ \mu}id)(\Gamma(a))))$$
\end{enumerate}
\end{prop}

\begin{proof}
By the previous lemma, we can compute the scalar product of i) in
the following way:
$$
\begin{aligned}
&\ \quad((\lambda_{\xi_2'}^{\beta,\alpha})^*(1\surl{\ _{\beta}
  \otimes_{\alpha}}_{\ \ N}U_{\mathcal H})\sigma_{2\mu^o}(1\surl{\ _{\alpha}
  \otimes_{\hat{\beta}}}_{\ \
  N^o}\sigma_{\mu})(1\surl{\ _{\alpha}
  \otimes_{\beta}}_{\ \ N^o}U_H)\lambda_{\xi_1}^{\alpha,\beta}(\xi'_1\surl{\ _{\alpha}
  \otimes_{\hat{\beta}}}_{\ \ \mu^o}\eta_1)|\xi_2\surl{\ _{\beta}
  \otimes_{\alpha}}_{\ \ \mu}\eta_2)\\
&=(U_{\mathcal
H}\lambda_{\xi_1}^{\alpha,\hat{\beta}}(\lambda_{\xi_2'}^{\beta,\alpha})^*U_H
(\xi'_1\surl{\ _{\alpha}
  \otimes_{\hat{\beta}}}_{\ \ \mu^o}\eta_1)|\xi_2\surl{\ _{\beta}
  \otimes_{\alpha}}_{\ \ \mu}\eta_2)\\
&=((\lambda_{\xi_2}^{\beta,\alpha})^*U_{\mathcal H}(\xi_1\surl{\
_{\alpha}
  \otimes_{\hat{\beta}}}_{\ \
  \mu^o}(\omega_{\xi'_1,\xi'_2}*id)(U_H)\eta_1 |\eta_2)\\
&=((\omega_{\xi_1,\xi_2}*id)(U_{\mathcal
H})(\omega_{\xi'_1,\xi'_2}*id)(U_H)\eta_1|\eta_2)
\end{aligned}$$

Also, the second assertion comes from the previous lemma and
\ref{rap}. Let's first assume that $\xi_1\in D(_{\alpha}{\mathcal
H},\mu)$. Then, we compute the vector in demand:
$$
\begin{aligned}
&\ \quad(\lambda_{\xi_2'}^{\beta,\alpha})^*(1\surl{\ _{\beta}
\otimes_{\alpha}}_{\ \ N}U_{\mathcal H})\sigma_{2\mu^o}(1\surl{\
_{\alpha}\otimes_{\hat{\beta}}}_{\ \ N^o}\sigma_{\mu})(1\surl{\
_{\alpha}\otimes_{\beta}}_{\ \
N^o}U_H)\lambda_{\xi_1}^{\alpha,\beta}(\xi'_1\surl{\ _{\alpha}
\otimes_{\hat{\beta}}}_{\ \ \mu^o}\Lambda_{\Phi}(a))\\
&=U_{\mathcal
H}\lambda_{\xi_1}^{\alpha,\hat{\beta}}(\lambda_{\xi_2'}^{\beta,\alpha})^*U_H
(\xi'_1\surl{\ _{\alpha}\otimes_{\hat{\beta}}}_{\ \
\mu^o}\Lambda_{\Phi}(a))\\
&=U_{\mathcal H}(\xi_1\surl{\ _{\alpha}\otimes_{\hat{\beta}}}_{\ \
N^o}\Lambda_{\Phi}((\omega_{\xi'_1,\xi'_2}\surl{\ _{\beta}
\star_{\alpha}}_{\ \mu}id)(\Gamma(a))))
\end{aligned}$$

So, we obtain the expected equality for all $\xi_1\in
D(_{\alpha}{\mathcal H},\mu)$. Since the two expressions are
continuous in $\xi_1$, density of $D(_{\alpha}{\mathcal H},\mu)$
in ${\mathcal H}$ implies that the equality is still true for all
$\xi_1\in {\mathcal H}$.
\end{proof}

\begin{prop}\label{igamma}
For all $a,b\in {\mathcal N}_{\Phi}\cap {\mathcal N}_{T_L}$, we
have:
$$(id\surl{\
_{\beta}\star_{\alpha}}_{\
\mu}\Gamma)(\Gamma(a))\rho_{J_{\Phi}\Lambda_{\Phi}(b)}^{\beta,\alpha}$$
$$=(1\surl{\ _{\beta}\otimes_{\alpha}}_{\ \ N}(1\surl{\
_{\beta}\otimes_{\alpha}}_{\ N}J_{\Phi}bJ_{\Phi})U_{\mathcal
H})\sigma_{2\mu^o}(1\!\!\surl{\
_{\alpha}\otimes_{\hat{\beta}}}_{\ \ N^o}\sigma_{\mu})(1\surl{\
_{\alpha}\otimes_{\beta}}_{\ \ N^o}U_H)(\sigma_{\mu}\surl{\
_{\alpha}\otimes_{\hat{\beta}}}_{\ \
N^o}1)\rho_{\Lambda_{\Phi}(a)}^{\alpha,\hat{\beta}}$$
\end{prop}

\begin{proof}
Let $\xi_1\in {\mathcal H}$ and $\xi'_1,\xi'_2\in
D(H_{\beta},\mu^o)$. We compose the second term of the equality on
the left by $(\lambda_{\xi_2'}^{\beta,\alpha})^*$ and we get:
$$(1\!\surl{\ _{\beta}\otimes_{\alpha}}_{\ N}\! J_{\Phi}bJ_{\Phi})
(\lambda_{\xi_2'}^{\beta,\alpha})^* (1\!\surl{\ _{\beta}
\otimes_{\alpha}}_{\ \ N}U_{\mathcal
H})\sigma_{2\mu^o}(1\!\surl{\ _{\alpha} \otimes_{\hat{\beta}}}_{\
\ N^o}\sigma_{\mu})(1\surl{\ _{\alpha}\otimes_{\beta}}_{\ \
N^o}U_H)(\sigma_{\mu}\surl{\ _{\alpha}\otimes_{\hat{\beta}}}_{\ \
N^o}\! 1)\rho_{\Lambda_{\Phi}(a)}^{\alpha,\hat{\beta}}$$ which we
evaluate on $\xi'_1\!\surl{\ _{\beta}\otimes_{\alpha}}_{\ \
\mu}\xi_1$, to get, by the previous proposition and
\ref{raccourci}:
$$
\begin{aligned}
&\ \quad(1\surl{\ _{\beta}\otimes_{\alpha}}_{\
N}J_{\Phi}bJ_{\Phi})U_{\mathcal H}(\xi_1\surl{\
_{\alpha}\otimes_{\hat{\beta}}}_{\ \
N^o}\Lambda_{\Phi}((\omega_{\xi'_1,\xi'_2}\surl{\ _{\beta}
\star_{\alpha}}_{\ \mu}id)(\Gamma(a))))\\
&=\Gamma((\omega_{\xi'_1,\xi'_2}\surl{\
_{\beta}\star_{\alpha}}_{\ \
\mu}id)(\Gamma(a)))\rho_{J_{\Phi}\Lambda_{\Phi}(b)}^{\beta,\alpha}\xi_1\\
&=(\lambda_{\xi_2'}^{\beta,\alpha})^*(id\surl{\
_{\beta}\star_{\alpha}}_{\
\mu}\Gamma)(\Gamma(a))\rho_{J_{\Phi}\Lambda_{\Phi}(b)}^{\beta,\alpha}(\xi'_1\surl{\
_{\beta}\otimes_{\alpha}}_{\ \ \mu}\xi_1)
\end{aligned}$$
So, the proposition holds.
\end{proof}

\begin{lemm}
For all $X\in M\surl{\ _{\beta}\star_{\alpha}}_{\ N}M\subset
(1\surl{\ _{\beta}\otimes_{\alpha}}_{\ N}\hat{\beta}(N))'$, we
have:
$$(\Gamma\surl{\ _{\beta}\star_{\alpha}}_{\
N}id)(X)=(U_H\surl{\ _{\beta}\otimes_{\alpha}}_{\ \ N}1)(1\surl{\
_{\alpha}\otimes_{\hat{\beta}}}_{\ \ N^o}X)(U_H^*\surl{\
_{\beta}\otimes_{\alpha}}_{\ \ N}1)$$
\end{lemm}

\begin{proof}
By \ref{implementation}, $\Gamma$ is implemented by $U_H$ so that
we easily deduce the lemma.
\end{proof}

\begin{prop}\label{gammai}
For all $a,b\in {\mathcal N}_{\Phi}\cap {\mathcal N}_{T_L}$, we
have:
$$(\Gamma\surl{\ _{\beta}\star_{\alpha}}_{\ N}id)(\Gamma(a))
\rho_{J_{\Phi}\Lambda_{\Phi}(b)}^{\beta,\alpha}$$
$$=(1\surl{\
_{\beta}\otimes_{\alpha}}_{\ N} 1\surl{\
_{\beta}\otimes_{\alpha}}_{\ N}J_{\Phi}bJ_{\Phi})(U_H\surl{\
_{\beta}\otimes_{\alpha}}_{\ N}1)(1\surl{\
_{\alpha}\otimes_{\hat{\beta}}}_{\ \ N^o}W^*)(U_H^*\surl{\
_{\alpha}\otimes_{\hat{\beta}}}_{\ \
N^o}1)\rho_{\Lambda_{\Phi}(a)}^{\alpha,\hat{\beta}}$$
\end{prop}

\begin{proof}
By the previous lemma and \ref{raccourci}, we can compute:
$$
\begin{aligned}
&\ \quad(1\surl{\ _{\beta}\otimes_{\alpha}}_{\ N} 1\surl{\
_{\beta}\otimes_{\alpha}}_{\ N}J_{\Phi}bJ_{\Phi})(U_H\surl{\
_{\beta}\otimes_{\alpha}}_{\ N}1)(1\surl{\
_{\alpha}\otimes_{\hat{\beta}}}_{\ \ N^o}W^*)(U_H^*\surl{\
_{\alpha}\otimes_{\hat{\beta}}}_{\ \
N^o}1)\rho_{\Lambda_{\Phi}(a)}^{\alpha,\hat{\beta}}\\
&=(U_H\surl{\ _{\beta}\otimes_{\alpha}}_{\ N}1)(1\surl{\
_{\alpha}\otimes_{\hat{\beta}}}_{\ \ N^o} 1\surl{\
_{\beta}\otimes_{\alpha}}_{\ N}J_{\Phi}bJ_{\Phi})(1\surl{\
_{\alpha}\otimes_{\hat{\beta}}}_{\ \
N^o}W^*)\rho_{\Lambda_{\Phi}(a)}^{\alpha,\hat{\beta}}U_H^*\\
&=(U_H\surl{\ _{\beta}\otimes_{\alpha}}_{\ N}1)(1\surl{\
_{\alpha}\otimes_{\hat{\beta}}}_{\ \ N^o}(1\surl{\
_{\beta}\otimes_{\alpha}}_{\
N}J_{\Phi}bJ_{\Phi})W^*\rho_{\Lambda_{\Phi}(a)}^{\alpha,\hat{\beta}})U_H^*\\
&=(U_H\surl{\ _{\beta}\otimes_{\alpha}}_{\ N}1)(1\surl{\
_{\alpha}\otimes_{\hat{\beta}}}_{\ \
N^o}\Gamma(a)\rho_{J_{\Phi}\Lambda_{\Phi}(b)}^{\beta,\alpha})U_H^*\\
&=(U_H\surl{\ _{\beta}\otimes_{\alpha}}_{\ N}1)(1\surl{\
_{\alpha}\otimes_{\hat{\beta}}}_{\ \ N^o}\Gamma(a))(U_H^*\surl{\
_{\alpha}\otimes_{\hat{\beta}}}_{\ \
N^o}1)\rho_{J_{\Phi}\Lambda_{\Phi}(b)}^{\beta,\alpha}=(\Gamma\surl{\
_{\beta}\star_{\alpha}}_{\ N}id)(\Gamma(a))
\rho_{J_{\Phi}\Lambda_{\Phi}(b)}^{\beta,\alpha}
\end{aligned}$$
\end{proof}

\begin{coro}The following relation is satisfied:
$$(U_H^*\!\!\surl{\ _{\alpha}
\otimes_{\hat{\beta}}}_{\ \ N^o} 1)(\sigma_{\mu^o}\!\! \surl{\
_{\alpha}\otimes_{\hat{\beta}}}_{\ \ N^o} 1)(1\!\! \surl{\
_{\alpha}\otimes_{\hat{\beta}}}_{\ \
N^o}U_H^*)\sigma_{2\mu}(1\!\! \surl{\ _{\beta}
\otimes_{\alpha}}_{\ N} \sigma_{\mu^o})(1\!\! \surl{\ _{\beta}
\otimes_{\alpha}}_{\ N}W)$$
$$=(1\!\!\surl{\ _{\alpha}\otimes_{\hat{\beta}}}_{\ \ N^o}W)(U_H^*\!\!\surl{\ _{\beta}
\otimes_{\alpha}}_{\ N} 1)$$
\end{coro}

\begin{proof}
We put together \ref{igamma} (with ${\mathcal H}=H_{\Phi}$) and
\ref{gammai} thanks to the co-product relation. We get:
$$(1 \surl{\ _{\beta} \otimes_{\alpha}}_{\ N} W^*)\sigma_{2\mu^o}(1
\surl{\ _{\alpha} \otimes_{\hat{\beta}}}_{\ \ N^o} \sigma_{\mu})(1
\surl{\ _{\alpha}
  \otimes_{\beta}}_{\ \ N^o}U_H)$$
$$=(U_H\surl{\ _{\beta}
\otimes_{\alpha}}_{\ N} 1)(1\surl{\ _{\alpha}
\otimes_{\hat{\beta}}}_{\ \ N^o} W^*)(U_H^*\surl{\
  _{\alpha} \otimes_{\hat{\beta}}}_{\ \ N^o} 1)(\sigma_{\mu^o}
  \surl{\ _{\alpha}
  \otimes_{\hat{\beta}}}_{\ \ N^o} 1)$$
Take adjoint and we are.
\end{proof}

\begin{theo}
$W$ is a pseudo-multiplicative unitary over $N$ w.r.t
$\alpha,\hat{\beta},\beta $.
\end{theo}

\begin{proof}
$W$ is a unitary from $H_{\Phi}\surl{\ _{\beta}
\otimes_{\alpha}}_{\ \mu}H_{\Phi}$ onto $H_{\Phi}\surl{\ _{\alpha}
\otimes_{\hat{\beta}}}_{\ \ \mu^o}H_{\Phi}$ which satisfies the
four required commutation relations. The previous corollary, with
$H=H_{\Phi}$, finishes the proof.
\end{proof}

Similar results hold for the right version:
\begin{theo}
If $W'=U'_{H_{\Psi}}$, then the following relation makes sense
and holds:
$$(W'\surl{\ _{\beta}
\otimes_{\alpha}}_{\ N} 1)(\sigma_{\mu}\surl{\ _{\beta}
\otimes_{\alpha}}_{\ N} 1)(1\surl{\ _{\beta} \otimes_{\alpha}}_{\
N}U'_H)\sigma_{2\mu^o}(1\surl{\ _{\hat{\alpha}}
\otimes_{\beta}}_{\ \ N^o}\sigma_{\mu})(1\surl{\ _{\hat{\alpha}}
\otimes_{\beta}}_{\ \ N^o}U'_H)$$
$$=(1\surl{\ _{\beta} \otimes_{\alpha}}_{\ N}U'_H)(W'\surl{\
  _{\hat{\alpha}} \otimes_{\beta}}_{\ \ N^o} 1)$$
If $H=H_{\Psi}$, then $W'$ is a pseudo-multiplicative unitary
over $N^o$ w.r.t $\beta,\alpha, \hat{\alpha}$.
\end{theo}

\begin{proof}
For example, it is sufficient to apply the previous results with
the opposite Hopf bimodule.
\end{proof}

\subsection{Von Neumann algebra generated by right leg
of the fundamental unitary}

\begin{defi}
We call $A(U'_H)$ (resp. ${\mathcal A}(U'_H)$) the weak closure
in ${\mathcal L}(H)$ of the vector space (resp. von Neumann
algebra) generated by $(\omega_{v,w}* id)(U'_H)$ with $v \in
D(_{\hat{\alpha}}(H_{\Psi}),\mu)$ and $w \in
D((H_{\Psi})_{\beta},\mu^o)$.
\end{defi}

\begin{prop}\label{situa}
$A(U'_H)$ is a non-degenerate involutive algebra i.e
$A(U'_H)={\mathcal A}(U'_H)$ such that:
$$\alpha(N)\cup\beta(N)\subseteq A(U'_H)={\mathcal A}(U'_H)\subseteq
M\subseteq\hat{\alpha}(N)'$$ Moreover, we have:
$$x\in {\mathcal A}(U'_H)'\cap\mathcal{L}(H)\Longleftrightarrow
U'_H(1 \surl{\ _{\hat{\alpha}} \otimes_{\beta}}_{\ \ N^o}
x)=(1\surl{\ _{\beta} \otimes_{\alpha}}_{\ N} x)U'_H$$
\end{prop}

In fact, we will see later that $A(U'_H)={\mathcal A}(U'_H)=M$.

\begin{proof}
The second and third points are obtained in \cite{EV} (theorem
6.1). As far as the first point is concerned, it comes from
\cite{E2} (proposition 3.6) and \ref{switch} which proves that
$A(U'_H)$ is involutive.
\end{proof}

To summarize the results of this section, we state the following
theorem:

\begin{theo}
Let $(N,M,\alpha,\beta,\Gamma)$ be a Hopf bimodule, $T_L$ (resp.
$T_R$) be a left (resp. right) invariant n.s.f operator-valued
weight. Then, for all n.s.f weight $\mu$ on $N$, if
$\Phi=\mu\circ\alpha^{-1}\circ T_L$, then the application:
$$v \surl{\ _{\alpha} \otimes_{\hat{\beta}}}_{\ \ \mu^o}
\Lambda_{\Phi}(a)\mapsto\sum_{i
    \in I} \xi_{i} \surl{\ _{\beta} \otimes_{\alpha}}_{\ \mu}
  \Lambda_{\Phi} ((\omega_{v,\xi_i} \surl{\
      _{\beta} \star_{\alpha}}_{\ \mu} id)(\Gamma(a)))$$
for all $v\in D((H_{\Phi})_{\beta},\mu^o)$, $a\in {\mathcal
N}_{T_L}\cap {\mathcal N}_{\Phi}$, $(N^o,\mu^o)$-basis
$(\xi_i)_{i\in I}$ of $(H_{\Phi})_{\beta}$ and where
$\hat{\beta}(n)=J_{\Phi}\alpha(n^*)J_{\Phi}$, extends to a
unitary $W$, the adjoint of which $W^*$ is a pseudo-multiplicative
unitary over $N$ w.r.t $\alpha,\hat{\beta},\beta$ from
$H_{\Phi}\surl{\ _{\alpha} \otimes_{\hat{\beta}}}_{\ \
\mu^o}H_{\Phi}$ onto $H_{\Phi}\surl{\ _{\beta}
\otimes_{\alpha}}_{\ \mu}H_{\Phi}$. Moreover, for all $m\in M$,
we have: $$\Gamma(m)=W^*(1\surl{\ _{\alpha}
\otimes_{\hat{\beta}}}_{\ \ N^o}m)W$$

Also, we have similar results from $T_R$.
\end{theo}

\section{Antipode}

Until the end, we introduce a new natural hypothesis which gives
a link between the right (resp. left) invariant operator-valued
weight and the (resp. anti-) representation of the basis. Then we
construct a closed antipode with polar decomposition which leads
to a co-involution and a one-parameter group of automorphisms of
$M$ called scaling group.

\subsection{Measured quantum groupoid's definition} \label{pop}
\begin{defi}
We say that a n.s.f operator-valued weight $T_L$ from $M$ to
$\alpha(N)$ is {\bf $\beta$-adapted} if there exists a n.s.f
weight $\nu_L$ on $N$ such that:
$$\sigma_t^{T_L}(\beta(n))=\beta(\sigma_{-t}^{\nu_L}(n))$$ for all
$n\in N$ and $t\in\mathbb{R}$. We also say that $T_L$ is
$\beta$-adapted w.r.t $\nu_L$.

We say that a n.s.f operator-valued weight $T_R$ from $M$ to
$\beta(N)$ is {\bf $\alpha$-adapted} if there exists a n.s.f
weight $\nu_R$ on $N$ such that:
$$\sigma_t^{T_R}(\alpha(n))=\alpha(\sigma_t^{\nu_R}(n))$$ for all
$n\in N$ and $t\in\mathbb{R}$. We also say that $T_R$ is
$\alpha$-adapted w.r.t $\nu_R$.
\end{defi}

\begin{defi}
A Hopf bimodule $(N,M,\alpha,\beta,\Gamma)$ with left (resp.
right) invariant n.s.f operator-valued weight $T_L$ (resp. $T_R$)
from $M$ to $\alpha(N)$ (resp. $\beta(N)$) is said to be a {\bf
measured quantum groupoid} if there exists a n.s.f weight $\nu$
on $N$ such that $T_L$ is $\beta$-adapted w.r.t $\nu$ and $T_R$
is $\alpha$-adapted w.r.t $\nu$. Then, we denote by
$(N,M,\alpha,\beta,\Gamma,\nu,T_L,T_R)$ the measured quantum
groupoid and we say that $\nu$ is {\bf quasi-invariant}.
\end{defi}

\begin{rema}
If a n.s.f operator-valued weight $T_L$ from $M$ to $\alpha(N)$
is $\beta$-adapted w.r.t $\nu$ and if $R$ is a co-involution of
$M$, then the n.s.f operator-valued weight $R\circ T_L\circ R$
from $M$ to $\beta(N)$ is $\alpha$-adapted w.r.t the same weight
$\nu$.
\end{rema}

\begin{lemm}\label{timpo}
If $\mu$ is a n.s.f weight on $N$ and if an operator-valued weight
$T_L$ is $\beta$-adapted w.r.t $\nu$, then there exists an
operator-valued weight $S^{\mu}$ from $M$ to $\beta(N)$, which is
$\alpha$-adapted w.r.t $\mu$ such that $\mu\circ\alpha^{-1}\circ
T_L=\nu\circ\beta^{-1}\circ S^{\mu}$. Also, if $\chi$ is a n.s.f
weight on $N$ and if an operator-valued weight $T_R$ is
$\alpha$-adapted w.r.t $\nu$, then there exists an
operator-valued weight $S_{\chi}$ from $M$ to $\alpha(N)$ normal,
which is $\beta$-adapted w.r.t $\chi$ such that
$\chi\circ\beta^{-1}\circ T_R=\nu\circ\beta^{-1}\circ S_{\chi}$.
\end{lemm}

\begin{proof}
For all $n\in N$ and $t\in\mathbb{R}$, we have
$\sigma_t^{\mu\circ\alpha^{-1}\circ
T_L}(\beta(n))=\sigma_t^{\nu\circ\beta^{-1}}(\beta(n))$. By
Haagerup's theorem, we obtain the existence of $S^{\mu}$ which is
clearly adapted. The second part of the lemma is very similar.
\end{proof}

Let $(N,M,\alpha,\beta,\Gamma,\nu,T_L,T_R)$ be a measured quantum
groupoid. Then the opposite measured quantum groupoid is
$(N^o,M,\beta,\alpha,\varsigma_N\circ\Gamma,\nu^o,T_R,T_L)$. We
put:
$$\Phi=\nu\circ\alpha^{-1}\circ T_L \quad\text{ and }\quad \Psi=\nu\circ\beta^{-1}\circ
T_R$$ We also put $S^{\nu}=S_L$ and $S_{\nu}=S_R$. By \ref{prem}
and \ref{evi}, we have: $$\Lambda_{\Phi}({\mathcal
T}_{\Phi,S_L})\subseteq J_{\Phi}\Lambda_{\Phi}({\mathcal
N}_{\Phi}\cap {\mathcal N}_{S_L})\subseteq
D((H_{\Phi})_{\beta},\nu^o)$$ and we have
$R^{\beta,\nu^o}(J_{\Phi}\Lambda_{\Phi}(a))=J_{\Phi}\Lambda_{S_L}(a)J_{\nu}$
for all $a\in {\mathcal N}_{\Phi}\cap {\mathcal N}_{S_L}$.

\subsection{The operator $G$}

We construct now an closed unbounded operator on $H_{\Phi}$ with
polar decomposition which gives needed elements to construct the
antipode. We have the following lemmas:

\begin{lemm}\label{jrel}
For all $\lambda \in \mathbb C$, $x\in {\mathcal
D}(\sigma_{i\lambda}^{\nu})$ and $\xi,\xi' \in
\Lambda_{\Phi}({\mathcal T}_{\Phi,T_L})$, we have:

\begin{equation}\label{rela}
\begin{aligned}
\alpha(x)\Delta_{\Phi}^{\lambda}
 &\subseteq \Delta_{\Phi}^{\lambda}
\alpha(\sigma_{i\lambda}^{\nu}(x)) \\
R^{\alpha,\nu}(\Delta_{\Phi}^{\lambda}\xi)\Delta_{\nu}^{\lambda}
&\subseteq \Delta_{\Phi}^{\lambda}R^{\alpha,\nu}(\xi) \\
\text{and }
\sigma_{i\lambda}^{\nu}(<\Delta_{\Phi}^{\lambda}\xi,\xi'>_{\alpha,\nu})&=<\xi,\Delta_{\Phi}^{\overline{\lambda}}\xi'>_{\alpha,\nu}
\end{aligned}
\end{equation}

and:

\begin{equation}\label{rela2}
\begin{aligned}
\hat{\beta}(x)\Delta_{\Phi}^{\lambda}
 &\subseteq \Delta_{\Phi}^{\lambda}
\hat{\beta}(\sigma_{i\lambda}^{\nu}(x)) \\
R^{\hat{\beta},\nu^o}(\Delta_{\Phi}^{\lambda}\xi)\Delta_{\nu}^{\lambda}
&\subseteq \Delta_{\Phi}^{\lambda}R^{\hat{\beta},\nu^o}(\xi) \\
\text{and }
\sigma_{i\lambda}^{\nu}(<\Delta_{\Phi}^{\lambda}\xi,\xi'>_{\hat{\beta},\nu^o})&=<\xi,\Delta_{\Phi}^{\overline{\lambda}}\xi'>_{\hat{\beta},\nu^o}.
\end{aligned}
\end{equation}
\end{lemm}

\begin{proof}
Straightforward.
\end{proof}

\begin{lemm}\cite{S3}\label{legitime}
We can define, for all $\lambda \in \mathbb C$, a closed operator
$\Delta_{\Phi}^{\lambda} \surl{\ _{\alpha}
\otimes_{\hat{\beta}}}_{\ \ \nu^o} \Delta_{\Phi}^{\lambda}$ which
naturally acts on elementary tensor products.
\end{lemm}

\begin{proof}
Let $\lambda \in \mathbb C$. We define a linear operator
$\Delta_{\lambda}$ on the algebraic tensor product
$\Lambda_{\Phi}({\mathcal T}_{\Phi,S_L})\odot {\mathcal
D}(\Delta_{\Phi}^{\lambda}) \subseteq D(_{\alpha}H_{\Phi},\nu)
\odot H_{\Phi}$ by the formula:
$$\Delta_{\lambda}(\xi \odot \eta) =\Delta_{\Phi}^{\lambda}
\xi \odot \Delta_{\Phi}^{\lambda}\eta$$ For all $\xi,\xi'\in
\Lambda_{\Phi}({\mathcal T}_{\Phi,S_L})$, $\eta\in {\mathcal
D}(\Delta_{\Phi}^{\lambda})$ and $\eta'\in {\mathcal
D}(\Delta_{\Phi}^{\overline{\lambda}})$, relations \eqref{rela}
imply:
$$
\begin{aligned}
(\Delta_{\lambda}(\xi\odot\eta)|\xi'\odot\eta') &=
(\hat{\beta}(<\Delta_{\Phi}^{\lambda}\xi,\xi'>_{\alpha,\nu})
\Delta_{\Phi}^{\lambda}\eta|\eta')\\
&=(\Delta_{\Phi}^{\lambda}\hat{\beta}(\sigma_{i\lambda}^{\nu}
(<\Delta_{\Phi}^{\lambda}\xi,\xi'>_{\alpha,\nu}))\eta|\eta')\\
&=(\hat{\beta}(<\xi,\Delta_{\Phi}^{\overline{\lambda}}\xi'>_{\alpha,\nu})
\eta|\Delta_{\Phi}^{\overline{\lambda}}\eta')\\
&=(\xi\surl{\ _{\alpha} \otimes_{\hat{\beta}}}_{\ \ \nu^o}\eta
|\Delta_{\Phi}^{\overline{\lambda}} \xi' \surl{\ _{\alpha}
\otimes_{\hat{\beta}}}_{\ \nu}
\Delta_{\Phi}^{\overline{\lambda}}\eta')
\end{aligned}$$

So, for all $\lambda\in \mathbb{C}$, $\Delta_{\lambda}$ go
through quotient to give a densely defined operator
$\tilde{\Delta}_{\lambda}$ on $H_{\Phi}\surl{\ _{\alpha}
\otimes_{\hat{\beta}}}_{\ \ \nu^o} H_{\Phi}$. Moreover, previous
equalities prove that
$\tilde{\Delta}_{\overline{\lambda}}\subseteq
\tilde{\Delta}_{\lambda}^*$. So, $\tilde{\Delta}_{\lambda}^*$ is
densely defined which means $\tilde{\Delta}_{\lambda}$ is
closable. The operator, we look for, is this closure. \cite{S3}
gives full ideas about the subject.
\end{proof}

Since, for all $x\in N$, we have
$J_{\Phi}\alpha(x)=\hat{\beta}(x^*)J_{\Phi}$, by \cite{S1}, we can
define a unitary anti-linear operator:
$$J_{\Phi}\surl{\ _{\alpha} \otimes_{\hat{\beta}}}_{\ \ \nu^o} J_{\Phi}:
H_{\Phi}\surl{\ _{\alpha} \otimes_{\hat{\beta}}}_{\ \
\nu^o}H_{\Phi} \rightarrow H_{\Phi}\surl{\ _{\hat{\beta}}
\otimes_{\alpha}}_{\ \nu} H_{\Phi}$$ such that the adjoint is
$J_{\Phi} \surl{\ _{\hat{\beta}} \otimes_{\alpha}}_{\ \nu}
J_{\Phi}$. Also, by composition, it is possible to define a
natural closed anti-linear operator:
$$S_{\Phi} \surl{\ _{\alpha}
\otimes_{\hat{\beta}}}_{\ \ \nu^o} S_{\Phi}: H_{\Phi}\surl{\
_{\alpha} \otimes_{\hat{\beta}}}_{\ \ \nu^o} H_{\Phi}\rightarrow
H_{\Phi}\surl{\ _{\hat{\beta}} \otimes_{\alpha}}_{\
\nu}H_{\Phi}$$ In the same way, if $F_{\Phi}=S_{\Phi}^*$, then it
is possible to define a natural closed anti-linear operator:
$F_{\Phi}\surl{\ _{\hat{\beta}}\otimes_{\alpha}}_{\ \nu}F_{\Phi}
: H_{\Phi}\surl{\ _{\hat{\beta}}\otimes_{\alpha}}_{\ \nu}
H_{\Phi}\rightarrow H_{\Phi}\surl{\
_{\alpha}\otimes_{\hat{\beta}}}_{\ \ \nu^o}H_{\Phi}$ and we have:
$$(S_{\Phi} \surl{\ _{\alpha} \otimes_{\hat{\beta}}}_{\ \ \nu^o}
S_{\Phi})^*=F_{\Phi}\surl{\ _{\hat{\beta}}\otimes_{\alpha}}_{\
\nu}F_{\Phi}$$

\begin{lemm}\label{chia}
For all $c\in ({\mathcal N}_{\Phi}\cap {\mathcal
N}_{T_L})^*({\mathcal N}_{\Psi}\cap {\mathcal N}_{T_R})$, $e\in
{\mathcal N}_{\Phi}\cap {\mathcal N}_{T_L}$ and all net
$(e_k)_{k\in K}$ of elements of ${\mathcal N}_{\Psi}\cap {\mathcal
N}_{T_R}$ weakly converging to $1$, then
$(\lambda_{J_{\Psi}\Lambda_{\Psi}(e_k)}^{\beta,\alpha})^*(1\surl{\
_{\beta} \otimes_{\alpha}}_{\
N}J_{\Phi}eJ_{\Phi})U_{H_{\Psi}}\rho_{\Lambda_{\Phi}(c^*)}^{\alpha,\hat{\beta}}
$ converges to
$(\lambda_{\Lambda_{\Psi}(c)}^{\hat{\alpha},\beta})^*U'^*_{H_{\Phi}}
\rho_{J_{\Phi}\Lambda_{\Phi}(e)}^{\beta,\alpha}$ in the weak
topology.
\end{lemm}

\begin{proof}
By \ref{raccourci}, we have, for all $k\in K$:

$$
\begin{aligned}
&\ \quad
(\lambda_{J_{\Psi}\Lambda_{\Psi}(e_k)}^{\beta,\alpha})^*(1\surl{\
_{\beta} \otimes_{\alpha}}_{\
\nu}J_{\Phi}eJ_{\Phi})U_{H_{\Psi}}\rho_{\Lambda_{\Phi}(c^*)}^{\alpha,\hat{\beta}}
\\
&=(\lambda_{J_{\Psi}\Lambda_{\Psi}(e_k)}^{\beta,\alpha})^*\Gamma(c^*)
\rho_{J_{\Phi}\Lambda_{\Phi}(e)}^{\beta,\alpha}=\left(
\Gamma(c)\lambda_{J_{\Psi}\Lambda_{\Psi}(e_k)}^{\beta,\alpha}\right)^*
\rho_{J_{\Phi}\Lambda_{\Phi}(e)}^{\beta,\alpha}\\
&=\left( (J_{\Psi}e_kJ_{\Psi}\surl{\ _{\beta} \otimes_{\alpha}}_{\
N}1)U'_{H_{\Phi}}\lambda_{\Lambda_{\Psi}(c)}^{\hat{\alpha},\beta}\right)^*
\rho_{J_{\Phi}\Lambda_{\Phi}(e)}^{\beta,\alpha}\\
&=(\lambda_{\Lambda_{\Psi}(c)}^{\hat{\alpha},\beta})^*U'^*_{H_{\Phi}}(J_{\Psi}e_k^*J_{\Psi}\surl{\
_{\beta} \otimes_{\alpha}}_{\ N}1)
\rho_{J_{\Phi}\Lambda_{\Phi}(e)}^{\beta,\alpha}=(\lambda_{\Lambda_{\Psi}(c)}^{\hat{\alpha},\beta})^*U'^*_{H_{\Phi}}
\rho_{J_{\Phi}\Lambda_{\Phi}(e)}^{\beta,\alpha}J_{\Psi}e_k^*J_{\Psi}
\end{aligned}$$
This computation implies the lemma.
\end{proof}

\begin{lemm}\label{dchia}
If $c\in ({\mathcal N}_{\Phi}\cap {\mathcal N}_{T_L})^*({\mathcal
N}_{\Psi}\cap {\mathcal N}_{T_R})$, $e\in {\mathcal N}_{\Phi}\cap
{\mathcal N}_{T_L}$, $\eta\in H_{\Psi}$, $v\in H_{\Phi}$ and a
net $(e_k)_{k\in K}$ of ${\mathcal N}_{\Psi}\cap {\mathcal
N}_{T_R}$ converges weakly to $1$, then the net:
$$((U_{H_{\Psi}}(\eta\surl{\ _{\alpha} \otimes_{\hat{\beta}}}_{\ \
\nu^o} \Lambda_{\Phi}(c^*))|J_{\Psi}\Lambda_{\Psi}(e_k)\surl{\
_{\beta} \otimes_{\alpha}}_{\ \nu} J_{\Phi}e^*J_{\Phi}v))_{k\in
K}$$ converges to $(
\eta|(\rho_{J_{\Phi}\Lambda_{\Phi}(e)})^*U'_{H_{\Phi}}(\Lambda_{\Psi}(c)\surl{\
_{\hat{\alpha}} \otimes_{\beta}}_{\ \ \nu^o}v))$.
\end{lemm}

\begin{proof}
It's a re-formulation of the previous lemma.
\end{proof}

\begin{prop}
Let $(\eta_i)_{i\in I}$ be a $(N,\nu)$-basis of $\ _{\alpha}H$,
$\Xi\in H_{\Psi}\surl{\ _{\beta} \otimes_{\alpha}}_{\ \nu} H$,
$u\in D(_{\alpha}H,\nu)$, $c\in ({\mathcal N}_{\Phi}\cap {\mathcal
N}_{T_L})^*({\mathcal N}_{\Psi}\cap {\mathcal N}_{T_R})$, $h\in
{\mathcal N}_{\Phi}\cap {\mathcal N}_{T_L}$ and $e$ be an element
of ${\mathcal N}_{\Phi}\cap {\mathcal N}_{T_L}\cap {\mathcal
N}_{\Phi}^*\cap {\mathcal N}_{T_L}^*$. Then, we have:

$$\lim_k\sum_{i\in I} (\eta_i\surl{\ _{\alpha}
\otimes_{\hat{\beta}}}_{\ \ \nu^o}
h^*(\lambda^{\beta,\alpha}_{J_{\Phi}\Lambda_{\Phi}(e_k)})^*U_{H_{\Psi}}(
(\rho_{\eta_i}^{\beta,\alpha})^*\Xi\surl{\ _{\alpha}
\otimes_{\hat{\beta}}}_{\ \ \nu^o}\Lambda_{\Phi}(c^*))| u\surl{\
_{\alpha} \otimes_{\hat{\beta}}}_{\ \
\nu^o}J_{\Phi}\Lambda_{\Phi}(e^*))$$ exists and is equal to
$((\rho^{\beta,\alpha}_u)^*\Xi|
(\rho^{\beta,\alpha}_{J_{\Phi}\Lambda_{\Phi}(e)})^*
U'_{H_{\Psi}}(\Lambda_{\Psi}(c)\surl{\ _{\hat{\alpha}}
\otimes_{\beta}}_{\ \ \nu^o}\Lambda_{\Phi}(h)))$.
\end{prop}

\begin{proof}
By \ref{base} and \ref{comm}, we can compute,  for all $i\in I$
and $k\in K$:

$$\begin{aligned}
&\quad\ (\eta_i\surl{\ _{\alpha} \otimes_{\hat{\beta}}}_{\ \
\nu^o}
h^*(\lambda^{\beta,\alpha}_{J_{\Phi}\Lambda_{\Phi}(e_k)})^*U_{H_{\Psi}}(
(\rho_{\eta_i}^{\beta,\alpha})^*\Xi\surl{\ _{\alpha}
\otimes_{\hat{\beta}}}_{\ \ \nu^o}\Lambda_{\Phi}(c^*))| u\surl{\
_{\alpha} \otimes_{\hat{\beta}}}_{\ \
\nu^o}J_{\Phi}\Lambda_{\Phi}(e^*))\\
&=\!\!(\hat{\beta}(<\eta_i,u>_{\alpha,\nu})
(\lambda^{\beta,\alpha}_{J_{\Phi}\Lambda_{\Phi}(e_k)})^*U_{H_{\Psi}}(
(\rho_{\eta_i}^{\beta,\alpha})^*\Xi\surl{\ _{\alpha}
\otimes_{\hat{\beta}}}_{\ \
\nu^o}\Lambda_{\Phi}(c^*))|gJ_{\Phi}\Lambda_{\Phi}(e^*))\\
&=\!\!((\lambda^{\beta,\alpha}_{J_{\Phi}\Lambda_{\Phi}(e_k)}\!)^*\!(1\!\!\surl{\
_{\beta} \otimes_{\alpha}}_{\
\nu}\!\!\hat{\beta}\!(<\eta_i,u>_{\alpha,\nu}\!)U_{H_{\Psi}}\!\!(
(\rho_{\eta_i}^{\beta,\alpha})^*\Xi\!\!\surl{\ _{\alpha}
\otimes_{\hat{\beta}}}_{\ \
\nu^o}\!\Lambda_{\Phi}(c^*)\!)|J_{\Phi}e^*J_{\Phi}\!\Lambda_{\Phi}(h)))\\
&=\!\!((\lambda^{\beta,\alpha}_{J_{\Phi}\Lambda_{\Phi}(e_k)})^*U_{H_{\Psi}}(
\beta(<\eta_i,u>_{\alpha,\nu})(\rho_{\eta_i}^{\beta,\alpha})^*
\Xi\surl{\ _{\alpha} \otimes_{\hat{\beta}}}_{\ \
\nu^o}\Lambda_{\Phi}(c^*))|J_{\Phi}e^*J_{\Phi}\Lambda_{\Phi}(h))
\end{aligned}$$
Take the sum over $i$ to obtain:
$$\sum_{i\in I}(\eta_i\surl{\ _{\alpha} \otimes_{\hat{\beta}}}_{\ \
\nu^o}
h^*(\lambda^{\beta,\alpha}_{J_{\Phi}\Lambda_{\Phi}(e_k)})^*U_{H_{\Psi}}(
(\rho_{\eta_i}^{\beta,\alpha})^*\Xi\surl{\ _{\alpha}
\otimes_{\hat{\beta}}}_{\ \ \nu^o}\Lambda_{\Phi}(c^*))| u\surl{\
_{\alpha} \otimes_{\hat{\beta}}}_{\ \
\nu^o}J_{\Phi}\Lambda_{\Phi}(e^*))$$
$$=(U_{H_{\Psi}}((\rho_u^{\beta,\alpha})^* \Xi\surl{\
_{\alpha} \otimes_{\hat{\beta}}}_{\ \
\nu^o}\Lambda_{\Phi}(c^*))|J_{\Phi}\Lambda_{\Phi}(e_k)\surl{\
_{\beta} \otimes_{\alpha}}_{\
\nu}J_{\Phi}e^*J_{\Phi}\Lambda_{\Phi}(h))$$ so that lemma
\ref{dchia} implies:
$$\lim_k\sum_{i\in I}(\eta_i\surl{\ _{\alpha} \otimes_{\hat{\beta}}}_{\ \
\nu^o}
h^*(\lambda^{\beta,\alpha}_{J_{\Phi}\Lambda_{\Phi}(e_k)})^*U_{H_{\Psi}}(
(\rho_{\eta_i}^{\beta,\alpha})^*\Xi\surl{\ _{\alpha}
\otimes_{\hat{\beta}}}_{\ \ \nu^o}\Lambda_{\Phi}(c^*))| u\surl{\
_{\alpha} \otimes_{\hat{\beta}}}_{\ \
\nu^o}J_{\Phi}\Lambda_{\Phi}(e^*))$$
$$=((\rho^{\beta,\alpha}_u)^*\Xi|
(\rho^{\beta,\alpha}_{J_{\Phi}\Lambda_{\Phi}(e)})^*
U'_{H_{\Phi}}(\Lambda_{\Psi}(c)\surl{\ _{\hat{\alpha}}
\otimes_{\beta}}_{\ \ \nu^o}\Lambda_{\Phi}(h)))\\
$$
\end{proof}

\begin{prop}\label{cx}
For all $a,c\in ({\mathcal N}_{\Phi}\cap {\mathcal
N}_{T_L})^*({\mathcal N}_{\Psi}\cap {\mathcal N}_{T_R})$, $b,d\in
{\mathcal T}_{\Psi,T_R}$ and $g,h\in {\mathcal T}_{\Phi,S_L}$,
the following vector:
$$U_{H_{\Phi}}^*\Gamma(g^*)(\Lambda_{\Phi}(h)\surl{\ _{\beta}
\otimes_{\alpha}}_{\ \nu}(\lambda^{\beta,\alpha}_{\Lambda_{\Psi}
(\sigma_{-i}^{\Psi}(b^*))})^*U_{H_{\Psi}}(\Lambda_{\Psi}(a)
\surl{\ _{\alpha} \otimes_{\hat{\beta}}}_{\ \
\nu^o}\Lambda_{\Phi}((cd)^*)))$$ belongs to ${\mathcal
D}(S_{\Phi}\surl{\ _{\alpha} \otimes_{\hat{\beta}}}_{\ \ \nu^o}
S_{\Phi})$ and the value of $\sigma_{\nu}(S_{\Phi}\surl{\
_{\alpha} \otimes_{\hat{\beta}}}_{\ \ \nu^o} S_{\Phi})$ on this
vector is equal to:
$$U_{H_{\Phi}}^*\Gamma(h^*)(\Lambda_{\Phi}(g)\surl{\ _{\beta}
\otimes_{\alpha}}_{\
\nu}(\lambda^{\beta,\alpha}_{\Lambda_{\Psi}(\sigma_{-i}^{\Psi}(d^*))})^*
U_{H_{\Psi}}(\Lambda_{\Psi}(c) \surl{\ _{\alpha}
\otimes_{\hat{\beta}}}_{\ \ \nu^o}\Lambda_{\Phi}((ab)^*)))$$
\end{prop}

\begin{proof}
For the proof, let denote by
$\Xi_1=U'_{H_{\Phi}}(\Lambda_{\Psi}(ab)\surl{\ _{\hat{\alpha}}
\otimes_{\beta}}_{\ \ \nu^o}\Lambda_{\Phi}(h))$ and by
$\Xi_2=U'_{H_{\Phi}}(\Lambda_{\Psi}(cd)\surl{\ _{\hat{\alpha}}
\otimes_{\beta}}_{\ \ \nu^o}\Lambda_{\Phi}(g))$. Then, for all
$e,f\in {\mathcal N}_{T_L}\cap {\mathcal N}_{\Phi}\cap{\mathcal
N}_{T_L}^*\cap {\mathcal N}_{\Phi}^*$, the scalar product of
$F_{\Phi}J_{\Phi}\Lambda_{\Phi}(e^*)\!\!\surl{\ _{\alpha}
\otimes_{\hat{\beta}}}_{\ \
\nu^o}\!F_{\Phi}J_{\Phi}\Lambda_{\Phi}(f)$ by:
$$U_{H_{\Phi}}^*\Gamma(g^*)(\Lambda_{\Phi}(h)\!\!\surl{\
_{\beta} \otimes_{\alpha}}_{\
\nu}(\lambda^{\beta,\alpha}_{\Lambda_{\Psi}
(\sigma_{-i}^{\Psi}(b^*))})^* U_{H_{\Psi}}(\Lambda_{\Psi}(a)
\!\!\surl{\ _{\alpha} \otimes_{\hat{\beta}}}_{\ \
\nu^o}\!\Lambda_{\Phi}((cd)^*)))$$ is equal to the scalar product
of $J_{\Phi}\Lambda_{\Phi}(e)\surl{\ _{\alpha}
\otimes_{\hat{\beta}}}_{\ \ \nu^o}J_{\Phi}\Lambda_{\Phi}(f^*)$ by:
$$U_{H_{\Phi}}^*\Gamma(g^*)(\Lambda_{\Phi}(h)\surl{\
_{\beta} \otimes_{\alpha}}_{\
\nu}(\lambda^{\beta,\alpha}_{\Lambda_{\Psi}
(\sigma_{-i}^{\Psi}(b^*))})^* U_{H_{\Psi}}(\Lambda_{\Psi}(a)
\surl{\ _{\alpha} \otimes_{\hat{\beta}}}_{\ \
\nu^o}\Lambda_{\Phi}((cd)^*)))$$ By \ref{precau}, this scalar
product is equal to the limit over $k$ of the sum over $i$ of:
$$(J_{\Phi}\Lambda_{\Phi}(e)\!\!\surl{\ _{\alpha}
\otimes_{\hat{\beta}}}_{\ \
\nu^o}\!J_{\Phi}\Lambda_{\Phi}(f^*)|\eta_i\!\!\surl{\ _{\alpha}
\otimes_{\hat{\beta}}}_{\ \ \nu^o}\!\!
g^*(\lambda^{\beta,\alpha}_{J_{\Psi}\Lambda_{\Psi}(e_k)})^*U_{H_{\Psi}}\!(
(\rho_{\eta_i}^{\beta,\alpha})^*\Xi_1\!\!\surl{\ _{\alpha}
\otimes_{\hat{\beta}}}_{\ \ \nu^o}\!\!\Lambda_{\Phi}((cd)^*)))$$
By the previous proposition applied with $\Xi=\Xi_1$, we get the
symmetric expression:
$$((\rho_{J_{\Phi}\Lambda_{\Phi}(f)}^{\beta,\alpha})^*\Xi_2|
(\rho_{J_{\Phi}\Lambda_{\Phi}(e)}^{\beta,\alpha})^*\Xi_1)$$ so
that, again by the previous proposition applied, this time, with
$\Xi=\Xi_2$ we obtain the limit over $k$ of the sum over $i$ of:
$$(\eta_i\!\!\surl{\ _{\alpha} \otimes_{\hat{\beta}}}_{\ \ \nu^o}\!
h^*(\lambda^{\beta,\alpha}_{J_{\Psi}\Lambda_{\Psi}(e_k)})^*U_{H_{\Psi}}(
(\rho_{\eta_i}^{\beta,\alpha})^*\Xi_2\!\!\surl{\ _{\alpha}
\otimes_{\hat{\beta}}}_{\ \ \nu^o}\!\Lambda_{\Phi}((ab)^*))|
J_{\Phi}\Lambda_{\Phi}(f)\!\!\surl{\ _{\alpha}
\otimes_{\hat{\beta}}}_{\ \ \nu^o}\!
J_{\Phi}\Lambda_{\Phi}(e^*))$$

This last expression is equal to the scalar product of:
$$U_{H_{\Phi}}^*\Gamma(h^*)(\Lambda_{\Phi}(g)\surl{\ _{\beta}
\otimes_{\alpha}}_{\
\nu}(\lambda^{\beta,\alpha}_{\Lambda_{\Psi}(\sigma_{-i}^{\Psi}(d^*))})^*
U_{H_{\Psi}}(\Lambda_{\Psi}(c) \surl{\ _{\alpha}
\otimes_{\hat{\beta}}}_{\ \ \nu^o}\Lambda_{\Phi}((ab)^*)))$$ by
$J_{\Phi}\Lambda_{\Phi}(f)\surl{\ _{\alpha}
\otimes_{\hat{\beta}}}_{\ \ \nu^o}J_{\Phi}\Lambda_{\Phi}(e^*)$
and to the scalar product of:
$$\sigma_{\nu^o}U_{H_{\Phi}}^*\Gamma(h^*)(\Lambda_{\Phi}(g)\surl{\
_{\beta} \otimes_{\alpha}}_{\
\nu}(\lambda^{\beta,\alpha}_{\Lambda_{\Psi}(\sigma_{-i}^{\Psi}(d^*))})^*
U_{H_{\Psi}}(\Lambda_{\Psi}(c) \surl{\ _{\alpha}
\otimes_{\hat{\beta}}}_{\ \ \nu^o}\Lambda_{\Phi}((ab)^*)))$$ by
$J_{\Phi}\Lambda_{\Phi}(e^*)\surl{\ _{\hat{\beta}}
\otimes_{\alpha}}_{\ \ \nu}J_{\Phi}\Lambda_{\Phi}(f)$. Since the
linear span of $J_{\Phi}\Lambda_{\Phi}(e^*)\surl{\ _{\hat{\beta}}
\otimes_{\alpha}}_{\ \ \nu}J_{\Phi}\Lambda_{\Phi}(f)$ where
$e,f\in {\mathcal N}_{T_L}\cap {\mathcal N}_{\Phi}\cap{\mathcal
N}_{T_L}^*\cap {\mathcal N}_{\Phi}^*$ is a core of
$F_{\Phi}\surl{\ _{\hat{\beta}}\otimes_{\alpha}}_{\
\nu}F_{\Phi}$, we get that:
$$U_{H_{\Phi}}^*\Gamma(g^*)(\Lambda_{\Phi}(h)\surl{\ _{\beta}
\otimes_{\alpha}}_{\ \nu}(\lambda^{\beta,\alpha}_{\Lambda_{\Psi}
(\sigma_{-i}^{\Psi}(b^*))})^*U_{H_{\Psi}}^*(\Lambda_{\Psi}(a)
\surl{\ _{\alpha} \otimes_{\hat{\beta}}}_{\ \
\nu^o}\Lambda_{\Phi}((cd)^*)))$$ belongs to ${\mathcal
D}(S_{\Phi}\surl{\ _{\alpha} \otimes_{\hat{\beta}}}_{\ \ \nu^o}
S_{\Phi})$ and the value of $S_{\Phi}\surl{\ _{\alpha}
\otimes_{\hat{\beta}}}_{\ \ \nu^o} S_{\Phi}$ on this vector is:
$$\sigma_{\nu^o}U_{H_{\Phi}}^*\Gamma(h^*)(\Lambda_{\Phi}(g)\surl{\ _{\beta}
\otimes_{\alpha}}_{\
\nu}(\lambda^{\beta,\alpha}_{\Lambda_{\Psi}(\sigma_{-i}^{\Psi}(d^*))})^*
U_{H_{\Psi}}(\Lambda_{\Psi}(c) \surl{\ _{\alpha}
\otimes_{\hat{\beta}}}_{\ \ \nu^o}\Lambda_{\Phi}((ab)^*)))$$
\end{proof}

\begin{prop}
There exists a closed densely defined anti-linear operator $G$ on
$H_{\Phi}$ such that the linear span of:
$$(\lambda^{\beta,\alpha}_{\Lambda_{\Psi}(\sigma_{-i}^{\Psi}(b^*))})^*
U_{H_{\Psi}}(\Lambda_{\Psi}(a) \surl{\ _{\alpha}
\otimes_{\hat{\beta}}}_{\ \ \nu^o}\Lambda_{\Phi}((cd)^*))$$ with $
a,c\in ({\mathcal N}_{\Phi}\cap {\mathcal N}_{T_L})^*({\mathcal
N}_{\Psi}\cap {\mathcal N}_{T_R}), b,d\in {\mathcal T}_{\Psi,T_R}
$, is a core of $G$ and we have:

$$G(\lambda^{\beta,\alpha}_{\Lambda_{\Psi}(\sigma_{-i}^{\Psi}(b^*))})^*
U_{H_{\Psi}}(\Lambda_{\Psi}(a) \surl{\ _{\alpha}
\otimes_{\hat{\beta}}}_{\ \ \nu^o}\Lambda_{\Phi}((cd)^*))$$
$$=(\lambda^{\beta,\alpha}_{\Lambda_{\Psi}(\sigma_{-i}^{\Psi}(d^*))})^*
U_{H_{\Psi}}(\Lambda_{\Psi}(c) \surl{\ _{\alpha}
\otimes_{\hat{\beta}}}_{\ \ \nu^o}\Lambda_{\Phi}((ab)^*))$$
Moreover, $G{\mathcal D}(G)={\mathcal D}(G)$ and $G^2=
id_{|{\mathcal D}(G)}$. \label{ferme}
\end{prop}

\begin{proof}
For all $n\in \mathbb N$, let $k_n \in\mathbb N$,
$a(n,l),c(n,l)\in ({\mathcal N}_{\Phi}\cap {\mathcal
N}_{T_L})^*({\mathcal N}_{\Psi}\cap {\mathcal N}_{T_R})$ and
$b(n,l),d(n,l)\in {\mathcal T}_{\Psi,T_R}$ and let $w\in H_{\Phi}$
such that:
$$v_n=\sum_{l=1}^{k_n} (\lambda^{\beta,\alpha}_{\Lambda_{\Psi}
(\sigma_{-i}^{\Psi}(b(n,l)^*))})^* U_{H_{\Psi}}
(\Lambda_{\Psi}(a(n,l)) \surl{\ _{\alpha}
\otimes_{\hat{\beta}}}_{\ \
\nu^o}\Lambda_{\Phi}((c(n,l)d(n,l))^*))\rightarrow 0$$
$$w_n=\sum_{l=1}^{k_n} (\lambda^{\beta,\alpha}_{\Lambda_{\Psi}
(\sigma_{-i}^{\Psi}(d(n,l)^*))})^*U_{H_{\Psi}}(\Lambda_{\Psi}(c(n,l))
\surl{\ _{\alpha} \otimes_{\hat{\beta}}}_{\ \
\nu^o}\Lambda_{\Phi}((a(n,l)b(n,l))^*))\rightarrow w$$ We have
$U_{H_{\Phi}}^*\Gamma(g^*)(\Lambda_{\Phi}(h)\surl{\ _{\beta}
\otimes_{\alpha}}_{\ \nu}v_n)\in {\mathcal D}(S_{\Phi}\surl{\
_{\alpha} \otimes_{\hat{\beta}}}_{\ \ \nu^o} S_{\Phi})$ for all
$g,h\in {\mathcal T}_{\Phi,S_L}$ and $n\in\mathbb{N}$ by the
previous proposition. Moreover, we have:
$$\sigma_{\nu}(S_{\Phi}\surl{\ _{\alpha}
\otimes_{\hat{\beta}}}_{\ \ \nu^o} S_{\Phi})
U_{H_{\Phi}}^*\Gamma(g^*)(\Lambda_{\Phi}(h)\surl{\ _{\beta}
\otimes_{\alpha}}_{\
\nu}v_n)=U_{H_{\Phi}}^*\Gamma(h^*)(\Lambda_{\Phi}(g)\surl{\
_{\beta} \otimes_{\alpha}}_{\ \nu}w_n)$$

Since $\Lambda_{\Phi}(g)$ and $\Lambda_{\Phi}(h)$ belongs to
$D((H_{\Phi})_{\beta},\nu^o)$, we obtain:
$$\sigma_{\nu}(S_{\Phi}\surl{\ _{\alpha}
\otimes_{\hat{\beta}}}_{\ \ \nu^o} S_{\Phi})
U_{H_{\Phi}}^*\Gamma(g^*)\lambda_{\Lambda_{\Phi}(h)}^{\beta,\alpha}v_n
=U_{H_{\Phi}}^*\Gamma(h^*)\lambda_{\Lambda_{\Phi}(g)}^{\beta,\alpha}w_n$$
The closure of $S_{\Phi}\surl{\ _{\alpha}
\otimes_{\hat{\beta}}}_{\ \ \nu^o} S_{\Phi}$ implies that
$U_{H_{\Phi}}^*\Gamma(h^*)\lambda_{\Lambda_{\Phi}(g)}^{\beta,\alpha}w=0$.
So, apply $U_{H_{\Phi}}$, to get
$\Gamma(h^*)\lambda_{\Lambda_{\Phi}(g)}^{\beta,\alpha}w=0$. Now,
${\mathcal T}_{\Phi,S_L}$ is dense in $M$ that's why
$\lambda_{\Lambda_{\Phi}(g)}^{\beta,\alpha}w=0$ for all $g\in
{\mathcal T}_{\Phi,S_L}$. Then, by \ref{evi}, we have:
$$||\lambda_{\Lambda_{\Phi}(g)}^{\beta,\alpha}w||^2
=(\alpha(<\Lambda_{\Phi}(g),\Lambda_{\Phi}(g)>_{\beta,\nu^o})w|w)\\
=(S_L(\sigma_{i/2}^{\Phi}(g)\sigma_{-i/2}^{\Phi}(g^*))w|w)$$ By
density of ${\mathcal T}_{\Phi,S_L}$, we obtain $||w||^2=0$ i.e
$w=0$. Consequently, the formula given in the proposition for $G$
gives rise to a closable densely defined well-defined operator on
$H_{\Phi}$. So the required operator is the closure of the
previous one.
\end{proof}

Thanks to polar decomposition of the closed operator $G$, we can
give the following definitions:

\begin{defi}\label{Ndefi}
We denote by $D$ the strictly positive operator $G^*G$ on
$H_{\Phi}$ (that means positive, self-adjoint and injective) and
by $I$ the anti-unitary operator on $H_{\Phi}$ such that
$G=ID^{1/2}$.
\end{defi}

Since $G$ is involutive, we have $I=I^*$, $I^2=1$ and
$IDI=D^{-1}$.

\subsection{A fundamental commutation relation}
In this section, we establish a commutation relation between $G$
and the elements $(\omega_{v,w}*id)(U'_{H_{\Phi}})$. We recall
that $W'=U'_{H_{\Psi}}$. We begin by two lemmas borrowed from
\cite{E2}.

\begin{lemm}
Let $\xi_i$ be a $(N^o,\nu^o)$-basis of $(H_{\Psi})_{\beta}$. For
all $w'\in D(_{\hat{\alpha}}H_{\Psi},\nu)$ and $w\in H_{\Psi}$,
we have: $$W'(w'\surl{\ _{\hat{\alpha}} \otimes_{\beta}}_{\ \
\nu^o} w)=\sum_{i}\xi_i \surl{\ _{\beta} \otimes_{\alpha}}_{\ \nu}
(\omega_{w',\xi_i} * id)(W')w$$ If we put
$\delta_i=(\omega_{w',\xi_i}
* id)(W')w$, then
$\alpha(<\xi_i,\xi_i>_{\beta,\nu^o})\delta_i=\delta_i$. Moreover,
if $w\in D(_{\hat{\alpha}}(H_{\Psi}),\nu)$, then $\delta_i\in
D(_{\hat{\alpha}}(H_{\Psi}),\nu)$.

For all $v,v'\in D((H_{\Psi})_{\beta},\nu^o)$ and $i\in I$, there
exists \mbox{$\zeta_i\in D((H_{\Psi})_{\beta},\nu^o)$} such that
$\alpha(<\xi_i,\xi_i>_{\beta,\nu^o})\zeta_i=\zeta_i$ and:
$$W'(v'\surl{\ _{\hat{\alpha}} \otimes_{\beta}}_{\ \ \nu^o}
v)=\sum_{i}\xi_i \surl{\ _{\beta} \otimes_{\alpha}}_{\ \nu}
\zeta_i$$
\end{lemm}

\begin{proof}
Lemma 3.4 of \cite{E2}.
\end{proof}

\begin{rema}
If $v,v'\in \Lambda_{\Psi}({\mathcal T}_{\Psi,T_R}) \subseteq
D(_{\hat{\alpha}}H, \nu) \cap D(H_{\beta},\nu^o)$, then, with
notations of the previous lemma, we have $\zeta_i\in
D(_{\hat{\alpha}}H, \nu) \cap D(H_{\beta},\nu^o)$.
\end{rema}

\begin{lemm}\label{alg}
Let $v,v'\in D(H_{\beta},\nu^o)$ and $w,w'\in D(_{\hat{\alpha}}H,
\nu)$. With notations of the previous lemma, we have:
$$(\omega_{v,w}*id)(U_H'{}^{\hspace{-.15cm}*})(\omega_{v',w'}*
id)(U_H'{}^{\hspace{-.15cm}*})=\sum_i(\omega_{\zeta_i,\delta_i} *
id)(U_H'{}^{\hspace{-.15cm}*})$$ in the norm convergence (and
also in the weak convergence).
\end{lemm}

\begin{proof}
Proposition 3.6 of \cite{E2}.
\end{proof}

\begin{lemm}
Let $a,c$ belonging to $({\mathcal N}_{\Phi}\cap {\mathcal
N}_{T_L})^*({\mathcal N}_{\Psi}\cap {\mathcal N}_{T_R})$. For all
$b,d,a',b',c',d'\in {\mathcal T}_{\Psi,T_R}$, the value of
$(\lambda_{\Lambda_{\Psi}(\sigma_{-i}^{\Psi}(b'^*))}
^{\beta,\alpha})^*U_{H_{\Psi}}$ on the sum over $i$ of:
$$\Lambda_{\Psi}((\omega_{\Lambda_{\Psi}(ab),\xi_i}\!*id)
(W')a')\!\!\surl{\ _{\alpha} \otimes_{\hat{\beta}}}_{\ \
\nu^o}\!\Lambda_{\Phi}(\!(c'd')^*
(\omega_{\xi_i,\Lambda_{\Psi}(cd)}*id)(W'^*))$$ is equal to:
$$(\omega_{\Lambda_{\Psi}(a'b'),\Lambda_{\Psi}(c'd')}*id)(U_{H_{\Phi}}'{}^{\hspace{-.3cm}*}\,)
(\lambda_{\Lambda_{\Psi}(\sigma_{-i}^{\Psi}(b^*))}^{\beta,\alpha})^*
U_{H_{\Psi}}(\Lambda_{\Psi}(a)\surl{\ _{\alpha}
\otimes_{\hat{\beta}}}_{\ \ \nu^o}\Lambda_{\Phi}((cd)^*))$$
\end{lemm}

\begin{proof}
First, let's suppose that $a\in {\mathcal T}_{\Psi,T_R}$. By
\ref{rap} and \ref{corres}, we have:
$$
\begin{aligned}
&\ \quad
(\omega_{\Lambda_{\Psi}(a'b'),\Lambda_{\Psi}(c'd')}*id)(U_{H_{\Phi}}'{}^{\hspace{-.3cm}*}\,)
(\lambda_{\Lambda_{\Psi}(\sigma_{-i}^{\Psi}(b^*))}^{\beta,\alpha})^*
U_{H_{\Psi}}(\Lambda_{\Psi}(a)\surl{\ _{\alpha}
\otimes_{\hat{\beta}}}_{\ \ \nu^o}\Lambda_{\Phi}((cd)^*))\\
&=(\omega_{\Lambda_{\Psi}(a'b'),\Lambda_{\Psi}(c'd')}*id)(U_{H_{\Phi}}'{}^{\hspace{-.3cm}*}\,)
\Lambda_{\Phi}((\omega_{\Lambda_{\Psi}(a),
\Lambda_{\Psi}(\sigma_{-i}^{\Psi}(b^*))}\surl{\ _{\beta}
\otimes_{\alpha}}_{\ \nu}id)(\Gamma((cd)^*)))\\
&=(\omega_{\Lambda_{\Psi}(a'b'),\Lambda_{\Psi}(c'd')}*id)(U_{H_{\Phi}}'{}^{\hspace{-.3cm}*}\,)
\Lambda_{\Phi}((\omega_{\Lambda_{\Psi}(ab),\Lambda_{\Psi}(cd)}*id)(U_{H_{\Phi}}'{}^{\hspace{-.3cm}*}\,))
\end{aligned}$$
By \ref{alg} and the closure of $\Lambda_{\Phi}$, this expression
is equal to the sum over $i\in I$ of:
$$\Lambda_{\Phi}((\omega_{(\omega_{\Lambda_{\Psi}(ab),\xi_i}*id)
(W')\Lambda_{\Psi}(a'b'),(\omega_{\Lambda_{\Psi}(cd),\xi_i}*id)(W')
\Lambda_{\Psi}(c'd')}*id)(U_{H_{\Phi}}'{}^{\hspace{-.3cm}*}\,))$$
Again, \ref{rap} and \ref{corres}, we obtain the sum over $i\in
I$ of the value of
$(\lambda_{\Lambda_{\Psi}(\sigma_{-i}^{\Psi}(b'^*))}
^{\beta,\alpha})^*U_{H_{\Psi}}$ on:
$$\Lambda_{\Psi}((\omega_{\Lambda_{\Psi}(ab),\xi_i}*id)
(W')a')\!\!\surl{\ _{\alpha} \otimes_{\hat{\beta}}}_{\ \
\nu^o}\Lambda_{\Phi}((c'd')^*(\omega_{\xi_i,\Lambda_{\Psi}(cd)}*id)(W'^*))$$
A density argument finishes the proof.
\end{proof}

\begin{prop}
If $v,w \in \Lambda_{\Psi}({\mathcal T}_{\Psi,T_R}^2) \subseteq
D(_{\hat{\alpha}}(H_{\Psi}), \nu) \cap
D((H_{\Psi})_{\beta},\nu^o)$, then we have:
\begin{align} (\omega_{v,w}* id)(U_{H_{\Phi}}'{}^{\hspace{-.3cm}*}\,)G & \subseteq
G(\omega_{w,v}* id)(U_{H_{\Phi}}'{}^{\hspace{-.3cm}*}\,) \label{inclu1}\\
\text{and} \ (\omega_{v,w}* id)(U'_{H_{\Phi}})G^* & \subseteq
G^*(\omega_{v,w}* id)(U'_{H_{\Phi}})\label{inclu2}
\end{align}

\end{prop}

\begin{proof}
Let $a,c\in ({\mathcal N}_{\Phi}\cap {\mathcal
N}_{T_L})^*({\mathcal N}_{\Psi}\cap {\mathcal N}_{T_R})$ and
$b,d,a',b',c',d'\in {\mathcal T}_{\Psi,T_R}$. By definition of
$G$, we have:
$$(\lambda_{\Lambda_{\Psi}(\sigma_{-i}^{\Psi}(d^*))}^{\beta,\alpha})^*
U_{H_{\Psi}}(\Lambda_{\Psi}(c)\surl{\ _{\alpha}
\otimes_{\hat{\beta}}}_{\ \ \nu^o}\Lambda_{\Phi}((ab)^*))\in
{\mathcal D}(G)$$ and:
$$(\omega_{\Lambda_{\Psi}(a'b'),\Lambda_{\Psi}(c'd')}*id)(U_{H_{\Phi}}'{}^{\hspace{-.3cm}*}\,)
G(\lambda_{\Lambda_{\Psi}(\sigma_{-i}^{\Psi}(d^*))}^{\beta,\alpha})^*
U_{H_{\Psi}}(\Lambda_{\Psi}(c)\surl{\ _{\alpha}
\otimes_{\hat{\beta}}}_{\ \ \nu^o}\Lambda_{\Phi}((ab)^*))$$
$$=(\omega_{\Lambda_{\Psi}(a'b'),\Lambda_{\Psi}(c'd')}*id)(U_{H_{\Phi}}'{}^{\hspace{-.3cm}*}\,)
(\lambda_{\Lambda_{\Psi}(\sigma_{-i}^{\Psi}(b^*))}^{\beta,\alpha})^*
U_{H_{\Psi}}(\Lambda_{\Psi}(a)\surl{\ _{\alpha}
\otimes_{\hat{\beta}}}_{\ \ \nu^o}\Lambda_{\Phi}((cd)^*))$$ By the
previous lemma, this is the sum over $i\in I$ of
$G(\lambda_{\Lambda_{\Psi}(\sigma_{-i}^{\Psi}(d'^*))}
^{\beta,\alpha})^*U_{H_{\Psi}}$ on:
$$\Lambda_{\Psi}((\omega_{\Lambda_{\Psi}(cd),\xi_i}*id)
(W')c')\!\!\surl{\ _{\alpha} \otimes_{\hat{\beta}}}_{\ \
\nu^o}\!\!\Lambda_{\Phi}((a'b')^*(\omega_{\xi_i,\Lambda_{\Psi}(ab)}*id)(W'^*))$$
Now, $G$ is a closed operator, so we deduce that the sum over
$i\in I$ of $(\lambda_{\Lambda_{\Psi}(\sigma_{-i}^{\Psi}(d'^*))}
^{\beta,\alpha})^*U_{H_{\Psi}}$ on:
$$\Lambda_{\Psi}((\omega_{\Lambda_{\Psi}(cd),\xi_i}*id)
(W')c')\!\!\surl{\ _{\alpha} \otimes_{\hat{\beta}}}_{\ \
\nu^o}\Lambda_{\Phi}((a'b')^*
(\omega_{\xi_i,\Lambda_{\Psi}(ab)}*id)(W'^*))$$ belongs to $
{\mathcal D}(G)$ and by the previous lemma, we obtain:
$$
\begin{aligned}
&\
\quad(\omega_{\Lambda_{\Psi}(a'b'),\Lambda_{\Psi}(c'd')}*id)(W'^*)
G(\lambda_{\Lambda_{\Psi}(\sigma_{-i}^{\Psi}(d^*))}^{\beta,\alpha})^*
U_{H_{\Psi}}(\Lambda_{\Psi}(c)\surl{\ _{\alpha}
\otimes_{\hat{\beta}}}_{\ \ \nu^o}\Lambda_{\Phi}((ab)^*))\\
&=G(\omega_{\Lambda_{\Psi}(c'd'),\Lambda_{\Psi}(a'b')}*id)(U_{H_{\Phi}}'{}^{\hspace{-.3cm}*}\,)
(\lambda_{\Lambda_{\Psi}(\sigma_{-i}^{\Psi}(d^*))}^{\beta,\alpha})^*
U_{H_{\Psi}}(\Lambda_{\Psi}(c)\surl{\ _{\alpha}
\otimes_{\hat{\beta}}}_{\ \ \nu^o}\Lambda_{\Phi}((ab)^*))
\end{aligned}$$
Now the linear span:
$$(\lambda^{\beta,\alpha}_{\Lambda_{\Psi}(\sigma_{-i}^{\Psi}(b^*))})^*
U_{H_{\Psi}}(\Lambda_{\Psi}(a) \surl{\ _{\alpha}
\otimes_{\hat{\beta}}}_{\ \ \nu^o}\Lambda_{\Phi}((cd)^*))$$ with
$a,c\in ({\mathcal N}_{\Phi}\cap {\mathcal N}_{T_L})^*({\mathcal
N}_{\Psi}\cap {\mathcal N}_{T_R}), b,d\in {\mathcal T}_{\Psi,T_R}
\}$, is a core for $G$ that's why the first inclusion holds. The
second one is the adjoint of the first one.
\end{proof}

\begin{coro}\label{coN}
For all $v,w \in \Lambda_{\Psi}({\mathcal T}_{\Psi,T_R}^2)$, we
have:
$$(\omega_{v,w}*id)(U'_{H_{\Phi}})D \subseteq
D(\omega_{\Delta_{\Psi}^{-1}v,\Delta_{\Psi}w}*
id)(U'_{H_{\Phi}})$$ where $D=G^*G$ is defined in \ref{Ndefi}.
\end{coro}

\begin{proof}
We have:
$$
\begin{aligned}
(\omega_{w,v}* id)(U'_{H_{\Phi}})G &=
(\omega_{S_{\Psi}w,\Delta_{\Psi}S_{\Psi}v} *
id)(U_{H_{\Phi}}'{}^{\hspace{-.3cm}*}\hspace{.1cm})G &&\quad
\text{by lemma
\ref{switch}} \\
& \subseteq G(\omega_{\Delta_{\Psi}S_{\Psi}v,S_{\Psi}w} \star
id)(U_{H_{\Phi}}'{}^{\hspace{-.3cm}*}\hspace{.1cm}) && \quad
\text{by inclusion\  \eqref{inclu1}} \\
&= G(\omega_{\Delta_{\Psi}^{-1}v,\Delta_{\Psi}w} *
id)(U_{H_{\Phi}}'{}^{\hspace{-.3cm}*}\hspace{.1cm}) && \quad
\text{by lemma \ref{switch}}
\end{aligned}$$

In the same way, we can finish the proof:

$$
\begin{aligned}
(\omega_{v,w}* id)(U'_{H_{\Phi}})D &= (\omega_{v,w}*
id)(U'_{H_{\Phi}})G^* &&
\quad \text{by definition \ref{Ndefi}} \\
& \subseteq G^*(\omega_{w,v}* id)(U'_{H_{\Phi}})G &&\quad
\text{by inclusion \
\eqref{inclu2}} \\
& \subseteq G^*G(\omega_{\Delta_{\Psi}^{-1}v,\Delta_{\Psi}w}*
id)(U'_{H_{\Phi}}) &&\\
&= D(\omega_{\Delta_{\Psi}^{-1}v,\Delta_{\Psi}w} *
id)(U'_{H_{\Phi}}) &&\quad \text{by definition \ref{Ndefi}}.
\end{aligned}$$

\end{proof}

\subsection{Scaling group}
In this section, we give a sense and we prove the following
commutation relation $U'_{H_{\Phi}}(\Delta_{\Psi} \surl{\
_{\hat{\alpha}} \otimes_{\beta}}_{\ \ \nu^o} D)=(\Delta_{\Psi}
\surl{\ _{\beta} \otimes_{\alpha}}_{\ \nu} D)U'_{H_{\Phi}}$ so as
to construct the scaling group $\tau$.

\begin{lemm}
For all $\lambda \in \mathbb C$ and $x$ analytic w.r.t $\nu$, we
have:
$$\alpha(x)D^{\lambda} \subseteq
D^{\lambda}\alpha(\sigma_{-i\lambda}^{\nu}(x))\text{ and }
\beta(x)D^{\lambda} \subseteq
D^{\lambda}\beta(\sigma_{-i\lambda}^{\nu}(x))$$\label{commuNa}
\end{lemm}

\begin{proof}
For all $a,c \in ({\mathcal N}_{\Phi}\cap {\mathcal
N}_{T_L})^*({\mathcal N}_{\Psi}\cap {\mathcal N}_{T_R})$, $b,d\in
{\mathcal T}_{\Psi,T_R}$ and $x$ analytic w.r.t $\nu$, we have by
\ref{base} and \ref{comm}:
$$
\begin{aligned}
&\ \quad
\beta(x)G(\lambda^{\beta,\alpha}_{\Lambda_{\Psi}(\sigma_{-i}^{\Psi}(b^*))})^*
U_{\Psi}(\Lambda_{\Psi}(a) \surl{\ _{\alpha}
\otimes_{\hat{\beta}}}_{\ \ \nu^o}\Lambda_{\Phi}((cd)^*))\\
&=\beta(x)(\lambda^{\beta,\alpha}_{\Lambda_{\Psi}(\sigma_{-i}^{\Psi}(d^*))})^*
U_{\Psi}(\Lambda_{\Psi}(c) \surl{\ _{\alpha}
\otimes_{\hat{\beta}}}_{\ \ \nu^o}\Lambda_{\Phi}((ab)^*))\\
&=(\lambda^{\beta,\alpha}_{\Lambda_{\Psi}(\sigma_{-i}^{\Psi}(d^*))})^*(1\surl{\
_{\beta} \otimes_{\alpha}}_{\
\nu}\beta(x))U_{\Psi}(\Lambda_{\Psi}(c) \surl{\ _{\alpha}
\otimes_{\hat{\beta}}}_{\ \ \nu^o}\Lambda_{\Phi}((ab)^*))\\
&=(\lambda^{\beta,\alpha}_{\Lambda_{\Psi}(\sigma_{-i}^{\Psi}(d^*))})^*U_{\Psi}(\Lambda_{\Psi}(c)
\surl{\ _{\alpha} \otimes_{\hat{\beta}}}_{\ \
\nu^o}\Lambda_{\Phi}(\beta(x)b^*a^*))\\
&=G(\lambda^{\beta,\alpha}_{\Lambda_{\Psi}(\beta(\sigma_i^{\nu}(x))\sigma_{-i}^{\Psi}(b^*))})^*
U_{\Psi}(\Lambda_{\Psi}(a) \surl{\ _{\alpha}
\otimes_{\hat{\beta}}}_{\ \ \nu^o}\Lambda_{\Phi}((cd)^*))\\
&=G\alpha(\sigma_{-i/2}^{\nu}(x^*))(\lambda^{\beta,\alpha}_{\Lambda_{\Psi}(\sigma_{-i}^{\Psi}(b^*))})^*
U_{\Psi}(\Lambda_{\Psi}(a) \surl{\ _{\alpha}
\otimes_{\hat{\beta}}}_{\ \ \nu^o}\Lambda_{\Phi}((cd)^*))
\end{aligned}$$

Now, the linear span of:
$$(\lambda^{\beta,\alpha}_{\Lambda_{\Psi}(\sigma_{-i}^{\Psi}(b^*))})^*
U_{\Psi}(\Lambda_{\Psi}(a) \surl{\ _{\alpha}
\otimes_{\hat{\beta}}}_{\ \ \nu^o}\Lambda_{\Phi}((cd)^*))$$ where
$a,c \in ({\mathcal N}_{\Phi}\cap {\mathcal N}_{T_L})^*({\mathcal
N}_{\Psi}\cap {\mathcal N}_{T_R}), b,d\in {\mathcal T}_{\Psi,T_R}
\}$, is a core for $G$, so that we have:
$$\beta(x)G \subseteq G\alpha(\sigma_{-i/2}^{\nu}(x^*))$$
Take adjoint to obtain $\alpha(x)G^* \subseteq
G^*\beta(\sigma_{i/2}^{\nu}(x^*))$. So, we conclude by:
$$\alpha(x)D=\alpha(x)G^*G \subseteq
G^*\beta(\sigma_{i/2}^{\nu}(x^*))G \subseteq
D\alpha(\sigma_{-i}^{\nu}(x))$$ The second part of the lemma can
be proved in a very similar way.
\end{proof}

We now state two lemmas analogous to relations (\ref{rela}) and
(\ref{rela2}) for $\Psi$ and we justify the existence of natural
operators:

\begin{lemm}
For all $\lambda\in\mathbb C$, $x\in {\mathcal
D}(\sigma_{-i\lambda}^{\nu})$ and $\xi,\xi' \in
\Lambda_{\Psi}({\mathcal T}_{\Psi,T_R})$, we have:

\begin{equation}\label{relat1}
\begin{aligned}
\beta(x)\Delta_{\Psi}^{\lambda}
 &\subseteq \Delta_{\Psi}^{\lambda}
\beta(\sigma_{-i\lambda}^{\nu}(x)) \\
R^{\beta,\nu^o}(\Delta_{\Psi}^{\lambda}\xi)\Delta_{\nu}^{-\lambda}
&\subseteq \Delta_{\Psi}^{\lambda}R^{\beta,\nu^o}(\xi) \\
\text{and }
\sigma_{-i\lambda}^{\nu}(<\Delta_{\Psi}^{\lambda}\xi,\xi'>_{\beta,\nu^o})&
=<\xi,\Delta_{\Psi}^{\overline{\lambda}}\xi'>_{\beta,\nu^o}
\end{aligned}
\end{equation}

and:

\begin{equation}\label{relat2}
\begin{aligned}
\hat{\alpha}(x)\Delta_{\Psi}^{\lambda}
 &\subseteq \Delta_{\Psi}^{\lambda}
\hat{\alpha}(\sigma_{-i\lambda}^{\nu}(x)) \\
R^{\hat{\alpha},\nu^o}(\Delta_{\Psi}^{\lambda}\xi)\Delta_{\nu}^{-\lambda}
&\subseteq \Delta_{\Psi}^{\lambda}R^{\hat{\alpha},\nu^o}(\xi) \\
\text{and }
\sigma_{-i\lambda}^{\nu}(<\Delta_{\Psi}^{\lambda}\xi,\xi'>_{\hat{\alpha},\nu^o})
&=<\xi,\Delta_{\Psi}^{\overline{\lambda}}\xi'>_{\hat{\alpha},\nu^o}
\end{aligned}
\end{equation}

\end{lemm}

\begin{proof}
It is sufficient to apply \ref{jrel} to the opposite measured
quantum groupoid for example.
\end{proof}

\begin{lemm}\label{legitime2}
We can define, for all $\lambda \in \mathbb C$, a closed linear
operator $\Delta_{\Psi}^{\lambda} \surl{\ _{\beta}
\otimes_{\alpha}}_{\ \nu} D^{\lambda}$ which naturally acts on
elementary tensor products.
\end{lemm}

\begin{proof}
The proof is very similar to \ref{legitime}.
\end{proof}

With relations \eqref{relat2} in hand, it's also possible to
define a closed linear operator $\Delta_{\Psi}^{\lambda} \surl{\
_{\hat{\alpha}} \otimes_{\beta}}_{\ \ \nu^o} D^{\lambda}$ on
$H_{\Psi} \surl{\ _{\hat{\alpha}} \otimes_{\beta}}_{\ \ \nu^o}
H_{\Phi}$.

\begin{prop}
The following relation holds:
\begin{equation}\label{fond}
U'_{H_{\Phi}}(\Delta_{\Psi} \surl{\ _{\hat{\alpha}}
\otimes_{\beta}}_{\ \ \nu^o} D)=(\Delta_{\Psi} \surl{\ _{\beta}
\otimes_{\alpha}}_{\ \nu} D)U'_{H_{\Phi}}
\end{equation}
\end{prop}

\begin{proof}
By \ref{coN}, we have, for all $v,w \in \Lambda_{\Psi}({\mathcal
T}_{\Psi,T_R})$ and $v',w' \in {\mathcal D}(D)$:

$$
\begin{aligned}
(U'_{H_{\Phi}}(v\surl{\ _{\hat{\alpha}} \otimes_{\beta}}_{\ \
\nu^o} v')|\Delta_{\Psi}w\surl{\ _{\beta} \otimes_{\alpha}}_{\
\nu} Dw')
&=((\omega_{v,\Delta_{\Psi}w}* id)(U'_{H_{\Phi}})v'|Dw') \\
&=(D(\omega_{\Delta_{\Psi}^{-1}(\Delta_{\Psi}v),\Delta_{\Psi}w}* id)(U'_{H_{\Phi}})v'|w') \\
&= ((\omega_{\Delta_{\Psi}v,w}*id)(U'_{H_{\Phi}})Dv'|w')\\
&=(U'_{H_{\Phi}}(\Delta_{\Psi}v\surl{\ _{\hat{\alpha}}
\otimes_{\beta}}_{\ \ \nu^o} Dv')|w \surl{\ _{\beta}
\otimes_{\alpha}}_{\ \nu} w')
\end{aligned}$$

By definition, we know that $\Lambda_{\Psi}({\mathcal
T}_{\Psi,T_R})\odot {\mathcal D}(D)$ is a core for
$\Delta_{\Psi}\surl{\ _{\beta} \otimes_{\alpha}}_{\ \nu}D$ so,
for all $u \in {\mathcal D}(\Delta_{\Psi}\surl{\ _{\beta}
\otimes_{\alpha}}_{\ \nu}D)$, we have:
$$(U'_{H_{\Phi}}(v\surl{\ _{\hat{\alpha}}
\otimes_{\beta}}_{\ \ \nu^o} v')|(\Delta_{\Psi}\surl{\ _{\beta}
\otimes_{\alpha}}_{\ \nu}
D)u)=(U'_{H_{\Phi}}(\Delta_{\Psi}v\surl{\ _{\hat{\alpha}}
\otimes_{\beta}}_{\ \ \nu^o} Dv')|u)$$ Since $\Delta_{\Psi}\surl{\
_{\beta} \otimes_{\alpha}}_{\ \nu}D$ is self-adjoint, we get:
$$(\Delta_{\Psi}\surl{\ _{\beta} \otimes_{\alpha}}_{\ \nu}D)U'_{H_{\Phi}}(v\surl{\
_{\hat{\alpha}} \otimes_{\beta}}_{\ \ \nu^o}
v')=U'_{H_{\Phi}}(\Delta_{\Psi}v\surl{\ _{\hat{\alpha}}
\otimes_{\beta}}_{\ \ \nu^o}Dv')$$ Finally, since
$\Lambda_{\Psi}({\mathcal T}_{\Psi,T_R})\odot {\mathcal D}(D)$ is
a core for $\Delta_{\Psi}\surl{\ _{\hat{\alpha}}
\otimes_{\beta}}_{\ \ \nu^o}D$ and by closeness of
$\Delta_{\Psi}\surl{\ _{\beta} \otimes_{\alpha}}_{\ \nu}D$, we
deduce that:
$$U'_{H_{\Phi}}(\Delta_{\Psi}\surl{\ _{\hat{\alpha}}
\otimes_{\beta}}_{\ \ \nu^o}D) \subseteq (\Delta_{\Psi}\surl{\
_{\beta} \otimes_{\alpha}}_{\ \nu} D)U'_{H_{\Phi}}$$ Because of
unitarity of $U'_{H_{\Phi}}$, we get that $(\Delta_{\Psi}\surl{\
_{\hat{\alpha}} \otimes_{\beta}}_{\ \ \nu^o}
D)U_{H_{\Phi}}'{}^{\hspace{-.3cm}*} \subseteq
U_{H_{\Phi}}'{}^{\hspace{-.3cm}*}(\Delta_{\Psi}\surl{\ _{\beta}
\otimes_{\alpha}}_{\ \nu}D)$ and by taking the adjoint, we get
the reverse inclusion:
$$(\Delta_{\Psi}\surl{\ _{\beta} \otimes_{\alpha}}_{\ \nu}D)U'_{H_{\Phi}}
\subseteq U'_{H_{\Phi}}(\Delta_{\Psi}\surl{\ _{\hat{\alpha}}
\otimes_{\beta}}_{\ \ \nu^o}D)$$
\end{proof}

We know begin the construction of the scaling group $\tau$
strictly speaking. We also prove a theorem which state that
$A(U'_H)=M$ and generalize proposition 1.5 of \cite{KV2}.

\begin{defi}
We denote by $M_R$ the weakly closed linear span of: $$\{
(\omega\surl{\ _{\beta} \star_{\alpha}}_{\ \nu} id)(\Gamma(x))\ |
\ x\in M,\ \omega \in M^+_* \text{ s.t } \exists k\in {\mathbb
R}^+\!\!,\ \omega\circ\beta \leq k\nu\}$$ Also, we denote by
$M_L$ the weakly closed linear span of:
$$\{ (id\surl{\
_{\beta} \star_{\alpha}}_{\ \nu} \omega)(\Gamma(x))\ | \ x\!\in\!
M,\ \omega\! \in\! M^+_* \text{ s.t } \exists k\in {\mathbb
R}^+\!\!, \ \omega\circ\alpha \leq k\nu\}$$
\end{defi}

By \ref{corres} and \ref{situa}, $M_R$ is equal to the von
Neumann subalgebra $A(U'_H)$ of $M$. Also, $M_L$ is a von Neumann
subalgebra of $M$. Moreover, we know $\alpha(N)\subseteq M_R$ and
$\beta(N) \subseteq M_L$, so that $M_L \surl{\ _{\beta}
\star_{\alpha}}_{\ \nu} M_R$ makes sense. Also, we have, for all
$m\in M$:

\begin{equation}\label{appart}
\Gamma(m)\in M_L \surl{\ _{\beta} \star_{\alpha}}_{\ N} M_R
\end{equation}

\begin{lemm}\label{tau}
There exists a unique strongly continuous one-parameter group
$\tau$ of automorphisms of $M_R$ such that
$\tau_t(x)=D^{-it}xD^{it}$ for all $t \in \mathbb R$ and $x\in
M_R$.
\end{lemm}

\begin{proof}
By commutation relation \eqref{fond}, for all $t\in \mathbb R$ and
$v,w \in \Lambda_{\Psi}({\mathcal T}_{\Psi,T_R})$, we get that:
$$D^{-it}(\omega_{v,w} * id)(U'_{H_{\Phi}})D^{it}
=(\omega_{\Delta_{\Psi}^{-it}v,\Delta_{\Psi}^{it}w}*
id)(U'_{H_{\Phi}})$$ Consequently, we obtain
$D^{-it}M_RD^{it}=M_R$ which is the only point to show.
\end{proof}

\begin{lemm}\label{mieux}
We have $\tau_t(\alpha(n))=\alpha(\sigma_t^{\nu}(n))$ for all
$n\in N$ and $t\in\mathbb{R}$.
\end{lemm}

\begin{proof}
Straightforward by lemma \ref{commuNa}.
\end{proof}

\begin{lemm}\label{tau2}
It is possible to define a normal *-automorphism $\sigma_t^{\Psi}
\surl{\ _{\beta} \star_{\alpha}}_{\ N} \tau_{-t}$ of $M \surl{\
_{\beta} \otimes_{\alpha}}_{\ N} M_R$ which naturally acts for all
$t\in\mathbb R$.
\end{lemm}

\begin{proof}
By the previous lemma and relations \eqref{relat1}, we have, for
all $n\in N$ and $t\in\mathbb{R}$:
$$\tau_t(\alpha(n))=\alpha(\sigma_t^{\nu}(n))\quad\text{ and }\quad
\sigma_t^{\Psi}(\beta(n))=\beta(\sigma_{-t}^{\nu}(n))$$ so that
it is possible to define a morphism:
$$\sigma_t^{\Psi} \surl{\ _{\beta} \star_{\alpha}}_{\ N} \tau_{-t} \ \text{:} \
M \surl{\ _{\beta} \star_{\alpha}}_{\ N} M_R \rightarrow M
\surl{\ _{\beta\circ\sigma_{-t}^{\nu}}
\star_{\alpha\circ\sigma_{-t}^{\nu}}}_{\ N} M_R$$ Then, it is
sufficient to prove that the range is equal to the domain. For all
$\xi \in D(H_{\beta},\nu^o)$ and $y\in {\mathcal N}_{\nu}$, we
compute:
$$
\begin{aligned}
\beta(\sigma_{-t}^{\nu}(y^*))\xi &
=\beta(\sigma_{-t}^{\nu}(y)^*)\xi=R^{\beta,\nu^o}(\xi)J_{\nu}\Lambda_{\nu}(\sigma_{-t}^{\nu}(y))\\
&=R^{\beta,\nu^o}(\xi)J_{\nu}\Delta_{\nu}^{-it}\Lambda_{\nu}(y)=R^{\beta,\nu^o}(\xi)\Delta_{\nu}^{-it}J_{\nu}\Lambda_{\nu}(y)
\end{aligned}$$
and we get, for all $\xi \in D(H_{\beta},\nu^o)$:
$$\xi \in D(H_{\beta\circ\sigma_{-t}^{\nu}},\nu^o) \ \text{and} \
R^{\beta\circ\sigma_{-t}^{\nu},\nu^o}(\xi)=R^{\beta,\nu^o}(\xi)\Delta_{\nu}^{-it}$$
To conclude, we show that scalar products on $H \odot H$ used to
define $H \surl{\ _{\beta} \otimes_{\alpha}}_{\ \nu}H$ and
$H\surl{\ _{\beta\circ\sigma_{-t}^{\nu}}
\otimes_{\alpha\circ\sigma_{-t}^{\nu}}}_{\ \nu} H$ are equal. For
all $\xi,\xi' \in D(H_{\beta},\nu^o)$ and $\eta,\eta' \in H$, we
have:
$$
\begin{aligned}
(\xi \surl{\ _{\beta\circ\sigma_{-t}^{\nu}}
\otimes_{\alpha\circ\sigma_{-t}^{\nu}}}_{\ \nu} \eta|\xi' \surl{\
_{\beta\circ\sigma_{-t}^{\nu}}
\otimes_{\alpha\circ\sigma_{-t}^{\nu}}}_{\ \nu} \eta')
&=(\alpha(\sigma_{-t}^{\nu}(<\xi,\xi'>_{\beta\circ\sigma_{-t}^{\nu},\nu^o}))\eta|\eta')\\
&=(\alpha(\sigma_{-t}^{\nu}(\Delta_{\nu}^{it}<\xi,\xi'>_{\beta,\nu^o}\Delta_{\nu}^{-it}))\eta|\eta')\\
&=(\alpha(<\xi,\xi'>_{\beta,\nu^o})\eta|\eta')\\
&=(\xi\surl{\ _{\beta} \otimes_{\alpha}}_{\ \nu}\xi'|\eta \surl{\
_{\beta} \otimes_{\alpha}}_{\ \nu} \eta')
\end{aligned}$$
\end{proof}

\begin{prop}\label{egal}
We have $(\sigma_t^{\Psi} \surl{\ _{\beta} \star_{\alpha}}_{\ N}
\tau_{-t})\circ\Gamma=\Gamma \circ \sigma_t^{\Psi}$ for all $t\in
\mathbb R$.
\end{prop}

\begin{proof}
By relation \eqref{appart}, the formula makes sense ($\tau$ is
just defined on $M_R$). By relation \eqref{fond}, we can compute
for all $m\in M$ and $t\in\mathbb R$:
$$
\begin{aligned}
(\sigma_t^{\Psi} \surl{\ _{\beta} \star_{\alpha}}_{\ \nu}
\tau_{-t})\circ\Gamma(m)&=(\Delta_{\Psi}^{it} \surl{\ _{\beta}
\otimes_{\alpha}}_{\ \nu}D^{it})\Gamma(m)(\Delta_{\Psi}^{-it}
\surl{\ _{\beta} \otimes_{\alpha}}_{\ \nu}
D^{-it})\\
&=(\Delta_{\Psi}^{it} \surl{\ _{\beta} \otimes_{\alpha}}_{\
\nu}D^{it})U'_{H_{\Phi}}(m\surl{\ _{\hat{\alpha}}
\otimes_{\beta}}_{\ \ \nu^o} \!\!
1)U_{H_{\Phi}}'{}^{\hspace{-.3cm}*} (\Delta_{\Psi}^{-it} \surl{\
_{\beta} \otimes_{\alpha}}_{\ \nu}D^{-it})\\
&=U'_{H_{\Phi}}(\Delta_{\Psi}^{it} \surl{\ _{\hat{\alpha}}
\otimes_{\beta}}_{\ \nu^o}D^{it})(m\surl{\ _{\hat{\alpha}}
\otimes_{\beta}}_{\ \ \nu^o} \!\! 1)(\Delta_{\Psi}^{-it} \surl{\
_{\hat{\alpha}} \otimes_{\beta}}_{\ \nu^o}
D^{-it})U_{H_{\Phi}}'{}^{\hspace{-.3cm}*}\\
&=U'_{H_{\Phi}}(\sigma_t^{\Psi}(m)\surl{\ _{\hat{\alpha}}
\otimes_{\beta}}_{\ \nu^o}
1)U_{H_{\Phi}}'{}^{\hspace{-.3cm}*}=\Gamma(\sigma_t^{\Psi}(m))
\end{aligned}$$

\end{proof}

We are now able to prove that we can re-construct $M$ thanks to
the fundamental unitary.

\begin{theo}\label{densevn}
If $<F>^{-\textsc{w}}$ is the weakly closed linear span of $F$ in
$M$, then following vector spaces:
$$\begin{aligned}
M_R&=<(\omega\surl{\ _{\beta} \star_{\alpha}}_{\ \nu}
id)(\Gamma(m))\ | \ m\in M,\omega\in M^+_*, k\in {\mathbb
R}^+\text{ s.t } \omega\circ\beta \leq k\nu>^{-\textsc{w}}\\
A(U'_H)&=<(\omega_{v,w}* id)(U'_H) \ |\ v \in
D(_{\hat{\alpha}}(H_{\Psi}),\mu), w \in
D((H_{\Psi})_{\beta},\mu^o)>^{-\textsc{w}}\\
M_L&=<(id\surl{\ _{\beta}\star_{\alpha}}_{\ \nu}
\omega)(\Gamma(m))\ | \ m\in M, \omega\in M^+_*, k\in {\mathbb
R}^+ \ \text{s.t } \omega\circ\alpha \leq k\nu>^{-\textsc{w}}\\
A(U_H)&=<(id*\omega_{v,w})(U_H) \ |\ v \in
D((H_{\Psi}),\mu^o)_{\hat{\beta}}, w \in
D(_{\alpha}(H_{\Psi}),\mu)>^{-\textsc{w}}
\end{aligned}$$
are equal to the whole von Neumann algebra $M$.
\end{theo}

\begin{proof}
We have already noticed that $M_R=A(U'_H)$ and $M_L=A(U_H)$.
Then, we get inspired by \cite{KV2}. By \ref{mieux}, we have
$\tau_t(\alpha(n))=\alpha(\sigma_t^{\nu}(n))$ so:
$$M_L=<(id\surl{\ _{\beta} \star_{\alpha}}_{\ \nu}
\omega\circ\tau_t)(\Gamma(m))\ | \ m\in M,\omega \in (M_R)^+_* ,
k\in {\mathbb R}^+ \ \text{s.t } \omega\circ\alpha \leq
k\nu>^{-\textsc{w}}$$ By \ref{egal}, we have
$\sigma_t^{\Psi}((id\surl{\ _{\beta} \star_{\alpha}}_{\ \nu}
\omega)\Gamma(m))=(id\surl{\ _{\beta} \star_{\alpha}}_{\ \nu}
\omega\circ\tau_t)\Gamma(\sigma_t^{\Psi}(m))$ that's why
$\sigma_t^{\Psi}(M_L)=M_L$ for all $t\in\mathbb R$. On the other
hand, by \ref{semi}, restriction of $\Psi$ to $M_L$ is
semi-finite. By Takesaki's theorem (\cite{St}, theorem 10.1),
there exists a unique normal and faithful conditional expectation
$E$ from $M$ to $M_L$ such that $\Psi(m)=\Psi(E(m))$ for all $m\in
M^+$. Moreover, if $P$ is the orthogonal projection on the
closure of $\Lambda_{\Psi}({\mathcal N}_{\Psi} \cap M_L)$ then
$E(m)P=PmP$.

So the range of $P$ contains $\Lambda_{\Psi}((id \surl{\ _{\beta}
\star_{\alpha}}_{\ \nu} \omega )\Gamma(x))$ for all $\omega$ and
$x\in {\mathcal N}_{\Psi}$. By right version of \ref{dense}
implies that $P=1$ so that $E$ is the identity and $M=M_L$. If we
apply the previous result to the opposite measured quantum
groupoid, then we get that $M=M_R$.
\end{proof}

\begin{coro}
There exists a unique strongly continuous one-parameter group
$\tau$ of automorphisms of $M$ such that, for all $t \in \mathbb
R$, $m\in M$ and $n\in N$:
$$\tau_t(m)=D^{-it}mD^{it},\ \tau_t(\alpha(n))=\alpha(\sigma_t^{\nu}(n)) \text{
and } \tau_t(\beta(n))=\beta(\sigma_t^{\nu}(n))$$
\end{coro}

\begin{proof}
Straightforward from the previous theorem and \ref{tau}. First
property comes from \ref{mieux} and the second one from
\ref{commuNa}.
\end{proof}

\begin{defi}
The one-parameter group $\tau$ is called {\bf scaling group}.
\end{defi}

\begin{lemm}
It is possible to define normal *-automorphisms $\tau_t \surl{\
_{\beta} \star_{\alpha}}_{\ N} \tau_t$ and $\tau_t \surl{\
_{\beta} \star_{\alpha}}_{\ N} \sigma_t^{\Phi}$ of $M \surl{\
_{\beta} \otimes_{\alpha}}_{\ N}M$ for all $t\in\mathbb{R}$.
\end{lemm}

\begin{proof}
The proof is very similar to \ref{tau2}.
\end{proof}

\begin{prop}
We have $\Gamma \circ \tau_t=(\tau_t \surl{\ _{\beta}
\star_{\alpha}}_{\ N}\tau_t)\circ\Gamma$ for all $t\in \mathbb R$.
\end{prop}

\begin{proof}
By \ref{egal} and co-product relation, we have for all
$t\in\mathbb R$:
$$\begin{aligned}
(id \surl{\ _{\beta} \star_{\alpha}}_{\ \nu}
\Gamma)(\sigma_t^{\Psi} \surl{\ _{\beta} \star_{\alpha}}_{\ \nu}
\tau_{-t})\circ\Gamma &=(id \surl{\ _{\beta} \star_{\alpha}}_{\
\nu} \Gamma)\Gamma \circ \sigma_t^{\Psi}\\
&=(\Gamma \surl{\ _{\beta} \star_{\alpha}}_{\ \nu} id)\Gamma
\circ \sigma_t^{\Psi}=(\Gamma \circ\sigma_t^{\Psi}\surl{\
_{\beta} \star_{\alpha}}_{\
\nu} \tau_{-t})\Gamma\\
&=(\sigma_t^{\Psi}\surl{\ _{\beta} \star_{\alpha}}_{\ \nu}
\tau_{-t}\surl{\ _{\beta} \star_{\alpha}}_{\ \nu}
\tau_{-t})(\Gamma\surl{\ _{\beta} \star_{\alpha}}_{\ \nu}
id)\Gamma\\
&=(\sigma_t^{\Psi}\surl{\ _{\beta}\star_{\alpha}}_{\ \nu}
[(\tau_{-t}\surl{\ _{\beta}\star_{\alpha}}_{\ \nu}
\tau_{-t})\circ\Gamma])\circ\Gamma
\end{aligned}$$
Consequently, for all $m\in M$, $\omega\in M_*^+$,
$k\in\mathbb{R}^+$ such that $\omega\circ\beta\leq k\nu$, we have:

$$
\begin{aligned}
\Gamma\circ\tau_{-t}\circ((\omega\circ\sigma_t^{\Psi})\surl{\
_{\beta} \star_{\alpha}}_{\ \nu} id)\Gamma&=(\omega\surl{\
_{\beta} \star_{\alpha}}_{\ \nu} id\surl{\ _{\beta}
\star_{\alpha}}_{\ \nu} id)(\sigma_t^{\Psi}\!\!\surl{\ _{\beta}
\star_{\alpha}}_{\ \nu}\!(\Gamma\circ\tau_{-t}))\circ\Gamma \\
&=(\omega\surl{\ _{\beta}\star_{\alpha}}_{\ \nu} id\surl{\
_{\beta} \star_{\alpha}}_{\ \nu} id)(\sigma_t^{\Psi}\!\!\surl{\
_{\beta} \star_{\alpha}}_{\ \nu}\![(\tau_{-t}\surl{\
_{\beta}\star_{\alpha}}_{\ \nu}\tau_{-t})\circ\Gamma])\\
&=[(\tau_{-t}\surl{\ _{\beta}\star_{\alpha}}_{\ \nu}
\tau_{-t})\circ\Gamma]\circ((\omega\circ\sigma_t^{\Psi})\surl{\
_{\beta} \star_{\alpha}}_{\ \nu} id)\Gamma
\end{aligned}$$

The theorem \ref{densevn} allows us to conclude.
\end{proof}

\begin{prop}\label{clef2}
For all $x\in M\cap\alpha(N)'$, we have $\Gamma(x)=1 \!\surl{\
_{\beta} \otimes_{\alpha}}_{\ N} x \Leftrightarrow x\in\beta(N)$.
Also, for all $x\in M\cap\beta(N)'$, we have $\Gamma(x)=x\!
\surl{\ _{\beta} \otimes_{\alpha}}_{\ N} 1 \Leftrightarrow
x\in\alpha(N)$.
\end{prop}

\begin{proof}
Let $x\in M\cap\alpha(N)'$ such that $\Gamma(x)=1\surl{\ _{\beta}
\otimes_{\alpha}}_{\ N} x$. For all $n \in \mathbb N$, we define
in the strong topology:
$$x_n=\frac{n}{\sqrt{\pi}}\int\!
exp(-n^2t^2)\sigma_t^{\Psi}(x)\ dt \quad \text{analytic w.r.t }
\sigma^{\Psi},$$ and:
$$y_n=\frac{n}{\sqrt{\pi}}\int\!
exp(-n^2t^2)\tau_{-t}(x)\ dt \quad \text{belongs to }
\alpha(N)'.$$ By \ref{egal}, we have $\Gamma(x_n)=1\!\surl{\
_{\beta} \otimes_{\alpha}}_{\ N}y_n$. If $d \in ({\mathcal
M}_{\Psi} \cap {\mathcal M}_{T_R})^+$, then, for all $n\in\mathbb
N$, we have $dx_n \in {\mathcal M}_{\Psi} \cap {\mathcal
M}_{T_R}$. Let $\omega \in M^+_*$ and $k\in {\mathbb R}^+$ such
that $\omega\circ\alpha\leq k\nu$. By right invariance, we get:

$$
\begin{aligned}
\omega\circ T_R(dx_n)&=\omega((\Psi \surl{\ _{\beta}
\star_{\alpha}}_{\ \nu} id)(\Gamma(dx_n)))\\
&=\Psi((id \surl{\ _{\beta} \star_{\alpha}}_{\ \nu}
\omega)(\Gamma(dx_n)))=\Psi((id \surl{\ _{\beta}
\star_{\alpha}}_{\ \nu} (y_n\omega))(\Gamma(d)))\\
&=\omega((\Psi \surl{\ _{\beta} \star_{\alpha}}_{\ \nu}
id)(\Gamma(d))y_n)=\omega(T_R(d)y_n)
\end{aligned}$$

Take the limit over $n\in\mathbb{N}$ to obtain $T_R(dx)=T_R(d)x$
for all $d\in {\mathcal M}_{\Psi} \cap {\mathcal M}_{T_R}$ and, by
semi-finiteness of $T_R$, we conclude that $x$ belongs to
$\beta(N)$. Reverse inclusion comes from axioms. If we apply this
result to the opposite measured quantum groupoid, then we get the
second point.
\end{proof}

\subsection{The antipode and its polar decomposition}
We now approach definition of the antipode.

\begin{lemm}\label{ense}
We have $(\omega_{v,w} * id)(U'_{H_{\Phi}})D^{\lambda} \subset
D^{\lambda}(\omega_{\Delta_{\Psi}^{-\lambda}v,\Delta_{\Psi}^{\lambda}w}
* id)(U'_{H_{\Phi}})$ for all $\lambda \in \mathbb C$ and $v,w\in
\Lambda_{\Psi}({\mathcal T}_{\Psi,T_R})$.
\end{lemm}

\begin{proof}
Straightforward from relation \eqref{fond}.
\end{proof}

\begin{prop}\label{proi}
If $I$ is the unitary part of the polar decomposition of $G$,
then, for all $v,w\in D((H_{\Psi})_{\beta},\nu^o)$, we have:
$$I(\omega_{J_{\Psi}w,v}*id)(U_{H_{\Phi}}'{}^{\hspace{-.3cm}*}\hspace{.1cm})I
=(\omega_{J_{\Psi}v,w}*id)(U'_{H_{\Phi}})$$
\end{prop}

\begin{proof}
We have $(\omega_{v,w} * id)(U'_{H_{\Phi}})D^{1/2}\subseteq
D^{1/2}(\omega_{\Delta_{\Psi}^{-1/2}v,\Delta_{\Psi}^{1/2}w}\star
id)(U'_{H_{\Phi}})$ for all $v,w\in\Lambda_{\Psi}({\mathcal
T}_{\Psi,T_R})$ by the previous lemma. On the other hand, by
inclusion \eqref{inclu2}, we have:
$$
(\omega_{v,w} * id)(U'_{H_{\Phi}})D^{1/2}= (\omega_{v,w} \star
id)(U'_{H_{\Phi}})G^*I\subseteq D^{1/2}I(\omega_{w,v} *
id)(U'_{H_{\Phi}})I$$ So $I(\omega_{w,v}*id)(U'_{H_{\Phi}})I=
(\omega_{\Delta_{\Psi}^{-1/2}v,\Delta_{\Psi}^{1/2}w} *
id)(U'_{H_{\Phi}})$ and, by \ref{switch}, we have:
$$I(\omega_{w,v}*id)(U_{H_{\Phi}}'{}^{\hspace{-.3cm}*}\hspace{.1cm})I
=(\omega_{\Delta_{\Psi}^{1/2}w,\Delta_{\Psi}^{-1/2}v} * id)(U_{H_{\Phi}}'{}^{\hspace{-.3cm}*}\,) \\
=(\omega_{J_{\Psi}v,J_{\Psi}w}* id)(U'_{H_{\Phi}})$$
\end{proof}

\begin{coro}
There exists a *-anti-automorphism $R$ of $M$ defined by
$R(m)=Im^*I$ such that $R^2=id$. (We recall that $I$ denotes the
unitary part of the polar decomposition of $G$).
\end{coro}

\begin{proof}
Straightforward from the previous proposition and theorem
\ref{densevn}.
\end{proof}

\begin{defi}
The unique *-anti-automorphism $R$ of $M$ such that $R(m)=Im^*I$,
where $I$ denotes the unitary part of the polar decomposition of
$G$, is called \textbf{unitary antipode}.
\end{defi}

\begin{defi}\label{antipode}
The application $S=R\tau_{-i/2}$ is called {\bf antipode}.
\end{defi}

The next proposition states elementary properties of the
antipode. Straightforward proofs are omitted.

\begin{prop}\label{numero}The antipode $S$ satisfies:
\begin{center}
\begin{minipage}{11cm}
\begin{enumerate}[i)]
\item for all $t\in\mathbb R$, we have $\tau_t\circ R=R\circ\tau_t$ and $\tau_t\circ S=S\circ\tau_t$
\item $SR=RS$ and $S^2=\tau_{-i}$
\item $S$ is densely defined and has dense range
\item $S$ is injective and $S^{-1}=R\tau_{i/2}=\tau_{i/2}R$
\item for all $x\in {\mathcal D}(S)$, $S(x^*)\in {\mathcal D}(S)$
and $S(S(x)^*)^*=x$ \label{test}
\end{enumerate}
\end{minipage}
\end{center}
\end{prop}

\subsection{Characterization of the antipode}

In \ref{antipode}, we define the antipode by giving its polar
decomposition. However, we have to verify that $S$ is what it
should be.

\subsubsection{Usual characterization of the antipode.}

\begin{prop}\label{espoir}
For all $v,w\in \Lambda_{\Psi}({\mathcal T}_{\Psi,T_R})$,
$(\omega_{w,v} * id)(U'_{H_{\Phi}})$ belongs to ${\mathcal D}(S)$
and we have: $$S((\omega_{w,v}*id)(U'_{H_{\Phi}}))=(\omega_{w,v} *
id)(U_{H_{\Phi}}'{}^{\hspace{-.3cm}*}\hspace{.1cm})$$ Moreover,
the linear span of $(\omega_{v,w} * id)(U'_{H_{\Phi}})$, where
$v,w\in \Lambda_{\Psi}({\mathcal T}_{\Psi,T_R})$, is a core for
$S$.
\end{prop}

\begin{proof}
By \ref{ense}, we have $(\omega_{w,v}\star id)(U'_{H_{\Phi}}) \in
{\mathcal D}(\tau_{-i/2})={\mathcal D}(S)$ and:
$$
\begin{aligned}
S((\omega_{w,v} * id)(U'_{H_{\Phi}})) &=
R((\omega_{\Delta_{\Psi}^{-1/2}w,\Delta_{\Psi}^{1/2}v}*id)(U'_{H_{\Phi}})) \\
&= (\omega_{S_{\Psi}v,\Delta_{\Psi}S_{\Psi}w}*
id)(U'_{H_{\Phi}}) &&\text{by proposition \ref{proi},} \\
&= (\omega_{w,v} * id)(U_{H_{\Phi}}'{}^{\hspace{-.3cm}*}\,)
&&\text{by lemma \ref{switch}}.
\end{aligned}$$

The involved subspace of $M$ is included in ${\mathcal
D}(\tau_{-i/2})$ by \ref{ense}, weakly dense in $M$ by theorem
\ref{densevn} and $\tau$-invariant by \ref{tau} which finishes the
proof.
\end{proof}

\begin{coro}\label{gammas}
For $a,b,c,d\in {\mathcal T}_{\Psi,T_R}$,
$(\omega_{\Lambda_{\Psi}(a),\Lambda_{\Psi}(b)} \surl{\ _{\beta}
\star_{\alpha}}_{\ \nu} id)(\Gamma(cd))$ belongs to ${\mathcal
D}(S)$ and we have:
$$S((\omega_{\Lambda_{\Psi}(a),\Lambda_{\Psi}(b)} \surl{\ _{\beta}
\star_{\alpha}}_{\
\nu}id)(\Gamma(cd)))=(\omega_{\Lambda_{\Psi}(c),\Lambda_{\Psi}(\sigma_{-i}^{\Psi}(d^*))}
\surl{\ _{\beta} \star_{\alpha}}_{\ \nu}
id)(\Gamma(\sigma_{i}^{\Psi}(a)b^*))$$
\end{coro}

\begin{proof}
By \ref{corres}, we know that:
$$(\omega_{\Lambda_{\Psi}(a),\Lambda_{\Psi}(b)}
\surl{\ _{\beta} \star_{\alpha}}_{\ \nu}
id)(\Gamma(cd))=(\omega_{\Lambda_{\Psi}(cd),\Lambda_{\Psi}
(b\sigma_{-i}^{\Psi}(a^*))}*id)(U'_{H_{\Phi}})$$ which belongs to
${\mathcal D}(S)$. Then, by \ref{corres} and \ref{switch}, we
have:
$$
\begin{aligned}
S((\omega_{\Lambda_{\Psi}(a),\Lambda_{\Psi}(b)} \surl{\ _{\beta}
\star_{\alpha}}_{\ \nu} id)(\Gamma(cd)))&=
S((\omega_{\Lambda_{\Psi}(cd),\Lambda_{\Psi}
(b\sigma_{-i}^{\Psi}(a^*))}*id)(W'))\\
&=(\omega_{\Lambda_{\Psi}(cd),\Lambda_{\Psi}
(b\sigma_{-i}^{\Psi}(a^*))}*id)(W'^*)\\
&=(\omega_{\Lambda_{\Psi}(\sigma_{i}^{\Psi}(a)b^*),\Lambda_{\Psi}
(\sigma_{-i}^{\Psi}(d^*c^*))}*id)(W')\\
&=(\omega_{\Lambda_{\Psi}(c),
\Lambda_{\Psi}(\sigma_{-i}^{\Psi}(d^*))} \surl{\ _{\beta}
\star_{\alpha}}_{\ \nu} id)(\Gamma(\sigma_{i}^{\Psi}(a)b^*))
\end{aligned}$$
\end{proof}

\subsubsection{The co-involution $R$}

In this section, we give a new expression of $R$ and we show that
it is a co-involution of the measured quantum groupoid.

\begin{prop}\label{defR}
For all $a,b\in {\mathcal N}_{\Psi}\cap {\mathcal N}_{T_R}$, we
have:
$$R((\omega_{J_{\Psi}\Lambda_{\Psi}(a)}\surl{\ _{\beta}\star_{\alpha}}_{\ \
\nu}id)(\Gamma(b^*b)))=(\omega_{J_{\Psi}\Lambda_{\Psi}(b)}\surl{\
_{\beta} \star_{\alpha}}_{\ \ \nu}id)(\Gamma(a^*a))$$
\end{prop}

\begin{proof}
The proposition comes from the following computation:
$$
\begin{aligned}
&\ \quad
R((\omega_{J_{\Psi}\Lambda_{\Psi}(a),J_{\Psi}\Lambda_{\Psi}(a)}\surl{\
_{\beta} \star_{\alpha}}_{\ \nu}id)(\Gamma(b^*b)))&\\
&=R((\omega_{\Lambda_{\Psi}(b^*b),J_{\Psi}\Lambda_{\Psi}(a^*a)}*id)(U'_{H_{\Phi}})
&\text{by corollary \ref{corres},}\\
&=(\omega_{\Lambda_{\Psi}(a^*a),J_{\Psi}\Lambda_{\Psi}(b^*b)}*id)(U'_{H_{\Phi}})
&\text{by definition of } R,\\
&=(\omega_{J_{\Psi}\Lambda_{\Psi}(b),J_{\Psi}\Lambda_{\Psi}(b)}\surl{\
_{\beta} \star_{\alpha}}_{\ \nu}id)(\Gamma(a^*a))
&\text{by corollary \ref{corres}.}\\
\end{aligned}$$

\end{proof}

\begin{rema}
We notice that $R$ is $T_L$-independent.
\end{rema}

\begin{prop}
We have $I\alpha(n^*)=\beta(n)I$ for all $n\in N$ and
$R\circ\alpha=\beta$.
\end{prop}

\begin{proof}
By \ref{commuNa}, we have, for all $x \in {\mathcal
T}_{\Psi,T_R}$:
$$\beta(x)GD^{-1/2}\subseteq G\alpha(\sigma_{-i/2}((x^*))
\subseteq GD^{-1/2}\alpha(x^*)\subseteq I\alpha(x^*)$$ and, on
the other hand, $\beta(x)GD^{-1/2}\subseteq \beta(x)I$ so that
$I\alpha(x^*)=\beta(x)I$. The result holds by normality of
$\alpha$ and $\beta$.
\end{proof}

By \cite{S2}, there exists a unitary and anti-linear operator $I
\surl{\ _{\beta} \otimes_{\alpha}}_{\ \nu}I$ from $H\!\!\!
\surl{\ _{\beta} \otimes_{\alpha}}_{\ \nu}\!\! H$ onto
$H\!\!\!\surl{\ _{\alpha} \otimes_{\beta}}_{\ \ \nu^o}\!\! H$,
the adjoint of which is $I\surl{\ _{\alpha} \otimes_{\beta}}_{\ \
\nu^o}I$. Also, there exists an anti-isomorphism $R\surl{\
_{\beta} \star_{\alpha}}_{\ N} R$ from $M\surl{\ _{\beta}
\star_{\alpha}}_{\ N}M$ onto $M\surl{\ _{\alpha} \star_{\beta}}_{\
N^o}M$ and, by definition of $R$, we have, for all $X\in M\surl{\
_{\beta} \star_{\alpha}}_{\ N}M$:
$$(R\surl{\ _{\beta} \star_{\alpha}}_{\ N} R)(X)=(I \surl{\
_{\beta} \otimes_{\alpha}}_{\ \nu}I)X^*(I\surl{\ _{\alpha}
\otimes_{\beta}}_{\ \ \nu^o} I)$$ We underline the fact that, if
$\omega\in M_*^+$, then $\omega\circ R \in M_*^+$ and, if there
exists $k\in\mathbb{R}^+$ such that $\omega\circ\alpha\leq k\nu$,
then $\omega\circ R\circ\beta\leq k\nu$. Also, if $\theta\in
M_*^+$ and $k'\in{\mathbb R}^+$ are such that
$\theta\circ\beta\leq k'\nu$, then $\theta\circ R\circ\alpha\leq
k\nu$. Then, we have $\omega R\surl{\ _{\beta} \star_{\alpha}}_{\
\nu} \theta R=(\omega \surl{\ _{\alpha} \star_{\beta}}_{\ \
\nu^o}\theta)\circ(R\surl{\ _{\beta} \star_{\alpha}}_{\ \nu} R)$.

\begin{lemm}
For all $a,x\in {\mathcal N}_{T_R}\cap {\mathcal N}_{\Psi}$,
$\omega\in M^+_*$ and $k\in\mathbb{R}^+$ such that
$\omega\circ\alpha\leq k\nu$, we have: $$\omega\circ
R((\omega_{J_{\Psi}\Lambda_{\Psi}(a)}\surl{\ _{\beta}
\star_{\alpha}}_{\ \nu}id)(\Gamma(x)))=(\Lambda_{\Psi}((id\surl{\
_{\beta} \star_{\alpha}}_{\
\nu}\omega)(\Gamma(a^*a)))|J_{\Psi}\Lambda_{\Psi}(x))$$
\end{lemm}

\begin{proof}
Let $b\in {\mathcal N}_{T_R}\cap {\mathcal N}_{\Psi}$. By
\ref{defR}, we can compute:
$$\begin{aligned}
\omega\circ R((\omega_{J_{\Psi}\Lambda_{\Psi}(a)}\surl{\ _{\beta}
\star_{\alpha}}_{\
\nu}id)(\Gamma(b^*b)))&=\omega((\omega_{J_{\Psi}\Lambda_{\Psi}(b)}\surl{\
_{\beta} \star_{\alpha}}_{\ \nu}id)(\Gamma(a^*a)))\\
&=((id\surl{\ _{\beta} \star_{\alpha}}_{\
\nu}\omega)(\Gamma(a^*a))J_{\Psi}\Lambda_{\Psi}(b)|J_{\Psi}\Lambda_{\Psi}(b))\\
&=(J_{\Psi}bJ_{\Psi}\Lambda_{\Psi}((id\!\!\surl{\ _{\beta}
\star_{\alpha}}_{\ \nu}\omega)(\Gamma(a^*a)))|J_{\Psi}\Lambda_{\Psi}(b))\\
&=(\Lambda_{\Psi}((id\surl{\ _{\beta} \star_{\alpha}}_{\
\nu}\omega)(\Gamma(a^*a)))|J_{\Psi}\Lambda_{\Psi}(b^*b))
\end{aligned}$$

Linearity and normality of the expressions imply the lemma.
\end{proof}

\begin{prop}
We have $\varsigma_{N^o}\circ(R\surl{\ _{\beta} \star_{\alpha}}_{\
N} R)\circ\Gamma=\Gamma \circ R$.
\end{prop}

\begin{proof}
Let $a,b\in {\mathcal N}_{T_R}\cap {\mathcal N}_{\Psi}$,
$\omega,\theta\in M^+_*$ and $k,k'\in\mathbb{R}^+$ such that
$\omega\circ\alpha\leq k\nu$ and $\theta\circ\beta\leq k'\nu$.
Then, we can compute by \ref{defR} and the previous lemma:
$$
\begin{aligned}
&\ \quad(\theta\surl{\ _{\beta} \star_{\alpha}}_{\
\nu}\omega)(\Gamma\circ
R((\omega_{J_{\Psi}\Lambda_{\Psi}(a)}\surl{\ _{\beta}
\star_{\alpha}}_{\ \nu}id)(\Gamma(b^*b))))\\
&=(\theta\surl{\ _{\beta} \star_{\alpha}}_{\
\nu}\omega)(\Gamma((\omega_{J_{\Psi}\Lambda_{\Psi}(b)}\surl{\
_{\beta}\star_{\alpha}}_{\ \nu}id)(\Gamma(a^*a))))\\
&=(\omega_{J_{\Psi}\Lambda_{\Psi}(b)}\surl{\
_{\beta}\star_{\alpha}}_{\ \nu}\theta\surl{\ _{\beta}
\star_{\alpha}}_{\ \nu}\omega)(id\surl{\
_{\beta}\star_{\alpha}}_{\ N}\Gamma)(\Gamma(a^*a))\\
&=(\omega_{J_{\Psi}\Lambda_{\Psi}(b)}\surl{\
_{\beta}\star_{\alpha}}_{\ \nu}\theta\surl{\ _{\beta}
\star_{\alpha}}_{\ \nu}\omega)(\Gamma\surl{\
_{\beta}\star_{\alpha}}_{\ N}id)(\Gamma(a^*a))\\
&=(\omega_{J_{\Psi}\Lambda_{\Psi}(b)}\surl{\
_{\beta}\star_{\alpha}}_{\ \nu}\theta)[\Gamma((id\surl{\ _{\beta}
\star_{\alpha}}_{\ \nu}\omega)(\Gamma(a^*a)))]\\
&=(\Lambda_{\Psi}((id\surl{\ _{\beta} \star_{\alpha}}_{\
\nu}\theta\circ
R)(\Gamma(b^*b)))|J_{\Psi}\Lambda_{\Psi}((id\surl{\ _{\beta}
\star_{\alpha}}_{\ \nu}\omega)(\Gamma(a^*a))))
\end{aligned}$$
Observe the symmetry of the last expression and use it to proceed
towards the computation:
$$
\begin{aligned}
&\ \quad (\Lambda_{\Psi}((id\surl{\ _{\beta} \star_{\alpha}}_{\
\nu}\omega)(\Gamma(a^*a)))|J_{\Psi}\Lambda_{\Psi}((id\surl{\
_{\beta} \star_{\alpha}}_{\ \nu}\theta\circ
R)(\Gamma(b^*b))))\\
&=(\omega_{J_{\Psi}\Lambda_{\Psi}(a)}\surl{\
_{\beta}\star_{\alpha}}_{\ \nu}\omega\circ R)[\Gamma((id\surl{\
_{\beta} \star_{\alpha}}_{\ \nu}\theta\circ R)(\Gamma(b^*b)))]\\
&=(\omega_{J_{\Psi}\Lambda_{\Psi}(a)}\surl{\
_{\beta}\star_{\alpha}}_{\ \nu}\omega\circ R\surl{\ _{\beta}
\star_{\alpha}}_{\ \nu}\theta\circ R)(\Gamma\surl{\
_{\beta}\star_{\alpha}}_{\ N}id)(\Gamma(b^*b))\\
&=(\omega_{J_{\Psi}\Lambda_{\Psi}(a)}\surl{\
_{\beta}\star_{\alpha}}_{\ \nu}\omega\circ R\surl{\ _{\beta}
\star_{\alpha}}_{\ \nu}\theta\circ R)(id\surl{\
_{\beta}\star_{\alpha}}_{\ N}\Gamma)(\Gamma(b^*b))\\
&=(\omega\circ R\surl{\ _{\beta} \star_{\alpha}}_{\
\nu}\theta\circ
R)(\Gamma((\omega_{J_{\Psi}\Lambda_{\Psi}(a)}\surl{\ _{\beta}
\star_{\alpha}}_{\ \nu}id)(\Gamma(b^*b))))\\
&=(\omega\surl{\ _{\alpha} \star_{\beta}}_{\
\nu^o}\theta)(R\surl{\ _{\beta} \star_{\alpha}}_{\
N}R)(\Gamma((\omega_{J_{\Psi}\Lambda_{\Psi}(a)}\surl{\ _{\beta}
\star_{\alpha}}_{\ \nu}id)(\Gamma(b^*b))))\\
&=(\theta\surl{\ _{\beta} \star_{\alpha}}_{\
\nu}\omega)\varsigma_{N^o}(R\surl{\ _{\beta} \star_{\alpha}}_{\
N}R)(\Gamma((\omega_{J_{\Psi}\Lambda_{\Psi}(a)}\surl{\ _{\beta}
\star_{\alpha}}_{\ \nu}id)(\Gamma(b^*b))))
\end{aligned}$$

Theorem \ref{densevn} easily implies the result.
\end{proof}

\subsubsection{Left strong invariance w.r.t the antipode.}

In this section, $T'$ denotes a left invariant n.s.f weight from
$M$ to $\alpha(N)$. We put $\Phi'=\nu\circ\alpha^{-1}\circ T'$,
$J_{\Phi'}$ the anti-linear operator and $\Delta_{\Phi'}$ the
modular operator which come from Tomita's theory of $\Phi'$,
$\sigma^{\Phi'}$ its modular group and $V=(U_{T'})_{H_{\Phi}}^*$
i.e the fundamental unitary associated with $T'$. The next
proposition is the left strong invariance w.r.t $S$.

\begin{prop}\label{invafort}
Elements $(id*\omega_{v,w})(V)$ belong to the domain of $S$ for
all $v,w\in\Lambda_{\Phi'}({\mathcal T}_{\Phi',T'})$ and we have
$S((id*\omega_{v,w})(V))=(id*\omega_{v,w})(V^*)$.
\end{prop}

\begin{proof}
By \ref{corres}, we have $(id*\omega)(V)=(\omega\circ
R*id)(U'_{H_{\Phi}})$ for all $\omega$. If
$\overline{\omega}(x)=\overline{\omega(x^*)}$, then, by
\ref{espoir}, we have:
$$\begin{aligned}
S((id*\omega)(V))=S((\omega\circ R*id)(U'_{H_{\Phi}}))
&=(\omega\circ R*id)(U_{H_{\Phi}}'{}^{\hspace{-.3cm}*}\,)\\
&=[(\overline{\omega}\circ R*id)(U'_{H_{\Phi}})]^*\\
&=[(id*\overline{\omega})(V)]^*=(id*\omega)(V^*)
\end{aligned}$$

\end{proof}

\begin{lemm}
For all $v\in {\mathcal D}(D^{1/2})$ and $w\in {\mathcal
D}(D^{1/2})$, we have:
$$(\omega_{v,w}*id)(V)^*=(\omega_{ID^{-1/2}v,ID^{1/2}w}*id)(V)$$
\end{lemm}

\begin{proof}
We have $(id*\omega_{w',v'})(V)\in {\mathcal D}(S)={\mathcal
D}(\tau_{-i/2})$ for all $v',w'$ belonging to
$\Lambda_{\Phi'}({\mathcal T}_{\Phi',T'})$ by \ref{invafort} and,
since $\tau$ is implemented by $D^{-1}$, we have:

$$
\begin{aligned}
(id*\omega_{w',v'})(V)D^{1/2} &\subseteq
D^{1/2}\tau_{-i/2}((id*\omega_{w',v'})(V)) \\
&=D^{1/2}R(S((id*\omega_{w',v'})(V))) \\
&=D^{1/2}I[(id*\omega_{w',v'})(V^*)]^*I \\
&=D^{1/2}I(id*\omega_{v',w'})(V)I.
\end{aligned}$$

Then, for all $v\in {\mathcal D}(D^{1/2})$ and $w\in {\mathcal
D}(D^{1/2})$, we have:

$$
\begin{aligned}
((\omega_{ID^{-1/2}v,ID^{1/2}w} * id)(V)w'|v') &=
((id*\omega_{w',v'})(V)
D^{1/2}Iv|D^{-1/2}Iw) \\
&=(D^{1/2}I(id*\omega_{v',w'})(V)v|D^{-1/2}Iw) \\
&=(w|(id*\omega_{v',w'})v) \\
&=((\omega_{v,w} * id)(V)^*w',v')
\end{aligned}$$

Then, the proposition holds.
\end{proof}

\begin{prop}\label{clef}
The following relations are satisfied:
\begin{center}
\begin{minipage}{10cm}
\begin{enumerate}[i)]
\item $(I \surl{\ _{\alpha} \otimes_{\epsilon}}_{\ \ N^o} J_{\Phi'})V
=V^*(I \surl{\ _{\beta} \otimes_{\alpha}}_{\ N} J_{\Phi'})$;
\item $(D^{-1}\surl{\ _{\alpha} \otimes_{\epsilon}}_{\ \ \nu^o} \Delta_{\Phi'})V
=V(D^{-1}\surl{\ _{\beta} \otimes_{\alpha}}_{\ \nu}
\Delta_{\Phi'})$;
\item $(\tau_t \surl{\ _{\beta} \star_{\alpha}}_{\ N}
\sigma^{\Phi'}_t)\circ \Gamma=\Gamma\circ\sigma^{\Phi'}_t$ for all
$t\in\mathbb R$ .
\end{enumerate}
\end{minipage}
\end{center}
where $\epsilon(n)=J_{\Phi'}\alpha(n^*)J_{\Phi'}$ for all $n\in
N$.
\end{prop}

\begin{proof}
We denote by $S_{\Phi'}$ the operator of Tomita's theory
associated with $\Phi'$ and defined as the closed operator on
$H_{\Phi'}$ such that $\Lambda_{\Phi'}({\mathcal N}_{\Phi'} \cap
{\mathcal N}_{\Phi'}^*)$ is a core for $S_{\Phi'}$ and
$S_{\Phi'}\Lambda_{\Phi'}(x)=\Lambda_{\Phi'}(x^*)$ for all $x\in
{\mathcal N}_{\Phi'} \cap {\mathcal N}_{\Phi'}^*$. Then, by
definition, we have $\Delta_{\Phi'}=S_{\Phi'}^*S_{\Phi'}$ and
$S_{\Phi'}=J_{\Phi'}\Delta_{\Phi'}^{1/2}$. Moreover, for all $m\in
M$ and $t\in\mathbb R$, we have
$\sigma^{\Phi'}_t(m)=\Delta_{\Phi'}^{it}m\Delta_{\Phi'}^{-it}$.

First of all, we verify these relations make sense. We have to
prove some commutation relations (\ref{legitime}, \ref{legitime2}
and \cite{S3}). We can write for all $n\in {\mathcal T}_{\nu}$ and
$y\in {\mathcal N}_{\Phi'} \cap {\mathcal N}_{\Phi'}^*$:
$$S_{\Phi'}\alpha(n)\Lambda_{\Phi'}(y)=S_{\Phi'}\Lambda_{\Phi'}(\alpha(n)y)$$
$$=\Lambda_{\Phi'}(y^*\alpha(n^*))=
\hat{\alpha}(\sigma_{-i/2}^{\nu}(n))S_{\Phi'}\Lambda_{\Phi'}(y)$$
so $\hat{\alpha}(\sigma_{-i/2}^{\nu}(n))S_{\Phi'} \subseteq
S_{\Phi'}\alpha(n)$ and by adjoint $\alpha(n)S_{\Phi'}^* \subseteq
S_{\Phi'}^*\hat{\alpha}(\sigma_{i/2}^{\nu}(n)$. Then:
$$\alpha(n)\Delta_{\Phi'}=\alpha(n)S_{\Phi'}^*S_{\Phi'} \subseteq
S_{\Phi'}^*\hat{\alpha}(\sigma_{i/2}^{\nu}(n)S_{\Phi'} \subseteq
\Delta_{\Phi'}\alpha(\sigma_i^{\nu}(n))$$ Since $\beta(n)D^{-1}
\subseteq D^{-1}\beta(\sigma_i^{\nu}(n))$, the second relation
makes sense. On an other hand, we know that
$I\beta(n)=\alpha(n^*)I$ and ${\mathcal
J}\alpha(n)=\epsilon(n^*)J_{\Phi'}$ to terms of the first
relation. Finally, for all $t\in\mathbb R$, we have:
$$\tau_t \circ\beta=\beta\circ\sigma_t^{\nu} \quad\text{ and }
\sigma^{\Phi'}_t(\alpha(n))=\Delta_{\Phi'}^{it}\alpha(n)\Delta_{\Phi'}^{-it}=
\alpha(\sigma_t^{\nu}(n))$$ which finishes verifications.

Let $v,w \in \Lambda_{\Phi}({\mathcal T}_{\Phi,S_L})$. By
\ref{prem}, we know that $(\omega_{v,w} \surl{\ _{\beta}
\star_{\alpha}}_{\ \nu} id)(\Gamma(y))$ belongs to ${\mathcal
N}_{T'}\cap {\mathcal N}_{\Phi'} \cap {\mathcal N}_{T'}^*\cap
{\mathcal N}_{\Phi'}^*$ for all $y \in {\mathcal N}_{T'}\cap
{\mathcal N}_{\Phi'} \cap {\mathcal N}_{T'}^*\cap {\mathcal
N}_{\Phi'}^*$. By \ref{lienGV}, we can write $(\omega_{v,w} \star
id)(V^*)\Lambda_{\Phi'}(y)=\Lambda_{\Phi'}((\omega_{v,w} \surl{\
_{\beta} \star_{\alpha}}_{\ \nu} id)(\Gamma(y)))$ so that
$(\omega_{v,w} * id)(V^*)\Lambda_{\Phi'}(y)$ belongs to
${\mathcal D}(S_{\Phi'})$. Then, we compute:
$$
\begin{aligned}
S_{\Phi'}(\omega_{v,w}*id)(V^*)\Lambda_{\Phi'}(y)&=
S_{\Phi'}\Lambda_{\Phi'}((\omega_{v,w} \surl{\ _{\beta}
\star_{\alpha}}_{\ \nu}
id)(\Gamma(y)))\\
&=\Lambda_{\Phi'}((\omega_{w,v} \surl{\ _{\beta}
\star_{\alpha}}_{\ \nu}
id)(\Gamma(y^*)))\\
&=(\omega_{w,v}*id)(V^*)\Lambda_{\Phi'}(y^*)\\
&=(\omega_{w,v}*id)(V^*)S_{\Phi'} \Lambda_{\Phi'}(y)
\end{aligned}$$ Since $\Lambda_{\Phi'}({\mathcal
N}_{T'}\cap {\mathcal N}_{\Phi'} \cap {\mathcal N}_{T'}^*\cap
{\mathcal N}_{\Phi'}^*)$ is a core for $S_{\Phi'}$, this implies:
\begin{equation}\label{parallele1}
(\omega_{w,v} * id)(V^*)S_{\Phi'} \subseteq S_{\Phi'}(\omega_{v,w}
* id)(V^*)
\end{equation}
Take adjoint so as to get:
\begin{equation}\label{parallele2}
(\omega_{w,v} * id)(V)S_{\Phi'}^* \subseteq
S_{\Phi'}^*(\omega_{v,w} * id)(V)
\end{equation}
Then, we deduce by the previous lemma:

$$
\begin{aligned}
(\omega_{v,w}* id)(V)\Delta_{\Phi'} &= (\omega_{v,w}*
id)(V)S_{\Phi'}^*S_{\Phi'}
\\
&\subseteq S_{\Phi'}^*(\omega_{v,w} * id)(V)S_{\Phi'} \\
&=S_{\Phi'}^*[(\omega_{ID^{-1/2}w,ID^{1/2}v} * id)(V)]^*S_{\Phi'}
\end{aligned}$$\\
Then by inclusion \eqref{parallele1} and the previous lemma, we
have:

$$
\begin{aligned}
(\omega_{v,w}* id)(V)\Delta_{\Phi'} &\subseteq
S_{\Phi'}^*S_{\Phi'}[(\omega_{ID^{1/2}v,ID^{-1/2}w} * id)(V)]^*\\
&=\Delta_{\Phi'} (\omega_{D^{1/2}IID^{1/2}v,D^{-1/2}IID^{-1/2}w}*id)(V)\\
&=\Delta_{\Phi'} (\omega_{Dv,N^{-1}w} * id)(V)
\end{aligned}$$
Consequently, like relation \eqref{fond}, we easily deduce that:
$$(D^{-1}\surl{\ _{\alpha} \otimes_{\epsilon}}_{\ \ \nu^o}
\Delta_{\Phi'})V =V(D^{-1}\surl{\ _{\beta} \otimes_{\alpha}}_{\
\nu} \Delta_{\Phi'})$$

Let's prove the first relation. By inclusion \eqref{parallele1},
for all $v\in {\mathcal D}(N^{-1/2})$ and $w\in {\mathcal
D}(D^{1/2})$, we have:

\begin{equation}\label{comenc}
\begin{aligned}
J_{\Phi'}(\omega_{w,v}* id)(V^*)J_{\Phi'}\Delta_{\Phi'}^{1/2}
&=J_{\Phi'}(\omega_{w,v}* id)(V^*)S_{\Phi'} \\
&\subseteq J_{\Phi'}S_{\Phi'}(\omega_{v,w}* id)(V^*)\\
&=\Delta_{\Phi'}^{1/2}(\omega_{v,w}* id)(V^*)
\end{aligned}
\end{equation}

For all $p,q\in {\mathcal D}(\Delta_{\Phi'}^{1/2})$, we have by
ii):

$$
\begin{aligned}
((\omega_{v,w}* id)(V^*)p,\Delta_{\Phi'}^{1/2}q) &= (V^*(v\surl{\
_{\alpha} \otimes_{\epsilon}}_{\ \ \nu^o} p)|w\surl{\ _{\beta}
\otimes_{\alpha}}_{\ \nu} \Delta_{\Phi'}^{1/2}q) \\
&=(V^*(v \surl{\ _{\alpha} \otimes_{\epsilon}}_{\ \
\nu^o}p)|D^{-1/2}(D^{1/2}w) \surl{\
_{\beta} \otimes_{\alpha}}_{\ \nu} \Delta_{\Phi'}^{1/2}q) \\
&=((D^{-1/2} \surl{\ _{\beta} \otimes_{\alpha}}_{\ \nu}
\Delta_{\Phi'}^{1/2})V^*(v\surl{\ _{\alpha}
\otimes_{\epsilon}}_{\ \ \nu^o} p)|D^{1/2}w\surl{\ _{\beta}
\otimes_{\alpha}}_{\ \nu} q) \\
&=(V^*(D^{-1/2}v\surl{\ _{\alpha} \otimes_{\epsilon}}_{\ \ \nu^o}
\Delta_{\Phi'}^{1/2}p)|D^{1/2}w\surl{\ _{\beta}
\otimes_{\alpha}}_{\ \nu}
q) \\
&=((\omega_{D^{-1/2}v,D^{1/2}w} *
id)(V^*)\Delta_{\Phi'}^{1/2}p|q).
\end{aligned}$$

Since $\Delta_{\Phi'}^{1/2}$ is self-adjoint, we get:
$$(\omega_{D^{-1/2}v,D^{1/2}w} * id)(V^*)\Delta_{\Phi'}^{1/2} \subseteq
\Delta_{\Phi'}^{1/2}(\omega_{v,w}* id)(V^*)$$ Also, by the
previous lemma, we have:
$$
\begin{aligned}
(\omega_{D^{-1/2}v,D^{1/2}w} *
id)(V^*)&=(\omega_{D^{1/2}w,D^{-1/2}v} * id) (V)^*\\
&=(\omega_{Iw,Iv}* id)(V)
\end{aligned}$$ That's why
$(\omega_{Iw,Iv}* id)(V)\Delta_{\Phi'}^{1/2}\subseteq
\Delta_{\Phi'}^{1/2}(\omega_{v,w}* id)(V^*)$. Since
$\Delta_{\Phi'}^{1/2}$ has dense range, this last inclusion and
\eqref{comenc} imply that:
$$(\omega_{Iw,Iv}*id)(V)=J_{\Phi'}(\omega_{v,w}*
id)(V^*)J_{\Phi'}$$ Then, we can compute:
$$
\begin{aligned}
&\quad((I\surl{\ _{\beta} \otimes_{\alpha}}_{\
\nu}J_{\Phi'})V^*(I \surl{\ _{\beta} \otimes_{\alpha}}_{\
\nu}J_{\Phi'})(v\surl{\ _{\beta} \otimes_{\alpha}}_{\ \nu}
q)|w\surl{\ _{\alpha}
\otimes_{\epsilon}}_{\ \ \nu^o} q)\\
&=(V(Iw \surl{\ _{\beta} \otimes_{\alpha}}_{\
\nu}J_{\Phi'}q)|Iv\surl{\ _{\alpha} \otimes_{\epsilon}}_{\ \
\nu^o}J_{\Phi'}p)\\
&=((\omega_{Iw,Iv}\star id)(V)J_{\Phi'}q|J_{\Phi'}p)
=(J_{\Phi'}(\omega_{w,v}\star id)(V^*)q|J_{\Phi'}p)\\
&=((\omega_{v,w}* id)(V)p|q)=(V(v\surl{\ _{\beta}
\otimes_{\alpha}}_{\ \nu} q)|w\surl{\ _{\alpha}
\otimes_{\epsilon}}_{\ \ \nu^o} q)
\end{aligned}$$
so that the first relation is proved. We end the proof by the
last equality. We know that $\Gamma$ is implemented by $V$,
$\sigma^{\Phi'}$ by $\Delta_{\Phi'}$ and $\tau$  by $D$ so that
the relation comes from $(D^{-1}\surl{\ _{\alpha}
\otimes_{\epsilon}}_{\ \nu} \Delta_{\Phi'})V =V(D^{-1}\surl{\
_{\beta} \otimes_{\alpha}}_{\ \nu} \Delta_{\Phi'})$ like
\ref{egal}.
\end{proof}

If we take $T'=T_L$ then $V=W^*$, $J_{\Phi'}=J_{\Phi}$ and
$\Delta_{\Phi'}=\Delta_{\Phi}$ so that we have the following
propositions:

\begin{prop}
For all $v,w\in \Lambda_{\Phi}({\mathcal T}_{\Phi,S_L})$,
$(id*\omega_{v,w})(W)$ belongs to ${\mathcal D}(S)$ and:
$$S((id*\omega_{v,w})(W))=(id*\omega_{v,w})(W^*)$$
\end{prop}

\begin{prop}
We have $(\omega_{v,w}*id)(W^*)^*=(\omega_{ID^{-1/2}v,ID^{1/2}w}
*id)(W^*)$ for all $v\in {\mathcal D}(D^{1/2})$ and $w\in
{\mathcal D}(D^{1/2})$.
\end{prop}

\begin{prop}\label{besoin}
The following relations are satisfied:
\begin{center}
\begin{minipage}{10cm}
\begin{enumerate}[i)]
\item $(I \surl{\ _{\alpha} \otimes_{\hat{\beta}}}_{\ \ N^o} J_{\Phi})W^*
=W(I \surl{\ _{\beta} \otimes_{\alpha}}_{\ N} J_{\Phi})$ ;
\item $(D^{-1}\surl{\ _{\beta} \otimes_{\alpha}}_{\ \nu}
\Delta_{\Phi})W^* =W^*(D^{-1}\surl{\ _{\beta}
\otimes_{\alpha}}_{\ \nu} \Delta_{\Phi})$ ;
\item $(\tau_t \surl{\ _{\beta} \star_{\alpha}}_{\ N}
\sigma^{\Phi}_t)\circ \Gamma=\Gamma\circ\sigma^{\Phi}_t$ for all
$t\in\mathbb R$.
\end{enumerate}
\end{minipage}
\end{center}
\end{prop}

We summarize the results of this section in the two following
theorems:

\begin{theo}
Let $(N,M,\alpha,\beta,\Gamma,\nu,T_L,T_R)$ be a measured quantum
groupoid and $W$ the pseudo-multiplicative unitary associated
with. Then the closed linear span of $(id*\omega_{v,w})(W)$ for
all $v\in D(_{\alpha}H_{\Phi},\nu)$ and $w\in
D((H_{\Phi})_{\hat{\beta}},\nu^o)$ is equal to the whole von
Neumann algebra $M$.
\end{theo}

\begin{theo}\label{invforte}
Let $(N,M,\alpha,\beta,\Gamma,\nu,T_L,T_R)$ be a measured quantum
groupoid and $W$ the pseudo-multiplicative associated with. If we
put $\Phi=\nu\circ\alpha^{-1}\circ T_L$, then there exists an
unbounded antipode $S$ which satisfies:
\begin{enumerate}[i)]
\item for all $x\in {\mathcal D}(S)$, $S(x)^*\in {\mathcal D}(S)$
and $S(S(x)^*)^*=x$
\item for all $v,w\in \Lambda_{\Phi}({\mathcal T}_{\Phi,S_L})$, $(id*\omega_{v,w})(W)$
belongs to ${\mathcal D}(S)$ and:
$$S((id*\omega_{v,w})(W))=(id*\omega_{v,w})(W^*)$$
\end{enumerate}
$S$ has the following polar decomposition $S=R\tau_{i/2}$, where
$R$ is a co-involution of $M$ satisfying $R^2=id$,
$R\circ\alpha=\beta$ and $\varsigma_{N^o}\circ(R\surl{\ _{\beta}
\star_{\alpha}}_{\ N}R)\circ\Gamma=\Gamma\circ R$, and where
$\tau$, the so-called scaling group, is a one-parameter group of
automorphisms such that
$\tau_t\circ\alpha=\alpha\circ\sigma_t^{\nu}$,
$\tau_t\circ\beta=\beta\circ\sigma_t^{\nu}$ satisfying
$\Gamma\circ\tau_t=(\tau_t\surl{\ _{\beta} \star_{\alpha}}_{\
N}\tau_t)\circ\Gamma$ for all $t\in\mathbb{R}$. $S,R$ and $\tau$
are independent of $T_L$ and of $T_R$.

Moreover, $R\circ T_L\circ R$ is a n.s.f operator-valued weight
which is right invariant and $\alpha$-adapted w.r.t $\nu$.
\end{theo}

\section{Uniqueness, modulus and scaling operator}

In this section, the quasi-invariant weight $\nu$ is fixed. We
establish uniqueness of invariant operator-valued weight which is
adapted w.r.t $\nu$ up to a strictly positive element affiliated
with the center of the basis. We construct a modulus and a scaling
operator which link the left invariant operator-valued weight
$T_L$ and the right invariant operator-valued weight $R\circ
T_L\circ R$. Their properties imply that the fundamental
pseudo-multiplicative unitary satisfies a condition similar to
Woronowicz's manageability. Also, we study conditions so that a
n.s.f weight $\nu'$ on $N$ is quasi-invariant.

\subsection{Commutation relations}

In this section, we denote by $T'$ an other n.s.f operator-valued
weight from $M$ to $\alpha(N)$ which is left invariant and
$\beta$-adapted w.r.t $\nu$. We put
$\Phi'=\nu\circ\alpha^{-1}\circ T'$.

\begin{lemm}\label{prep}
If we put $\kappa_t=\sigma_t^{\Phi'}\circ\tau_{-t}$ for all
$t\in\mathbb{R}$, then $\Psi$ is $\kappa$-invariant. Also, if we
put $\kappa'_t=\sigma_t^{\Psi}\circ\tau_t$ for all
$t\in\mathbb{R}$, then $\Phi$ is $\kappa'$-invariant.
\end{lemm}

\begin{proof}
We know $\kappa_t\circ\alpha=\alpha$ for all $t\in\mathbb{R}$ and
we have:

$$
\begin{aligned}
\Gamma\circ\kappa_t=\Gamma\circ\sigma_t^{\Phi'}\circ\tau_{-t}&=(\tau_t\surl{\
_{\beta} \star_{\alpha}}_{\
N}\sigma_t^{\Phi'})\Gamma\circ\tau_{-t}\\
&=(\tau_t\tau_{-t}\surl{\ _{\beta} \star_{\alpha}}_{\
N}\sigma_t^{\Phi'}\tau_{-t})\Gamma=(id\surl{\ _{\beta}
\star_{\alpha}}_{\ N}\kappa_t)\Gamma
\end{aligned}$$

By right invariance of $T_R$, we deduce, for all $a\in {\mathcal
M}_{T_R}^+$:

$$T_R\circ\kappa_t(a)=(\Psi\surl{\ _{\beta} \star_{\alpha}}_{\ \nu}id)\Gamma(\kappa_t(a))
=\kappa_t((\Psi\surl{\ _{\beta} \star_{\alpha}}_{\
\nu}id)\Gamma(a)) =\kappa_t\circ T_R(a)$$

Since $T'$ is $\beta$-adapted w.r.t $\nu$, we get for all $a\in
{\mathcal M}_{T_R}\cap {\mathcal M}_{\Psi}^+$ and
$t\in\mathbb{R}$:
$$\begin{aligned}
\Psi\circ\kappa_t(a)=\nu\circ\beta^{-1}\circ
T_R\circ\kappa_t(a)&=\nu\circ\beta^{-1}\circ\kappa_t\circ
T_R(a)\\
&=\nu\circ\beta^{-1}\circ\sigma_t^{\Phi'}\circ\tau_{-t}\circ
T_R(a)\\
&=\nu\circ\sigma_{-t}^{\nu}\circ\beta^{-1}\circ\tau_{-t}\circ
T_R(a)\\
&=\nu\circ\sigma_{-2t}^{\nu}\circ T_R(a)=\nu\circ\beta^{-1}\circ
T_R(a)=\Psi(a)
\end{aligned}$$

The proof of the second part is very similar.
\end{proof}

\begin{prop}
$\sigma^{\Phi'}$ and $\tau$ (resp. $\sigma^{\Psi}$ and $\tau$)
commute each other.
\end{prop}

\begin{proof}
Since $\Psi$ is $\kappa$-invariant, we know that
$\sigma^{\Psi}_s\circ\sigma_t^{\Phi'}\circ\tau_{-t}
=\sigma_t^{\Phi'}\circ\tau_{-t}\circ\sigma_s^{\Psi}$, for all
$s,t\in\mathbb{R}$ so that:
$$\begin{aligned}
(id\surl{\ _{\beta} \star_{\alpha}}_{\
N}\kappa_t)\Gamma=\Gamma\circ\kappa_t
&=\Gamma\circ\sigma^{\Psi}_{-s}\circ\kappa_t\circ\sigma^{\Psi}_s
=(\sigma^{\Psi}_{-s}\surl{\ _{\beta} \star_{\alpha}}_{\
N}\tau_s)\circ\Gamma\circ\kappa_t\circ\sigma^{\Psi}_s\\
&=(\sigma^{\Psi}_{-s}\!\!\surl{\ _{\beta} \star_{\alpha}}_{\
N}\tau_s\circ\kappa_t)\circ\Gamma\circ\sigma^{\Psi}_s=(id\!\!\surl{\
_{\beta} \star_{\alpha}}_{\
N}\tau_s\circ\kappa_t\circ\tau_{-s})\circ\Gamma
\end{aligned}$$

So, for all $a\in M$, $\omega\in M_*^+$ and $k\in\mathbb{R}^+$
such that $\omega\circ\beta\leq k\nu$, we get:
$$\sigma_t^{\Phi'}\circ\tau_{-t}((\omega\surl{\ _{\beta} \star_{\alpha}}_{\ \nu}id)\Gamma(a))=
\tau_s\circ\sigma_t^{\Phi'}\circ\tau_{-t}\circ\tau_{-s}((\omega\surl{\
_{\beta} \star_{\alpha}}_{\ \nu}id)\Gamma(a))$$ and by theorem
\ref{densevn}, we easily obtain commutation between
$\sigma^{\Phi'}$ and $\tau$. A similar reasoning from $\kappa'$
implies the second part.
\end{proof}

\begin{coro}\label{mod}
$\sigma^{\Phi'}$ and $\sigma^{\Psi}$ commute each other.
\end{coro}

\begin{proof}
By the previous proposition, we compute, for all
$s,t\in\mathbb{R}$:
$$
\begin{aligned}
\Gamma\circ\sigma^{\Phi'}_s\circ\sigma^{\Psi}_t=(\tau_s \surl{\
_{\beta} \star_{\alpha}}_{\ N} \sigma^{\Phi'}_s)\circ
\Gamma\circ\sigma^{\Psi}_t &=(\tau_s\sigma_t^{\Psi} \surl{\
_{\beta} \star_{\alpha}}_{\ N}
\sigma^{\Phi'}_s\tau_{-t})\circ\Gamma\\
&=(\sigma_t^{\Psi}\tau_s \surl{\ _{\beta} \star_{\alpha}}_{\ N}
\tau_{-t}\sigma^{\Phi'}_s)\circ\Gamma\\
&=(\sigma_t^{\Psi}\surl{\ _{\beta} \star_{\alpha}}_{\ N}
\tau_{-t})\circ\Gamma\circ\sigma^{\Phi'}_s=\Gamma\circ\sigma^{\Psi}_t\circ\sigma^{\Phi'}_s\\
\end{aligned}$$

Since $\Gamma$ is injective, we have done.
\end{proof}

\subsection{First result about uniqueness of invariant operator-valued weight}

In this section, we choose to work with left invariant
operator-valued weights, but it is clear that we have similar
results for right invariant operator-valued weights. Let $T_1$ and
$T_2$ be two n.s.f left invariant operator-valued weights from $M$
to $\alpha(N)$ such that $T_1\leq T_2$. For all $i\in\{1,2\}$, we
put $\Phi_i=\nu\circ\alpha^{-1}\circ T_i$ and
$\hat{\beta}_i(n)=J_{\Phi_i}\alpha(n^*)J_{\Phi_i}$.

We define, as we have done for $U_H$, an isometry $(U_2)_H$ by
the following formula:
$$(U_2)_H(v\surl{\ _{\alpha} \otimes_{\hat{\beta}_2}}_{\
\nu^o}\Lambda_{\Phi_2}(a))=\sum_{i\in I}\xi_i \surl{\ _{\beta}
\otimes_{\alpha}}_{\
\nu}\Lambda_{\Phi_2}((\omega_{v,\xi_i}\surl{\ _{\beta}
\star_{\alpha}}_{\ \nu}id)(\Gamma(a)))$$ for all $v\in
D(H_{\beta},\nu^o)$ and $a\in{\mathcal N}_{\Phi_2}\cap {\mathcal
N}_{T_2}$. Then, we know that $(U_2)_H$ is unitary and
$\Gamma(m)=(U_2)_H(1\surl{\ _{\alpha} \otimes_{\hat{\beta}_2}}_{\
\ N^o}m)(U_2)_H^*$ for all $m\in M$.

Since $T_1\leq T_2$, there exists $F\in
\mathcal{L}(H_{\Phi_2},H_{\Phi_1})$ such that, for all $x\in
{\mathcal N}_{\Phi_2}\cap {\mathcal N}_{T_2}$, we have
$F\Lambda_{\Phi_2}(x)=\Lambda_{\Phi_1}(x)$. It is easy to verify
that, for all $n\in N$, we have $F\hat{\beta}_2(n)=
\hat{\beta}_1(n)F$. If we put $P=F^*F$,then $P$ belongs to
$M'\cap \hat{\beta}_2(N)'$ and $J_{\Phi_2}PJ_{\Phi_2}$ belongs to
$M\cap\alpha(N)'$.

\begin{lemm}
We have $\Gamma(J_{\Phi_2}PJ_{\Phi_2})=1\surl{\ _{\beta}
\otimes_{\alpha}}_{\ N}J_{\Phi_2}PJ_{\Phi_2}$.
\end{lemm}

\begin{proof}
We have, for all $v,w\in D(H_{\beta},\nu^o)$ and $a,b\in{\mathcal
N}_{\Phi_2}\cap {\mathcal N}_{T_2}$:

$$
\begin{aligned}
&\quad((1\surl{\ _{\beta} \otimes_{\alpha}}_{\
N}P)(U_2)_H(v\surl{\ _{\alpha} \otimes_{\hat{\beta}_2}}_{\
\nu^o}\Lambda_{\Phi_2}(a))|(U_2)_H(w\surl{\ _{\alpha}
\otimes_{\hat{\beta}_2}}_{\ \nu^o}\Lambda_{\Phi_2}(b)))\\
&=((U_1)_H(v\surl{\ _{\alpha} \otimes_{\hat{\beta}_1}}_{\
\nu^o}\Lambda_{\Phi_1}(a))|(U_1)_H(w\surl{\ _{\alpha}
\otimes_{\hat{\beta}_1}}_{\ \nu^o}\Lambda_{\Phi_1}(b)))
\end{aligned}$$
where $(U_1)_H$ is defined in the same way as $(U_2)_H$. The two
expressions are continuous in $v$ and $w$, so by density of
$D(H_{\beta},\nu^o)$ in $H$, we get, for all $v,w\in H$ and
$a,b\in{\mathcal N}_{\Phi_2}\cap {\mathcal N}_{T_2}$:

$$
\begin{aligned}
&\quad((1\surl{\ _{\beta} \otimes_{\alpha}}_{\
N}P)(U_2)_H(v\surl{\ _{\alpha} \otimes_{\hat{\beta}_2}}_{\
\nu^o}\Lambda_{\Phi_2}(a))|(U_2)_H(w\surl{\ _{\alpha}
\otimes_{\hat{\beta}_2}}_{\ \nu^o}\Lambda_{\Phi_2}(b)))\\
&=((U_1)_H(v\surl{\ _{\alpha} \otimes_{\hat{\beta}_1}}_{\
\nu^o}\Lambda_{\Phi_1}(a))|(U_1)_H(w\surl{\ _{\alpha}
\otimes_{\hat{\beta}_1}}_{\ \nu^o}\Lambda_{\Phi_1}(b)))\\
&=(v\surl{\ _{\alpha} \otimes_{\hat{\beta}_1}}_{\
\nu^o}\Lambda_{\Phi_1}(a)|w\surl{\ _{\alpha}
\otimes_{\hat{\beta}_1}}_{\ \nu^o}\Lambda_{\Phi_1}(b))\\
&=((1\surl{\ _{\alpha} \otimes_{\hat{\beta}_2}}_{\
N^o}P)(v\surl{\ _{\alpha} \otimes_{\hat{\beta}_2}}_{\
\nu^o}\Lambda_{\Phi_2}(a))|w\surl{\ _{\alpha}
\otimes_{\hat{\beta}_2}}_{\ \nu^o}\Lambda_{\Phi_2}(b))
\end{aligned}$$
so that $(U_2)_H^*(1\surl{\ _{\beta} \otimes_{\alpha}}_{\
N}P)(U_2)_H=1\surl{\ _{\alpha} \otimes_{\hat{\beta}_2}}_{\ N^o}P$.
In particular, if $H=H_{\Phi}$, then by \ref{clef} we get
$(U_2)_H(1\surl{\ _{\alpha} \otimes_{\hat{\beta}_2}}_{\
N^o}J_{\Phi_2}PJ_{\Phi_2})(U_2)_H^*=1\surl{\ _{\beta}
\otimes_{\alpha}}_{\ N}J_{\Phi_2}PJ_{\Phi_2}$. Finally, since
$J_{\Phi_2}PJ_{\Phi_2}\in M$, we have
$\Gamma(J_{\Phi_2}PJ_{\Phi_2})=1\surl{\ _{\beta}
\otimes_{\alpha}}_{\ N}J_{\Phi_2}PJ_{\Phi_2}$.
\end{proof}

\begin{prop}
If $T_1$ and $T_2$ are n.s.f left invariant weights from $M$ to
$\alpha(N)$ such that $T_1\leq T_2$, then there exists an
injective $p\in N$ such that $0\leq p\leq 1$ and, for all $x,y\in
{\mathcal N}_{\Phi_2}\cap {\mathcal N}_{T_2}$, we have
$(\Lambda_{\Phi_1}(x)|\Lambda_{\Phi_1}(y))=(J_{\Phi_2}\beta(p)J_{\Phi_2}\Lambda_{\Phi_2}(x)|
\Lambda_{\Phi_2}(y))$.
\end{prop}

\begin{proof}
Straightforward from the previous lemma and \ref{clef2}.
\end{proof}

\begin{prop}\label{uni}
Let $T_1$ and $T_2$ be n.s.f left invariant weights from $M$ to
$\alpha(N)$ such that $T_1\leq T_2$. If $T_1$ and $T_2$ are
$\beta$-adapted w.r.t $\nu$, then there exists an injective $p\in
Z(N)$ such that $0\leq p\leq 1$ and $\Phi_1=(\Phi_2)_{\beta(p)}$
in the sense of \cite{St}. (We recall that
$\Phi_i=\nu\circ\alpha^{-1}\circ T_i$ for $i=1,2$).
\end{prop}

\begin{proof}
Since $T_1$ and $T_2$ are $\beta$-adapted w.r.t $\nu$,
$\sigma^{T_1}$ and $\sigma^{T_2}$ are equal on $\beta(N)$. For all
$x,y\in {\mathcal N}_{\Phi_2}\cap {\mathcal N}_{T_2}$ and $n\in
{\mathcal T}_{\nu}$, we compute:
$$
\begin{aligned}
(PJ_{\Phi_2}\beta(n)J_{\Phi_2}\Lambda_{\Phi_2}(x)|\Lambda_{\Phi_2}(y))
&=(P\Lambda_{\Phi_2}(x\sigma_{-i/2}^{T_2}(\beta(n))|\Lambda_{\Phi_2}(y))\\
&=(\Lambda_{\Phi_1}(x\sigma_{-i/2}^{T_1}(\beta(n))|\Lambda_{\Phi_1}(y))\\
&=\Phi_1(y^*x\sigma_{-i/2}^{T_1}(\beta(n)))
\end{aligned}$$

Then K.M.S conditions, applied for the n.s.f $\Phi_1$ on $M$,
imply that the last expression is equal to:
$$
\begin{aligned}
\Phi_1(\sigma_{i/2}^{T_1}(\beta(n))y^*x)
&=(\Lambda_{\Phi_1}(x)|\Lambda_{\Phi_1}(y\sigma_{-i/2}^{T_1}(\beta(n^*))))\\
&=(P\Lambda_{\Phi_2}(x)|\Lambda_{\Phi_2}\!(y\sigma_{-i/2}^{T_2}(\beta(n^*))))\\
&=(J_{\Phi_2}\beta(n)J_{\Phi_2}P\Lambda_{\Phi_2}(x)|\Lambda_{\Phi_2}(y))
\end{aligned}$$
That's why $P\in (J_{\Phi_2}\beta(N)J_{\Phi_2})'$ and,
consequently, by the previous proposition and injectivity of
$\beta$, we get $p\in Z(N)$. We know that $0\leq P\leq 1$ and $P$
is injectif with dense range, so the same is for $p$. Finally, by
the previous proposition and \cite{St} (proposition 3.13), we get
that $\beta(p)$ coincides with the analytic continuation in $-i$
of the cocycle $[D\Phi_1:D\Phi_2]$. Then, we have:
$$[D\Phi_1:D\Phi_2]_t=\beta(p)^{it}$$ for all $t\in\mathbb{R}$.
Since $p\in Z(N)$ and $T_2$ is $\beta$-adapted w.r.t $\nu$, we
have $\beta(p)$ belongs to the centralizer of $\Phi_2$ and
$\Phi_1=(\Phi_2)_{\beta(p)}$.
\end{proof}

\begin{lemm}\label{formicru}
Let $T$ and $T'$ be n.s.f operator-valued weights which are
$\beta$-adapted w.r.t $\nu$. If $T+T'$ is semi-finite, then $T+T'$
is $\beta$-adapted w.r.t $\nu$.
\end{lemm}

\begin{proof}
Since $T$ and $T'$ are $\beta$-adapted w.r.t $\nu$, by
\ref{timpo}, there exists n.s.f operator-valued weights $S$ from
$M$ to $\beta(N)$ such that:
$$\nu\circ\alpha^{-1}\circ T=\nu\circ\beta^{-1}\circ S\quad\text{
and }\quad \nu\circ\alpha^{-1}\circ T'=\nu\circ\beta^{-1}\circ
S'$$ Consequently $\nu\circ\alpha^{-1}\circ
(T+T')=\nu\circ\beta^{-1}\circ (S+S')$. This weight is
semi-finite, since $T+T'$ is. Then $S+S'$ is also semi-finite. We
deduce, for all $n\in N$ and $t\in\mathbb{R}$:
$$
\begin{aligned}
\sigma_t^{T+T'}(\beta(n))=\sigma_t^{\nu\circ\alpha^{-1}\circ
(T+T')}(\beta(n)) &=\sigma_t^{\nu\circ\beta^{-1}\circ
(S+S')}(\beta(n))\\
&=\sigma_t^{\nu\circ\beta^{-1}}(\beta(n))=\beta(\sigma_{-t}^{\nu}(n))
\end{aligned}$$
\end{proof}

We recall a technical lemma of \cite{Ku1}:

\begin{lemm}
If $\phi$ and $\eta$ are n.s.f weights on $M$ and if there exists
a strictly positive operator $\lambda$ which is affiliated with
$(M^{\phi})^+$ satisfying:
$$||\Lambda_{\eta}(\sigma_t^{\phi}(x))||=||\lambda^{\frac{t}{2}}\Lambda_{\eta}(x)||$$
for all $x\in {\mathcal N}_{\eta}$ and  $t\in\mathbb{R}$, then
${\mathcal N}_{\eta}\cap {\mathcal N}_{\phi}$ is a core for both
$\Lambda_{\eta}$ and $\Lambda_{\phi}$.
\end{lemm}

\begin{proof}
We can define unitary $T_t$ such that
$T_t\Lambda_{\eta}(a)=\lambda^{-t/2}\Lambda_{\eta}(\sigma_t^{\phi}(a))$
for all $t\in\mathbb{R}$. Moreover, there exists a strictly
positive operator $T$ such that $T^{it}=T_t$ for all
$t\in\mathbb{R}$. The end of the proof is similar to \cite{KV3}
(proposition 1.14).
\end{proof}

\begin{prop}\label{struc}
Let $T_1$ be a n.s.f left invariant operator-valued weight, which
is $\beta$-adapted w.r.t $\nu$, such that there exists a strictly
positive operator $\lambda$ which is affiliated to $(M^{\Phi})^+$
satisfying
$||\Lambda_1(\sigma_t^{\Phi}(x))||=||\lambda^{\frac{t}{2}}\Lambda_1(x)||$
for all $x\in {\mathcal N}_{\Phi_1}$ and  $t\in\mathbb{R}$. Then,
there exists a strictly positive operator $q$ which is affiliated
to the center of $N$ such that $\Phi_1=(\Phi)_{\beta(q)}$.
\end{prop}

\begin{proof}
We put $T_2=T_L+T_1$. Since
$||\Lambda_1(\sigma_t^{\Phi}(x))||=||\lambda^{\frac{t}{2}}\Lambda_1(x)||$
for all $x\in {\mathcal N}_{\Phi_1}$ and $t\in\mathbb{R}$, the
left invariant operator-valued weight $T_2$ is n.s.f. So, by
\ref{formicru}, $T_2$ is $\beta$-adapted w.r.t $\nu$. Finally,
since $T_1\leq T_2$ and $T_L\leq T_2$, by \ref{uni}, there exists
an injective $p\in N$ between $0$ and $1$ such that
$\Phi_1=(\Phi_2)_{\beta(p)}$ and $\Phi=(\Phi_2)_{\beta(1-p)}$. By
\cite{St}, we have:
$$[D\Phi_1:D\Phi_2]_t=\beta(p)^{it} \text{ and } [D\Phi:D\Phi_2]_t=\beta(1-p)^{it}$$ Then, we have, for all
$t\in\mathbb{R}$:
$$[D\Phi_1:D\Phi]_t=[D\Phi_1:D\Phi_2]_t[D\Phi_2
:D\Phi]_t=\beta(\frac{p}{1-p})^{it}$$ that's why
$q=\frac{p}{1-p}$ is the suitable element.
\end{proof}

\subsection{Modulus and scaling operator}

From now, we study the following measured quantum group
$(N,M,\alpha,\beta,\Gamma,\nu,T_L,R\circ T_L\circ R)$ so that we
look at $\Phi=\nu\circ\alpha^{-1}\circ T_L$ and $\Phi\circ R$.
Then, we recall $\sigma_t^{\Phi\circ
R}=R\circ\sigma_{-t}^{\Phi}\circ R$ for all $t\in\mathbb{R}$.

By \ref{mod}, we know that modular groups associated with $\Phi$
and $\Phi\circ R$ commute each other and, by \cite{Vae}
(proposition 2.5), there exist a strictly positive operator
$\delta$ affiliated with $M$ and a strictly positive operator
$\lambda$ affiliated to the center of $M$ such that, for all
$t\in\mathbb{R}$, we have $[D\Phi\circ R
:D\Phi]_t=\lambda^{\frac{1}{2}it^2}\delta^{it}$. Modular groups of
$\Phi$ and $\Phi\circ R$ are linked by $\sigma_t^{\Phi\circ
R}(m)=\delta^{it}\sigma_t^{\Phi}(m)\delta^{-it}$ for all
$t\in\mathbb{R}$ and $m\in M$.

\begin{defi}
We call \textbf{scaling operator} the strictly positive operator
$\lambda$ affiliated to $Z(M)$ and \textbf{modulus} the strictly
positive operator $\delta$ affiliated to $M$ such that, for all
$t\in\mathbb{R}$, we have:
$$[D\Phi\circ R:D\Phi]_t=\lambda^{\frac{1}{2}it^2}\delta^{it}$$
\end{defi}

In this section, we establish properties of scaling operator and
modulus e.g compatibility of these objects with the Hopf bimodule
structure.

\begin{prop}\label{deltadot}
The scaling operator does not depend on the quasi-invariant weight
but just on the modular group associated with. If $\dot{\delta}$
is the class of $\delta$ for the equivalent relation
$\delta_1\sim\delta_2$ if, and only if there exists a strictly
positive operator $h$ affiliated to $Z(N)$ such that
$\delta_2^{it}=\beta(h^{it})\delta_1^{it}\alpha(h^{-it})$, then
$\dot{\delta}$ does not depend on the quasi-invariant weight but
just on the modular group associated with.
\end{prop}

\begin{proof}
If $\nu'$ is a n.s.f weight on $N$ such that
$\sigma^{\nu'}=\sigma^{\nu}$, then there exists a strictly
positive $h$ affiliated to $Z(N)$ such that $\nu'=\nu_h$. We just
have to compute:
$$
\begin{aligned}
&\ \quad [D\nu'\circ\alpha^{-1}\circ T_L\circ
R:D\nu'\circ\alpha^{-1}\circ
T_L]_t\\
&=[D\nu_h\circ\alpha^{-1}\circ T_L\circ R:D\Phi\circ
R]_t[D\Phi\circ R:D\Phi]_t[D\Phi:D\nu_h\circ\alpha^{-1}\circ T_L]_t\\
&=\beta([D\nu_h:D\nu]^*_{-t})\lambda^{\frac{1}{2}it^2}\delta^{it}\alpha([D\nu:D\nu_h]_t)=\lambda^{\frac{1}{2}it^2}\beta(h^{it})\delta^{it}\alpha(h^{-it})
\end{aligned}$$
\end{proof}

\begin{lemm}\label{prep1}
For all $s,t\in\mathbb{R}$, we have
$[D\Phi\circ\sigma_s^{\Phi\circ R}:D\Phi]_t=\lambda^{ist}$.
\end{lemm}

\begin{proof}
The computation of the cocycle is straightforward:
$$
\begin{aligned}
\ [D\Phi\circ\sigma_s^{\Phi\circ R}:D\Phi]_t
&=[D\Phi\circ\sigma_s^{\Phi\circ R}:D\Phi\circ
R\circ\sigma_s^{\Phi\circ R}]_t
[D\Phi\circ R:D\Phi]_t\\
&=\sigma_{-s}^{\Phi\circ R}([D\Phi:D\Phi\circ R]_t)[D\Phi\circ R:D\Phi]_t\\
&=\delta^{-is}\sigma_{-s}^{\Phi}(\lambda^{-\frac{it^2}{2}}\delta^{-it})
\delta^{is}\lambda^{\frac{it^2}{2}}\delta^{it}\\
&=\delta^{-is}\lambda^{-\frac{it^2}{2}}\lambda^{ist}
\delta^{-it}\delta^{is}\lambda^{\frac{it^2}{2}}\delta^{it}=\lambda^{ist}
\end{aligned}$$
\end{proof}

\begin{prop}
We have $R(\lambda)=\lambda$, $R(\delta)=\delta^{-1}$ and
$\tau_t(\delta)=\delta$, $\tau_t(\lambda)=\lambda$ for all
$t\in\mathbb{R}$.
\end{prop}

\begin{proof}
Relations between $R$, $\lambda$ and $\delta$ come from uniqueness
of Radon-Nikodym cocycle decomposition. By \ref{prep}, we have
$\Phi\circ\tau_{-s}=\Phi\circ\sigma_s^{\Phi\circ R}$ for all
$s,t\in\mathbb{R}$, so:
$$\tau_s([D\Phi\circ R:D\Phi]_t)=[D\Phi\circ
R\circ\tau_{-s}:D\Phi\circ\tau_{-s}]_t=[D\Phi\circ\sigma_s^{\Phi\circ
R}\circ R:D\Phi\circ\sigma_s^{\Phi\circ R}]_t$$ Consequently, by
the previous lemma, we get:
$$
\begin{aligned}
&\ \quad\tau_s([D\Phi\circ
R:D\Phi]_t)\\
&=[D\Phi\circ\sigma_s^{\Phi\circ R}\circ
R:D\Phi\circ R]_t[D\Phi\circ R:D\Phi]_t[D\Phi:D\Phi\circ\sigma_s^{\Phi\circ R}]_t\\
&=R([D\Phi\circ\sigma_s^{\Phi\circ R}:D\Phi]^*_{-t})
[D\Phi\circ R:D\Phi]_t[D\Phi\circ\sigma_s^{\Phi\circ R}:D\Phi]_t^*\\
&=R(\lambda^{ist})\lambda^{-\frac{it^2}{2}}\delta^{it}\lambda^{-ist}
=\lambda^{-\frac{it^2}{2}}\delta^{it}
\end{aligned}$$
\end{proof}

\begin{coro}
The modulus $\delta$ is affiliated with
$M\cap\alpha(N)'\cap\beta(N)'$.
\end{coro}

\begin{proof}
Since $\Phi=\nu\circ\beta^{-1}\circ S_L$, we have:
$$\lambda^{\frac{it^2}{2}}\delta^{it}=[D\Phi\circ
R:D\Phi]_t=[DR\circ T_L\circ R:DS_L]_t$$ which belongs to
$M\cap\beta(N)'$. Since $\lambda$ is affiliated with $Z(M)$, we
get that $\delta$ is affiliated with $M\cap\beta(N)'$. Finally,
since $R(\delta)=\delta$, we obtain that $\delta$ is affiliated
with $M\cap\alpha(N)'\cap\beta(N)'$.
\end{proof}

\begin{lemm}
For all $t\in\mathbb{R}$, $\tau_{-t}\circ T_L\circ\tau_t$ is a
n.s.f left invariant operator-valued weight from $M$ to
$\alpha(N)$. Moreover, $\tau_{-t}\circ T_L\circ\tau_t$ is
$\beta$-adapted for $\nu$.
\end{lemm}

\begin{proof}
For all $t\in\mathbb{R}$, we have
$\nu\circ\alpha^{-1}\circ\tau_{-t}\circ
T_L\circ\tau_t=\Phi\circ\tau_t$. Then:
$$
\begin{aligned}
(id\surl{\ _{\beta} \star_{\alpha}}_{\
\nu}\nu\circ\alpha^{-1}\circ\tau_{-t}\circ
T_L\circ\tau_t)\circ\Gamma&=(id\surl{\ _{\beta} \star_{\alpha}}_{\
\nu}\Phi\circ\tau_t)\circ\Gamma\\
&=\tau_{-t}\circ(id\surl{\ _{\beta} \star_{\alpha}}_{\
\nu}\Phi)\circ\Gamma\circ\tau_t=\tau_{-t}\circ T_L\circ\tau_t
\end{aligned}$$

On the other hand, for all $s,t\in\mathbb{R}$ and $n\in N$, we
have:
$$
\begin{aligned}
\sigma_s^{\Phi\circ\tau_t}(\beta(n))=\tau_{-t}\circ\sigma_s^{\Phi}\circ\tau_t(\beta(n))
&=\tau_{-t}\circ\sigma_s^{\Phi}(\beta(\sigma_{-t}^{\nu}(n)))\\
&=\tau_{-t}(\beta(\sigma_{-(s+t)}^{\nu}(n)))=\beta(\sigma_{-s}^{\nu}(n))
\end{aligned}$$
\end{proof}

\begin{prop}
There exists a strictly positive operator $q$ affiliated with
$Z(N)$ such that the scaling operator
$\lambda=\alpha(q)=\beta(q)$. In particular, $\lambda$ is
affiliated with $Z(M)\cap\alpha(N)\cap\beta(N)$.
\end{prop}

\begin{proof}
By the previous lemma, $\tau_s\circ T_L\circ\tau_{-s}$ is left
invariant and $\beta$-adapted w.r.t $\nu$. Moreover, since
$\sigma^{\Phi}$ and $\tau$ commute, $\Phi\circ\tau_{-s}$ is
$\sigma^{\Phi}$-invariant. That's why, we are in \ref{struc}
conditions so that we get a strictly positive operator $q_s$
affiliated with $Z(N)$ such that
$[D\Phi\circ\tau_{-s}:D\Phi]_t=\beta(q_s)^{it}$. On the other
hand, by \ref{prep1}, we have $[D\Phi\circ\sigma_s^{\Phi\circ
R}:D\Phi]_t=\lambda^{ist}$. By \ref{prep}, we have
$\Phi\circ\tau_{-s}=\Phi\circ\sigma_s^{\Phi\circ R}$, so we
obtain that $\lambda^{ist}=\beta(q_s)^{it}$ for all
$s,t\in\mathbb{R}$. We easily deduce that there exists a strictly
positive operator $q$ affiliated with $Z(N)$ such that
$\lambda=\beta(q)$. Finally, since $R(\lambda)=\lambda$, we also
have $\lambda=\alpha(q)$.
\end{proof}

\begin{lemm}
We have, for all $a,b\in {\mathcal N}_{\Phi}\cap {\mathcal
N}_{T_L}$ and  $t\in\mathbb{R}$:
$$\omega_{J_{\Phi}\Lambda_{\Phi}(\lambda^{\frac{t}{2}}\tau_t(a))}
=\omega_{J_{\Phi}\Lambda_{\Phi}(a)}\circ\tau_{-t}\text{ and }
\omega_{J_{\Phi}\Lambda_{\Phi}(b)}\circ\sigma_t^{\Phi\circ R}=
\omega_{J_{\Phi}\Lambda_{\Phi}(\lambda^{\frac{t}{2}}\sigma_{-t}^{\Phi\circ
R}(b))}$$
\end{lemm}

\begin{proof}
$\tau$ is implemented by $P$ that's why the first relation holds.
By \cite{Vae} (proposition 2.4), we know that $\Delta_{\Phi\circ
R}=J_{\Phi}\delta J_{\Phi}\delta\Delta_{\Phi}$ so that we can
compute:
$$
\begin{aligned}
&\ \quad (\sigma_t^{\Phi\circ
R}(x)J_{\Phi}\Lambda_{\Phi}(b)|J_{\Phi}\Lambda_{\Phi}(b))\\
&=(x\Delta_{\Phi\circ R}^{-it}J_{\Phi}\Lambda_{\Phi}(b)|
\Delta_{\Phi\circ R}^{-it}J_{\Phi}\Lambda_{\Phi}(b))\\
&=(xJ_{\Phi}\delta^{it}\!\!J_{\Phi}\delta^{-it}\!\Delta_{\Phi}^{-it}\!\!J_{\Phi}\Lambda_{\Phi}(b)|
J_{\Phi}\delta^{it}\!J_{\Phi}\delta^{-it}\!\Delta_{\Phi}^{-it}\!J_{\Phi}\Lambda_{\Phi}(b))\\
&=(x\delta^{-it}J_{\Phi}\Lambda_{\Phi}(\sigma_{-t}^{\Phi}(b))|
\delta^{-it}J_{\Phi}\Lambda_{\Phi}(\sigma_{-t}^{\Phi}(b)))\\
&=(xJ_{\Phi}\Lambda_{\Phi}(\lambda^{\frac{t}{2}}\sigma_{-t}^{\Phi}(b)\delta^{it})|
J_{\Phi}\Lambda_{\Phi}(\lambda^{\frac{t}{2}}\sigma_{-t}^{\Phi}(b)\delta^{it}))\\
&=(xJ_{\Phi}\Lambda_{\Phi}(\lambda^{\frac{t}{2}}\delta^{-it}\sigma_{-t}^{\Phi}(b)\delta^{it})|
J_{\Phi}\Lambda_{\Phi}(\lambda^{\frac{t}{2}}\delta^{-it}\sigma_{-t}^{\Phi}(b)\delta^{it}))\\
&=(xJ_{\Phi}\Lambda_{\Phi}(\lambda^{\frac{t}{2}}\sigma_{-t}^{\Phi\circ
R}(b))|
J_{\Phi}\Lambda_{\Phi}(\lambda^{\frac{t}{2}}\sigma_{-t}^{\Phi\circ
R}(b)))
\end{aligned}$$

\end{proof}

\begin{prop}
We have $\Gamma\circ\tau_t=(\sigma_t^{\Phi}\surl{\ _{\beta}
  \star_{\alpha}}_{\
  N}\sigma_{-t}^{\Phi\circ R})\circ\Gamma$ for all $t\in\mathbb{R}$.
\end{prop}

\begin{proof}
For all $a,b\in {\mathcal N}_{\Phi}\cap {\mathcal N}_{T_L}$ and
$t\in\mathbb{R}$, we compute:
$$
\begin{aligned}
&\ \quad(id\surl{\ _{\beta}\star_{\alpha}}_{\
\nu}\omega_{J_{\Phi}\Lambda_{\Phi}(b)})[(\sigma^{\Phi}_{-t}\surl{\
_{\beta}\star_{\alpha}}_{\ N}\sigma^{\Phi\circ
R}_t)\circ\Gamma\circ\tau_t(a^*a)]\\
&=\sigma^{\Phi}_{-t}[(id\surl{\ _{\beta}\star_{\alpha}}_{\
\nu}\omega_{J_{\Phi}\Lambda_{\Phi}(b)}\circ\sigma_t^{\Phi\circ
R})(\Gamma\circ\tau_t(a^*a))]
\end{aligned}$$ By the previous
lemma, this last expression is equal to:
$$
\begin{aligned}
&\ \quad\sigma^{\Phi}_{-t}[(id\surl{\ _{\beta}\star_{\alpha}}_{\
\nu}\omega_{J_{\Phi}\Lambda_{\Phi}(\lambda^{\frac{t}{2}}\sigma_{-t}^{\Phi\circ
R}(b))})(\Gamma\circ\tau_t(a^*a))]\\
&=\sigma^{\Phi}_{-t}\circ R[(id\surl{\ _{\beta}\star_{\alpha}}_{\
\nu}\omega_{J_{\Phi}\Lambda_{\Phi}(\tau_t(a))}
)(\Gamma(\lambda^t\sigma_{-t}^{\Phi\circ R}(b^*b)))]\\
&=R\circ\sigma_t^{\Phi\circ R}[(id\surl{\
_{\beta}\star_{\alpha}}_{\
\nu}\omega_{J_{\Phi}\Lambda_{\Phi}(\lambda^{\frac{t}{2}}\tau_t(a))}
)(\Gamma\circ\sigma_{-t}^{\Phi\circ R}(b^*b))]
\end{aligned}$$ Again, by the previous lemma, this last expression
is equal to:
$$
\begin{aligned}
&\ \quad R\circ\sigma_t^{\Phi\circ R}[(id\surl{\
_{\beta}\star_{\alpha}}_{\
\nu}\omega_{J_{\Phi}\Lambda_{\Phi}(a)}\circ\tau_{-t}
(\Gamma\circ\sigma_{-t}^{\Phi\circ R}(b^*b))]\\
&=R[(id\surl{\ _{\beta}\star_{\alpha}}_{\
\nu}\omega_{J_{\Phi}\Lambda_{\Phi}(a)}
)(\Gamma(b^*b))]=(id\surl{\ _{\beta}\star_{\alpha}}_{\
\nu}\omega_{J_{\Phi}\Lambda_{\Phi}(b)} )(\Gamma(a^*a))
\end{aligned}$$ So, we conclude that $(\sigma^{\Phi}_{-t}\surl{\
_{\beta}\star_{\alpha}}_{\ N}\sigma^{\Phi\circ
R}_t)\circ\Gamma\circ\tau_t=\Gamma$ for all $t\in\mathbb{R}$.
\end{proof}

Since $\delta$ is affiliated with $M\cap\alpha(N)'\cap\beta(N)'$,
we can define an operator $\delta^{it}\surl{\ _{\beta}
\otimes_{\alpha}}_{\ N}\delta^{it}$ which belongs to
$(M\cap\beta(N)')\surl{\ _{\beta}\otimes_{\alpha}}_{\
N}(M\cap\alpha(N)')\subset M\surl{\ _{\beta}\star_{\alpha}}_{\
N}M$ for all $t\in\mathbb{R}$. Now, we prove that $\delta$ is a
group-like element i.e $\Gamma(\delta)=\delta\surl{\
_{\beta}\otimes_{\alpha}}_{\ N}\delta$.

\begin{lemm}
For all $s,t\in\mathbb{R}$, $\Gamma(\delta^{is})$ and
$\delta^{it}\surl{\ _{\beta}\otimes_{\alpha}}_{\ N}\delta^{it}$
commute each other.
\end{lemm}

\begin{proof}
For all $t\in\mathbb{R}$, we have:
$$
\begin{aligned}
(\sigma_{-t}^{\Phi}\circ\sigma_t^{\Phi\circ R}\surl{\ _{\beta}
\star_{\alpha}}_{\ N}\sigma_{-t}^{\Phi}\circ\sigma_t^{\Phi\circ
R})\circ\Gamma&=(\sigma_{-t}^{\Phi}\surl{\
_{\beta}\star_{\alpha}}_{\
N}\sigma_{-t}^{\Phi}\circ\sigma_t^{\Phi\circ
R}\circ\tau_t)\circ\Gamma\circ\sigma_t^{\Phi\circ R}\\
&=(\sigma_{-t}^{\Phi}\circ\tau_t\surl{\ _{\beta}
\star_{\alpha}}_{\ N}\sigma_t^{\Phi\circ R}\circ\tau_t)\circ
\Gamma\circ\sigma_{-t}^{\Phi}\circ\sigma_t^{\Phi\circ R}\\
&=(\sigma_{-t}^{\Phi}\surl{\ _{\beta}\star_{\alpha}}_{\
N}\sigma_t^{\Phi\circ
R})\circ\Gamma\circ\tau_t\sigma_{-t}^{\Phi}\circ\sigma_t^{\Phi\circ
R}\\
&=\Gamma\circ\sigma_{-t}^{\Phi}\circ\sigma_t^{\Phi\circ R}
\end{aligned}$$ We know that
$\sigma_{-t}^{\Phi}\sigma_t^{\Phi\circ
R}(m)=\delta^{it}m\delta^{-it}$ for all $m\in M$, that's why we
get: $$(\delta^{it}\surl{\ _{\beta} \otimes_{\alpha}}_{\
N}\delta^{it})\Gamma(m)(\delta^{-it}\surl{\
_{\beta}\otimes_{\alpha}}_{\
N}\delta^{-it})=\Gamma(\delta^{it}m\delta^{-it})$$ In particular,
for all $s\in\mathbb{R}$, we have:
$$(\delta^{it}\surl{\ _{\beta}
  \otimes_{\alpha}}_{\ N}\delta^{it})\Gamma(\delta^{is})(\delta^{-it}\surl{\ _{\beta}
  \otimes_{\alpha}}_{\
  N}\delta^{-it})=\Gamma(\delta^{it}\delta^{is}\delta^{-it})=\Gamma(\delta^{is})$$
\end{proof}

We recall \cite{KV1} (result 7.6) in the following lemma:

\begin{lemm}
There exists a subspace $C$ of $M$ such that, for all $c\in C$:
\begin{enumerate}[i)]
\item $c\in {\mathcal T}_{\Phi\circ R}$;
\item $\delta^{-1/2}\sigma^{\Phi\circ R}_{-i/2}(c^*)$ is bounded and belongs to
${\mathcal D}(\sigma_{i/2}^{\Phi\circ R})\cap {\mathcal
N}_{\Phi\circ R}$;
\item $\delta^{-1/2}c^*$ is bounded and belongs to ${\mathcal N}_{\Phi\circ R}^*$;
\item $\sigma_{i/2}^{\Phi\circ R}
(\delta^{-1/2}\sigma^{\Phi\circ
R}_{-i/2}(c^*))=\lambda^{-i/4}\delta^{-1/2}c^*$;
\item $\Lambda_{\Phi\circ R}(C)$ is a core for $\delta^{-1/2}$;
\item $c\in {\mathcal N}_{\Phi}$,
$\Phi(c^*c) =\Phi\circ R((\delta^{-1/2}\sigma_{-i/2}^{\Phi\circ
R}(c^*))^*\delta^{-1/2}\sigma_{-i/2}^{\Phi\circ R}(c^*))$ ;
\item $T_L(c^*c)=S_R(\delta^{-1/2}c^*c\delta^{-1/2})$
\end{enumerate}
(We recall that $S_R$ satisfies $\nu\circ\beta^{-1}\circ
T_R=\Phi\circ R=\nu\circ\alpha^{-1}\circ S_R$).
\end{lemm}

\begin{prop}
If $S_R$ is the unique n.s.f operator-valued weight which
satisfies $\nu\circ\beta^{-1}\circ T_R=\Phi\circ
R=\nu\circ\alpha^{-1}\circ S_R$, then we have: $$\Phi\circ
R((\omega_v\surl{\ _{\beta} \star_{\alpha}}_{\
\nu}id)(\Gamma(a)))=(S_R(a)\delta^{-1/2}v|\delta^{-1/2}v)$$ for
all $v\in {\mathcal D}(\delta^{-1/2})\cap D((H_{\Phi\circ
R})_{\beta},\nu^o)$ and $a\in {\mathcal N}_{\Phi\circ R}\cap
{\mathcal N}_{S_R}$.
\end{prop}

\begin{proof}
Let $c\in C$ and $d\in {\mathcal N}_{\Phi\circ R}\cap {\mathcal
N}_{T_R}$. By left invariance of $T_L$, we have:
$$
\begin{aligned}
\Phi\circ R((\omega_{J_{\Phi\circ R}\Lambda_{\Phi\circ
R}(c)}\surl{\ _{\beta} \star_{\alpha}}_{\ \nu}id)(\Gamma(d^*d)))
&=\Phi((\omega_{J_{\Phi\circ R}\Lambda_{\Phi\circ R}(d)}\surl{\
_{\beta} \star_{\alpha}}_{\
\nu}id)(\Gamma(c^*c)))\\
&=(T_L(c^*c)J_{\Phi\circ R}\Lambda_{\Phi\circ R}(d)|J_{\Phi\circ
R}\Lambda_{\Phi\circ R}(d))
\end{aligned}$$
By properties of elements of $C$, this last expression is equal
to:
$$(S_R(\delta^{-1/2}c^*c\delta^{-1/2})J_{\Phi\circ R}\Lambda_{\Phi\circ R}(d)|
J_{\Phi\circ R}\Lambda_{\Phi\circ R}(d))$$ If we denote by
$\varepsilon$ the anti-representation of $N$ such that
$\varepsilon(n)=J_{\Phi\circ R}\alpha(n^*)J_{\Phi\circ R}$ for
all $n\in N$, then the expression is equal to:
$$
\begin{aligned}
&\ \quad ||J_{\Phi\circ R}\Lambda_{\Phi\circ
R}(c\delta^{-1/2})\surl{\ _{\alpha}
\otimes_{\varepsilon}}_{\ \nu^o}\Lambda_{\Phi\circ R}(d)||^2\\
&=(S_R(d^*d)
J_{\Phi\circ R}\Lambda_{\Phi\circ R}(c\delta^{-1/2})|J_{\Phi\circ R}\Lambda_{\Phi\circ R}(c\delta^{-1/2}))\\
&=(S_R(d^*d)\Lambda_{\Phi\circ R}(\sigma_{-i/2}^{\Phi\circ
R}(\delta^{-1/2}c^*))| \Lambda_{\Phi\circ
R}(\sigma_{-i/2}^{\Phi\circ R}(\delta^{-1/2}c^*)))
\end{aligned}$$
Then, properties of elements of $C$ allow to finish the
computation to get:
$$
\begin{aligned}
&\ \quad\Phi\circ R((\omega_{J_{\Phi\circ R}\Lambda_{\Phi\circ
R}(c)}\surl{\ _{\beta} \star_{\alpha}}_{\ \nu}id)(\Gamma(d^*d)))\\
&=(S_R(d^*d)\delta^{-1/2}J_{\Phi\circ R}\Lambda_{\Phi\circ
R}(c)|\delta^{-1/2}J_{\Phi\circ R}\Lambda_{\Phi\circ R}(c))
\end{aligned}$$

By continuity, the proposition holds.
\end{proof}

\begin{coro}
For all $v\in {\mathcal D}(\delta^{-1/2})\cap D((H_{\Phi\circ
R})_{\beta},\nu^o)$ and for all element $a\in {\mathcal
N}_{\Phi\circ R}\cap {\mathcal N}_{S_R}$, we have:
$$S_R((\omega_v\surl{\
_{\beta} \star_{\alpha}}_{\
\nu}id)(\Gamma(a)))=\alpha(<S_R(a)\delta^{-1/2}v,\delta^{-1/2}v>_{\beta,\nu^o})$$
\end{coro}

\begin{proof}
Straightforward by the formula
$\nu\circ\alpha^{-1}(<x\xi,\eta>_{\beta,\nu^o})=(x\xi|\eta)$ for
all $x\in\beta(N)'$ and $\xi,\eta\in D(H_{\beta},\nu^o)$, and the
previous proposition.
\end{proof}

We put $\Gamma^{(2)}$ the *-homomorphism $(\Gamma\surl{\ _{\beta}
\star_{\alpha}}_{\ N}id)\circ\Gamma=(id\surl{\ _{\beta}
\star_{\alpha}}_{\ N}\Gamma)\circ\Gamma$ from $M$ to $M\surl{\
_{\beta} \otimes_{\alpha}}_{\ N}M\surl{\ _{\beta}
\otimes_{\alpha}}_{\ N}M$.

\begin{lemm}
If $v\in {\mathcal D}(\delta^{-1/2}\surl{\ _{\beta}
\otimes_{\alpha}}_{\ N}\delta^{-1/2})\cap D((H_{\Phi\circ
R}\surl{\ _{\beta} \otimes_{\alpha}}_{\ \nu}H_{\Phi\circ
R})_{\beta},\nu^o)$ and $a\in {\mathcal M}_{\Phi\circ R}\cap
{\mathcal N}_{S_R}$, then we have:
$$\begin{aligned}
&\ \quad\Phi\circ R((\omega_v\surl{\ _{\beta}\star_{\alpha}}_{\
\nu}id)(\Gamma^{(2)}(a)))\\
&=((S_R(a)\surl{\ _{\beta} \otimes_{\alpha}}_{\
N}1)(\delta^{-1/2}\surl{\ _{\beta} \otimes_{\alpha}}_{\
N}\delta^{-1/2})v|(\delta^{-1/2}\surl{\ _{\beta}
\otimes_{\alpha}}_{\ N}\delta^{-1/2})v)
\end{aligned}$$
\end{lemm}

\begin{proof}
Let $\xi,\eta\in {\mathcal D}(\delta^{-1/2}))\cap D((H_{\Phi\circ
R})_{\beta},\nu^o)$. By the previous proposition and its
corollary, we have:
$$
\begin{aligned}
&\ \quad\Phi\circ R((\omega_{\xi\surl{\ _{\beta}
\otimes_{\alpha}}_{\ \nu}\eta}\surl{\ _{\beta} \star_{\alpha}}_{\
\nu}id)(\Gamma^{(2)}(a)))\\
&=\Phi\circ R((\omega_{\eta}\surl{\ _{\beta} \star_{\alpha}}_{\
\nu}id)(\Gamma((\omega_{\xi}\surl{\
_{\beta}\star_{\alpha}}_{\ \nu}id)(\Gamma(a)))))\\
&=(S_R((\omega_{\xi}\surl{\ _{\beta} \star_{\alpha}}_{\
\nu}id)(\Gamma(a)))\delta^{-1/2}\eta|\delta^{-1/2}\eta)\\
&=(\alpha(<S_R(a)\delta^{-1/2}\xi,\delta^{-1/2}\xi>_{\beta,\nu^o}
\delta^{-1/2}\eta|\delta^{-1/2}\eta)
\end{aligned}$$
which is equal to: $$((S_R(a)\surl{\ _{\beta} \otimes_{\alpha}}_{\
N}1)(\delta^{-1/2}\xi\surl{\ _{\beta} \otimes_{\alpha}}_{\
\nu}\delta^{-1/2}\eta)|\delta^{-1/2}\xi\surl{\ _{\beta}
\otimes_{\alpha}}_{\ N}\delta^{-1/2}\eta)$$ Since ${\mathcal
D}(\delta^{-1/2}))\cap D((H_{\Phi\circ R})_{\beta},\nu^o)\surl{\
_{\beta} \otimes_{\alpha}}_{\ \nu}{\mathcal
D}(\delta^{-1/2}))\cap D((H_{\Phi\circ R})_{\beta},\nu^o)$ is a
core for $\delta^{-1/2}\surl{\ _{\beta}\otimes_{\alpha}}_{\
N}\delta^{-1/2}$, the lemma is proved.
\end{proof}

\begin{lemm}We have:
$$\Phi\circ R((\omega_v\surl{\ _{\beta} \star_{\alpha}}_{\
\nu}id)(\Gamma^{(2)}(a)))=((S_R(a)\surl{\ _{\beta}
\otimes_{\alpha}}_{\
N}1)(\Gamma(\delta^{-1/2})v|(\Gamma(\delta^{-1/2})v)$$ for all
$v\in {\mathcal D}(\Gamma(\delta^{-1/2}))\cap D((H_{\Phi\circ
R}\surl{\ _{\beta} \otimes_{\alpha}}_{\ \nu}H_{\Phi\circ
R})_{\beta},\nu^o)$ and $a\in {\mathcal N}_{\Phi\circ R}\cap
{\mathcal N}_{S_R}$.
\end{lemm}

\begin{proof}
Let $(\eta_i)_{i\in I}$ be a $D(H_{\beta},\nu^o)$-basis of
$H_{\beta}$. We put
$w_i=(\lambda_{\eta_i}^{\alpha,\hat{\beta}})^*Wv$ for all $i\in
I$ and we have, for all $m\in M$:
$$
\begin{aligned}
(\Gamma(m)v|v)&=((1\surl{\ _{\alpha} \otimes_{\hat{\beta}}}_{\
N^o}m)Wv|Wv)\\
&=\sum_{i\in i}((\lambda_{\eta_i}^{\alpha,\hat{\beta}})^*(1\surl{\
_{\alpha} \otimes_{\hat{\beta}}}_{\
N^o}m)Wv|(\lambda_{\eta_i}^{\alpha,\hat{\beta}})^*Wv)=\sum_{i\in
I}(mw_i|w_i)
\end{aligned}$$

Since $v\in {\mathcal D}(\Gamma(\delta^{-1/2}))$ and
$\Gamma(x)=W^*(1\surl{\ _{\alpha} \otimes_{\hat{\beta}}}_{\
N^o}x)W$, we notice that $Wv$ belongs to ${\mathcal D}(1\surl{\
_{\alpha} \otimes_{\hat{\beta}}}_{\ N^o}\delta^{-1/2})$ and $w_i$
belongs to ${\mathcal D}(\delta^{-1/2})$. Then, by the previous
proposition and normality of $\Phi\circ R$, we have:
$$\begin{aligned}
\Phi\circ R((\omega_v\surl{\ _{\beta} \star_{\alpha}}_{\
\nu}id)(\Gamma^{(2)}(a)))&=\Phi\circ R((\omega_v\circ\Gamma\surl{\
_{\beta} \star_{\alpha}}_{\ \nu}id)(\Gamma(a)))\\
&=\sum_{i\in I}\Phi\circ R((\omega_{w_i}\surl{\
_{\beta} \star_{\alpha}}_{\ \nu}id)(\Gamma(a)))\\
&=\sum_{i\in I}(S_R(a)\delta^{-1/2}w_i|\delta^{-1/2}w_i)\\
&=((S_R(a)\surl{\ _{\beta} \otimes_{\alpha}}_{\
N}1)(\Gamma(\delta^{-1/2})v|(\Gamma(\delta^{-1/2})v)
\end{aligned}$$

By continuity, the proposition holds.
\end{proof}

\begin{prop}
We have $\Gamma(\delta)=\delta\surl{\ _{\beta}
\otimes_{\alpha}}_{\ N}\delta$.
\end{prop}

\begin{proof}
For all element $v$ in the intersection of ${\mathcal
D}(\delta^{-1/2}\surl{\ _{\beta} \otimes_{\alpha}}_{\
N}\delta^{-1/2})$, ${\mathcal D}(\Gamma(\delta^{-1/2}))$ and
$D((H_{\Phi\circ R}\surl{\ _{\beta} \otimes_{\alpha}}_{\
\nu}H_{\Phi\circ R})_{\beta},\nu^o)$, we have:
$$||(\delta^{-1/2}\surl{\ _{\beta} \otimes_{\alpha}}_{\
N}\delta^{-1/2})v||^2=||\Gamma(\delta^{-1/2})v||^2$$ But, these
operators commute each other so that there exists a subspace
which is a core for both $\delta^{-1/2}\surl{\ _{\beta}
\otimes_{\alpha}}_{\ N}\delta^{-1/2}$ and
$\Gamma(\delta^{-1/2})$. By the previous equality, we get that
the two operators have the same domain and consequently they are
equal.
\end{proof}

\subsection{Uniqueness of left invariant operator-valued weight}

\begin{theo}
If $T'$ a n.s.f left invariant operator-valued weight which is
$\beta$-adapted w.r.t $\nu$, then there exists a strictly
positive operator $h$ affiliated with $Z(N)$ such that, for all
$t\in\mathbb{R}$, we have:
$$\nu\circ\alpha^{-1}\circ T'=(\nu\circ\alpha^{-1}\circ
T_L)_{\beta(h)}\text{ and } [DT':DT_L]_t=\beta(h^{it})$$
\end{theo}

\begin{proof}
We put $\Phi'=\nu\circ\alpha^{-1}\circ T'$. By \ref{clef}, we
have for all $s\in\mathbb{R}$:
$$\Gamma\circ\sigma_{-s}^{\Phi}\circ\sigma_s^{\Phi'}
=(\tau_{-s}\surl{\ _{\beta}\star_{\alpha}}_{\
N}\sigma_{-s}^{\Phi})\circ\Gamma\circ\sigma_s^{\Phi'}= (id\surl{\
_{\beta} \star_{\alpha}}_{\
N}\sigma_{-s}^{\Phi}\circ\sigma_s^{\Phi'})\circ\Gamma$$

By right invariance $T_R$, we have for all $a\in {\mathcal
M}_{T_R}^+$:
$$
\begin{aligned}
T_R(\sigma_{-s}^{\Phi}\circ\sigma_s^{\Phi'}(a))&=(\Phi\circ
R\surl{\ _{\beta} \star_{\alpha}}_{\
\nu}id)(\Gamma(\sigma_{-s}^{\Phi}\circ\sigma_s^{\Phi'}(a)))\\
&=\sigma_{-s}^{\Phi}\circ\sigma_s^{\Phi'}((\Phi\circ R\surl{\
_{\beta} \star_{\alpha}}_{\ \nu}id)\Gamma(a))
=\sigma_{-s}^{\Phi}\circ\sigma_s^{\Phi'}(T_R(a))
\end{aligned}$$ Since $T$ and $T'$ are $\alpha$-adapted w.r.t $\nu$,
we get that $\Phi\circ R$ is
$\sigma_{-s}^{\Phi}\circ\sigma_s^{\Phi'}$-invariant and, so
$\sigma_t^{\Phi\circ R}$ and
$\sigma_{-s}^{\Phi}\circ\sigma_s^{\Phi'}$ commute each other. But
$\sigma^{\Phi\circ R}$ and $\sigma^{\Phi}$ commute each other
that's why $\sigma^{\Phi\circ R}$ and $\sigma^{\Phi'}$ also
commute each other. For all $s,t\in\mathbb{R}$, we have:
$$\Gamma(\sigma_t^{\Phi'}(\delta^{is})\delta^{-is})=
(\tau_t\surl{\ _{\beta} \star_{\alpha}}_{\
N}\sigma_t^{\Phi'})(\Gamma(\delta^{is}))(\delta^{-is}\surl{\
_{\beta} \otimes_{\alpha}}_{\ N}\delta^{-is})=1\surl{\ _{\beta}
\otimes_{\alpha}}_{\ N}\sigma_t^{\Phi'}(\delta^{is})\delta^{-is}$$
Consequently $\sigma_t^{\Phi'}(\delta^{is})\delta^{-is}$ belongs
to $\beta(N)$. Since $T'$ is $\alpha$-adapted w.r.t $\nu$,
$\beta(N)$ is $\sigma^{T'}$-invariant and $M\cap\beta(N)'$ is
$\sigma^{\Phi'}$-invariant so that, in fact,
$\sigma_t^{\Phi'}(\delta^{is})\delta^{-is}$ belongs to
$\beta(Z(N))$ and we easily get that there exists a strictly
positive operator $k$ affiliated with $Z(N)$ such that
$\sigma_t^{\Phi'}(\delta^{is})=\beta(k^{ist})\delta^{is}$. Then,
we have:
$$
\begin{aligned}
\sigma_s^{\Phi'}\circ\sigma_t^{\Phi}(m)
=\sigma_s^{\Phi'}(\delta^{-it}\sigma_t^{\Phi\circ
R}(m)\delta^{it})&=\beta(k^{-ist})\delta^{-it}\sigma_s^{\Phi'}\circ\sigma_t^{\Phi\circ
R}(m)\delta^{it}\beta(k^{ist})\\
&=\beta(k^{-ist})\sigma_t^{\Phi}\circ\sigma_s^{\Phi'}(m)\beta(k^{ist})
\end{aligned}$$ Since $T_L$ is $\beta$-adapted w.r.t $\nu$,
$\beta(k)$ is affiliated to the centralizer of $\sigma^T$. Apply
$\Phi$ to the previous formula and get:
$$\begin{aligned}
\Phi\circ\sigma_s^{\Phi'}\circ\sigma_s^{\Phi}(m^*m)
&=\Phi(\beta(k^{-ist})\sigma_t^{\Phi}\circ\sigma_s^{\Phi'}(m^*m)\beta(k^{ist}))\\
&=\Phi(\sigma_t^{\Phi}\circ\sigma_s^{\Phi'}(m^*m))=\Phi\circ\sigma_s^{\Phi'}(m^*m)
\end{aligned}$$
So, by \ref{struc} and left invariance $\sigma_{-s}^{\Phi'}\circ
T_L\circ\sigma_s^{\Phi'}$, there exists a strictly positive
operator $q_s$ affiliated with $Z(N)$ such that
$\Phi\circ\sigma_s^{\Phi'}=\Phi_{\beta(q_s)}$. By usual
arguments, we deduce that there exists a strictly positive $q$
affiliated to $Z(N)$ such that
$\Phi\circ\sigma_s^{\Phi'}=\Phi_{\beta(q^{-s})}$ and
$[D\Phi\circ\sigma_s^{\Phi'}:D\Phi]_s=\beta(q^{-s})$. Then, again
by \ref{struc}, there exists a strictly positive operator $h$
affiliated to $Z(N)$ such that $\Phi=\Phi_{\beta(h)}$ avec
$[DT':DT_L]_t=\beta(h^{it})$.
\end{proof}

Also, we have a similar result for right invariant
operator-valued weight.

\begin{coro}
There exists a strictly positive operator $h$ affiliated with
$Z(N)$ such that:
$$T_R=(R\circ T_L\circ R)_{\alpha(h)}$$
\end{coro}

\begin{proof}
$T_R$ and $R\circ T_L\circ R$ satisfy hypothesis of the previous
theorem.
\end{proof}

\subsection{Manageability of the fundamental unitary}

In this section, we prove that the fundamental unitary satisfies
a proposition similar to Woronowicz's manageability of \cite{W}.

\begin{lemm}
There exists a strictly positive operator $P$ on $H_{\Phi}$ such
that, for all $x\in {\mathcal N}_{\Phi}$ and  $t\in\mathbb{R}$, we
have
$P^{it}\Lambda_{\Phi}(x)=\lambda^{\frac{t}{2}}\Lambda_{\Phi}(\tau_t(x))$.
\end{lemm}

\begin{proof}
Since $\Phi\circ R=\Phi_{\delta}$, by \cite{Vae} (5.3), we have:
$$\Lambda_{\Phi}(\sigma_t^{\Phi\circ
R}(x))=\delta^{it}J_{\Phi}
\lambda^{\frac{t}{2}}\delta^{it}J_{\Phi}\Delta^{it}_{\Phi}\Lambda_{\Phi}(x)$$
and since $\lambda$ is affiliated with $Z(M)$, we get
$||\Lambda_{\Phi}(\sigma_t^{\Phi\circ
R}(x))||=||\lambda^{\frac{t}{2}}\Lambda_{\Phi}(x)||$ for all
$x\in {\mathcal N}_{\Phi}\cap {\mathcal N}_{T_L}$ and
$t\in\mathbb{R}$. But, we know that $\Phi$ is $\sigma_t^{\Phi\circ
R}\circ\tau_t$-invariant, so
$||\Lambda_{\Phi}(x)||=||\lambda^{\frac{t}{2}}\Lambda_{\Phi}(\tau_t(x))||$.
Then, there exists $P_t$ on $H_{\Phi}$ such that:
$$P_t\Lambda_{\Phi}(x)=\lambda^{\frac{t}{2}}\Lambda_{\Phi}(\tau_t(x))$$
for all $x\in {\mathcal N}_{\Phi}\cap {\mathcal N}_{T_L}$ and
$t\in\mathbb{R}$. For all $s,t\in\mathbb{R}$, we verify that
$P_sP_t=P_{st}$ thanks to relation $\tau_t(\lambda)=\lambda$ and
the existence of $P$ follows.
\end{proof}

\begin{defi}
We call \textbf{manageable operator} the strictly positive
operator $P$ on $H_{\Phi}$ such that
$P^{it}\Lambda_{\Phi}(x)=\lambda^{\frac{t}{2}}\Lambda_{\Phi}(\tau_t(x))$,
for all $x\in {\mathcal N}_{\Phi}$ and  $t\in\mathbb{R}$.
\end{defi}

\begin{prop}
For all $m\in M$, $n\in N$ and $t\in\mathbb{R}$, we have:
$$
\begin{aligned}
&P^{it}mP^{-it}=\tau_t(m)
&\quad P^{it}\alpha(n)P^{-it}=\alpha(\sigma^{\nu}_t(n))\\
&P^{it}\beta(n)P^{-it}=\beta(\sigma^{\nu}_t(n)) &\quad
P^{it}\hat{\beta}(n)P^{-it}=\hat{\beta}(\sigma^{\nu}_t(n))
\end{aligned}$$
\end{prop}

\begin{proof}
Straightforward.
\end{proof}

\noindent Then, we can define operators $P^{it}\surl{\
_{\beta}\otimes_{\alpha}}_{\ \nu}P^{it}$ on $H_{\Phi}\surl{\
_{\beta} \otimes_{\alpha}}_{\ \nu}H_{\Phi}$ and $P^{it}\surl{\
_{\alpha}\otimes_{\hat{\beta}}}_{\ \nu^o}P^{it}$ on
$H_{\Phi}\surl{\ _{\alpha}\otimes_{\hat{\beta}}}_{\
\nu^o}H_{\Phi}$ for all $t\in\mathbb{R}$.

\begin{theo}
$W$ satisfies a manageability relation. More exactly, we have:
$$\begin{aligned}
&\ \quad (\sigma_{\nu}W^*\sigma_{\nu}(q\surl{\ _{\hat{\beta}}
\otimes_{\alpha}}_{\ \nu}v)|p\surl{\ _{\alpha}\otimes_{\beta}}_{\
\nu^o}w)\\
&=(\sigma_{\nu^o}W\sigma_{\nu^o}(J_{\Phi}p\surl{\
_{\alpha}\otimes_{\beta}}_{\ \nu^o}P^{-1/2}v)| J_{\Phi}q\surl{\
_{\hat{\beta}} \otimes_{\alpha}}_{\ \nu}P^{1/2}w)
\end{aligned}$$
for all $v\in {\mathcal D}(P^{-\frac{1}{2}})$, $w\in {\mathcal
D}(P^{\frac{1}{2}})$ and $p,q\in D(_{\alpha}H_{\Phi},\nu)\cap
D((H_{\Phi})_{\hat{\beta}},\nu^o)$. Moreover, we have
$W(P^{it}\surl{\ _{\beta}\otimes_{\alpha}}_{\
\nu}P^{it})=(P^{it}\surl{\ _{\alpha}\otimes_{\hat{\beta}}}_{\
\nu^o}P^{it})W$ for all $t\in\mathbb{R}$. \label{mania}
\end{theo}

\begin{proof}
Let $p,q\in D(_{\alpha}H_{\Phi},\nu)\cap
D((H_{\Phi})_{\hat{\beta}},\nu^o)$. For all $v\in {\mathcal
D}(D^{1/2})$ and $w\in {\mathcal D}(D^{-1/2})$, we know that:
$$(I(id*\omega_{q,p})(W)Iv|w)=((id*\omega_{p,q})(W)D^{1/2}v|D^{-1/2}w)$$

Since $(id*\omega_{p,q})(W)\in {\mathcal D}(S)={\mathcal
D}(\tau_{-i/2})$ and $\tau$ is implemented by $D^{-1}$, we have
$\tau_{-i/2}((id*\omega_{p,q})(W))=I(id*\omega_{q,p})(W)I$. But
$\tau$ is also implemented by $P$, so that:
$$(I(id*\omega_{q,p})(W)Iv|w)=((id*\omega_{p,q})(W)P^{1/2}v|P^{-1/2}w)$$
for all $v\in {\mathcal D}(P^{1/2})$ and $w\in {\mathcal
D}(P^{-1/2})$. By \ref{besoin}, we rewrite the formula:
$$\begin{aligned}
&\ \quad (\sigma_{\nu}W^*\sigma_{\nu}(q\surl{\ _{\hat{\beta}}
\otimes_{\alpha}}_{\ \nu}v)|p\surl{\ _{\alpha}\otimes_{\beta}}_{\
\nu^o}w)\\
&=(\sigma_{\nu^o}W\sigma_{\nu^o}(J_{\Phi}p\surl{\
_{\alpha}\otimes_{\beta}}_{\ \nu^o}P^{-1/2}v)| J_{\Phi}q\surl{\
_{\hat{\beta}} \otimes_{\alpha}}_{\ \nu}P^{1/2}w)
\end{aligned}$$
Now, we have to prove $W^*(P^{it}\surl{\
_{\alpha}\otimes_{\hat{\beta}}}_{\ \nu^o}P^{it})= (P^{it}\surl{\
_{\beta}\otimes_{\alpha}}_{\ \nu}P^{it})W^*$ for all
$t\in\mathbb{R}$. First of all, because of the commutation
relation between $P$ and $\beta$, $D((H_{\Phi})_{\beta},\nu^o)$
is $P^{it}$-invariant and if $(\xi_i)_{i\in I}$ is a
$(N^o,\nu^o)$-basis of $(H_{\Phi})_{\beta}$, then
$(P^{it}\xi_i)_{i\in I}$ is also. Let $v\in
D((H_{\Phi})_{\beta},\nu^o)$ and $a\in {\mathcal N}_{T_L}\cap
{\mathcal N}_{\Phi}$. We compute:

$$
\begin{aligned}
&\ \quad(P^{it}\surl{\ _{\beta}\otimes_{\alpha}}_{\ \nu}P^{it})
W^*(v\surl{\ _{\alpha}\otimes_{\hat{\beta}}}_{\
\nu^o}\Lambda_{\Phi}(a))\\
&=\sum_{i\in I}P^{it}\xi_i\surl{\ _{\beta} \otimes_{\alpha}}_{\
\nu}\lambda^{t/2}\Lambda_{\Phi}(\tau_t((\omega_{v,\xi_i}
\surl{\ _{\beta} \star_{\alpha}}_{\ \nu}id)(\Gamma(a))))\\
&=\sum_{i\in I}P^{it}\xi_i\surl{\ _{\beta} \otimes_{\alpha}}_{\
\nu}\Lambda_{\Phi}((\omega_{P^{it}v,P^{it}\xi_i}
\surl{\ _{\beta} \star_{\alpha}}_{\ \nu}id)(\Gamma(\lambda^{t/2}\tau_t(a))))\\
&=W^*(P^{it}v\surl{\ _{\alpha}\otimes_{\hat{\beta}}}_{\
\nu^o}\lambda^{t/2} \Lambda_{\Phi}(\tau_t(a)))=W^*(P^{it}\surl{\
_{\alpha}\otimes_{\hat{\beta}}}_{\ \nu^o}P^{it})(v\surl{\
_{\alpha}\otimes_{\hat{\beta}}}_{\ \nu^o}\Lambda_{\Phi}(a))
\end{aligned}$$

\end{proof}

Following \cite{E2} (definition 4.1), we define the notion of
weakly regular pseudo-multiplicative unitary.

\begin{defi}
A pseudo-multiplicative unitary ${\mathcal W}$ w.r.t
$\alpha,\beta,\hat{\beta}$ is said to be \textbf{weakly regular}
if the weakly closed linear span of
$(\lambda_v^{\alpha,\beta})^*\mathcal{W}\rho_w^{\hat{\beta},\alpha}$
where $v,w$ belongs to $D(_{\alpha}H,\nu)$ is equal to
$\alpha(N)'$.
\end{defi}

\begin{prop}\label{reg}
The operator $\widehat{W}=\sigma_{\nu}W^*\sigma_{\nu}$ from
$H_{\Phi}\surl{\ _{\hat{\beta}}\otimes_{\alpha}}_{\ \nu}H_{\Phi}$
to $H_{\Phi}\surl{\ _{\alpha}\otimes_{\beta}}_{\ \nu^o}H_{\Phi}$
is a pseudo-multiplicative unitary over $N$ w.r.t
$\alpha,\beta,\hat{\beta}$ which is weakly regular in the sense of
\cite{E2} (definition 4.1).
\end{prop}

\begin{proof}
By \cite{EV}, we know that $\widehat{W}$ is a
pseudo-multiplicative unitary. We also know that
$<(\lambda_v^{\alpha,\beta})^*\widehat{W}\rho_w^{\hat{\beta},\alpha}>
^{-\textsc{w}}\subset\alpha(N)'$. For all $v\in {\mathcal
D}(P^{-\frac{1}{2}})$, $w\in {\mathcal D}(P^{\frac{1}{2}})$ and
$p,q\in D(_{\alpha}H_{\Phi},\nu)\cap
D((H_{\Phi})_{\hat{\beta}},\nu^o)$, we have, by theorem
\ref{mania}:
$$((\lambda_p^{\alpha,\beta})^*\widehat{W}\rho_v^{\hat{\beta},\alpha}q|w)
=(\sigma_{\nu^o}W\sigma_{\nu^o}(J_{\Phi}p\surl{\
_{\alpha}\otimes_{\beta}}_{\ \nu^o}P^{-1/2}v)| J_{\Phi}q\surl{\
_{\hat{\beta}} \otimes_{\alpha}}_{\ \nu}P^{1/2}w)$$ and on the
other hand:
$$
\begin{aligned}
(R^{\alpha,\nu}(v)R^{\alpha,\nu}(p)^*q|w)
&=(R^{\alpha,\nu}(v)J_{\nu}R^{\hat{\beta},\nu^o}(J_{\Phi}p)^*J_{\Phi}q|w)\\
&=(R^{\alpha,\nu}(v)J_{\nu}\Lambda_{\nu}(<J_{\Phi}q,J_{\Phi}p>_{\hat{\beta},\nu^o_L})|w)\\
&=(P^{-1/2}R^{\alpha,\nu}(v)J_{\nu}\Lambda_{\nu}(<J_{\Phi}q,J_{\Phi}p>_{\hat{\beta},\nu^o_L})|P^{1/2}w)\\
&=(R^{\alpha,\nu}(P^{-1/2}v)\Delta_{\nu}^{-1/2}J_{\nu}\Lambda_{\nu}(<J_{\Phi}q,J_{\Phi}p>_{\hat{\beta},\nu^o_L})|P^{1/2}w)\\
&=(R^{\alpha,\nu}(P^{-1/2}v)\Lambda_{\nu}(<J_{\Phi}p,J_{\Phi}q>_{\hat{\beta},\nu^o_L})|P^{1/2}w)\\
&=(\alpha(<J_{\Phi}p,J_{\Phi}q>_{\hat{\beta},\nu^o_L})P^{-1/2}v|P^{1/2}w)\\
&=(J_{\Phi}p\surl{\ _{\hat{\beta}} \otimes_{\alpha}}_{\
\nu}P^{-1/2}v| J_{\Phi}q\surl{\ _{\hat{\beta}}
\otimes_{\alpha}}_{\ \nu}P^{1/2}w)
\end{aligned}$$
There exists $\Xi\in H_{\Phi}\surl{\
_{\hat{\beta}}\otimes_{\alpha}}_{\ \nu}H_{\Phi}$ such that
$\sigma_{\nu^o}W\sigma_{\nu^o}\Xi=J_{\Phi}p\surl{\ _{\hat{\beta}}
\otimes_{\alpha}}_{\ \nu}P^{-1/2}v$ since $W$ is onto. By
definition, there exists a net
$(\sum_{k=1}^{n(i)}J_{\Phi}p_k^i\surl{\
_{\alpha}\otimes_{\beta}}_{\ \nu^o}P^{-1/2}v_k^i)_{i\in I}$ which
converges to $\Xi$. So, the net
$((\sum_{k=1}^{n(i)}(\lambda_{p_k^i}^{\alpha,\beta})^*
\widehat{W}\rho_{v_k^i}^{\hat{\beta},\alpha}q|w))_{i\in I}$
converges to:
$$
\begin{aligned}
(\sigma_{\nu^o}W\sigma_{\nu^o}\Xi|J_{\Phi}q\surl{\ _{\hat{\beta}}
\otimes_{\alpha}}_{\ \nu}P^{1/2}w)&=(J_{\Phi}p\surl{\
_{\hat{\beta}} \otimes_{\alpha}}_{\ \nu}P^{-1/2}v|
J_{\Phi}q\surl{\ _{\hat{\beta}} \otimes_{\alpha}}_{\
\nu}P^{1/2}w)\\
&=(R^{\alpha,\nu}(v)R^{\alpha,\nu}(p)^*q|w)
\end{aligned}$$

Then, we obtain $\alpha(N)'=
<R^{\alpha,\nu}(v)R^{\alpha,\nu}(p)^*>^{-\textsc{w}}\subset
<(\omega_{v,p}*id)(\widehat{W}\sigma_{\nu^o})>^{-\textsc{w}}$.
\end{proof}

\begin{coro}\label{alsjdu}
If $\hat{M}$ denote the weak closed linear span of
$(\omega_{\xi,\eta}*id)(W)$ where $\xi\in
D((H_{\Phi})_{\beta},\nu^o)$ and $\eta\in
D(_{\alpha}H_{\Phi},\nu)$, then $\hat{M}$ is a von Neumann
algebra.
\end{coro}

\begin{proof}
Comes from weak regularity of $\hat{W}$ and \cite{E2} (proposition
3.2).
\end{proof}

\subsection{Changing the quasi-invariant weight}

Let $\nu'$ be a n.s.f weight on $N$ such that there exist
strictly positive operator $h$ and $k$ affiliated with $N$
strongly commuting and $[D\nu':D\nu]_t=k^{\frac{it^2}{2}}h^{it}$
for all $t\in\mathbb{R}$. By \cite{Vae} (proposition 5.1), it is
equivalent to $\sigma_t^{\nu}(h^{is})=k^{ist}h^{is}$ for all
$s,t\in\mathbb{R}$ and $\nu'=\nu_h$ in the sense of \cite{Vae}.
This hypothesis is satisfied, in particular, if $\sigma^{\nu}$ and
$\sigma^{\nu'}$ commute each other. In this cas, $k$ is
affiliated with $Z(N)$.

\begin{prop}
There exists a n.s.f operator-valued weight $T_L'$ from $M$ to
$\alpha(N)$ which is $\beta$-adapted w.r.t $\nu'$ such that, for
all $t\in\mathbb{R}$, we have:
$$[DT_L':DT_L]_t=\beta(k^{\frac{-it^2}{2}}h^{it})$$
\end{prop}

\begin{proof}
By \ref{timpo}, there exists a n.s.f operator-valued weight $S_L$
from $M$ to $\beta(N)$ such that $\nu\circ\alpha^{-1}\circ
T_L=\nu\circ\beta^{-1}\circ S_L$ so that $S_L$ is $\alpha$-adapted
w.r.t $\nu$. Then, again by \ref{timpo},there exists a n.s.f
operator-valued weight $T_L'$ from $M$ to $\alpha(N)$ such that
$\nu'\circ\beta^{-1}\circ S=\nu\circ\alpha^{-1}\circ T_L'$ so
that $T_L'$ is $\beta$-adapted w.r.t $\nu'$. Then, we compute the
Radon-Nikodym cocycle for all $t\in\mathbb{R}$:
$$
\begin{aligned}
\ [DT_L':DT_L]_t &=[D\nu\circ\alpha^{-1}\circ
T_L':D\nu\circ\alpha^{-1}\circ
T_L]_t\\
&=[D\nu'\circ\beta^{-1}\circ
S:D\nu\circ\beta^{-1}\circ S]_t\\
&=\beta([D\nu':D\nu]_{-t}^*)=\beta(k^{\frac{-it^2}{2}}h^{it})
\end{aligned}$$

\end{proof}

\begin{coro}
We have:
$$\nu\circ\alpha^{-1}\circ
T_L'=(\nu\circ\alpha^{-1}\circ T_L)_{\beta(h)}\quad\text{ and
}\quad\nu'\circ\alpha^{-1}\circ T_L'=(\nu\circ\alpha^{-1}\circ
T_L)_{\alpha(h)\beta(h)}$$
\end{coro}

\begin{proof}
Come from \cite{Vae} (proposition 5.1) and the following equality,
for all $t\in\mathbb{R}$, $[D\nu'\circ\alpha^{-1}\circ
T_L':D\nu\circ\alpha^{-1}\circ
T_L]_t=\alpha(k^{\frac{it^2}{2}})\beta(k^{\frac{-it^2}{2}})\alpha(h^{it})\beta(h^{it})$.
\end{proof}

\begin{prop}
$T_L'$ is left invariant.
\end{prop}

\begin{proof}
Let $a\in {\mathcal M}_{T_L'}^+$. By left invariance of $T_L$, we
have:
$$
\begin{aligned}
(id\surl{\ _{\beta} \star_{\alpha}}_{\
\nu'}\nu'\circ\alpha^{-1}\circ T_L')(\Gamma(a))&=(id\surl{\
_{\beta} \star_{\alpha}}_{\ \nu}\nu\circ\alpha^{-1}\circ
T_L')(\Gamma(a))\\
&=(id\surl{\ _{\beta} \star_{\alpha}}_{\
\nu}(\nu\circ\alpha^{-1}\circ
T_L)_{\beta(h)})(\Gamma(a))\\
&=(id\surl{\ _{\beta} \star_{\alpha}}_{\
\nu}\nu\circ\alpha^{-1}\circ
T_L)(\Gamma(\beta(h^{1/2})a\beta(h^{1/2})))\\
&=T_L(\beta(h^{1/2})a\beta(h^{1/2}))=T'(a)
\end{aligned}$$

\end{proof}

We state the right version of these results:

\begin{prop}
There exists a n.s.f right invariant operator-valued weight $T_R'$
which is $\alpha$-adapted w.r.t $\nu'$ such that, for all
$t\in\mathbb{R}$, we have:
$$[DT_R':DT_R]_t=\alpha(k^{\frac{it^2}{2}}h^{it})$$
Moreover, we have:
$$\nu\circ\beta^{-1}\circ
T_R'=(\nu\circ\beta^{-1}\circ T_R)_{\alpha(h)}\quad\text{ and
}\quad\nu'\circ\beta^{-1}\circ T_R'=(\nu\circ\beta^{-1}\circ
T_R)_{\alpha(h)\beta(h)}$$
\end{prop}

\begin{lemm}
The application $I_{\nu}^{\nu'}$ defined by the following formula:
$$I_{\nu}^{\nu'}(\xi\surl{\ _{\beta} \otimes_{\alpha}}_{\
\nu}\eta)=\beta(k^{-i/8})\xi\surl{\ _{\beta} \otimes_{\alpha}}_{\
\nu'}\alpha(h^{1/2})\eta$$ for all $\xi\in H$ and  $\eta\in
D(_{\alpha}H,\nu)\cap {\mathcal D}(\alpha(h^{1/2}))$, is an
isomorphism of $\beta(N)'-\alpha(N)'^o$-bimodules from $H\surl{\
_{\beta} \otimes_{\alpha}}_{\ \nu}H$ onto $H\surl{\ _{\beta}
\otimes_{\alpha}}_{\ \nu'}H$.
\end{lemm}

\begin{proof}
For all $x\in {\mathcal N}_{\nu'}$, we have:
$$\alpha(x)\alpha(h^{1/2})\eta=\alpha(xh^{1/2})\eta
=R^{\alpha,\nu}(\eta)\Lambda_{\nu}(xh^{1/2})=R^{\alpha,\nu}(\eta)\Lambda_{\nu'}(x)$$
so that $\alpha(h^{1/2})\eta\in D(_{\alpha}H,\nu)$ and
$R^{\alpha,\nu'}(\alpha(h^{1/2})\eta)=R^{\alpha,\nu}(\eta)$.
Also, we recall that
$J_{\nu'}=J_{\nu}k^{-i/8}J_{\nu}k^{i/8}J_{\nu}$ by \cite{Vae}
(proposition 2.5). Then, we have:
$$
\begin{aligned}
&\quad\ (\beta(k^{-i/8})\xi_1\surl{\ _{\beta} \otimes_{\alpha}}_{\
\nu}\alpha(h^{1/2})\eta_1|\beta(k^{-i/8})\xi_2\surl{\ _{\beta}
\otimes_{\alpha}}_{\ \nu}\alpha(h^{1/2})\eta_2)\\
&=(\beta(J_{\nu'}<\alpha(h^{1/2})\eta_1,\alpha(h^{1/2})\eta_2>_{\alpha,\nu'}^*J_{\nu'})
\beta(k^{-i/8})\xi_1|\beta(k^{-i/8})\xi_2)\\
&=(\beta(k^{-i/8}J_{\nu}k^{-i/8}J_{\nu}k^{i/8}J_{\nu}<\eta_1,\eta_2>_{\alpha,\nu}^*
J_{\nu}k^{-i/8}J_{\nu}k^{i/8}J_{\nu}k^{i/8}) \xi_1|\xi_2)\\
&=(\beta(J_{\nu}<\eta_1,\eta_2>_{\alpha,\nu}^*J_{\nu})
\xi_1|\xi_2)=(\xi_1\surl{\ _{\beta} \otimes_{\alpha}}_{\
\nu}\eta_1|\xi_2\surl{\ _{\beta} \otimes_{\alpha}}_{\ \nu}\eta_2)
\end{aligned}$$
\end{proof}

\begin{rema}
For all $\xi\in D(H_{\beta},\nu^o)$ and $\eta\in
D(_{\alpha}H,\nu)$, we have:
$$
\begin{aligned}
I_{\nu}^{\nu'}(\xi\surl{\ _{\beta} \otimes_{\alpha}}_{\
\nu}\eta)&=\beta(k^{-i/8})\xi\surl{\ _{\beta} \otimes_{\alpha}}_{\
\nu'}\alpha(h^{1/2})\eta=\xi\surl{\ _{\beta} \otimes_{\alpha}}_{\
\nu'}\alpha(k^{-i/8}h^{1/2})\eta\\
&=\beta(\sigma_{i/2}^{\nu}(h^{1/2}))\xi\surl{\
_{\beta}\otimes_{\alpha}}_{\
\nu'}\alpha(k^{-i/8})\eta=\beta(\sigma_{i/2}^{\nu}(k^{-i/8}h^{1/2}))\xi\surl{\ _{\beta} \otimes_{\alpha}}_{\ \nu'}\eta\\
&=\beta(k^{i/8})\xi\surl{\ _{\beta} \otimes_{\alpha}}_{\
\nu'}\alpha(\sigma_{-i/2}^{\nu}(h^{1/2}))\eta=\beta(k^{i/8}h^{1/2})\xi\surl{\
_{\beta} \otimes_{\alpha}}_{\ \nu'}\eta
\end{aligned}$$
\end{rema}

\begin{prop}
Let $(N,M,\alpha,\beta,\Gamma,\nu,T_L,T_R)$ be a measured quantum
groupoid. There exists a measured quantum groupoid
$(N,M,\alpha,\beta,\Gamma,\nu',T_L',T_R')$ fundamental objects of
which, $R'$, $\tau'$, $\lambda'$, $\delta'$ and $P'$, are
expressed, for all $t\in\mathbb{R}$, in the following way:
\begin{center}
\begin{minipage}{14cm}
\begin{enumerate}[i)]
\item $R'=R$, $\lambda'=\lambda$ and $\delta'=\delta$;
\item
$\tau_t'=Ad_{\alpha(k^{\frac{-it^2}{2}}h^{-it})\beta(k^{\frac{it^2}{2}}h^{it})}\circ\tau_t
=Ad_{\alpha([D\nu':D\nu]_t^*)\beta([D\nu':D\nu]_t)}\circ\tau_t$
\item $P'^{it}=\alpha(k^{\frac{it^2}{2}}h^{it})\beta(k^{\frac{-it^2}{2}}h^{-it})
J_{\Phi}\alpha(k^{\frac{it^2}{2}}h^{it})\beta(k^{\frac{-it^2}{2}}h^{-it})J_{\Phi}P^{it}$
\end{enumerate}
\end{minipage}
\end{center}
\end{prop}

\begin{proof}
The existence of $(N,M,\alpha,\beta,\Gamma,\nu',T_L',T_R')$ has
been already proved. We put $\Phi'=\nu'\circ\alpha^{-1}\circ
T_L'$ and $\Psi'=\nu'\circ\beta^{-1}\circ T_R'$. Let $x,y\in
{\mathcal N}_{T_R'}\cap {\mathcal N}_{\Psi'}$. By \cite{Vae}
(proposition 2.5), we have:
$$
\begin{aligned}
J_{\Psi'}\Lambda_{\Psi'}(x)&=J_{\Psi}\alpha(k^{-i/8})\beta(k^{i/8})
J_{\Psi}\alpha(k^{i/8})\beta(k^{-i/8})J_{\Psi}\Lambda_{\Psi}(x\alpha(h^{1/2})\beta(h^{1/2}))\\
\omega_{J_{\Psi'}\Lambda_{\Psi'}(x)}
&=\omega_{\alpha(k^{i/8})\beta(k^{-i/8})J_{\Psi}\Lambda_{\Psi}(x\alpha(h^{1/2})\beta(h^{1/2}))}
\end{aligned}$$
Then, we easily verify
$$\lambda^{\beta,\alpha,\nu'}
_{\alpha(k^{i/8})\beta(k^{-i/8})J_{\Psi}\Lambda_{\Psi}(x\alpha(h^{1/2})\beta(h^{1/2}))}=
I_{\nu}^{\nu'}\lambda^{\beta,\alpha,\nu}_
{\alpha(k^{i/8})J_{\Psi}\Lambda_{\Psi}(x\alpha(h^{1/2}))}$$ We
compute:
$$
\begin{aligned}
&\ \quad(\omega_{J_{\Psi'}\Lambda_{\Psi'}(x)}\surl{\ _{\beta}
\star_{\alpha}}_{\ \nu'}id)(\Gamma(y^*y))\\
&=(\omega_{\alpha(k^{i/8})\beta(k^{-i/8})J_{\Psi}\Lambda_{\Psi}(x\alpha(h^{1/2})\beta(h^{1/2}))}\surl{\
_{\beta} \star_{\alpha}}_{\ \nu'}id)(\Gamma(y^*y))\\
&=(\omega_{\alpha(k^{i/8})J_{\Psi}\Lambda_{\Psi}(x\alpha(h^{1/2}))}\surl{\
_{\beta} \star_{\alpha}}_{\
\nu}id)(\Gamma(y^*y))\\
&=(\omega_{J_{\Psi}\Lambda_{\Psi}(x\alpha(k^{-i/8}h^{1/2}))}\surl{\
_{\beta} \star_{\alpha}}_{\ \nu}id)(\Gamma(y^*y))
\end{aligned}$$
Apply $R$ to get:
$$
\begin{aligned}
&\ \quad R[(\omega_{J_{\Psi'}\Lambda_{\Psi'}(x)}\surl{\ _{\beta}
\star_{\alpha}}_{\ \nu'}id)(\Gamma(y^*y))]\\
&=(\omega_{J_{\Psi}\Lambda_{\Psi}(y)}\surl{\
_{\beta} \star_{\alpha}}_{\ \nu}\!id)(\Gamma(\alpha(k^{i/8}h^{1/2})x^*x\alpha(k^{-i/8}h^{1/2})))\\
&=(\omega_{\alpha(k^{-i/8}h^{1/2})J_{\Psi}\Lambda_{\Psi}(y)}\surl{\
_{\beta} \star_{\alpha}}_{\ \nu}id)(\Gamma(x^*x))\\
&=(\omega_{\alpha(k^{i/8})J_{\Psi}\Lambda_{\Psi}(y\alpha(h^{1/2}))}\surl{\
_{\beta} \star_{\alpha}}_{\ \nu}id)(\Gamma(x^*x))\\
&=(\omega_{J_{\Psi'}\Lambda_{\Psi'}(y)}\surl{\ _{\beta}
\star_{\alpha}}_{\
\nu'}id)(\Gamma(x^*x))=R'[(\omega_{J_{\Psi'}\Lambda_{\Psi'}(x)}\surl{\
_{\beta}
\star_{\alpha}}_{\ \nu'}id)(\Gamma(y^*y))]\\
\end{aligned}$$
so that $R=R'$. For all $a\in M$, $\xi\in D(H_{\beta},\nu'^o)$ and
$t\in\mathbb{R}$, we have:
$$
\begin{aligned}
&\quad\ \tau_t((\omega_{\xi}\surl{\ _{\beta} \star_{\alpha}}_{\
\nu'}id)(\Gamma(a)))\\
&=\tau_t(\alpha(k^{-i/8}h^{-1/2})(\omega_{\xi}\surl{\ _{\beta}
\star_{\alpha}}_{\
\nu}id)(\Gamma(a))\alpha(k^{i/8}h^{-1/2}))\\
&=\alpha(\sigma_t^{\nu}(k^{-i/8}h^{-1/2}))\tau_t((\omega_{\xi}\surl{\
_{\beta} \star_{\alpha}}_{\
\nu}id)(\Gamma(a)))\alpha(\sigma_t^{\nu}(k^{i/8}h^{-1/2}))\\
&=\alpha(k^{-t/2-i/8}h^{-1/2})(\omega_{\Delta_{\Psi}^{-it}\xi}\surl{\
_{\beta} \star_{\alpha}}_{\
\nu}id)(\Gamma(\sigma_{-t}^{\Psi}(a)))\alpha(k^{-t/2+i/8}h^{-1/2})
\end{aligned}$$
By \cite{Vae} (proposition 2.4 and corollaire 2.6), we know that:
$$
\begin{aligned}
&\ \quad(\omega_{\Delta_{\Psi}^{-it}\xi}\surl{\ _{\beta}
\star_{\alpha}}_{\ \nu}id)(\Gamma(\sigma_{-t}^{\Psi}(a)))\\
&=(\omega_{\alpha(k^{\frac{-it^2}{2}}h^{it})
\beta(k^{\frac{it^2}{2}}h^{it})\Delta_{\Psi'}^{-it}\xi}\surl{\
_{\beta} \star_{\alpha}}_{\
\nu}id)(\Gamma(Ad_{\alpha(k^{\frac{it^2}{2}}h^{-it})\beta(k^{\frac{-it^2}{2}}h^{-it})}\circ
\sigma_{-t}^{\Psi'}(a)))
\end{aligned}$$ so that:
$$
\begin{aligned}
&\quad\ \tau_t((\omega_{\xi}\surl{\ _{\beta} \star_{\alpha}}_{\
\nu'}id)(\Gamma(a)))\\
&=Ad_{\alpha(k^{-t/2+i/8}h^{-1/2})\beta(k^{\frac{it^2}{2}}h^{it})}\circ
(\omega_{\beta(k^{\frac{it^2}{2}}h^{it})\Delta_{\Psi'}^{-it}\xi}\surl{\
_{\beta} \star_{\alpha}}_{\
\nu}id)(\Gamma(\sigma_{-t}^{\Psi'}(a)))\\
&=\alpha(k^{\frac{-it^2}{2}}h^{-it})\beta(k^{\frac{it^2}{2}}h^{it})
(\omega_{\Delta_{\Psi'}^{-it}\xi}\surl{\ _{\beta}
\star_{\alpha}}_{\ \nu'}id)(\Gamma(\sigma_{-t}^{\Psi'}(a)))
\alpha(k^{\frac{it^2}{2}}h^{it})\beta(k^{\frac{-it^2}{2}}h^{-it})\\
&=\alpha(k^{\frac{-it^2}{2}}h^{-it})\beta(k^{\frac{it^2}{2}}h^{it})
\tau'_t((\omega_{\xi}\surl{\ _{\beta} \star_{\alpha}}_{\
\nu'}id)(\Gamma(a)))\alpha(k^{\frac{it^2}{2}}h^{it})\beta(k^{\frac{-it^2}{2}}h^{-it})
\end{aligned}$$
Consequently, we have:
$$\tau_t'(z)=\alpha(k^{\frac{it^2}{2}}h^{it})\beta(k^{\frac{-it^2}{2}}h^{-it})
\tau_t(z)\alpha(k^{\frac{-it^2}{2}}h^{-it})\beta(k^{\frac{it^2}{2}}h^{it})$$
for all $z\in M$ and  $t\in\mathbb{R}$. Now, we compute the
Radon-Nikodym cocycle:
$$
\begin{aligned}
&\quad\ [D\nu'\circ\alpha^{-1}\circ T'\circ
R:D\nu'\circ\alpha^{-1}\circ
T']_t\\
&=[D\nu'\alpha^{-1}T'R:D\nu\alpha^{-1}TR]_t[D\nu\alpha^{-1}
TR:D\nu\alpha^{-1}T]_t[D\nu\alpha^{-1}T:D\nu'\alpha^{-1}T']_t\\
&=\alpha([D\nu':D\nu]_t)\beta([D\nu':D\nu]_{-t}^*)\lambda^{\frac{it^2}{2}}\delta^{it}
\alpha([D\nu:D\nu']_t)\beta([D\nu:D\nu']_{-t}^*)
\end{aligned}$$
which is equal to $\lambda^{\frac{it^2}{2}}\delta^{it}$. Finally,
we express the manageable operator $P'$ in terms of $P$. We have,
for all $x\in {\mathcal N}_{T_L'}\cap {\mathcal N}_{\Phi'}$ and
$t\in\mathbb{R}$:
$$
\begin{aligned}
&\quad\ P'^{it}\Lambda_{\Phi'}(x)=\lambda'^{t/2}\Lambda_{\Phi'}(\tau'_t(x))\\
&=\lambda^{t/2}\Lambda_{\Phi}(\alpha(k^{\frac{it^2}{2}}h^{it})\beta(k^{\frac{-it^2}{2}}h^{-it})
\tau_t(x)\alpha(k^{\frac{-it^2}{2}}h^{-it})\beta(k^{\frac{it^2}{2}}h^{it})
\alpha(h^{1/2})\beta(h^{1/2}))
\end{aligned}$$
which is equal to the value of:
$$\lambda^{t/2}\alpha(k^{\frac{it^2}{2}}h^{it})\beta(k^{\frac{-it^2}{2}}h^{-it})
J_{\Phi}\alpha(k^{\frac{it^2}{2}}h^{it})\beta(k^{\frac{-it^2}{2}}h^{-it})
\alpha(k^{t/2})\beta(k^{t/2})J_{\Phi}$$ on
$\Lambda_{\Phi}(\tau_t(x)\alpha(h^{1/2})\beta(h^{1/2}))$ and the
value of:
$$\lambda^{t/2}\alpha(k^{\frac{it^2}{2}}h^{it})\beta(k^{\frac{-it^2}{2}}h^{-it})
J_{\Phi}\alpha(k^{\frac{it^2}{2}}h^{it})\beta(k^{\frac{-it^2}{2}}h^{-it})J_{\Phi}$$
on $\Lambda_{\Phi}(\tau_t(x\alpha(h^{1/2})\beta(h^{1/2})))$ which
is:
$$\alpha(k^{\frac{it^2}{2}}h^{it})\beta(k^{\frac{-it^2}{2}}h^{-it})
J_{\Phi}\alpha(k^{\frac{it^2}{2}}h^{it})\beta(k^{\frac{-it^2}{2}}h^{-it})J_{\Phi}
P^{it}\Lambda_{\Phi'}(x)$$
\end{proof}

Thanks to these formulas, we verify for example that
$\tau'_t(\alpha(n))=\alpha(\sigma_t^{\nu'}(n))$,
$\tau'_t(\beta(n))=\beta(\sigma_t^{\nu'}(n))$ and $\tau'$ is
implemented by $P'$.

\begin{prop}
Let $(N,M,\alpha,\beta,\Gamma,\nu,T_L,T_R)$ be measured quantum
groupoid and let $\tilde{T_L}$ be an other n.s.f left invariant
operator-valued weight which is $\beta$-adapted w.r.t $\nu$. Then
fundamental objects $\tilde{R}$, $\tilde{\tau}$,
$\tilde{\lambda}$, $\tilde{\delta}$ and $\tilde{P}$ of the
measured quantum groupoid
$(N,M,\alpha,\beta,\Gamma,\nu,\tilde{T_L},T_R)$ can be expressed
in the following way:
\begin{center}
\begin{minipage}{14cm}
\begin{enumerate}[i)]
\item $\tilde{R}=R$, $\tilde{\tau}=\tau$, $\tilde{\lambda}=\lambda$
and $\tilde{P}=P$
\item $\tilde{\delta}=\delta\alpha(h)\beta(h^{-1})$
where $h$ is affiliated with $Z(N)$ s.t.
$\tilde{T_L}=(T_L)_{\beta(h)}$
\end{enumerate}
\end{minipage}
\end{center}
\end{prop}

\begin{proof}
By uniqueness theorem, there exists a strictly positive operator
$h$ affiliated with $Z(N)$ such that
$\nu\circ\alpha^{-1}\circ\tilde{T_L}=(\nu\circ\alpha^{-1}\circ
T_L)_{\beta(h)}$ and, for all $t\in\mathbb{R}$, we have
$[D\tilde{T_L}:DT_L]_t=\beta(h^{it})$. We have already noticed
that $R$ and $\tau$ are independent w.r.t left invariant
operator-valued weight and $\beta$-adapted w.r.t $\nu$. We
compute then Radon-Nydodim cocycle:
$$
\begin{aligned}
&\ \quad
[D\nu\beta^{-1}R\tilde{T_L}R:D\nu\alpha^{-1}\tilde{T_L}]_t\\
&=[D\nu\beta^{-1}R\tilde{T_L}R:D\nu\beta^{-1}
RT_LR]_t[D\Psi:D\Phi]_t
[D\nu\alpha^{-1}T_L:D\nu\alpha^{-1}\tilde{T_L}]_t\\
&=R([D\tilde{T_L}:DT_L]^*_{-t})[D\Psi:D\Phi]_t
[DT_L:D\tilde{T_L}]_t\\
&=\alpha(h^{it})\lambda^{\frac{it^2}{2}}\delta^{it}\beta(h^{-it})
=\lambda^{\frac{it^2}{2}}\delta^{it}\alpha(h^{it})\beta(h^{-it})\\
\end{aligned}$$
Then, it remains to compute $\tilde{P}$. If, we put
$\tilde{\Phi}=\nu\circ\alpha^{-1}\circ\tilde{T_L}$, we have, for
all $t\in\mathbb{R}$ and  $x\in {\mathcal N}_{\tilde{T_L}}\cap
{\mathcal N}_{\tilde{\Phi}}$:
$$
\begin{aligned}
\tilde{P}^{it}\Lambda_{\tilde{\Phi}}(x)
=\tilde{\lambda}^{t/2}\Lambda_{\tilde{\Phi}}(\tilde{\tau}_t(x))
&=\lambda^{t/2}\Lambda_{\Phi}(\tau_t(x)\beta(h^{1/2}))
=\lambda^{t/2}\Lambda_{\Phi}(\tau_t(x\beta(h^{1/2}))\\
&=P^{it}\Lambda_{\Phi}(x\beta(h^{1/2}))\!=\!P^{it}\!\Lambda_{\tilde{\Phi}}(x)
\end{aligned}$$

\end{proof}

\begin{theo}
Let $(N,M,\alpha,\beta,\Gamma,\nu,T_L,T_R)$ and
$(N,M,\alpha,\beta,\Gamma,\nu',T_L',T_R')$ be measured quantum
groupoids such that there exist strictly positive operators $h$
and $k$ affiliated with $N$ which strongly commute and
$[D\nu':D\nu]_t=k^{\frac{it^2}{2}}h^{it}$ for all
$t\in\mathbb{R}$. For all $t\in\mathbb{R}$, fundamental objects
of the two structures are linked by:
\begin{center}
\begin{minipage}{14cm}
\begin{enumerate}[i)]
\item $R'=R$
\item
$\tau_t'=Ad_{\alpha(k^{\frac{-it^2}{2}}h^{-it})\beta(k^{\frac{it^2}{2}}h^{it})}\circ\tau_t
=Ad_{\alpha([D\nu':D\nu]_t^*)\beta([D\nu':D\nu]_t)}\circ\tau_t$
\item $\lambda'=\lambda$
\item $\dot{\delta'}=\dot{\delta}$ where $\dot{\delta}$ and $\dot{\delta}'$
have been defined in proposition \ref{deltadot}
\item $P'^{it}=\alpha(k^{\frac{it^2}{2}}h^{it})\beta(k^{\frac{-it^2}{2}}h^{-it})
J_{\Phi}\alpha(k^{\frac{it^2}{2}}h^{it})\beta(k^{\frac{-it^2}{2}}h^{-it})J_{\Phi}P^{it}$
\end{enumerate}
\end{minipage}
\end{center}
\end{theo}

\begin{proof}
We successively apply the two previous propositions.
\end{proof}

We state results of the section in the following theorems:

\begin{theo}
Let $(N,M,\alpha,\beta,\Gamma,\nu,T_L,T_R)$ be a measured quantum
groupoid. If $T'$ is a left invariant operator-valued weight
which is $\beta$-adapted w.r.t $\nu$, then there exists a strictly
positive operator $h$ affiliated with $Z(N)$ such that, for all
$t\in\mathbb{R}$:
$$\nu\circ\alpha^{-1}\circ T'=(\nu\circ\alpha^{-1}\circ
T_L)_{\beta(h)}$$ We have a similar result for the right
invariant operator-valued weights.
\end{theo}

\begin{theo}
Let $(N,M,\alpha,\beta,\Gamma,\nu,T_L,R\circ T_L\circ R)$ be a
measured quantum groupoid. Then there exists a strictly positive
operator $\delta$ affiliated with $M\cap\alpha(N)'\cap\beta(N)'$
called modulus and then there exists a strictly positive operator
$\lambda$ affiliated with $Z(M)\cap\alpha(N)\cap\beta(N)$ called
scaling operator such that $[D\nu\circ\alpha^{-1}\circ T_L\circ
R:D\nu\circ\alpha^{-1}\circ
T_L]_t=\lambda^{\frac{it^2}{2}}\delta^{it}$ for all
$t\in\mathbb{R}$.

Moreover, we have, for all $s,t\in\mathbb{R}$:
\begin{enumerate}[i)]
\item $\begin{aligned} &\ [D\nu\circ\alpha^{-1}\circ T_L\circ\tau_s
:D\nu\circ\alpha^{-1}\circ T_L]_t=\lambda^{-ist}\\
                       &\ [D\nu\circ\alpha^{-1}\circ T_L\circ R\circ\tau_s
:D\nu\circ\alpha^{-1}\circ T_L\circ R]_t=\lambda^{-ist}\\
                       &\ [D\nu\circ\alpha^{-1}\circ T_L\circ\sigma^{\nu\circ\alpha^{-1}\circ T_L\circ R}_s
:D\nu\circ\alpha^{-1}\circ T_L]_t=\lambda^{ist}\\
                       &\ [D\nu\circ\alpha^{-1}\circ T_L\circ R\circ\sigma^{\nu\circ\alpha^{-1}\circ T_L}_s:D\nu\circ\alpha^{-1}\circ T_L\circ R]_t=\lambda^{-ist}
       \end{aligned}$
\item $R(\lambda)=\lambda$, $R(\delta)=\delta^{-1}$,
$\tau_t(\delta)=\delta$ and $\tau_t(\lambda)=\lambda$ ;
\item $\delta$ is a group-like element i.e $\Gamma(\delta)
      =\delta\surl{\ _{\beta} \otimes_{\alpha}}_{\ N}\delta$.
\end{enumerate}
If $\nu'$ is a n.s.f weight on $N$ and $h$, $k$ are strictly
positive operators, affiliated with $N$, strongly commuting and
satisfying $[D\nu':D\nu]_t=k^{\frac{it^2}{2}}h^{it}$ for all
$t\in\mathbb{R}$, then there exists a n.s.f left invariant
operator-valued weight $\tilde{T_L}$ which is $\beta$-adapted
w.r.t $\nu'$. Moreover, if
$(N,M,\alpha,\beta,\Gamma,\nu',T_L',T_R')$ is an other measured
quantum groupoid, then, for all $t\in\mathbb{R}$, fundamental
objects are linked by:
\begin{center}
\begin{minipage}{14cm}
\begin{enumerate}[i)]
\item $R'=R$
\item
$\tau_t'=Ad_{\alpha(k^{\frac{-it^2}{2}}h^{-it})\beta(k^{\frac{it^2}{2}}h^{it})}\circ\tau_t
=Ad_{\alpha([D\nu':D\nu]_t^*)\beta([D\nu':D\nu]_t)}\circ\tau_t$
\item $\lambda'=\lambda$
\item $\dot{\delta'}=\dot{\delta}$ where $\dot{\delta}$ and $\dot{\delta}'$
have been defined in proposition \ref{deltadot}
\item $P'^{it}=\alpha(k^{\frac{it^2}{2}}h^{it})\beta(k^{\frac{-it^2}{2}}h^{-it})
J_{\Phi}\alpha(k^{\frac{it^2}{2}}h^{it})\beta(k^{\frac{-it^2}{2}}h^{-it})J_{\Phi}P^{it}$
\end{enumerate}
\end{minipage}
\end{center}
\end{theo}

\section{Examples}

\subsection{Groupoids}

\begin{defi}
A \textbf{groupoid} $G$ is a small category in which each
morphism $\gamma:x\rightarrow y$ is an isomorphism the inverse of
which is $\gamma^{-1}$. Let $G^{\{0\}}$ the set of objects of $G$
that we identify with $\{\gamma\in G |
\gamma\circ\gamma=\gamma\}$. For all $\gamma\in G$,
$\gamma:x\rightarrow y$, we denote
$x=\gamma^{-1}\gamma=s(\gamma)$ we call source object and
$y=\gamma\gamma^{-1}=r(\gamma)$ we call range object. If
$G^{\{2\}}$ is the set of pairs $(\gamma_1,\gamma_2)$ of $G$ such
that $s(\gamma_1)=r(\gamma_2)$, then composition of morphisms
makes sense in $G^{\{2\}}$.
\end{defi}

In \cite{R1}, J. Renault defines the structure of locally compact
groupoid $G$ with a Haar system $\{\lambda^u,u\in G^{\{0\}}\}$
and a quasi-invariant measure $\mu$ on $G^{\{0\}}$. We refer to
\cite{R1} for definitions and notations. We put
$\nu=\mu\circ\lambda$. We refer to \cite{C3} and \cite{ADR} for
discussions about transversal measures.

If $G$ is $\sigma$-compact, J.M Vallin constructs in \cite{V1} two
co-involutive Hopf bimodules on the same basis
$N=L^{\infty}(G^{\{0\}},\mu)$, following T. Yamanouchi's works in
\cite{Y}. The underlying von Neumann algebras are
$L^{\infty}(G,\nu)$ which acts by multiplication on $H=L^2(G,\nu)$
and $\mathcal{L}(G)$ generated by the left regular representation
$L$ of $G$.

We define a (resp. anti-) representation $\alpha$ (resp. $\beta$)
from $N$ in $L^{\infty}(G,\nu)$ such that, for all $f\in N$:
$$\alpha(f)=f\circ r\quad\text{ and }\quad\beta(f)=f\circ s$$

For all $i,j\in\{\alpha,\beta\}$, we define
$G^{\{2\}}_{i,j}\subset G\times G$ and a measure $\nu^2_{i,j}$
such that:
$$H\surl{\ _i\otimes_j}_{\ N}H\text{ is identified with }L^2(G^{\{2\}}_{i,j},\nu^2_{i,j})$$
For example, $G^{\{2\}}_{\beta,\alpha}$ is equal to $G^{\{2\}}$
and $\nu^2_{\beta,\alpha}$ to $\nu^2$. Then, we construct a
unitary $W_G$ from $H\surl{\ _{\alpha}\otimes_{\alpha}}_{\ \mu}H$
onto $H\surl{\ _{\beta}\otimes_{\alpha}}_{\ \mu}H$, defined for
all $\xi\in L^2(G^{\{2\}}_{\alpha,\alpha},\nu^2_{\alpha,\alpha})$
by:
$$W_G\xi(s,t)=\xi(s,st)$$ for $\nu^2$-almost all $(s,t)$ in $G^{\{2\}}$.

This leads to define co-products $\Gamma_G$ and
$\widehat{\Gamma_G}$ by formulas:
$$\Gamma_G(f)=W_G(1\surl{\ _{\alpha}
  \otimes_{\alpha}}_{\ N}f)W_G^*\quad\text{ and }\quad \widehat{\Gamma_G}(k)=W_G^*(k\surl{\ _{\beta}
  \otimes_{\alpha}}_{\ N}1)W_G$$
for all $f\in L^{\infty}(G,\nu)$ and $k\in\mathcal{L}(G)$, this
explicitly gives:
$$\Gamma_G(f)(s,t)=f(st)$$
for all $f\in L^{\infty}(G,\nu)$ and $\nu^2$-almost all $(s,t)$
in $G^{\{2\}}$,
$$\widehat{\Gamma_G}(L(h))\xi(x,y)=\int_Gh(s)\xi(s^{-1}x,s^{-1}y)d\lambda^{r(x)}(s)$$
for all $\xi\in
L^2(G^{\{2\}}_{\alpha,\alpha},\nu^2_{\alpha,\alpha})$, $h$ a
continuous function with compact support on $G$ and
$\nu^2_{\alpha,\alpha}$-almost all $(x,y)$ in
$G^{\{2\}}_{\alpha,\alpha}$. Moreover, we define two
co-involutions $j_G$ and $\widehat{j_G}$ by:
$$j_G(f)(x)=f(x^{-1})$$ for all $f\in L^{\infty}(G,\nu)$ and
almost all $x$,
$$\widehat{j_G}(g)=Jg^*J$$ for all $g\in\mathcal{L}(G)$ and where
$J$ is the involution $J\xi=\overline{\xi}$ for all $\xi\in
L^2(G)$. Finally, we define two n.s.f left invariant
operator-valued weights $P_G$ and $\widehat{P_G}$:
$$P_G(f)(y)=\int_Gf(x)d\lambda^{r(y)}(x)\quad\text{ and
}\quad\widehat{P_G}(L(f))=\alpha(f_{|G^{\{0\}}})$$ for all
continuous with compact support $f$ on $G$ $\nu$-almost all $y$ in
$G$.

\begin{theo}
Let $G$ be a $\sigma$-compact, locally compact groupoid with a
Haar system and a quasi-invariant measure $\mu$on units. Then:
$$
(L^{\infty}(G^{\{0\}},\mu),L^{\infty}(G,\nu),\alpha,\beta,\Gamma_G,\mu,P_G,j_GP_Gj_G)$$
is a commutative measured quantum groupoid and:
$$(L^{\infty}(G^{\{0\}},\mu),{\mathcal
L}(G),\alpha,\alpha,\widehat{\Gamma_G},\mu,\widehat{P_G},
\widehat{j_G}\widehat{P_G}\widehat{j_G})
$$ is a symmetric measured quantum groupoid. The unitary $V_G=W_G^*$
is the fundamental unitary of the commutative structure.
\end{theo}

\begin{proof}
By \cite{V1} (th. 3.2.7 and 3.3.7),
$(L^{\infty}(G^{\{0\}},\mu),L^{\infty}(G,\nu),\alpha,\beta,\Gamma_G)$
and $(L^{\infty}(G^{\{0\}},\mu),{\mathcal
L}(G),\alpha,\alpha,\widehat{\Gamma_G})$ are co-involutive Hopf
bimodules with left invariant operator-valued weights; to get
right invariants operator-valued weights, we consider $j_GP_Gj_G$
and $\widehat{j_G}\widehat{P_G}\widehat{j_G}$.

Since $L^{\infty}(G,\nu)$ is commutative, $P_G$ is adapted w.r.t
$\mu$ by \cite{V1} (theorem 3.3.4),
$\sigma_t^{\mu\circ\alpha^{-1}\circ\widehat{P_G}}$ fixes point by
point $\alpha(N)$ so that $\widehat{P_G}$ is adapted w.r.t $\mu$.

Finally, for all $e,f,g$ continuous functions with compact support
and almost all $(s,t)$ in $G^{\{2\}}$, we have, by
\ref{raccourci}:
$$
\begin{aligned}
(1\surl{\ _{\beta}\otimes_{\alpha}}_{\ N}JeJ)W_G(f\surl{\
_{\alpha}\otimes_{\alpha}}_{\ \mu}g)(s,t)&=\overline{e(t)}f(s)g(st)=\Gamma_G(g)(f\surl{\ _{\beta}\otimes_{\alpha}}_{\ \mu}\overline{e})(s,t)\\
&=(1\surl{\ _{\beta}\otimes_{\alpha}}_{\ N}JeJ)U_H(f\surl{\
_{\alpha}\otimes_{\alpha}}_{\ \mu}g)(s,t)
\end{aligned}$$
so that we get $U_H=W_G$.
\end{proof}

\begin{rema}
In the commutative structure, modular function
$\frac{d\nu^{-1}}{d\nu}$ and modulus coincide and the scaling
operator is trivial.
\end{rema}

We have a similar result for measured quantum groupoids in the
sense of Hahn (\cite{Ha1} and \cite{Ha2}):

\begin{theo}
From all measured groupoid $G$, we construct a commutative
measured quantum groupoid
$(L^{\infty}(G^{\{0\}},\mu),L^{\infty}(G,\nu),\alpha,\beta,\Gamma_G,\mu,P_G,j_GP_Gj_G)$
and a symmetric one $(L^{\infty}(G^{\{0\}},\mu),{\mathcal
L}(G),\alpha,\alpha,\widehat{\Gamma_G},\mu,\widehat{P_G},\widehat{j_G}\widehat{P_G}\widehat{j_G})$.
Objects are defined in a similar way as in the locally compact
case. The unitary $V_G$ is the fundamental unitary of the
commutative structure.
\end{theo}

\begin{proof}
Results come from \cite{Y} for the symmetric case. It is
sufficient to apply in this case, technics of \cite{V1} for the
commutative case and invariant operator-valued weights.
\end{proof}

\begin{conj}
If $(N,M,\alpha,\beta,\Gamma,\mu,T_L,T_R)$ is a measured quantum
groupoid such that $M$ is commutative, then there exists a
locally compact groupoid $G$ such that:
$$(N,M,\alpha,\beta,\Gamma,\mu,T_L,T_R)\simeq
(L^{\infty}(G^{\{0\}},\mu),L^{\infty}(G,\nu),\alpha,\beta,\Gamma_G,\mu,P_G,j_G\circ
P_G\circ j_G)$$
\end{conj}

\subsection{Finite quantum groupoids}

\begin{defi}(Weak Hopf C*-algebras \cite{BSz})
We call \textbf{weak Hopf C*-algebra } or finite quantum groupoid
all $(M,\Gamma,\kappa,\varepsilon)$ where $M$ is a finite
dimensional C*-algebra with a co-product $\Gamma:M\rightarrow
M\otimes M$, a co-unit $\varepsilon$ and an antipode
$\kappa:M\rightarrow M$ such that, for all $x,y\in M$:
\begin{center}
\begin{minipage}{12cm}
\begin{enumerate}[i)]
\item $\Gamma$ is a *-homomorphism (not necessary unital);
\item Unit and co-unit satisfy the following relation:
$$(\varepsilon\otimes\varepsilon)((x\otimes
1)\Gamma(1)(1\otimes y))=\varepsilon(xy)$$
\item $\kappa$ is an anti-homomorphism of algebra and co-algebra such that:
\begin{itemize}
\item $(\kappa\circ *)^2=\iota$
\item $(m(\kappa\otimes id)\otimes id)(\Gamma\otimes id)\Gamma
(x)=(1\otimes x)\Gamma(1)$.
\end{itemize}
where $m$ denote the product on $M$.
\end{enumerate}
\end{minipage}
\end{center}

\end{defi}

We recall some results \cite{NV1}, \cite{NV2} and \cite{BNSz}. If
$(M,\Gamma,\kappa,\varepsilon)$ is a weak Hopf C*-algebra. We call
co-unit range (resp. source) the application
$\varepsilon_t=m(id\otimes \kappa)\Gamma$ (resp.
$\varepsilon_s=m(\kappa\otimes id)\Gamma$). We have
$\kappa\circ\varepsilon_t=\varepsilon_s\circ\kappa$. There exists
a unique faithful positive linear form $h$, called normalized Haar
measure of $(M,\Gamma,\kappa,\varepsilon)$ which is
$\kappa$-invariant, such that $(id\otimes h)(\Gamma(1))=1$ and,
for all $x,y\in M$, we have:
$$(id\otimes h)((1\otimes
y)\Gamma(x))=\kappa((i\otimes h)(\Gamma(y)(1\otimes x)))$$
Moreover, $E^s_h=(h\otimes id)\Gamma$ (resp. $E^t_h=(id\otimes
h)\Gamma$) is a Haar conditional expectation to the source (resp.
range) Cartan subalgebra $\varepsilon_s(M)$ (resp. range
$\varepsilon_t(M)$) such that $h\circ E^s_h=h$ (resp. $h\circ
E^t_h=h$). Range and source Cartan subalgebras commute.

By \cite{V5} and \cite{Ni}, we can always assume that
$\kappa^2_{|\varepsilon_t(M)}=id$ thanks to a deformation.
\textbf{In the following, we assume that the condition holds.}

Since $h\circ\kappa=h$ and
$\kappa\varepsilon_t=\varepsilon_s\kappa$, we have
$h\circ\varepsilon_t=h\circ\varepsilon_s$.

\begin{theo}
Let $(M,\Gamma,\kappa,\varepsilon)$ be a weak Hopf C*-algebra,
$h$ its normalized Haar measure, $E^s_h$ (resp. $E^t_h$) its
source (resp. range) Haar conditional expectation and
$\varepsilon_t(M)$ its range Cartan subalgebra. We put
$N=\varepsilon_t(M)$, $\alpha=id_{|N}$, $\beta=\kappa_{|N}$,
$\tilde{\Gamma}$ the co-product $\Gamma$ viewed as an operator
which takes value in:
$$M\surl{\ _{\beta} \star_{\alpha}}_{\ N}M\simeq (M\otimes
M)_{\Gamma(1)}$$ and $\mu=h\circ\alpha=h\circ\beta$. Then
$(N,M,\alpha,\beta,\tilde{\Gamma},\mu,E^t_h,E^s_h)$ is a measured
quantum groupoid.
\end{theo}

\begin{proof}
$\alpha$ is a representation from $N$ in $M$ and, since
$\kappa^2_{|\varepsilon_t(M)}=id$, $\beta$ is a
anti-representation from $N$ in $M$. They commute each other
because Cartan subalgebras commute and
$\kappa\circ\varepsilon_t=\varepsilon_s\circ\kappa$. For all $n\in
N$, there exists $m\in M$ such that $n=\varepsilon_t(m)$. So, we
have:
$$\tilde{\Gamma}(\alpha(n))=\tilde{\Gamma}(\varepsilon_t(m))=\Gamma(1)(\varepsilon_t(m)\otimes 1)\Gamma(1)
=\alpha(n)\surl{\ _{\beta}\otimes_{\alpha}}_{\ N}1$$ Also, we have
$\tilde{\Gamma}(\beta(n))=1\surl{\ _{\beta} \otimes_{\alpha}}_{\
N}\beta(n)$ and $\tilde{\Gamma}$ is a co-product. Then
$(N,M,\alpha,\beta,\Gamma)$ is Hopf bimodule. Moreover, for all
$n\in N$ and $t\in\mathbb{R}$, we have:
$$
\begin{aligned} \sigma_t^{E_h^t}(\beta(n))=\sigma_t^{h\circ
E_h^t}(\beta(n))=\sigma_t^{h\circ
E_h^s}(\beta(n))&=\sigma_t^{h_{|\beta(N)}}(\beta(n))\\
&=\beta(\sigma_{-t}^{h_{|\beta(N)}\circ\beta}(n))
=\beta(\sigma_{-t}^{\mu}(n))
\end{aligned}$$ and $E_h^t$ is
$\beta$-adapted w.r.t $\mu$. Since $E_h^s=\kappa\circ
E_h^t\circ\kappa$, then $E_h^s$ is $\alpha$-adapted w.r.t $\mu$.
\end{proof}

\begin{theo}
Let $(N,M,\alpha,\beta,\Gamma,\nu,T_L,T_R)$ be a measured quantum
groupoid such that $M$ is finite dimensional. Then, there exist
$\tilde{\Gamma}$, $\kappa$ and $\varepsilon$ such that
$(M,\tilde{\Gamma},\kappa,\varepsilon)$ is a weak Hopf C*-algebra.
\end{theo}

\begin{proof}
By \ref{iden}, we identify via $I_{\beta,\alpha}^{\nu}$,
$L^2(M)\surl{\ _{\beta}\otimes_{\alpha}}_{\ N}L^2(M)$ with a
subspace of $L^2(M)\otimes L^2(M)$. We put
$\tilde{\Gamma}(x)=I_{\beta,\alpha}^{\nu}\Gamma(x)(I_{\beta,\alpha}^{\nu})^*$.
By \cite{V3} (definition 2.2.3), the fundamental
pseudo-multiplicative unitary becomes a multiplicative partial
isometry on $L^2(M)\otimes L^2(M)$ of basis
$(N,\alpha,\hat{\beta},\beta)$ by
$I=I_{\alpha,\hat{\beta}}^{\nu}W(I_{\beta,\alpha}^{\nu})^*$. $I$
is regular in the sense of \cite{V3} (definition 2.6.3) by
\ref{reg}. Moreover, if we put $H=L^2(M)$, then
$Tr_H(R(m))=Tr_H(m)$ for all $m\in M$ because $R$ is implemented
by an anti-unitary, so
$Tr_H\circ\beta=Tr_H\circ\alpha=Tr_H\circ\hat{\beta}$ and we
conclude by \cite{V3} (proposition 3.1.3).
\end{proof}

\begin{rema}
With notations of section \ref{iden}, $\kappa$ and $S$ are linked
by:
$$\kappa(x)=\alpha(n_o^{1/2}d^{1/2})\beta(n_o^{-1/2}d^{-1/2})S(x)
\alpha(n_o^{-1/2}d^{-1/2})\beta(n_o^{1/2}d^{1/2})$$
\end{rema}

\subsection{Quantum groups}

\begin{theo}
Measured quantum groupoids, basis $N$ on which is equal to
$\mathbb{C}$ are exactly locally compact quantum groups (in the
von Neumann setting) introduced by J. Kustermans and S. Vaes in
\cite{KV2}.
\end{theo}

\begin{proof}
In this case, the notion of relative tensor product is just usual
tensor product of Hilbert spaces, the notion of fibered product
is just tensor product of von Neumann algebras and the notion of
operator-valued weight is just weight.
\end{proof}

\subsection{Compact case}

In this section, we show that pseudo-multiplicative unitaries of
compact type in the sense of \cite{E2} correspond exactly to
measured quantum groupoids with a Haar conditional expectation.

\begin{defi}
Let $W$ be a pseudo-multiplicative unitary over $N$ w.r.t
$\alpha,\beta,\hat{\beta}$. Let $\nu$ be a n.s.f weight on $N$.
We say that $W$ is of \textbf{compact type} w.r.t $\nu$ if there
exists $\xi\in H$ such that:
\begin{center}
\begin{minipage}{13cm}
\begin{enumerate}[i)]
\item $\xi$ belongs to $D(H_{\hat{\beta}},\nu^o)\cap D(_{\alpha}H,\nu)\cap D(H_{\beta},\nu^o)$;
\item
$<\xi,\xi>_{\hat{\beta},\nu^o}=<\xi,\xi>_{\alpha,\nu}=<\xi,\xi>_{\beta,\nu^o}=1$
\item we have $W(\xi\surl{\ _{\hat{\beta}}\otimes_{\alpha}}_{\ \nu}\eta)
=\xi\surl{\ _{\alpha}\otimes_{\beta}}_{\ \nu^o}\eta$ for all
$\eta\in H$.
\end{enumerate}
\end{minipage}
\end{center}
In this case, $\xi$ is said to be \textbf{fixed and
bi-normalized}. We also say that $W$ is of \textbf{discrete type}
w.r.t $\nu$ if $\hat{W}$ is of compact type.
\end{defi}

By \cite{E2} (proposition 5.11), we recall that, if $W$ is of
compact type w.r.t $\nu$ and $\xi$ is a fixed and bi-normalized
vector, then $\nu$ shall be a faithful, normal, positive form on
$N$ which is equal to
$\omega_{\xi}\circ\alpha=\omega_{\xi}\circ\beta=\omega_{\xi}\circ\hat{\beta}$
and it is called \textbf{canonical form}.

\begin{prop}\label{compact}
Let $(N,M,\alpha,\beta,\Gamma)$ be a Hopf bimodule. Assume there
exist:
\begin{center}
\begin{minipage}{13cm}
\begin{enumerate}[i)]
\item a n.f left invariant conditional expectation from $E$
to $\alpha(N)$;
\item a n.f right invariant conditional expectation from $F$
to $\beta(N)$;
\item a n.f state $\nu$ on $N$ such that $\nu\circ\alpha^{-1}
\circ E=\nu\circ\beta^{-1}\circ F$.
\end{enumerate}
\end{minipage}
\end{center}

Then $(N,M,\alpha,\beta,\Gamma,\nu,E,F)$ is a measured quantum
groupoid. Moreover, if $R,\tau,\lambda$ and $\delta$ are
fondamental objects of the structure, then we have $F=R\circ
E\circ R$ and $\lambda=\delta=1$. Finally,
$\Lambda_{\nu\circ\alpha^{-1}\circ E}(1)$ is co-fixed and
bi-normalized, and the fundamental pseudo-multiplicative unitary
$W$ is weakly regular and of discrete type in sense of \cite{E2}
(paragraphe 5).
\end{prop}

\begin{proof}
For all $t\in\mathbb{R}$ and $n\in N$, we have:
$$\sigma_t^E(\beta(n))=\sigma_t^{\nu\circ\alpha^{-1}
\circ E}(\beta(n))=\sigma_t^{\nu\circ\beta^{-1} \circ
F}(\beta(n))=\beta(\sigma_{-t}^{\nu}(n))$$ Also, we have:
$$\sigma_t^F(\alpha(n))=\sigma_t^{\nu\circ\beta^{-1}
\circ F}(\alpha(n))=\sigma_t^{\nu\circ\alpha^{-1} \circ
E}(\alpha(n))=\alpha(\sigma_t^{\nu}(n))$$ so that
$(N,M,\alpha,\beta,\Gamma,\nu,E,F)$ is a measured quantum
groupoid. By definition, we have:
$$[D\nu\circ\alpha^{-1}
\circ E\circ R:D\nu\circ\alpha^{-1} \circ
E]_t=\lambda^{\frac{it^2}{2}}\delta^{it}$$ On the other hand,
since $\nu\circ\alpha^{-1}\circ E=\nu\circ\beta^{-1}\circ F$ and
by uniqueness, there exists a strictly positive element $h$
affiliated with $Z(N)$:
$$
\begin{aligned}
\ [D\nu\circ\alpha^{-1}\circ E\circ R:D\nu\circ\alpha^{-1}\circ
E]_t&=[DR\circ E\circ R:DF]_t=\alpha(h^{it})
\end{aligned}$$
We deduce that $\lambda=1$ and $\delta=\alpha(h)$, so
$\alpha(h^{-1})=\delta^{-1}=R(\delta)=\beta(h)$ and by \cite{E1}
(5.2), we get $h=1$.

We put $\Phi=\nu\circ\alpha^{-1}\circ E$. If $(\xi_i)_{i\in I}$
is a $(N^o,\nu^o)$-basis of $(H_{\Phi})_{\beta}$ then, for all
$v\in D(H_{\beta},\nu^o)$:
$$
\begin{aligned}
U_H(v \surl{\ _{\alpha} \otimes_{\hat{\beta}}}_{\ \ \nu^o}
\Lambda_{\Phi}(1))&=\sum_{i
    \in I} \xi_{i} \surl{\ _{\beta} \otimes_{\alpha}}_{\ \nu}
  \Lambda_{\Phi} ((\omega_{v,\xi_i} \surl{\
      _{\beta} \star_{\alpha}}_{\ \nu} id)(\Gamma(1)))\\
&=\sum_{i \in I} \xi_{i} \surl{\ _{\beta} \otimes_{\alpha}}_{\
\nu} \alpha(<v,\xi_i>_{\beta,\nu^o})\Lambda_{\Phi}(1)=v\surl{\
_{\beta} \otimes_{\alpha}}_{\ \nu}\Lambda_{\Phi}(1)
\end{aligned}$$
It is easy to see that $\Lambda_{\Phi}(1)$ belongs to
$D((H_{\Phi})_{\hat{\beta}},\nu^o)\cap D(_{\alpha}H_{\Phi},\nu)$
and satisfies
$<\Lambda_{\Phi}(1),\Lambda_{\Phi}(1)>_{\hat{\beta},\nu^o}=
<\Lambda_{\Phi}(1),\Lambda_{\Phi}(1)>_{\alpha,\nu}=1$ so that, by
continuity, we get $U_H(v \surl{\ _{\alpha}
\otimes_{\hat{\beta}}}_{\ \ \nu^o} \Lambda_{\Phi}(1))=v\surl{\
_{\beta} \otimes_{\alpha}}_{\ \nu}\Lambda_{\Phi}(1)$ for all $v\in
H$ i.e $\Lambda_{\Phi}(1)$ is co-fixed and bi-normalized. Since
$\nu\circ\alpha^{-1}\circ E=\Phi=\nu\circ\beta^{-1}\circ F$, we
have by \ref{prem}, for all $n\in \mathcal{N}_{\nu}$:
$$\beta(n^*)\Lambda_{\Phi}(1)=\beta(n^*)J_{\Phi}\Lambda_{\Phi}(1)
=J_{\Phi}\Lambda_F(1)\Lambda_{\nu}(n)$$ so that
$\Lambda_{\Phi}(1)$ is $\beta$-bounded w.r.t $\nu^o$ and
$R^{\beta,\nu^o}(\Lambda_{\Phi}(1))=J_{\Phi}\Lambda_F(1)J_{\nu}$.
Consequently, $\Lambda_{\Phi}(1)$ is bi-normalized and $W$ is of
discrete type.
\end{proof}

\begin{coro}
Let $W$ be a weakly regular pseudo-multiplicative unitary over $N$
w.r.t $\alpha,\beta,\hat{\beta}$ of compact type w.r.t the
canonical form $\nu$. If $\xi$ a fixed and bi-normalized vector,
we put:
\begin{center}
\begin{minipage}{13cm}
\begin{enumerate}[i)]
\item ${\mathcal A}$ the von Neumann algebra generated by right leg of $W$;
\item $\Gamma(x)=\sigma_{\nu^o}W(x\surl{\ _{\alpha}\otimes_{\beta}}_{\ \nu^o}1)
W^*\sigma_{\nu}$ for all $x\in {\mathcal A}$ ;
\item $E=(\omega_{\xi}\surl{\ _{\beta}\star_{\alpha}}_{\ \nu}id)\circ\Gamma$
and $F=(id\surl{\ _{\beta}\star_{\alpha}}_{\
\nu}\omega_{\xi})\circ\Gamma$.
\end{enumerate}
\end{minipage}
\end{center}
Then $(N,{\mathcal A},\alpha,\beta,\Gamma,\nu,E,F)$ is a measured
quantum groupoid. Moreover, if $R,\tau,\lambda$ and $\delta$ are
the fundamental objects of the structure, we have $F=R\circ E\circ
R$, $\lambda=\delta=1$ and the fundamental unitary is $\hat{W}$.
\end{coro}

\begin{proof}
By \cite{EV} (6.3), we know that $(N,{\mathcal
A},\alpha,\beta,\Gamma)$ is a Hopf bimodule. By \cite{E2} (theorem
6.6), $E$ is a n.f left invariant conditional expectation from
${\mathcal A}$ to $\alpha(N)$. By \cite{E2} (propositions 6.2 and
6.3), $F$ is a n.f right invariant conditional expectation from
${\mathcal A}$ to $\beta(N)$. Moreover, we clearly have
$\omega_{\xi}\circ E=\omega_{\xi}\circ F$ so that
$\nu\circ\alpha^{-1}\circ E=\nu\circ\beta^{-1}\circ F$. We are in
conditions of the previous proposition an we get that
$(N,{\mathcal A},\alpha,\beta,\Gamma,\nu,E,F)$ is a measured
quantum groupoid, $F=R\circ E\circ R$ and $\lambda=\delta=1$.
Finally, by \cite{E2} (corollaire 7.7), $\hat{W}$ is the
fundamental unitary. (More exactly, it is
$\sigma_{\nu^o}W_s^*\sigma_{\nu}$ where $W_s$ is the standard
form of $W$ in th sense of \cite{E2} (paragraphe 7)).
\end{proof}

The converse is also true and so we characterize the compact case:

\begin{coro}
Let $(N,M,\alpha,\beta,\Gamma)$ be a Hopf bimodule. We assume
there exist:
\begin{center}
\begin{minipage}{13cm}
\begin{enumerate}[i)]
\item a co-involution $R$;
\item a n.f left invariant conditional expectation from $E$
to $\alpha(N)$.
\end{enumerate}
\end{minipage}
\end{center}
Then there exists a n.f state $\nu$ on $N$ such that
$(N,M,\alpha,\beta,\Gamma,\nu,E,R\circ E\circ R)$ is a measured
quantum groupoid with trivial modulus and scaling operator and the
fundamental unitary of which is of discrete type w.r.t $\nu$.
\end{coro}

\begin{proof}
We put $F=R\circ E\circ R$ which is a n.f right invariant
conditional expectation from $M$ to $\beta(N)$. We also put:
$$\tilde{E}=E_{|\beta(N)}:\beta(N)\rightarrow\alpha(Z(N))\text{ and }
\tilde{F}=F_{|\alpha(N)}:\alpha(N)\rightarrow\beta(Z(N))$$ We
have, for all $m\in M$:
$$
\begin{aligned}
\tilde{F}E(m)\surl{\ _{\beta} \otimes_{\alpha}}_{\ N}1&=(F\surl{\
_{\beta} \star_{\alpha}}_{\ N}id)(E(m)\surl{\ _{\beta}
\otimes_{\alpha}}_{\ N}1)\\
&=(F\surl{\ _{\beta} \star_{\alpha}}_{\ N}id)(id\surl{\ _{\beta}
\star_{\alpha}}_{\ N}E)\Gamma(m)\\
&=(id\surl{\ _{\beta} \star_{\alpha}}_{\ N}E)(F\surl{\ _{\beta}
\star_{\alpha}}_{\ N}id)\Gamma(m)\\
&=(id\surl{\ _{\beta} \star_{\alpha}}_{\ N}E)(1\surl{\ _{\beta}
\otimes_{\alpha}}_{\ N}F(m))=1\surl{\ _{\beta}
\otimes_{\alpha}}_{\ N}\tilde{E}F(m)
\end{aligned}$$
so, if $\tilde{F}E(m)=\beta(n)$ for some $n\in Z(N)$, then
$\tilde{E}F(m)=\alpha(n)$. Moreover, we have:
$$
\begin{aligned}
\tilde{E}F(m)\surl{\ _{\beta} \otimes_{\alpha}}_{\
N}1&=EF(m)\surl{\ _{\beta} \otimes_{\alpha}}_{\ N}1=(id\surl{\
_{\beta} \star_{\alpha}}_{\ N}E)\Gamma(F(m))\\
&=(id\surl{\ _{\beta} \star_{\alpha}}_{\ N}E)(1\surl{\ _{\beta}
\otimes_{\alpha}}_{\ N}F(m))=1\surl{\ _{\beta}
\otimes_{\alpha}}_{\ N}\tilde{E}F(m)
\end{aligned}$$
so that $\alpha(n)=\beta(n)$. Consequently
$\tilde{E}F(m)=\tilde{F}E(m)$ and $EF=FE$ is a n.f conditional
expectation from $M$ to: $$\tilde{N}=\alpha(\{n\in
Z(N),\alpha(n)=\beta(n)\})=\beta(\{n\in
Z(N),\alpha(n)=\beta(n)\})$$ Also, we have $R_{|\tilde{N}}=id$.
So, if $\omega$ is a n.f state on $\tilde{N}$, we have
$\omega\circ\tilde{E}\circ\beta=\omega\circ\tilde{F}\circ\alpha$
and
$\nu=\omega\circ\tilde{E}\circ\beta=\omega\circ\tilde{F}\circ\alpha$
satisfies hypothesis of \ref{compact}: then, corollary holds.
\end{proof}

\begin{coro}
Let $(N,M,\alpha,\beta,\Gamma,\nu,T_L,T_R)$ be a measured quantum
groupoid such that $T_L$ is a conditional expectation. Then there
exists a n.f state $\nu'$ on $N$ such that
$\sigma^{\nu'}=\sigma^{\nu}$ and the fundamental unitary is of
discrete type w.r.t $\nu'$.
\end{coro}

\begin{proof}
Let $R$ be the co-involution. By the previous corollary, there
exists a n.f state $\nu'$ on $N$ such that
$(N,M,\alpha,\beta,\Gamma,\nu',T_L,R\circ T_L\circ R)$ is a
measured quantum groupoid. Since $T_L$ is $\beta$-adapted w.r.t
$\nu$ and $\nu'$, we have $\sigma^{\nu'}=\sigma^{\nu}$. We easily
verify that the fundamental unitary of the first structure
coincides with that of the last one which is of discrete type
w.r.t $\nu'$ by the previous corollary.
\end{proof}

\subsection{Depth $2$ inclusions}

Let $M_0\subseteq M_1$ be an inclusion of von Neumann algebras.
We call \textbf{basis construction} the following inclusions:

$$M_0\subseteq M_1\subseteq M_2=J_1M'_0J_1=End_{M_0^o}(L^2(M_1))$$

By iteration, we construct Jones' tower $M_0\subseteq M_1\subseteq
M_2\subseteq M_3\subseteq\cdots$

\begin{defi}If $M'_0\cap M_1\subseteq M'_0\cap M_2\subseteq M'_0\cap
M_3$ is a basis construction, then the inclusion is said to be
\textbf{of depth $2$}.
\end{defi}

Let $T_1$ be a n.s.f operator-valued weight from $M_1$ to $M_0$.
By Haagerup's construction \cite{St} (12.11) and \cite{EN} (10.1),
it is possible to define a canonical n.s.f operator-valued weight
$T_2$ from $M_2$ to $M_1$ such that, for all $x,y\in {\mathcal
N}_{T_1}$, we have:
$$T_2(\Lambda_{T_1}(x)\Lambda_{T_1}(y)^*)=xy^*$$ By iteration, we define,
for all $i\geq 1$, a n.s.f operator-valued weight $T_i$ from $M_i$
to $M_{i-1}$. If $\psi_0$ is n.s.f weight sur $M_0$, we put
$\psi_i=\psi_{i-1}\circ T_i$.

\begin{defi}
\cite{EN} (11.12), \cite{EV} (3.6). $T_1$ is said to be regular if
restrictions of $T_2$ to $M_0'\cap M_2$ and of $T_3$ to $M_1'\cap
M_3$ are semifinite.
\end{defi}

\begin{prop}\cite{EV} (3.2, 3.8, 3.10).
If $M_0\subset M_1$ is an inclusion with a regular n.s.f
operator-valued weight $T_1$ from $M_1$ to $M_0$, then there
exists a natural *-representation $\pi$ of $M_0'\cap M_3$ on
$L^2(M_0'\cap M_2)$ whose restriction to $M_0'\cap M_2$ is the
standard representation of $M_0'\cap M_2$. Moreover, the inclusion
is of depth $2$ if, and only if $\pi$ is faithful.
\end{prop}

The following theorem exhibits a structure of measured quantum
groupoid coming from inclusion of von Neumann algebras.

\begin{theo}
\cite{E3} (theorem 9.2) Let $M_0\subset M_1$ be a depth $2$
inclusion of von Neumann algebras with a regular n.s.f
operator-valued weight $T_1$ from $M_1$ to $M_0$. Let us assume
there exists a n.s.f weight $\chi$ on $M'_0\cap M_1$ such that
$\sigma_t^{\chi}=\sigma_t^{T_1}$ for all $t\in\mathbb{R}$.

Then, there exists $\Gamma$ such that $(M'_1\cap M_2,M'_1\cap
M_3,j_2,id,\Gamma)$ is a Hopf bimodule with a left invariant
operator-valued weight $T_{3|M'_1\cap M_3}$ which is
$j_2$-adapted and a co-involution $j_2$ where $j_2$ is the
canonical anti-isomorphism from $M_2$ onto $M'_2$ which sends $x$
to $J_2x^*J_2$ where $J_2$ is given by Tomita's theory on
$L^2(M_2)$. So, we get a structure of measured quantum groupoid
on $M'_1\cap M_3$.
\end{theo}

By \cite{EV}, theorem 7.3 and proposition 7.5, the previous
quantum groupoid acts on the von Neumann algebra $M_1$ such that
invariants are exactly elements of the von Neumann algebra $M_0$.

\begin{rema}
If $M_0$ and $M_1$ are semi-finite, then $\sigma^{T_1}$ is
interior. Moreover, if, $M_0'\cap M_1$ is semi-finite, then
$\sigma^{T_1}$ becomes the modular group of a n.s.f weight $\chi$
on $M_0'\cap M_1$. Consequently, all operator-valued weights are
adapted in sense of \cite{E1}. Then, we can also put a structure
of measured quantum groupoid on $M'_0\cap M_2$ such that the left
operator-valued weight $T_{2|M'_0\cap M_2}$ is invariant and
$j_1$-adapted where $j_1$ comes from Tomita's theory. This
situation will be developed in more details in a forthcoming
article about duality.
\end{rema}

\subsection{Quantum space quantum groupoid}

Let $M$ be a von Neumann algebra. $M$ acts on
$H=L^2(M)=L^2_{\nu}(M)$ where $\nu$ is a n.s.f weight on $M$. We
denote by $M'$, (resp. $Z(M)'$) the commutant of $M$ (resp.
$Z(M)$) in $\mathcal{L}(L^2(M))$. Let $tr$ be a n.s.f trace on
$Z(M)$. $M'\surl{\star}_{Z(M)}M=M'\surl{\otimes}_{Z(M)}M$ acts on
$L^2(M)\surl{\otimes}_{tr}L^2(M)$. There exists a n.s.f
operator-valued weight $T$ from $M$ to $Z(M)$ such that
$\nu=tr\circ T$.

Let $\alpha$ (resp. $\beta$) be the (resp. anti-) representation
of $M$ to $M'\surl{\otimes}_{Z(M)}M$ such that
$\alpha(m)=1\surl{\otimes}_{Z(M)}m$ (resp.
$\beta(m)=j(m)\surl{\otimes}_{Z(M)}1$) where
$j(x)=J_{\nu}x^*J_{\nu}$ for all $x\in \mathcal{L}(L^2_{\nu}(M))$.

\begin{prop}\label{renv2}
The following formula:
$$
\begin{aligned}
I: [L^2(M)\surl{\otimes}_{tr}L^2(M)]\surl{\ _{\beta}
  \otimes_{\alpha}}_{\ \nu}[L^2(M)\surl{\otimes}_{tr}L^2(M)]
  &\rightarrow L^2(M)\surl{\otimes}_{tr}L^2(M)
  \surl{\otimes}_{tr}L^2(M)\\
[\Lambda_{\nu}(y)\surl{\otimes}_{tr}\eta]\surl{\ _{\beta}
  \otimes_{\alpha}}_{\ \nu}\Xi&\mapsto \alpha(y)\Xi\surl{\otimes}_{tr}\eta
\end{aligned}$$
for all $\eta\in L^2(M),\Xi\in L^2(M)\surl{\otimes}_{tr}L^2(M)$
and $y\in M$, defines a canonical isomorphism such that we have
$I([m\surl{\otimes}_{Z(M)}z]\surl{\ _{\beta}\otimes_{\alpha}}_{\
\nu}Z)=(\alpha(M)Z\surl{\otimes}_{Z(M)}z)I$, for all $m\in M$,
$z\in Z(M)'$ and $Z\in \mathcal{L}(L^2(M))\surl{\star}_{Z(M)}M'$.
\end{prop}

\begin{proof}
Straightforward.
\end{proof}
We identify $(M'\surl{\otimes}_{Z(M)}M)\surl{\
_{\beta}\star_{\alpha}}_{\ M}(M'\surl{\otimes}_{Z(M)}M)$ with
$M'\surl{\otimes}_{Z(M)}Z(M)\surl{\otimes}_{Z(M)}M$ and so with
$M'\surl{\otimes}_{Z(M)}M$. We define a normal *-homomorphism
$\Gamma$ by:
$$
\begin{aligned}
M'\surl{\otimes}_{Z(M)}M &\rightarrow (M'\surl{\otimes}_{Z(M)}M)
\surl{\ _{\beta}\star_{\alpha}}_{\ \nu} (M'\surl{\otimes}_{Z(M)}M)\\
n\surl{\otimes}_{Z(M)}m&\mapsto
I^*(n\surl{\otimes}_{Z(M)}1\surl{\otimes}_{Z(M)}m)I=[1\surl{\otimes}_{Z(M)}m]
\surl{\ _{\beta}\otimes_{\alpha}}_{\ \nu}[n\surl{\otimes}_{Z(M)}1]
\end{aligned}$$
$\Gamma$ is, in fact, the identity trough the previous
isomorphism.

\begin{theo}
If we put $T_R=id\surl{\star}_{Z(M)}T$ and
$R=\varsigma_{Z(M)}\circ(j\surl{\otimes}_{Z(M)}j)$, then
$(M,M'\surl{\otimes}_{Z(M)}M,\alpha,\beta,\Gamma,\nu,R\circ
T_R\circ R,T_R)$ becomes a measured quantum groupoid w.r.t $\nu$
called \textbf{quantum space quantum groupoid}.
\end{theo}

\begin{proof}
By definition, $\Gamma$ is a morphism of Hopf bimodule. We have
to prove co-product relation. For all $m\in M$ and $n\in M'$, we
have:
$$
\begin{aligned}
(\Gamma\surl{\ _{\beta}\star_{\alpha}}_{\
\nu}id)\circ\Gamma(n\surl{\otimes}_{Z(M)}m )
&=[1\surl{\otimes}_{Z(M)}m]\surl{\ _{\beta}\otimes_{\alpha}}_{\
\nu}[1\surl{\otimes}_{Z(M)}1]
\surl{\ _{\beta}\otimes_{\alpha}}_{\ \nu}[n\surl{\otimes}_{Z(M)}1]\\
&=(id\surl{\ _{\beta}\star_{\alpha}}_{\
\nu}\Gamma)\circ\Gamma(n\surl{\otimes}_{Z(M)}m )
\end{aligned}$$
Now, we show that $T_R$ is right invariant and $\alpha$-adapted
w.r.t $\nu$. So, for all $m\in M,n\in M'$ and $\xi\in
D(_{\alpha}(L^2(M)\surl{\otimes}_{tr} L^2(M)),\nu^o)$, we put
$\Psi=\nu\circ\beta^{-1}\circ T_R$ and we compute:
$$
\begin{aligned}
\omega_{\xi}((\Psi\surl{\ _{\beta}\star_{\alpha}}_{\
\nu}id)\Gamma(n\surl{\otimes}_{Z(M)} m))&=\Psi((id\surl{\
_{\beta}\star_{\alpha}}_{\
\nu}\omega_{\xi})([1\surl{\otimes}_{Z(M)}m]\surl{\ _{\beta}
\otimes_{\alpha}}_{\ \nu}[n\surl{\otimes}_{Z(M)}1]))\\
&=\Psi([1\surl{\otimes}_{Z(M)}m]\beta(<[n\surl{\otimes}_{Z(M)}1]\xi,\xi>_{\alpha,\nu}))\\
&=\nu(<[n\surl{\otimes}_{Z(M)}T(m)]\xi,\xi>_{\alpha,\nu})\\
&=\omega_{\xi}(n\surl{\otimes}_{Z(M)}T(m))=\omega_{\xi}(T_R(n\surl{\otimes}_{Z(M)}m))
\end{aligned}$$ Finally, we have
for all $t\in\mathbb{R}$:
$$
\begin{aligned}
\sigma_t^{T_R}&=\sigma_t^{\nu'\surl{\star}_{Z(M)}\nu}\
_{|(M'\surl{\otimes}_{Z(M)} M)\cap
\beta(M)'}=\sigma_t^{\nu'\surl{\star}_{Z(M)}\nu}\
_{|(M'\surl{\otimes}_{Z(M)}
M)\cap (M\surl{\star}_{Z(M)}\mathcal{L}(L^2(M)))}\\
&=\sigma_t^{\nu'\surl{\star}_{Z(M)}\nu}\
_{|Z(M)\surl{\otimes}_{Z(M)}
M}=(id\surl{\otimes}_{Z(M)}\sigma_t^{\nu})\
_{|1\surl{\otimes}_{Z(M)}M}=1\surl{\otimes}_{Z(M)}\sigma_t^{\nu}
\end{aligned}$$ so that
$\sigma_t^{T_R}\circ\alpha(m)=1\surl{\otimes}_{Z(M)}\sigma_t^{\nu}(m)
=\alpha(\sigma_t^{\nu}(m))$ for all $t\in\mathbb{R}$ and $m\in
M$. Since it is easy to see that $R$ is a co-involution, we have
done.
\end{proof}

By \ref{raccourci}, we can compute the pseudo-multiplicative
unitary. Let first notice that $\Phi=\nu'\surl{\star}_{Z(M)}\nu
=\Psi$ so that $\lambda=\delta=1$ and:
$$\alpha=1\surl{\otimes}_{Z(M)}id, \hat{\alpha}=id\surl{\otimes}_{Z(M)}1,
\beta=j\surl{\otimes}_{Z(M)}1\text{ and }
\hat{\beta}=1\surl{\otimes}_{Z(M)}j$$ For example, we have
$D((H\surl{\otimes}_{tr} H)_{\hat{\beta},\nu^o})\supset
H\surl{\otimes}_{tr}
D(H_j,\nu^o)=H\surl{\otimes}_{tr}\Lambda_{\nu}({\mathcal
N}_{\nu})$ and for all $\eta\in H$ and $y\in {\mathcal N}_{\nu}$,
we have
$R^{\hat{\beta},\nu^o}(\eta\surl{\otimes}_{tr}\Lambda_{\nu}(y))=\lambda^{tr}_{\eta}
R^{j,\nu^o}(\Lambda_{\nu}(y))=\lambda^{tr}_{\eta}y$.

\begin{lemm}\label{simex2}
We have, for all $\eta\in H$ and $e\in {\mathcal N}_{\nu}$:
$$I\rho^{\beta,\alpha}_{\eta\surl{\otimes}_{tr}
J_{\nu}\Lambda_{\nu}(e)}=\lambda^{tr}_{\eta}J_{\nu}eJ_{\nu}
\surl{\otimes}_{Z(M)}1\text{ and }
I\lambda_{\Lambda_{\nu}(y)\otimes\eta}^{\beta,\alpha}=\rho_{\eta}^{tr}(1\surl{\otimes}_{Z(M)}
y)$$
\end{lemm}

\begin{proof}
Straightforward.
\end{proof}

\begin{prop}
We have, for all $\Xi\in H\surl{\otimes}_{tr}H,\eta\in H$ and
$m\in {\mathcal N}_{\nu}$:
$$W^*(\Xi\surl{\ _{\alpha} \otimes_{\hat{\beta}}}_{\
\nu^o}(\eta\surl{\otimes}_{tr}\Lambda_{\nu}(m)))
=I^*(\eta\surl{\otimes}_{tr}(1\surl{\otimes}_{Z(M)}m)\Xi)$$
\end{prop}

\begin{proof}
For all $m,e\in {\mathcal N}_{\nu}$ and $m',e'\in {\mathcal
N}_{\nu'}$, we have by the previous lemma:
$$
\begin{aligned}
I\Gamma(m'\surl{\otimes}_{Z(M)}m)\rho^{\beta,\alpha}_{J_{\nu'}
\Lambda_{\nu'}(e')\surl{\otimes}_{tr}
J_{\nu}\Lambda_{\nu}(e)}&=(m'\surl{\otimes}_{Z(M)}1\surl{\otimes}_{Z(M)}m)I
\rho^{\beta,\alpha}_{J_{\nu'}\Lambda_{\nu'}(e')\surl{\otimes}_{tr}
J_{\nu}\Lambda_{\nu}(e)}\\
&=(m'\otimes 1\otimes
m)\lambda^{tr}_{J_{\nu'}\Lambda_{\nu'}(e')}J_{\nu}eJ_{\nu}
\surl{\otimes}_{Z(M)}1\\
&=\lambda^{tr}_{J_{\nu'}e'J_{\nu'}\Lambda_{\nu'}(m')}J_{\nu}eJ_{\nu}\surl{\otimes}_{Z(M)}m
\end{aligned}$$
On the other hand, we have by \ref{renv2}:
$$
\begin{aligned}
&\ \quad I([1\surl{\otimes}_{Z(M)}1]\surl{\ _{\beta}
\otimes_{\alpha}}_{\ \nu} [J_{\nu'}e'J_{\nu'}\surl{\otimes}_{Z(M)}
J_{\nu}eJ_{\nu}])W^*
\rho^{\alpha,\hat{\beta}}_{\Lambda_{\nu'}(m')\surl{\otimes}_{tr}\Lambda_{\nu'}(m')}\\
&=(J_{\nu'}e'J_{\nu'}\surl{\otimes}_{Z(M)}J_{\nu}eJ_{\nu}\surl{\otimes}_{Z(M)}
1)IW^*\rho^{\alpha,\hat{\beta}}_{\Lambda_{\nu'}(m')\surl{\otimes}_{tr}\Lambda_{\nu'}(m')}
\end{aligned}$$
Then, by \ref{raccourci} and taking the limit over $e$ and $e'$
which go to $1$, we get for all $\Xi\in H\surl{\otimes}_{tr}H$:
$$W^*(\Xi\surl{\ _{\alpha} \otimes_{\hat{\beta}}}_{\
\nu^o}(\Lambda_{\nu'}(m')\surl{\otimes}_{tr}\Lambda_{\nu}(m)))
=I^*(\Lambda_{\nu'}(m')\surl{\otimes}_{tr}(1\surl{\otimes}_{Z(M)}m)\Xi)$$
Now, if $\Xi\in D(_{\alpha}(H\surl{\otimes}_{tr}H),\nu)$, by
continuity and density of $\Lambda_{\nu'}({\mathcal N}_{\nu'})$ we
have for all $\Xi\in D(_{\alpha}(H\surl{\otimes}_{tr}H),\nu)$:
$$W^*(\Xi\surl{\ _{\alpha} \otimes_{\hat{\beta}}}_{\
\nu^o}(\eta\surl{\otimes}_{tr}\Lambda_{\nu}(m)))
=I^*(\eta\surl{\otimes}_{tr}(1\surl{\otimes}_{Z(M)}m)\Xi)$$ Since
$\eta\surl{\otimes}_{tr}\Lambda_{\nu}(m)\in
D((H\surl{\otimes}_{tr} H)_{\hat{\beta},\nu^o})$, the relation
holds by continuity for all $\Xi\in H\surl{\otimes}_{tr}H$.
\end{proof}

\begin{rema}
If $\sigma_{tr}$ is the flip of $L^2(M)\surl{\otimes}_{tr}L^2(M)$,
then $\sigma_{tr}\circ\hat{\beta}=\beta\circ\sigma_{tr}$ and if
$I'=(1\surl{\otimes}_{Z(M)}\sigma_{tr})I(\sigma_{tr}\surl{\
_{\hat{\beta}}\otimes_{\alpha}}_{\ \nu} [1\surl{\otimes}_{Z(M)}
1])\sigma_{\nu^o}$, then $I'$ is the identification:

$$
\begin{aligned}
I': [L^2(M)\surl{\otimes}_{tr}L^2(M)]\surl{\ _{\alpha}
\otimes_{\hat{\beta}}}_{\ \nu^o}[L^2(M)\surl{\otimes}_{tr}
L^2(M)]&\rightarrow L^2(M)\surl{\otimes}_{tr}L^2(M)
\surl{\otimes}_{tr}L^2(M)\\
\Xi\surl{\ _{\beta}\otimes_{\alpha}}_{\
\nu}[\eta\surl{\otimes}_{tr}\Lambda_{\nu}(y)]&\mapsto\eta\surl{\otimes}_{tr}\alpha(y)\Xi
\end{aligned}$$
for all $\eta\in L^2(M),\Xi\in L^2(M)\surl{\otimes}_{tr}L^2(M)$
and $y\in M$. Consequently, by the previous proposition
$W^*=I^*I'$.
\end{rema}

\begin{coro}
We can reconstruct the von Neumann algebra thanks to $W$:
$$M'\surl{\otimes}_{Z(M)}
M=<(id*\omega_{\xi,\eta})(W^*)\,|\,\xi\in
D((H\surl{\otimes}_{tr}H)_{\hat{\beta}},\nu^o),\eta\in
D(_{\alpha}(H\surl{\otimes}_{tr}H),\nu)>^{-\textsc{w}}$$
\end{coro}

\begin{proof}
By \ref{appartenance}, we know that:
$$<(id*\omega_{\xi,\eta})(W^*)\,|\,\xi\in
D((H\surl{\otimes}_{tr}H)_{\hat{\beta}},\nu^o,\eta\in
D(_{\alpha}(H\surl{\otimes}_{tr}H ),\nu)>^{-\textsc{w}}\subset
M'\surl{\otimes}_{Z(M)}M$$ Let $\eta,\xi\in H$ and $m,e\in
{\mathcal N}_{\nu}$. Then, for all $\Xi_1,\Xi_2\in
H\surl{\otimes}_{tr}H$, we have by \ref{simex2}:
$$
\begin{aligned}
&\ \quad
((id*\omega_{\eta\surl{\otimes}_{tr}\Lambda_{\nu}(m),\xi\surl{\otimes}_{tr}
J_{\nu}\Lambda_\nu(e)})(W^*)\Xi_1|\Xi_2)\\
&=(W^*(\Xi_1\surl{\ _{\alpha} \otimes_{\hat{\beta}}}_{\ \nu^o}
[\eta\surl{\otimes}_{tr}\Lambda_{\nu}(m)])|\Xi_2\surl{\ _{\beta}
\otimes_{\alpha}}_{\ \nu}[\xi\surl{\otimes}_{tr}J_{\nu}\Lambda_\nu(e)])\\
&=(I^*(\eta\surl{\otimes}_{tr}(1\surl{\otimes}_{Z(M)}m)\Xi_1)|\Xi_2\surl{\
_{\beta}\otimes_{\alpha}}_{\ \nu}[\xi\surl{\otimes}_{tr}J_{\nu}\Lambda_\nu(e)])\\
&=(\eta\surl{\otimes}_{tr}(1\surl{\otimes}_{Z(M)}m)\Xi_1|\xi\surl{\otimes}_{tr}
(J_{\nu}eJ_{\nu}\surl{\otimes}_{Z(M)}1)\Xi_2)\\
&=((<\eta,\xi>_{tr}J_{\nu}e^*J_{\nu}\otimes m)\Xi_1|\Xi_2)
\end{aligned}$$
Consequently, we get the reverse inclusion thanks to the relation:
$$(id*\omega_{\eta\surl{\otimes}_{tr}\Lambda_{\nu}(m),\xi\surl{\otimes}_{tr}
J_{\nu}\Lambda_\nu(e)})(W^*)=<\eta,\xi>_{tr}J_{\nu}e^*J_{\nu}\surl{\otimes}_{Z(M)}
m$$
\end{proof}

Now, we compute $G$ so as to get the antipode.

\begin{prop}
If $F_{\nu}=S_{\nu}^*$ comes from Tomita's theory, then we have:
$$G=\sigma_{tr}\circ (F_{\nu}\surl{\otimes}_{tr}F_{\nu})$$
\end{prop}

\begin{proof}
Let $a=J_{\nu}a_1J_{\nu}\surl{\otimes}_{Z(M)}
a_2,b=J_{\nu}b_1J_{\nu}\surl{\otimes}_{Z(M)}
b_2,c=J_{\nu}c_1J_{\nu}\surl{\otimes}_{Z(M)} c_2$ and
$d=J_{\nu}d_1J_{\nu}\surl{\otimes}_{Z(M)} d_2$ be elements of
$M'\surl{\otimes}_{Z(M)} M$ analytic w.r.t
$\nu'\surl{\star}_{Z(M)}\nu$. Then, by \ref{simex2}, the value of
$(\lambda^{\beta,\alpha}_{\Lambda_{\nu}(\sigma_{i/2}^{\nu}(b_1))\surl{\otimes}_{tr}
\Lambda_{\nu}(\sigma_{-i}^{\nu}(b_2^*))})^*W^*$ on
$$[\Lambda_{\nu'}(J_{\nu}a_1J_{\nu})
\surl{\otimes}_{tr}\Lambda_{\nu}(a_2)]\surl{\ _{\alpha}
\otimes_{\hat{\beta}}}_{\
\nu^o}[\Lambda_{\nu'}(J_{\nu}d_1^*c_1^*J_{\nu})
\surl{\otimes}_{tr}\Lambda_{\nu}(d_2^*c_2^*)]$$ is equal to:
$$
\begin{aligned}
&\ \quad
(\lambda^{\beta,\alpha}_{\Lambda_{\nu}(\sigma_{i/2}^{\nu}(b_1))\surl{\otimes}_{tr}
\Lambda_{\nu}(\sigma_{-i}^{\nu}(b_2^*))})^*I^*(\Lambda_{\nu'}(J_{\nu}d_1^*c_1^*J_{\nu})
\surl{\otimes}_{tr}\Lambda_{\nu'}(J_{\nu}a_1J_{\nu})\surl{\otimes}_{tr}\Lambda_{\nu}(d_2^*c_2^*a_2))\\
&=\!\left[\rho^{tr}_{\Lambda_{\nu}(\sigma_{-i}^{\nu}(b_2^*))}(1\surl{\otimes}_{Z(M)}
\sigma_{i/2}^{\nu}(b_1))\right]^*\!\!\!\!(\Lambda_{\nu'}(J_{\nu}d_1^*c_1^*J_{\nu})
\surl{\otimes}_{tr}\Lambda_{\nu'}(J_{\nu}a_1J_{\nu})\surl{\otimes}_{tr}\Lambda_{\nu}(d_2^*c_2^*a_2))\\
&=<d_2^*c_2^*\Lambda_{\nu}(a_2),\Lambda_{\nu}(\sigma_{-i}^{\nu}(b_2^*))>_{tr}\,
\Lambda_{\nu'}(J_{\nu}d_1^*c_1^*J_{\nu})\surl{\otimes}_{tr}\sigma_{-i/2}^{\nu}(b_1^*)
\Lambda_{\nu'}(J_{\nu}a_1J_{\nu})\\
&=<\Lambda_{\nu}(a_2b_2),\Lambda_{\nu}(c_2d_2)>_{tr}\,
J_{\nu}\Lambda_{\nu}(d_1^*c_1^*)\surl{\otimes}_{tr}
J_{\nu}\Lambda_{\nu}(a_1b_1)
\end{aligned}$$
Consequently, by definition of $G$:
$$G\left[<\Lambda_{\nu}(a_2b_2),\Lambda_{\nu}(c_2d_2)>_{tr}\,
J_{\nu}\Lambda_{\nu}(d_1^*c_1^*)\surl{\otimes}_{tr}
J_{\nu}\Lambda_{\nu}(a_1b_1)\right]$$ is equal to the value of
$G(\lambda^{\beta,\alpha}_{\Lambda_{\nu}(\sigma_{i/2}^{\nu}(b_1))\surl{\otimes}_{tr}
\Lambda_{\nu}(\sigma_{-i}^{\nu}(b_2^*))})^*W^*$ on:
$$[\Lambda_{\nu'}(J_{\nu}a_1J_{\nu})
\surl{\otimes}_{tr}\Lambda_{\nu}(a_2)]\surl{\ _{\alpha}
\otimes_{\hat{\beta}}}_{\
\nu^o}[\Lambda_{\nu'}(J_{\nu}d_1^*c_1^*J_{\nu})
\surl{\otimes}_{tr}\Lambda_{\nu}(d_2^*c_2^*)]$$ which is equal to
the value of
$(\lambda^{\beta,\alpha}_{\Lambda_{\nu}(\sigma_{i/2}^{\nu}(d_1))\surl{\otimes}_{tr}
\Lambda_{\nu}(\sigma_{-i}^{\nu}(d_2^*))})^*W^*$ on:
$$[\Lambda_{\nu'}(J_{\nu}c_1J_{\nu})
\surl{\otimes}_{tr}\Lambda_{\nu}(c_2)]\surl{\ _{\alpha}
\otimes_{\hat{\beta}}}_{\
\nu^o}[\Lambda_{\nu'}(J_{\nu}b_1^*a_1^*J_{\nu})
\surl{\otimes}_{tr}\Lambda_{\nu}(b_2^*a_2^*)]$$ This last vector
is $<\Lambda_{\nu}(c_2d_2),\Lambda_{\nu}(a_2b_2)>_{tr}\,
J_{\nu}\Lambda_{\nu}(b_1^*a_1^*)\surl{\otimes}_{tr}
J_{\nu}\Lambda_{\nu}(c_1d_1)$. Since $G$ is closed, we get:
$$G\left[J_{\nu}\Lambda_{\nu}(d_1^*c_1^*)\surl{\otimes}_{tr}
J_{\nu}\Lambda_{\nu}(a_1b_1)\right]=\left[
J_{\nu}\Lambda_{\nu}(b_1^*a_1^*)\surl{\otimes}_{tr}
J_{\nu}\Lambda_{\nu}(c_1d_1)\right]$$ so that $G$ coincides with
$\sigma_{tr}(F_{\nu}\surl{\otimes}_{tr}F_{\nu})$.
\end{proof}

The polar decomposition of $G=ID^{1/2}$ is such that
$D=\Delta_{\nu}^{-1}\surl{\otimes}_{tr}\Delta_{\nu}^{-1}$ and
$I=\sigma_{tr}(J_{\nu}\surl{\otimes}_{tr} J_{\nu})$ so that the
scaling group is
$\tau_t=\sigma_{-t}^{\nu'}\surl{\star}_{Z(M)}\sigma_t^{\nu}$ for
all $t\in\mathbb{R}$ and the unitary antipode is
$R=\varsigma_{Z(M)}\circ (j\surl{\otimes}_{Z(M)}j)$. We also
notice that $\nu'\surl{\star}_{Z(M)}\nu$ is $\tau$-invariant.

\begin{rema}
If $M$ is the commutative von Neumann algebra $L^{\infty}(X)$,
then the structure coincides with the quantum space $X$.
\end{rema}

\subsection{Pairs quantum groupoid}

Let $M$ be a von Neumann algebra. $M$ acts on
$H=L^2(M)=L^2_{\nu}(M)$ where $\nu$ is a n.s.f weight on $M$. We
denote by $M'$ the commutant of $M$ in $\mathcal{L}(L^2(M))$.
$M'\otimes M$ acts on $L^2(M)\otimes L^2(M)$.

Let $\alpha$ (resp. $\beta$) be the (resp. anti-) representation
of $M$ to $M'\otimes M$ such that $\alpha(m)=1\otimes m$ (resp.
$\beta(m)=j(m)\otimes 1$) where $j(x)=J_{\nu}x^*J_{\nu}$ for all
$x\in \mathcal{L}(L^2_{\nu}(M))$.

\begin{prop}\label{renv}
The following formula:
$$
\begin{aligned}
I: [L^2(M)\otimes L^2(M)]\surl{\ _{\beta}
  \otimes_{\alpha}}_{\ \nu}[L^2(M)\otimes L^2(M)]&\rightarrow L^2(M)\otimes L^2(M)
  \otimes L^2(M)\\
[\Lambda_{\nu}(y)\otimes\eta]\surl{\ _{\beta}
  \otimes_{\alpha}}_{\ \nu}\Xi&\mapsto \alpha(y)\Xi\otimes\eta
\end{aligned}$$
for all $\eta\in L^2(M),\Xi\in L^2(M)\otimes L^2(M)$ and $y\in M$,
defines a canonical isomorphism such that we have $I([m\otimes
x]\surl{\ _{\beta}\otimes_{\alpha}}_{\ \nu}[y\otimes
n])=(y\otimes mn\otimes x)I$, for all $m\in M, n\in M'$ and
$x,y\in \mathcal{L}(L^2(M))$.
\end{prop}

\begin{proof}
Straightforward.
\end{proof}
Then, we can identify $(M'\otimes M)\surl{\
_{\beta}\star_{\alpha}}_{\ M}(M'\otimes M)$ with $M'\otimes
Z(M)\otimes M$. We define a normal *-homomorphism $\Gamma$ by:
$$
\begin{aligned}
M'\otimes M &\rightarrow (M'\otimes M)
\surl{\ _{\beta}\star_{\alpha}}_{\ \nu} (M'\otimes M)\\
n\otimes m&\mapsto I^*(n\otimes 1\otimes m)I=[1\otimes m] \surl{\
_{\beta}\otimes_{\alpha}}_{\ \nu}[n\otimes 1]
\end{aligned}$$

\begin{theo}
$(M,M'\surl{\otimes}_{Z(M)}M,\alpha,\beta,\Gamma,\nu,\nu'\otimes
id,id\otimes\nu)$ is a measured quantum groupoid w.r.t $\nu$
called \textbf{pairs quantum groupoid}.
\end{theo}

\begin{proof}
By definition, $\Gamma$ is a morphism of Hopf bimodule. We have
to prove co-product relation. For all $m\in M$ and $n\in M'$, we
have:
$$
\begin{aligned}
(\Gamma\surl{\ _{\beta}\star_{\alpha}}_{\
\nu}id)\circ\Gamma(n\otimes m )&=[1\otimes m]\surl{\
_{\beta}\otimes_{\alpha}}_{\ \nu}[1\otimes
1]\surl{\ _{\beta}\otimes_{\alpha}}_{\ \nu}[n\otimes 1]\\
&=(id\surl{\ _{\beta}\star_{\alpha}}_{\
\nu}\Gamma)\circ\Gamma(n\otimes m)
\end{aligned}$$
$R=\varsigma\circ (\beta_{\nu}\otimes\beta_{\nu})$, where
$\varsigma: M'\otimes M\rightarrow M\otimes M'$ is the flip, is a
co-involution so it is sufficient to show that $T_L=\nu'\otimes
id$ is left invariant and $\beta$-adapted w.r.t $\nu$. Let $m\in
M,n\in M'$ and $\xi\in D((L^2(M)\otimes L^2(M))_{\beta,\nu^o})$.
We put $\Phi=\nu\circ\alpha^{-1}\circ T_L$ and we compute:
$$
\begin{aligned}
\omega_{\xi}((id\surl{\ _{\beta}\star_{\alpha}}_{\
\nu}\Phi)\Gamma(n\otimes m))&=\Phi((\omega_{\xi}\surl{\
_{\beta}\star_{\alpha}}_{\ \nu}id)([1\otimes m]\surl{\ _{\beta}
  \otimes_{\alpha}}_{\ \nu}[n\otimes 1]))\\
&=\Phi([n\otimes 1]\alpha(<[1\otimes
m]\xi,\xi>_{\beta,\nu^o}))\\
&=\nu'(n)\nu(<[1\otimes
m]\xi,\xi>_{\beta,\nu^o})\\
&=\nu'(n)\omega_{\xi}(1\otimes m)=\omega_{\xi}(T_L(n\otimes m))
\end{aligned}$$ Finally, we prove that
$T_R=R\circ T_L\circ R=id\otimes \nu$ is $\alpha$-adapted w.r.t
$\nu$. For all $t\in\mathbb{R}$, we have:
$$\sigma_t^{T_R}=\sigma_t^{\nu'\otimes\nu}\ _{|(M'\otimes
M)\cap \beta(M)'}=\sigma_t^{\nu'\otimes\nu}\ _{|Z(M)\otimes
M}=id\otimes\sigma_t^{\nu}\ _{|Z(M)\otimes M}$$ so that we have
for all $t\in\mathbb{R}$ and $m\in M$:
$$\sigma_t^{T_R}\circ\alpha(m)=1\otimes\sigma_t^{\nu}(m)=\alpha(\sigma_t^{\nu}(m))$$
\end{proof}

\begin{rema}
If $M=L^{\infty}(X)$, we find the structure of pairs groupoid
$X\times X$.
\end{rema}

By \ref{raccourci}, we can compute the pseudo-multiplicative
unitary. Let first notice that $\Phi=\nu'\otimes\nu =\Psi$ so
that $\lambda=\delta=1$ and:
$$\alpha=1\otimes id, \hat{\alpha}=id\otimes 1, \beta=\beta_{\nu}\otimes 1
\text{ and } \hat{\beta}=1\otimes \beta_{\nu}$$ For example, we
have $D((H\otimes H)_{\hat{\beta},\nu^o})\supset H\otimes
D(H_{\beta_{\nu}},\nu^o)=H\otimes \Lambda_{\nu}({\mathcal
N}_{\nu})$ and for all $\eta\in H$ and $y\in {\mathcal N}_{\nu}$,
we have
$R^{\hat{\beta},\nu^o}(\eta\otimes\Lambda_{\nu}(y))=\lambda_{\eta}
R^{\beta_{\nu},\nu^o}(\Lambda_{\nu}(y))=\lambda_{\eta}y$.

\begin{lemm}\label{simex}
We have, for all $\eta\in H$ and $e\in {\mathcal N}_{\nu}$:
$$I\rho^{\beta,\alpha}_{\eta\otimes J_{\nu}\Lambda_{\nu}(e)}=
\lambda_{\eta}J_{\nu}eJ_{\nu} \otimes 1\text{ and }
I\lambda_{\Lambda_{\nu}(y)\otimes\eta}^{\beta,\alpha}=\rho_{\eta}(1\otimes
y)$$
\end{lemm}

\begin{proof}
Straightforward.
\end{proof}

\begin{prop}
We have, for all $\Xi\in H\otimes H,\eta\in H$ and $m\in {\mathcal
N}_{\nu}$:
$$W^*(\Xi\surl{\ _{\alpha} \otimes_{\hat{\beta}}}_{\
\nu^o}(\eta\otimes\Lambda_{\nu}(m))) =I^*(\eta\otimes(1\otimes
m)\Xi)$$
\end{prop}

\begin{proof}
For all $m,e\in {\mathcal N}_{\nu}$ and $m',e'\in {\mathcal
N}_{\nu'}$, we have by the previous lemma:
$$
\begin{aligned}
I\Gamma(m'\otimes
m)\rho^{\beta,\alpha}_{J_{\nu'}\Lambda_{\nu'}(e')\otimes
J_{\nu}\Lambda_{\nu}(e)}&=(m'\otimes 1\otimes m)I
\rho^{\beta,\alpha}_{J_{\nu'}\Lambda_{\nu'}(e')\otimes
J_{\nu}\Lambda_{\nu}(e)}\\
&=(m'\otimes 1\otimes
m)\lambda_{J_{\nu'}\Lambda_{\nu'}(e')}J_{\nu}eJ_{\nu}\otimes 1\\
&=\lambda_{J_{\nu'}e'J_{\nu'}\Lambda_{\nu'}(m')}
J_{\nu}eJ_{\nu}\otimes m
\end{aligned}$$
On the other hand, we have by \ref{renv}:
$$
\begin{aligned}
&\ \quad I([1\otimes 1]\surl{\ _{\beta} \otimes_{\alpha}}_{\ \nu}
[J_{\nu'}e'J_{\nu'}\otimes J_{\nu}eJ_{\nu}])W^*
\rho^{\alpha,\hat{\beta}}_{\Lambda_{\nu'}(m')\otimes\Lambda_{\nu'}(m')}\\
&=(J_{\nu'}e'J_{\nu'}\otimes J_{\nu}eJ_{\nu}\otimes 1)
IW^*\rho^{\alpha,\hat{\beta}}_{\Lambda_{\nu'}(m')\otimes\Lambda_{\nu'}(m')}
\end{aligned}$$
Then by \ref{raccourci} and taking the limit over $e$ and $e'$
which go to $1$, we get for all $\Xi\in H\otimes H$:
$$W^*(\Xi\surl{\ _{\alpha} \otimes_{\hat{\beta}}}_{\
\nu^o}(\Lambda_{\nu'}(m')\otimes\Lambda_{\nu}(m)))
=I^*(\Lambda_{\nu'}(m')\otimes(1\otimes m)\Xi)$$ Now, if $\Xi\in
D(_{\alpha}(H\otimes H),\nu)$, by continuity and density of
$\Lambda_{\nu'}({\mathcal N}_{\nu'})$, we have for all $\Xi\in
D(_{\alpha}(H\otimes H),\nu)$:
$$W^*(\Xi\surl{\ _{\alpha} \otimes_{\hat{\beta}}}_{\
\nu^o}(\eta\otimes\Lambda_{\nu}(m))) =I^*(\eta\otimes(1\otimes
m)\Xi)$$ Since $\eta\otimes\Lambda_{\nu}(m)\in D((H\otimes
H)_{\hat{\beta},\nu^o})$, the previous relation holds by
continuity for all $\Xi\in H\otimes H$.
\end{proof}

\begin{rema}
If $\sigma$ denotes the flip of $L^2(M)\otimes L^2(M)$, then
$\sigma\circ\hat{\beta}=\beta\circ\sigma$ and if
$I'=(1\otimes\sigma)I (\sigma\surl{\ _{\hat{\beta}}
\otimes_{\alpha}}_{\ \nu} [1\otimes 1])\sigma_{\nu^o}$, then $I'$
is the identification:
$$
\begin{aligned}
I': [L^2(M)\otimes L^2(M)]\surl{\ _{\alpha}
\otimes_{\hat{\beta}}}_{\ \nu^o}[L^2(M)\otimes
L^2(M)]&\rightarrow L^2(M)\otimes L^2(M)
  \otimes L^2(M)\\
\Xi\surl{\ _{\beta}
  \otimes_{\alpha}}_{\ \nu}[\eta\otimes\Lambda_{\nu}(y)]&\mapsto
  \eta\otimes\alpha(y)\Xi
\end{aligned}$$
for all $\eta\in L^2(M),\Xi\in L^2(M)\otimes L^2(M)$ and $y\in M$.
Consequently, by the previous proposition $W^*=I^*I'$.
\end{rema}

\begin{coro}
We can re-construct the underlying von Neumann algebra thanks to
$W$:
$$M'\otimes M=<(id*\omega_{\xi,\eta})(W^*)\,|\,\xi\in D((H\otimes
H)_{\hat{\beta}},\nu^o),\eta\in D(_{\alpha}(H\otimes
H),\nu)>^{-\textsc{w}}$$
\end{coro}

\begin{proof}
By \ref{appartenance}, we know that:
$$<(id*\omega_{\xi,\eta})(W^*)\,|\,\xi\in
D((H\otimes H)_{\hat{\beta}},\nu^o,\eta\in D(_{\alpha}(H\otimes H
),\nu)>^{-w}\subset M'\otimes M$$ Let $\eta,\xi\in H$ and $m,e\in
{\mathcal N}_{\nu}$. Then, for all $\Xi_1,\Xi_2\in H\otimes H$, we
have, by \ref{simex}:
$$
\begin{aligned}
&\ \quad((id*\omega_{\eta\otimes\Lambda_{\nu}(m),\xi\otimes
J_{\nu}\Lambda_\nu(e)})(W^*)\Xi_1|\Xi_2)\\
&=(W^*(\Xi_1\surl{\ _{\alpha} \otimes_{\hat{\beta}}}_{\ \nu^o}
\eta\otimes\Lambda_{\nu}(m))|\Xi_2\surl{\ _{\beta}
\otimes_{\alpha}}_{\ \nu}\xi\otimes J_{\nu}\Lambda_\nu(e))\\
&=(I^*(\eta\otimes(1\otimes m)\Xi_1)|\Xi_2\surl{\ _{\beta}
\otimes_{\alpha}}_{\ \nu}\xi\otimes J_{\nu}\Lambda_\nu(e))\\
&=(\eta\otimes(1\otimes m)\Xi_1|\xi\otimes (J_{\nu}eJ_{\nu}\otimes
1)\Xi_2)\\
&=(\eta|\xi)((J_{\nu}e^*J_{\nu}\otimes m)\Xi_1|\Xi_2)
\end{aligned}$$
Consequently, we get the reverse inclusion thanks to the relation:
\begin{equation}\label{equ1}
(id*\omega_{\eta\otimes\Lambda_{\nu}(m),\xi\otimes
J_{\nu}\Lambda_\nu(e)})(W^*)=(\eta|\xi)(J_{\nu}e^*J_{\nu}\otimes
m)
\end{equation}
\end{proof}

Now, we compute $G$ so as to get the antipode.

\begin{prop}
If $F_{\nu}=S_{\nu}^*$ comes from Tomita's theory, we have:
$$G=\sigma(F_{\nu}\otimes F_{\nu})$$
\end{prop}

\begin{proof}
For all $a=J_{\nu}a_1J_{\nu}\otimes a_2,b=J_{\nu}b_1J_{\nu}\otimes
b_2,c=J_{\nu}c_1J_{\nu}\otimes c_2$ and
$d=J_{\nu}d_1J_{\nu}\otimes d_2$ be analytic elements of
$M'\otimes M$ w.r.t $\nu'\otimes\nu$. Then, by \ref{simex}, we
have:
$$
\begin{aligned}
&\
\quad(\lambda^{\beta,\alpha}_{\Lambda_{\nu}(\sigma_{i/2}^{\nu}(b_1))\otimes
\Lambda_{\nu}(\sigma_{-i}^{\nu}(b_2^*))})^*W^*(\Lambda_{\nu'\otimes\nu}(a)\surl{\
_{\alpha} \otimes_{\hat{\beta}}}_{\
\nu^o}\Lambda_{\nu'\otimes\nu}((J_{\nu}d_1^*J_{\nu}\otimes
d^*_2)c^*))\\
&=(\lambda^{\beta,\alpha}_{\Lambda_{\nu}(\sigma_{i/2}^{\nu}(b_1))\otimes
\Lambda_{\nu}(\sigma_{-i}^{\nu}(b_2^*))})^*I^*(\Lambda_{\nu'}(J_{\nu}d_1^*c_1^*J_{\nu})
\otimes (1\otimes d_2^*c_2^*)\Lambda_{\nu'\otimes\nu}(a))\\
&=\left[\rho_{\Lambda_{\nu}(\sigma_{-i}^{\nu}(b_2^*))}(1\otimes
\sigma_{i/2}^{\nu}(b_1))\right]^*(\Lambda_{\nu'}(J_{\nu}d_1^*c_1^*J_{\nu})
\otimes (1\otimes d_2^*c_2^*)\Lambda_{\nu'\otimes\nu}(a))\\
&=(d_2^*c_2^*\Lambda_{\nu}(a_2)|\Lambda_{\nu}(\sigma_{-i}^{\nu}(b_2^*)))\,
\Lambda_{\nu'}(J_{\nu}d_1^*c_1^*J_{\nu})\otimes\sigma_{-i/2}^{\nu}(b_1^*)
\Lambda_{\nu'}(J_{\nu}a_1J_{\nu})\\
&=\nu(d_2^*c_2^*a_2b_2)\, J_{\nu}\Lambda_{\nu}(d_1^*c_1^*)\otimes
J_{\nu}\Lambda_{\nu}(a_1b_1)
\end{aligned}$$
Consequently, by definition of $G$, we have:
$$
\begin{aligned}
&\quad G\left[\nu(d_2^*c_2^*a_2b_2)\,
J_{\nu}\Lambda_{\nu}(d_1^*c_1^*)\otimes
J_{\nu}\Lambda_{\nu}(a_1b_1)\right]\\
&=G(\lambda^{\beta,\alpha}_{\Lambda_{\nu}(\sigma_{i/2}^{\nu}(b_1))\otimes
\Lambda_{\nu}(\sigma_{-i}^{\nu}(b_2^*))})^*W^*(\Lambda_{\nu'\otimes\nu}(a)\surl{\
_{\alpha} \otimes_{\hat{\beta}}}_{\
\nu^o}\Lambda_{\nu'\otimes\nu}((J_{\nu}d_1^*J_{\nu}\otimes
d^*_2)c^*))\\
&=(\lambda^{\beta,\alpha}_{\Lambda_{\nu}(\sigma_{i/2}^{\nu}(d_1))\otimes
\Lambda_{\nu}(\sigma_{-i}^{\nu}(d_2^*))})^*W^*(\Lambda_{\nu'\otimes\nu}(c)\surl{\
_{\alpha} \otimes_{\hat{\beta}}}_{\
\nu^o}\Lambda_{\nu'\otimes\nu}((J_{\nu}b_1^*J_{\nu}\otimes
b^*_2)a^*))\\
&=\nu(b_2^*a_2^*c_2d_2)\, J_{\nu}\Lambda_{\nu}(b_1^*a_1^*)\otimes
J_{\nu}\Lambda_{\nu}(c_1d_1)
\end{aligned}$$
Since $G$ is anti-linear, we get:
$$G\left[J_{\nu}\Lambda_{\nu}(d_1^*c_1^*)\otimes
J_{\nu}\Lambda_{\nu}(a_1b_1)\right]=\left[
J_{\nu}\Lambda_{\nu}(b_1^*a_1^*)\otimes
J_{\nu}\Lambda_{\nu}(c_1d_1)\right]$$ so that $G$ coincides with
$\sigma(F_{\nu}\otimes F_{\nu})$.
\end{proof}

The polar decomposition of $G=ID^{1/2}$ is such that
$D=\Delta_{\nu}^{-1}\otimes\Delta_{\nu}^{-1}$ and
$I=\Sigma(J_{\nu}\otimes J_{\nu})$ so that the scaling group is
$\tau_t=\sigma_{-t}^{\nu'}\otimes\sigma_t^{\nu}$ for all
$t\in\mathbb{R}$ and the unitary antipode is $R=\varsigma\circ
(\beta_{\nu}\otimes\beta_{\nu})$. We also notice that
$\nu'\otimes\nu$ is $\tau$-invariant.

\begin{coro}
We have ${\mathcal D}(S)={\mathcal D}(\sigma_{i/2}^{\nu'})\otimes
{\mathcal D}(\sigma_{-i/2}^{\nu})$ and we have
$$S(J_{\nu}eJ_{\nu}\otimes
m^*)=J_{\nu}\sigma_{i/2}^{\nu}(m)J_{\nu}\otimes\sigma_{-i/2}^{\nu}(e^*)$$
for all $e,m\in {\mathcal D}(\sigma_{i/2}^{\nu})$. Moreover
$(id*\omega_{\xi,\eta})(W)\in {\mathcal D}(S)$ and:
$$S((id*\omega_{\xi,\eta})(W))=(id*\omega_{\xi,\eta})(W^*)$$ for
all $\xi,\eta\in D(_{\alpha}(H\otimes H),\nu)\cap D((H\otimes
H)_{\hat{\beta}},\nu^o)$.
\end{coro}

\begin{proof}
The first part of the corollary is straightforward by what
precedes. Let $\zeta,\eta\in H$ and $e,m\in {\mathcal
D}(\sigma_{i/2}^{\nu})$. By \ref{equ1}, we have:
$$\begin{aligned}
S((id*\omega_{\zeta\otimes J_{\nu}\Lambda_{\nu}(e),\eta\otimes
\Lambda_{\nu}(m)})(W))&=S((\zeta|\eta)J_{\nu}eJ_{\nu}\otimes
m^*)\\
&=(\zeta|\eta)J\sigma_{i/2}^{\nu}(m)J\otimes\sigma_{-i/2}^{\nu}(e^*)\\
&=(id*\omega_{\zeta\otimes J_{\nu}\Lambda_{\nu}(e),\eta\otimes
\Lambda_{\nu}(m)})(W^*)
\end{aligned}$$
Since $S$ is closed, we can conclude.
\end{proof}

\subsection{Operations on measured quantum groupoids}

\subsubsection{Sum of measured quantum groupoids}\label{sgqm}

A union of groupoids is still a groupoid. We establish here a
similar result at the quantum level:

\begin{prop}
Let $(N_i,M_i,\alpha_i,\beta_i,\Gamma_i,\nu_i,T_L^i,T_R^i)_{i\in
I}$ be a family of measured quantum groupoids. Then, identifying
the von Neumann algebra $\bigoplus_{i\in I}M_i\surl{\
_{\beta}\star_{\alpha}}_{\ N_i}M_i$ with $\left(\bigoplus_{i\in
I}M_i\right)\surl{\ _{\beta}\star_{\alpha}}_{\ \bigoplus_{i\in
I}N_i}\left(\bigoplus_{i\in I}M_i\right)$, we get:
$$(\bigoplus_{i\in I}N_i,\bigoplus_{i\in I}M_i,\bigoplus_{i\in I}\alpha_i,\bigoplus_{i\in I}\beta_i,
\bigoplus_{i\in I}\Gamma_i,\bigoplus_{i\in I}\nu_i,\bigoplus_{i\in
I}T_L^i,\bigoplus_{i\in I}T_R^i)$$ a measured quantum groupoid
where operators act on the diagonal.
\end{prop}

\begin{proof}
Straightforward.
\end{proof}

In particular, the sum of two quantum groups with different
scaling constants (\cite{VaV} for examples) produce a measured
quantum groupoid with non scalar scaling operator.

\subsubsection{Tensor product of measured quantum groupoids}
Cartesian product of groups correspond to tensor product of
quantum groups. In the same way, we have:

\begin{prop}
Let $(N_i,M_i,\alpha_i,\beta_i,\Gamma_i,\nu_i,T_L^i,T_R^i)$ be
measured quantum groupoids for $i=1,2$. Then, if we identify
$(M_1\surl{\ _{\beta_1}\star_{\alpha_1}}_{\ N_1}M_1)\otimes
(M_2\surl{\ _{\beta_2}\star_{\alpha_2}}_{\ N_2}M_2)$ with the von
Neumann algebra $(M_1\otimes M_2)\surl{\
_{\beta_1\otimes\beta_2}\star_{\alpha_1\otimes\alpha_2}}_{\
N_1\otimes N_2}(M_1\otimes M_2)$, we have:
$$(N_1\otimes N_2,M_1\otimes M_2,\alpha_1\otimes\alpha_2,\beta_1\otimes\beta_2,
\Gamma_1\otimes\Gamma_2,\nu_1\otimes\nu_2,T_L^1\otimes
T_L^2,T_R^1\otimes T_R^2)$$ is a measured quantum groupoid.
\end{prop}

\begin{proof}
Straightforward.
\end{proof}

\subsubsection{Direct integrals of measured
quantum groupoids}

In this section, $X$ denote $\sigma$-compact, locally compact
space and $\mu$ a Borel measure on $X$. Theory of hilbertian
integrals is described in \cite{T}.

\begin{prop}
Let $(N_p,M_p,\alpha_p,\beta_p,\Gamma_p,\nu_p,T_L^p,T_R^p)_{p\in
X}$ be a family of measured quantum groupoids. If we identify the
von Neumann algebras $\int_X^{\oplus}\!M_p\surl{\
_{\beta}\star_{\alpha}}_{\ N_p}M_p\,d\mu(p)$ and
$\left(\int_X^{\oplus}\!M_p\,d\mu(p)\right)\surl{\
_{\beta}\star_{\alpha}}_{\
\int_X^{\oplus}\!N_p\,d\mu(p)}\left(\int_X^{\oplus}\!M_p\,d\mu(p)\right)$,
we have:
$$(\int_X^{\oplus}\!N_p\,d\mu(p),\int_X^{\oplus}\!M_p\,d\mu(p),\int_X^{\oplus}\!\alpha_p\,d\mu(p),\int_X^{\oplus}\!\beta_p\,d\mu(p),\cdots$$
$$\cdots\int_X^{\oplus}\!\Gamma_p\,d\mu(p),\int_X^{\oplus}\!\nu_p\,d\mu(p),\int_X^{\oplus}\!T_L^p\,d\mu(p),\int_X^{\oplus}\!T_R^p\,d\mu(p))$$
is a measured quantum groupoid.
\end{prop}

\begin{proof}
Left to the reader.
\end{proof}

\cite{Bl} gives examples. In this case, the basis is
$L^{\infty}(X)$ and $\alpha=\beta=\hat{\beta}$. The fundamental
unitary comes from a space onto the same space and then can be
viewed as a field of multiplicative unitaries.

\backmatter

\end{document}